\documentclass[12pt,a4paper,reqno]{amsart}
\usepackage[a4paper,top=3cm,bottom=3cm,inner=3cm,outer=3cm,includehead,includefoot]{geometry}
\usepackage{amsmath,amssymb,amsthm,bbm,bm,mathtools,calc,verbatim,tikz,url,hyperref,cite,fullpage,mathtools}
\usetikzlibrary{shapes.misc,calc,intersections,patterns,decorations.pathreplacing,backgrounds}
\hypersetup{colorlinks=true,citecolor=blue,linktoc=page}

\newtheorem{theorem}{Theorem}
\newtheorem{lemma}[theorem]{Lemma}
\newtheorem{prop}[theorem]{Proposition}
\newtheorem{claim}[theorem]{Claim}

\newtheorem{obs}[theorem]{Observation}
\theoremstyle{definition}
\newtheorem{definition}[theorem]{Definition}
\newtheorem{remark}[theorem]{Remark}

\newenvironment{clmproof}[1]{\begin{proof}[Proof of Claim~\ref{#1}]\let\qednow\qedsymbol\renewcommand{\qedsymbol}{}}{\; \qednow \end{proof}}

\numberwithin{theorem}{section}
\mathtoolsset{showonlyrefs}

\everydisplay{\thickmuskip=7mu}
\renewcommand{\leq}{\leqslant}
\renewcommand{\geq}{\geqslant}
\renewcommand{\le}{\leqslant}
\renewcommand{\ge}{\geqslant}
\renewcommand{\to}{\rightarrow}

\newcommand{\diam}{\operatorname{diam}}
\newcommand{\aff}{\operatorname{aff}}

\newcommand{\interior}{\operatorname{int}}
\newcommand{\ext}{\operatorname{ext}}
\newcommand{\fext}{\operatorname{ext^\to}}

\newcommand{\fret}{\operatorname{ret^\leftarrow}}
\newcommand{\clint}{\operatorname{\underline{int}}}

\newcommand{\sgn}{\operatorname{sgn}}

\let\eps\varepsilon

\def\A{\mathcal{A}}
\def\cB{\mathcal{B}}
\def\B{\mathcal{B}}

\def\C{\mathcal{C}}
\def\D{\mathcal{D}}
\def\E{\mathbb{E}}
\def\Ex{\mathbb{E}}
\def\cE{\mathcal{E}}
\def\F{\mathcal{F}}
\def\G{\mathcal{G}}
\def\H{\mathbb{H}}
\def\HH{\mathcal{H}}
\def\cH{\mathcal{H}}
\def\bH{\mathbf{H}}

\def\L{\mathbb{L}}
\def\cL{\mathcal{L}}
\def\N{\mathbb{N}}
\def\cN{\mathcal{N}}

\def\P{\mathbb{P}}
\def\Pr{\mathbb{P}}
\def\cP{\mathcal{P}}
\def\Q{\mathbb{Q}}
\def\QQ{\mathcal{Q}}
\def\R{\mathbb{R}}

\def\S{\mathcal{S}}
\def\SS{\mathbf{S}}
\def\T{\mathcal{T}}
\def\U{\mathcal{U}}

\def\W{\mathcal{W}}
\def\X{\mathcal{X}}
\def\Y{\mathcal{Y}}
\def\Z{\mathbb{Z}}

\def\x{\mathbf{x}}

\def\<{\langle}
\def\>{\rangle}

\def\0{\mathbf{0}}
\def\ih{\mathrm{IH}}
\def\iha{\mathrm{IH_a}}
\def\ihb{\mathrm{IH_b}}

\def\edge{\partial}

\def\Vor{\mathrm{Vor}}

\def\Cell{\mathrm{Cell}}

\def\Cbar{\overline{C}}
\def\Fbar{\overline{F}}

\def\Qbar{\overline{\QQ}}
\def\1{\mathbbm{1}}
\def\bfdelta{\boldsymbol\delta}

\tikzstyle{closed}=[draw, fill, circle, minimum size=0.1cm, inner sep=0pt]
\tikzstyle{open}=[draw, circle, minimum size=0.1cm, inner sep=0pt]

\title{The critical length for growing a droplet}

\author[P. Balister \and B. Bollob\'as \and R. Morris \and P.J. Smith]{Paul Balister \and B\'ela Bollob\'as \and Robert Morris \and Paul Smith}

\address{Mathematical Institute, University of Oxford, Radcliffe Observatory Quarter, Woodstock Road, Oxford, OX2 6GG, UK}\email{Paul.Balister@maths.ox.ac.uk}

\address{Department of Pure Mathematics and Mathematical Statistics, Wilberforce Road, Cambridge, CB3 0WA, UK, and Department of Mathematical Sciences, University of Memphis, Memphis, TN 38152, USA}
\email{b.bollobas@dpmms.cam.ac.uk}

\address{IMPA, Estrada Dona Castorina 110, Jardim Bot\^anico, Rio de Janeiro, 22460-320, Brazil}
\email{rob@impa.br}

\address{Clerkenwell, London}
\email{paulsmith@cantab.net}

\thanks{P.B.\ and B.B.\ were partially supported by NSF grant DMS~1855745, R.M.\ by FAPERJ (Proc.~E-26/200.977/2021) and CNPq (Proc.~303681/2020-9), and by ERC Starting Grant 680275 MALIG, and P.S. by ERC Starting Grant 676970, by ISF grants 1147/14 and 1207/15, and by a CNPq bolsa PDJ} 

\date{\today}

\begin{document}

\begin{abstract}
In many interacting particle systems, relaxation to equilibrium is thought to occur via the growth of `droplets', and it is a question of fundamental importance to determine the critical length at which such droplets appear. In this paper we construct a mechanism for the growth of droplets in an arbitrary finite-range monotone cellular automaton on a $d$-dimensional lattice. Our main application is an upper bound on the critical probability for percolation that is sharp up to a constant factor in the exponent. Our method also provides several crucial tools that we expect to have applications to other interacting particle systems, such as kinetically constrained spin models on $\Z^d$. 

This is one of three papers that together confirm the Universality Conjecture of Bollob\'as, Duminil-Copin, Morris and Smith. 
\end{abstract}

\maketitle

\tableofcontents

\section{Introduction}

Consider a collection of particles on a $d$-dimensional lattice, interacting over a finite range. Important examples include the Ising model of ferromagnetism~\cite{Hugo20,FV17,M99} and kinetically constrained models of the liquid-glass transition~\cite{CMRT2,GST,T22} (see Section~\ref{KCM:sec}). In many such settings, the dynamical relaxation of the system is thought to occur via the growth of so-called `droplets' (see, e.g.,~\cite{BH16,H04,Sch98}), and it is therefore an important and fundamental problem to understand the most likely way in which such droplets grow.

In this paper we construct a mechanism for the growth of a droplet in an arbitrary (finite-range) monotone cellular automaton in $d$ dimensions. These models are monotone versions of arbitrary finite-range Ising models, and results for specific automata have previously been applied to the study of the Ising model (see, e.g.,~\cite{CMnuc2,FSS,MGlauber}) and kinetically constrained spin models (see~\cite{HMT1,HMT2,MMT}). We shall use our mechanism to prove an upper bound on the critical length for percolation (infection of the entire vertex set) that is, in a certain sense, optimal. Together with the results of~\cite{BBMSlower,BBMSsub}, this confirms the so-called Universality Conjecture of Bollob\'as, Duminil-Copin, Morris and Smith~\cite{BDMS}. 

In a monotone cellular automaton, a set of `infected' sites grows according to a deterministic (and homogeneous) update rule. The study of specific monotone cellular automata was initiated by Chalupa, Leath and Reich~\cite{CLR} in 1979, who introduced the family of models now known as `$r$-neighbour bootstrap percolation' (see below). The early work on bootstrap percolation~\cite{Adler,ADE,AL,And93,AMS,vE,Mount,Mount92,Sch1,Sch3} was motivated by and closely related to contemporaneous developments in the study of metastability~\cite{DS1,NS91,NS92,KO,MOS91,Sch91,Sch4,SS98}. Over the past 25 years a great deal of progress has been made in our understanding of the $r$-neighbour process~\cite{BBM3d,BBDM,CC,CM,HM,Hol} and other related automata~\cite{BDMSDuarte,DH,DE,DEH}.

The study of completely general monotone cellular automata, however, was initiated\footnote{Large families of monotone cellular automata were studied earlier by Gravner and Griffeath~\cite{GG93,GG99} and by Duminil-Copin and Holroyd~\cite{DH}.} 
only recently, by Bollob\'as, Smith and Uzzell~\cite{BSU}, who introduced the following general family of models, called \emph{\,$\U$-bootstrap percolation}.

\begin{definition}\label{def:MCA}
Let\/ $\U = \big\{ X_1,\ldots,X_m \big\}$ be a finite collection of finite, non-empty subsets of $\Z^d \setminus \{\0\}$. Let $A \subset \Z_n^d$ be a set of initially \emph{infected} sites, set $A_0 = A$, and define
$$A_{t+1} = A_t \cup \big\{ x \in \Z_n^d \,:\, x + X \subset A_t \text{ for some } X \in \U \big\}$$
for each $t \ge 0$. We write $[A]_\U = \bigcup_{t\geq 0} A_t$ for the set of eventually infected sites.
\end{definition}

We call any such collection $\U$ an \emph{update family}, and each $X \in \U$ an \emph{update rule}. For example, the classical $r$-neighbour process is the $\U$-bootstrap process whose update family $\cN_r^d$ consists of the ${2d \choose r}$ subsets of size $r$ of the $2d$ nearest neighbours of the origin. 

Motivated by the applications to statistical physics mentioned above, we are interested in understanding the behaviour of this process when the initial set $A$ is chosen randomly. To make this precise, let us say that $A$ is a \emph{$p$-random} subset of $\Z_n^d$ if each $x \in \Z_n^d$ is included in $A$ independently with probability $p$, and write $\P_p$ for the associated product probability measure. We say that $A$ \emph{percolates} if $[A]_\U = \Z_n^d$. Our main aim is to determine the behaviour (for fixed $d \in \N$ and as $n \to \infty$) of the \emph{critical probability} 
\begin{equation}\label{def:pc}
p_c\big( \Z_n^d,\U \big) := \inf \Big\{ p \in (0,1] \,:\, \P_p\big( [A]_\U = \Z_n^d \big) \ge 1/2 \Big\}
\end{equation}
for each $d$-dimensional update family $\U$. Similarly, the \emph{critical length} of $\U$ is
$$L_c\big( \U,p \big) := \sup \Big\{ n \in \N \,:\, \P_p\big( [A]_\U = \Z_n^d \big) \le 1/2 \Big\},$$
and we are interested in the rate of growth of $L_c(\U,p)$ as $p \to 0$. 

For the $r$-neighbour process on $\Z^d$, it was first proved by Schonmann~\cite{Sch1} that\footnote{More precisely, Schonmann proved that for any $p > 0$, a $p$-random set percolates almost surely in the $r$-neighbour process on $\Z^d$, which is defined by replacing $\Z_n^d$ by $\Z^d$ in Definition~\ref{def:MCA}.}
$$p_c\big( \Z^d,\cN_r^d \big) = 0$$
for every $d \ge r \ge 1$. If $r > d$ then it is not hard to see that $p_c\big( \Z^d,\cN_r^d \big) = 1$, so Schonmann's theorem is, in one sense, best possible. In fact, a more careful analysis of his proof reveals a bound of the form 
$$p_c\big( \Z_n^d, \cN_r^d \big) \, \le \, \bigg( \frac{C(d,r)}{\log_{(r-1)} n} \bigg)^{d-r+1},$$
for some constant $C(d,r) > 0$, where $\log_{(r)}$ denotes an $r$-times iterated logarithm, so $\log_{(0)} n = n$ and $\log_{(r)}n = \log \log_{(r-1)} n$ for each $r \geq 1$. This bound was shown to be best possible up to a constant factor by Aizenman and Lebowitz~\cite{AL} (in the case $r = 2$), by Cerf and Cirillo~\cite{CC} (in the case $d = r = 3$) and by Cerf and Manzo~\cite{CM} (in general), and a sharp threshold was determined by Holroyd~\cite{Hol} (in the case $d = r = 2$) and by Balogh, Bollob\'as, Duminil-Copin and Morris~\cite{BBM3d,BBDM} (in general), who proved that 
$$p_c\big( \Z_n^d, \cN_r^d \big) \, = \, \bigg( \frac{\lambda(d,r) + o(1)}{\log_{(r-1)} n} \bigg)^{d-r+1}$$
as $n \to \infty$, where $\lambda(d,r) > 0$ is an explicit constant. Despite considerable interest in such models, the order of the critical probability has been determined for only a small number of additional three-dimensional update families~\cite{Blanq19,Blanq22,vEF}. 

In this paper we develop an analogue of Schonmann's method for an arbitrary finite-range model. In particular, this (together with the main result of~\cite{BBMSsub,HS22}) allows us to characterize the update families $\U$ for which $p_c( \Z^d,\U ) = 0$. It also allows us to prove an upper bound on $p_c( \Z_n^d,\U )$ that is (by the results of~\cite{BBMSlower}) not far from best possible, and implies that \emph{every}\, $\U$-bootstrap process resembles (in some weak sense) one of the $r$-neighbour models. We begin by stating a non-technical version of our main theorem; in order to do so, let us first recall some important definitions from~\cite{BBPS,BDMS,BSU}.

One of the key insights of Bollob\'as, Smith and Uzzell~\cite{BSU}, who were the first to study completely general update families in two dimensions, was that the typical behaviour of the $\U$-bootstrap process should, roughly speaking, be determined by its action on discrete half-spaces. Given an update family $\U$, we define the \emph{stable set} of $\U$ to be\footnote{We write $\SS^{d-1}$ to denote the unit sphere in $\R^d$, that is, $\SS^{d-1} = \{ x \in \R^d : \|x\|_2 = 1 \}$.}
\begin{equation}\label{def:S}
\S(\U) := \big\{ u \in \SS^{d-1} :\, [ \H_u ]_\U = \H_u \big\},
\end{equation}
where, for each $u \in \SS^{d-1}$, we write $\H_u := \{ x \in \Z^d : \< x,u \> < 0 \}$ for the discrete half-space with normal $u$. Note that $u \in \S(\U)$ if and only if $X \not\subset \H_u$ for each $X \in \U$. 

The following classification of $d$-dimensional update families was proposed by Bollob\'as, Smith and Uzzell~\cite{BSU} in the case $d = 2$, by Balister, Bollob\'as, Przykucki and Smith~\cite{BBPS} for subcritical families, and by Bollob\'as, Duminil-Copin, Morris and Smith~\cite{BDMS} in general. Given a sphere $\SS \subset \R^d$ of arbitrary dimension and a set $\T \subset \R^d$, we write $\interior_\SS(\T)$ for the interior of $\T$ in $\SS$ with respect to the topology induced by geodesic distance.

\begin{definition}\label{def:tri}
A $d$-dimensional update family $\U$ with stable set $\S = \S(\U)$ is:
\begin{itemize}
\item \emph{supercritical} if $H \cap \S = \emptyset$ for some open hemisphere $H \subset \SS^{d-1}$;\smallskip
\item \emph{critical} if there exists a hemisphere $H \subset \SS^{d-1}$ such that $\interior_{\SS^{d-1}}(H \cap \S) = \emptyset$ and if $H \cap \S \neq \emptyset$ for every open hemisphere $H \subset \SS^{d-1}$;\smallskip
\item \emph{subcritical} if $\interior_{\SS^{d-1}}(H \cap \S) \neq \emptyset$ for every hemisphere $H \subset \SS^{d-1}$.
\end{itemize}
\end{definition}

In this paper we prove the upper bounds in the following theorem; the matching lower bounds are proved in two companion papers~\cite{BBMSlower,BBMSsub} (see also~\cite{HS22}). 
The theorem confirms a conjecture of Bollob\'as, Duminil-Copin, Morris and Smith~\cite{BDMS}. 

\begin{theorem}\label{conj:universality}
Let\/ $\U$ be a $d$-dimensional update family.
\begin{itemize}
\item[$(a)$] If\/ $\U$ is supercritical then $p_c\big( \Z_n^d,\U \big) = n^{-\Theta(1)}$.\smallskip
\item[$(b)$] If\/ $\U$ is critical then there exists $r \in \{2,\dots,d\}$ such that
$$p_c\big( \Z_n^d,\U \big) = \bigg(\frac{1}{\log_{(r-1)} n}\bigg)^{\Theta(1)}.$$
\item[$(c)$] If\/ $\U$ is subcritical then $p_c\big( \Z^d,\U \big) > 0$.
\end{itemize}
\end{theorem}

We shall define the quantity $r = r(\U)$ explicitly in Section~\ref{sec:r}, and use it to state a more precise version of the theorem above (see Theorem~\ref{thm:universality}). Both theorems can equivalently be stated in terms of the critical length of\, $\U$, for which they say that $L_c(\U,p)$ is polynomial if $\U$ is supercritical; grows like a tower of exponentials of height $r-1$ if $\U$~is critical; and is infinite for all sufficiently small $p > 0$ if\, $\U$ is subcritical. 

Theorem~\ref{conj:universality} was originally conjectured in the case $d = 2$ by Bollob\'as, Smith and Uzzell~\cite{BSU}, who moreover proved their conjecture for supercritical and critical two-dimensional families. Balister, Bollob\'as, Przykucki and Smith~\cite{BBPS} then completed the proof of the conjecture when $d = 2$, by showing that the critical probability for subcritical two-dimensional families is bounded away from zero, and also conjectured that $p_c(\Z^d,\U) > 0$ if and only if\/ $\U$ is subcritical. More recently, Bollob\'as, Duminil-Copin, Morris and Smith~\cite{BDMS} determined the critical probability up to a constant factor for all critical two-dimensional families; see also~\cite{BDMSDuarte,DE,DH,DHar,DEH} for related results.

In higher dimensions the problem becomes much harder, and results have previously been obtained in only a few special cases: the $r$-neighbour process (see above), and a small number of other specific models~\cite{Blanq19,vEF,HMod}. For each of the update families studied previously, the upper bounds on $p_c(\Z_n^d,\U)$ have been relatively easy to prove, and the main challenge has been to prove corresponding lower bounds. For general models, however, this is no longer true, and both upper and lower bounds pose unique and distinct challenges. Indeed, in order to prove the upper bounds in Theorem~\ref{conj:universality} we need to develop a number of significant new tools and techniques in order to deal with arbitrary update families. In Section~\ref{sec:outline} we describe some of the main challenges that we face in proving the upper bound, and outline how we overcome them. 

One indication of the difficulty of the problem is that we are only able to bound the critical probability up to a constant factor in the exponent. It would be natural for the reader to wonder whether the upper bound we prove in this paper is likely to be sharp. In fact the exponents we obtain are very far from best possible, and we make no attempt to optimize them. This is partly to simplify the proof, but mainly because it turns out that the correct exponent is in fact uncomputable in general! Indeed, in a companion paper~\cite{BBMSuncomp} we show that for each $1 \le r < d$, an algorithm that determines for every update family $\U$ with $r(\U) = r$ whether the exponent is at most $2/3$ or at least $1$  would also solve the halting problem. It is a very interesting (and likely difficult) open problem to determine the exponent for $d$-dimensional update families with $r(\U) = d$ (which we expect to be computable). More generally, one might hope to determine the exponent for all $\U$ as a function of the (uncomputable) exponents for supercritical update families. 

In Sections~\ref{sec:r} and~\ref{sec:outline} we shall state a more precise version of our main theorem, and give an outline of its proof. The proof itself is given in Sections~\ref{sec:memphis}--\ref{proof:sec} and Appendices~\ref{polytope:app}--\ref{cover:app}. Before embarking on this journey, however, let us provide some further motivation for Theorem~\ref{conj:universality}, by discussing a potential application of the techniques introduced in this paper to the study of kinetically constrained models of the liquid-glass transition.

\subsection{Kinetically constrained models}\label{KCM:sec}

A glass is a disordered material that nevertheless behaves mechanically like a solid, and is formed by rapidly cooling a viscous liquid. Understanding this liquid-glass transition is an important open problem in condensed matter physics, see for example~\cite{ALBB,DS01}. Kinetically constrained models were introduced in the 1980s (see~\cite{FA}, or the reviews~\cite{GST,RS03}) in order to model the liquid-glass transition, and exhibit several key properties of super-cooled liquids near the glass transition point.

Perhaps the simplest way to understand a kinetically constrained model is as a biased random walk on the family of percolating sets\footnote{More precisely, this is true if the initial set of empty sites percolates for the $\U$-bootstrap process. It follows from Theorem~\ref{conj:universality} that at equilibrium (i.e., if the empty sites are $p$-random) then this holds almost surely for every $p > 0$ if and only if $\U$ is not subcritical.} in $\U$-bootstrap percolation on $\Z^d$. The state (either `empty' or `occupied') of a site $x \in \Z^d$ can update only if the set $x + X$ is entirely empty for some $X \in \U$; when this occurs, it updates at rate~$1$, becoming empty with probability $p$, and occupied with probability $1-p$. Well-studied examples of kinetically constrained models include the East model (see, e.g.,~\cite{CFM}), whose update family consists of the single set $\{-1\}$ when $d = 1$, and the $r$-facilitated Friedrickson--Andersen model, introduced in~\cite{FA}, which corresponds to $r$-neighbour bootstrap percolation.

The connection with bootstrap percolation was first observed in a seminal paper by Cancrini, Martinelli, Roberto and Toninelli~\cite{CMRT}, who were also the first to study general kinetically constrained models. In particular, they showed that the infection time 
$$\tau(\Z^d,\U) \, := \, \inf\big\{ t \ge 0 : \0 \text{ is empty at time $t$}\big\}$$
starting from equilibrium 
is almost surely finite if $p > p_c(\Z^d,\U)$. It therefore follows from Theorem~\ref{conj:universality} that $\tau(\Z^d,\U)$ is almost surely finite for all $p > 0$ if and only if\, $\U$ is not subcritical. Moreover, a lower bound on the mean infection time is given (up to a constant factor) by the median infection time of the origin in the corresponding $\U$-bootstrap process (see~\cite[Lemma~4.3]{MT}). However, due to the more complex (non-monotone) behaviour of kinetically constrained models, the scaling of the infection times in the two models are in general qualitatively different.

Over the past few years, there have been some dramatic advances in our understanding of kinetically constrained models in two dimensions, mirroring those in the study of $\U$-bootstrap processes. Motivated by the connection with $\U$-bootstrap percolation, a number of conjectures were made in~\cite{M17} regarding the rate of growth of the typical infection time as $p \to 0$ for critical update families in two dimensions and supercritical update families in $d$ dimensions. These have now all been either proved or disproved (see~\cite{HMT1,HMT2,MaMaT,MMT}), and in two dimensions the situation is now extremely well-understood, with the full universality picture determined (see~\cite{Har,HarMa,HMT3} for the most recent developments). In particular, the lower bound given by coupling with the $\U$-bootstrap process is sometimes sharp up to a constant factor in the exponent (for example, for the $2$-neighbour model, see~\cite{HMT3}), and sometimes not\footnote{More precisely, for many models the existence of certain `energy barriers' dominates the expected infection time. As a result of this, both the critical and supercritical families need to be partitioned into two different universality classes (and then further refined for logarithmic corrections, see~\cite{Har,HarMa}).} (for example, for the Duarte model, see~\cite{MaMaT}). 

In higher dimensions, it follows from Theorem~\ref{conj:universality} that
\begin{equation}\label{eq:KCMlower}
\Ex\big[ \tau(\Z^d,\U) \big] \ge \exp_{(r-1)}( p^{-c} ),
\end{equation}
for every critical $d$-dimensional update family $\U$, where $c = c(\U) > 0$ is a constant, $r = r(\U) \in \{2,\ldots,d\}$ is defined in Section~\ref{sec:r} (see Definition~\ref{def:r}), and we write $\exp_{(r)}$ for an $r$-times iterated exponential, so $\exp_{(0)} n = n$ and $\exp_{(r)}n = e^{\exp_{(r-1)} n}$ for each~$r \ge 1$. It seems reasonable to conjecture that~\eqref{eq:KCMlower} is sharp up to the value of the constant $c$ for all critical families (even though this is not true for supercritical families, see~\cite{MaMaT}), and we expect the techniques introduced in this paper to play a central role in the proof of this conjecture. Indeed, the general method for proving such upper bounds developed in the papers~\cite{CMRT,MT,MMT,HMT1,Har} relies on the existence of a `low-energy' mechanism to infect droplets of roughly the `critical' size, which is precisely what our method supplies in the setting of $\U$-bootstrap percolation. Nevertheless, we expect there to be significant technical challenges involved in making such an approach rigorous.

\section{The resistance of an update family}\label{sec:r}

The primary objective of this section is to define explicitly, in Definition~\ref{def:r}, the parameter $r = r(\U)$ for an arbitrary $d$-dimensional update family, which we shall refer to as the \emph{resistance} of\, $\U$. Having done so, we will state a refined version of Theorem~\ref{conj:universality} (see Theorem~\ref{thm:universality}), and then derive a few simple consequences of the definition of $r(\U)$; in~particular, we show in Lemma~\ref{lem:tri} that $r(\U) \in \{2,\ldots,d\}$ for every critical update family. The idea behind the definition 
is quite simple, but the details are somewhat technical, so the reader may find it helpful to have in mind pictures of the stable sets of the 2- and 3-neighbour models in three dimensions. For $\cN_2^3$, the stable set is the six points $\{\pm e_1, \pm e_2, \pm e_3\}$, and for $\cN_3^3$ the stable set is the union of the three great circles orthogonal to the standard basis vectors.

The first step is to introduce a family of objects called `$\SS$-stable sets', which are subsets of the sphere $\SS$, and generalize the notion of a stable set of an update family. We define these next, and then prove a simple property of such sets (Lemma~\ref{lem:inducedwelldef}), which essentially says that if one takes an $\SS$-stable set $\T$, and a point $u$ on the sphere $\SS$ in which $\T$ is embedded, then the sets\footnote{Here, and throughout the paper, we write $\| \cdot \|$ for the Euclidean norm on $\R^d$.}
$$\T \cap \big\{ v \in \SS : \| u - v \| = \eta \big\}$$
all `look the same' whenever $\eta$ is sufficiently small. This property will enable us to show that $r(\U)$ (see Definition~\ref{def:r}) is well-defined.

\begin{definition}\label{def:happy}
Let $\T\subset\SS^{d-1}$ and let $\SS$ be a sphere (of arbitrary size and dimension) embedded in $\SS^{d-1}$.
We say that $\T$ is \emph{$\SS$-stable} if there exists a finite collection $\cH_1,\ldots,\cH_m$ of finite families of closed hemispheres of $\SS$, such that
\begin{equation}\label{eq:happy}
\T \cap \SS = \bigcap_{i=1}^m \bigcup_{H \in \cH_i} H.
\end{equation}
\end{definition}

To see why this is a natural definition, let us observe that the stable set of a $d$-dimensional update family is $\SS^{d-1}$-stable.

\begin{lemma}\label{lem:hemispheres}
If\/ $\U$ is a $d$-dimensional update family, 
then
$$\S(\U) = \bigcap_{X \in \, \U} \bigcup_{x \in X} \big\{ u\in \SS^{d-1} : \< x,u \> \geq 0 \big\}.$$
In particular, $\S(\U)$ is\/ $\SS^{d-1}$-stable.
\end{lemma}

\begin{proof}
Simply note that a direction $u \in \SS^{d-1}$ is unstable for $\U$ (that is, $u \not\in \S(\U)$) if and only if there exists $X \in \U$ such that $\< x, u \> < 0$ for every $x \in X$.
\end{proof}

The function $r = r(\U)$ will depend on both local and global properties of the stable set $\S(\U)$. Given a sphere $\SS\subset\SS^{d-1}$, $u\in\SS$ and $\eta>0$, define the sub-sphere
\[
S_\eta(\SS,u) := \big\{ v\in \SS : \|u - v\| = \eta \big\}.
\]
We say that two subsets of $\R^d$ are \emph{equivalent} if one can be obtained from the other by a composition of translations, dilations and rotations. We remark that all of the equivalences used in this paper will be proved using either a rotation or a homothety\footnote{That is, a map of the form $M \colon x \mapsto \lambda x + a$, for some $\lambda > 0$ and $a \in \R^d$.}. We write $A\equiv B$ if $A$ is equivalent to $B$.

\begin{lemma}\label{lem:inducedwelldef}
Let $\SS \subset \SS^{d-1}$ be a sphere, let $\T \subset \SS^{d-1}$ be $\SS$-stable, and let $u \in \SS$. Then there exists $\eta_0 = \eta_0(\SS,\T,u) > 0$ such that, for all $0 < \eta < \eta_0$, 
$$\T \cap S_\eta(\SS,u) \equiv \T \cap S_{\eta_0}(\SS,u)$$
and $\T$ is $S_\eta(\SS,u)$-stable.
\end{lemma}

\begin{proof}
Since $\T$ is $\SS$-stable, it can be written as a finite number of intersections and unions of closed hemispheres in $\SS$. Hence, by choosing $\eta$ small enough, we can ensure that if $S_\eta(\SS,u)$ intersects the boundary of one of those hemispheres, then that boundary also passes through $u$. Now, let $\varphi \colon \R^d \to \R^d$ be the unique homothety such that $\varphi\big( S_\eta(\SS,u) \big) = S_{\eta_0}(\SS,u)$, and observe that
$$\varphi\big( \T \cap S_\eta(\SS,u) \big) = \T \cap S_{\eta_0}(\SS,u)$$
if $0 < \eta < \eta_0$ and $\eta_0$ is sufficiently small, by the comments above. Moreover, $\T$ is a $S_\eta(\SS,u)$-stable set, since each closed hemisphere that defines $\T$ either contains $S_\eta(\SS,u)$, avoids it, or intersects it in a hemisphere. 
\end{proof}

In the following definition (and throughout the paper) we shall use Lemma~\ref{lem:inducedwelldef} implicitly by writing $S_\eta(\SS,u)$ without specifying $\eta$. This should always be taken to mean that $\eta$ is smaller than the $\eta_0 = \eta_0(\SS,\T,u)$ of Lemma~\ref{lem:inducedwelldef} (for the $\T$ currently under consideration). 

For each $k \in \{0,\ldots,d - 1\}$, let us write $\C^k$ for the set of all $k$-dimensional spheres embedded in $\SS^{d-1}$, and for each $\SS \in \C^k$, let us write $\B(\SS)$ for the set of all $\SS$-stable sets. We are now ready to define the resistance of a $d$-dimensional update family. 

\begin{definition}\label{def:r}
For each $1 \leq k \leq d$ and each $\SS \in \C^{k-1}$, we define two functions
\begin{align*}
\rho^{k-1}(\SS; \, \cdot, \, \cdot) &\colon \B(\SS) \times \SS \to \{0,1,\dots,k\}\\
r^k(\SS; \, \cdot) &\colon \B(\SS) \to \{1,\dots,k+1\}
\end{align*}
inductively as follows. Set $r^0 \equiv 1$, and let $1 \leq k \leq d$, $\SS \in \C^{k-1}$, $\T \in \B(\SS)$ and $u \in \SS$. We define the \emph{induced resistance of $u$ with respect to $\T$ in $\SS$} to be
\begin{equation}\label{eq:rho}
\rho^{k-1}(\SS; \T, u) := 
\begin{cases} 
r^{k-1}\big( S_\eta(\SS,u); \T \big) \quad & \text{if } u \in \T, \\ 
0 & \text{otherwise}, 
\end{cases}
\end{equation}
and the \emph{resistance of $\T$ in $\SS$} to be
\begin{equation}\label{eq:r}
r^k(\SS; \T) := \min_H \, \max _{u \in H} \, \rho^{k-1}\big( \SS; \T, u \big) + 1,
\end{equation}
where the minimum is taken over all open hemispheres $H \subset \SS$.

Now, given a $d$-dimensional update family $\U$, define the \emph{resistance} of $\U$ to be
$$r(\U) := r^d\big( \SS^{d-1}; \S(\U) \big).$$
\end{definition}

In words, the resistance of\, $\U$ is determined by the largest induced resistance of a direction in the `easiest' hemisphere of $\SS^{d-1}$, and the induced resistance of a stable direction $u$ is given by the resistance of a small sphere centred at $u$. 

\pagebreak

The integer $r$ in Theorem~\ref{conj:universality} turns out to be equal to the resistance $r(\U)$: this is the precise version of the Universality Theorem.


\begin{theorem}\label{thm:universality}
Let $\U$ be a $d$-dimensional update family.
\begin{itemize}
\item[$(a)$] If\/ $\U$ is supercritical then $p_c(\Z_n^d,\U) = n^{-\Theta(1)}$.\smallskip
\item[$(b)$] If\/ $\U$ is critical then 
$$p_c(\Z_n^d,\U) = \bigg(\frac{1}{\log_{(r-1)} n}\bigg)^{\Theta(1)},$$
where $r = r(\U)$.\smallskip
\item[$(c)$] If\/ $\U$ is subcritical then $p_c(\Z^d,\U) > 0$.
\end{itemize}
\end{theorem}

Unsurprisingly, the definition of $r(\U)$ plays an extremely important role in the proof of Theorem~\ref{thm:universality}, so let us take some time to understand the various functions defined above. Roughly speaking, we think of the functions $\rho^{k-1}(\SS; \T, u)$ and $r^k(\SS; \T)$ as specifying, for certain `induced' processes (see Section~\ref{sec:chat-induced}) determined by $\SS$, $\T$ and $u$, which of the classical $r$-neighbour models they most closely resemble. These induced processes take place in sub-lattices of $\Z^d$, the superscript indicating the dimension of the sub-lattice, and can be thought of as the $\U$-bootstrap process `assisted' by one or more half-spaces. 

In order to develop some intuition, it is perhaps instructive to consider the classical $r$-neighbour model itself, for which the following holds.

\begin{obs}\label{obs:rneighbour}
For each $1 \le r \le d+1$, the update family $\cN_r^d$ has resistance $r$.  
\end{obs}

\begin{proof}
We prove the observation by induction on $d$; for $d = 1$ one easily checks that $\S(\cN_1^1) = \emptyset$ and $\S(\cN_2^1) = \SS^0$, so $\rho^0( \SS^0; \cN_1^1, u ) = 0$ and $\rho^0( \SS^0; \cN_2^1, u ) = 1$ for each $u \in \SS^0$, and hence
$$r(\cN_1^1) = 1 \qquad \text{and} \qquad r(\cN_2^1) = 2.$$ 

If $d \ge 2$, then let $\T := \S(\cN_r^d)$ be the stable set of the $r$-neighbour update family $\cN_r^d$, and observe that the elements of $\T$ are exactly those that have non-zero inner product with at most $r - 1$ of the standard basis vectors. It follows that 
$$\T \cap S_\eta(\SS^{d-1},u) \equiv \S(\cN_{r-1}^{d-1})$$ 
for every $u \in \{\pm e_1, \ldots, \pm e_d \}$, where $e_1,\dots,e_d$ are the standard basis vectors. Indeed, if $v \in S_\eta(\SS^{d-1},u)$, then $v \in \T$ if and only if $v$ has non-zero inner product with at most $r - 2$ of the standard basis vectors other than $u$. Hence, by the induction hypothesis,
\[
r^{d-1}\big( S_\eta(\SS^{d-1},u); \T \big) = r - 1,
\]
and therefore, since $u \in \T$, we have
\[
\rho^{d-1}\big( \SS^{d-1}; \T, u \big) = r - 1
\]
by Definition~\ref{def:r}.
Moreover, for any $u \in \SS^{d-1}$, the set $\T \cap S_\eta(\SS^{d-1},u)$ is equivalent to a subset of $\S(\cN_{r-1}^{d-1})$, since if $u$ has $i$ non-zero coordinates and $v \in S_\eta(\SS^{d-1},u)$, then $v \in \T$ if and only if $v$ has non-zero inner product with at most $r - i - 1$ of the remaining standard basis vectors. Hence, by the induction hypothesis, we have
\[
r^{d-1}\big( S_\eta(\SS^{d-1},u); \T \big) \leq r - 1,
\]
and therefore
\[
\rho^{d-1}\big( \SS^{d-1}; \T, u \big) \le r - 1
\]
for every $u \in \SS^{d-1}$, by Definition~\ref{def:r}. Since every open hemisphere in $\SS^{d-1}$ contains at least one element of the set $\{\pm e_1, \ldots, \pm e_d \}$, it follows from~\eqref{eq:rho} and~\eqref{eq:r} that $r(\cN_r^d) = r$, as claimed.
\end{proof}

We are now ready to consider general $d$-dimensional update families. First let us record in the following lemma the important fact that $2 \le r(\U) \le d$ for every critical family $\U$.

\begin{lemma}\label{lem:tri}
Let\/ $\U$ be a $d$-dimensional update family. The following hold:
\begin{itemize}
\item[$(a)$] $\U$ is supercritical if and only if\/ $r(\U) = 1$.\smallskip
\item[$(b)$] $\U$ is critical if and only if\/ $r(\U) \in \{2,\dots,d\}$.\smallskip
\item[$(c)$] $\U$ is subcritical if and only if\/ $r(\U) = d + 1$.
\end{itemize}
\end{lemma}

\begin{proof}
By Definition~\ref{def:tri}, the family $\U$ is supercritical if and only if $H \cap \S(\U) = \emptyset$ for some open hemisphere $H \subset \SS^{d-1}$, and by Definition~\ref{def:r}, 
$$\rho^{d-1} \big( \SS^{d-1}; \S(\U), u \big) \ge 1 \qquad \Leftrightarrow \qquad u \in \S(\U)$$ 
for every $u \in \SS^{d-1}$. It follows that $\U$ is supercritical if and only if there exists an open hemisphere $H \subset \SS^{d-1}$ such that $\rho^{d-1} \big( \SS^{d-1}; \S(\U), u \big) = 0$ for every $u \in H$, and such an $H$ exists if and only if $r(\U) = 1$, as required.

Next, suppose that $\U$ is subcritical, so, by Definition~\ref{def:tri}, for every open hemisphere $H \subset \SS^{d-1}$, the set $\interior_{\SS^{d-1}} \big( H \cap \S(\U) \big)$ is non-empty. This implies that $S_\eta(\SS^{d-1},u) \subset \S(\U)$ for some $u \in H \cap \S(\U)$, and therefore\footnote{Note that it follows from~\eqref{eq:rho} and~\eqref{eq:r} by induction on $k$ that if $\SS \subset \T$ for some $\SS \in \C^{k-1}$, with $1 \le k \le d$, then $r^k(\SS; \T) = k + 1$.}
$$r^{d-1}\big( S_\eta(\SS^{d-1},u); \S(\U) \big) = d.$$
Hence, for every open hemisphere $H \subset \SS^{d-1}$, there exists $u \in H$ with $\rho^{d-1}(\SS^{d-1}; \T, u) = d$, which implies that $r(\U) = d + 1$.

Since $1 \le r(\U) \le d+1$ for every $d$-dimensional update family $\U$ (by Definition~\ref{def:r}), and critical families are exactly those that are neither supercritical nor subcritical (by Definition~\ref{def:tri}), to complete the proof of the lemma it only remains to show that if $r(\U) = d + 1$, then $\interior_{\SS^{d-1}} \big( H \cap \S(\U) \big) \neq \emptyset$ for every hemisphere $H \subset \SS^{d-1}$. To show this, we shall prove the following more general statement by induction on $k$: For each $1 \le k \le d$, and each $\SS \in \C^{k-1}$, if $\T$ is an $\SS$-stable subset of $\SS^{d-1}$ and $r^k(\SS; \T) = k + 1$, then $\interior_\SS (H \cap \T) \neq \emptyset$ for every hemisphere $H \subset \SS$.

When $k = 1$, the condition $r^1(\SS; \T) = 2$ implies that $\SS \subset \T$, which implies the claim in this case since the topology induced by geodesic distance is the discrete topology for a 0-dimensional sphere. So let $\SS \in \C^{k-1}$, let $\T$ be $\SS$-stable with $r^k(\SS; \T) = k + 1$, and assume that the claim holds for smaller values of $k$. By~\eqref{eq:rho} and~\eqref{eq:r}, for every open hemisphere $H \subset \SS$, there exists an element $u \in H \cap \T$ such that
$$r^{k-1}\big( S_\eta(\SS,u); \T \big) = \rho^{k-1}( \SS; \T, u ) = k,$$
where we may assume that $\eta$ is small enough so that $S_\eta(\SS,u) \subset H$. By Lemma~\ref{lem:inducedwelldef} and the induction hypothesis, it follows that
$$\interior_{S_\eta(\SS,u)} (\T) \neq \emptyset,$$
for all sufficiently small values of $\eta$, and therefore (again using Lemma~\ref{lem:inducedwelldef}) the set
$$\T \cap \bigcup_{\eta' \le \eta} S_{\eta'}(\SS,u) \subset H$$
has non-empty interior in $\SS$, as required. This proves the induction step, and the induction hypothesis with $\SS = \SS^{d-1}$ and $\T = \S(\U)$ completes the proof of the lemma.
\end{proof}

\section{The main theorem and an outline of the proof}\label{sec:outline} 

In this paper we shall prove the following theorem, which implies the upper bounds in Theorem~\ref{thm:universality} for critical and supercritical update families. 

\begin{theorem}\label{thm:upper}
Let\/ $\U$ be a $d$-dimensional update family such that $r = r(\U) \leq d$. Then 
$$p_c\big( \Z_n^d,\U \big) \le \bigg( \frac{1}{\log_{(r-1)} n} \bigg)^{\eps(\U)},$$
for some $\eps(\U) > 0$ and all sufficiently large $n \in \N$.
\end{theorem}

It will also follow immediately from the proof of Theorem~\ref{thm:upper} (see Section~\ref{final:proof:sec}) that 
$$p_c\big( \Z^d,\U \big) = 0$$
for every non-subcritical $d$-dimensional update family $\U$. We remark that the bound given by Theorem~\ref{thm:upper} is in one sense best possible, since for each $d \ge r \ge 1$ and each $c > 0$, there exists a $d$-dimensional update family $\U$ with $r(\U) = r$ such that 
$$p_c\big( \Z_n^d,\U \big) \ge \big( \log_{(r-1)} n \big)^{-c}.$$
For example, it follows from standard techniques (see~\cite{CC,CM}) that this holds for the model obtained from the classical $r$-neighbour model by replacing each element $u$ of each update set by the set $\{u,2u,\ldots,ku\}$, for some $k > 1/c$. Moreover, as discussed in the introduction, one cannot expect to determine the optimal value of the constant $\eps(\U)$ for every update family $\U$, since it is shown in~\cite{BBMSuncomp} that this constant is in general uncomputable if $r < d$. We do not make any attempt to optimize the value of $\eps(\U)$ given by our proof; instead we have tried to simplify the argument wherever possible.  

In order to prove Theorem~\ref{thm:upper}, we need to overcome a number of significant technical challenges; in doing so, we shall develop a novel toolkit for studying the growth of droplets in interacting particle systems. The reader who is familiar with the bootstrap percolation literature may find the difficulty and complexity of the proof surprising, since upper bounds on $p_c( \Z_n^d,\U )$ are usually easier to prove than lower bounds. Indeed, to prove an upper bound one `only' has to find a single way in which to percolate, rather than deal with all possible ways. We therefore begin by explaining why the theorem should be true, describing the main challenges that we need to overcome in order to prove it, and outlining how we shall go about overcoming them. 

In order to get warmed up, let us first discuss why we should expect $r(\U)$ to control the probability of percolation in the $\U$-bootstrap process. For concreteness, let us consider the case $d = 3$ and $r(\U) = 3$, and let $H$ be an open hemisphere with $\rho^{2}\big( \SS^2; \S(\U), u \big) \le 2$ for every $u \in H$, whose existence is guaranteed (via~\eqref{eq:r}) by the assumption $r(\U) \le 3$. 

Roughly speaking, the idea is that for some constant $C > 0$, a `droplet'\footnote{Our droplets will always be polytopes, but the choice of the faces will be important and rather delicate, and will only be made later.} $D$ of size $\exp( p^{-C} ) = (\log n)^{o(1)}$ is likely to grow in the direction of the centre $w$ of $H$. This is because we expect the induced process (see Definition~\ref{def:induced}) 
on the face of $D$ in some direction $u \in \SS^2$ to resemble a two-dimensional process with resistance $\rho^{2}\big( \SS^2; \S(\U), u \big)$. Since the critical probability for percolation in such a process is poly-logarithmic in the size of the face (as was first proved by Bollob\'as, Smith and Uzzell~\cite{BSU}), we would expect growth to occur (with high probability) on any such face. 

The alert reader may have noticed that there are several serious problems with the sketch given above. First, what happens at the boundaries between faces of $D$? Second, how far can we grow in direction $w$ before we meet a face on which the induced process does not percolate? Third, the results of~\cite{BSU} (and, more to the point, the induction hypothesis) are valid on a torus, whereas the faces of $D$ are (non-toral) droplets. 

It is this third problem, which may seem like a trivial point at first sight, that actually appears to be fatal. To see this, consider a face $D_u$ of $D$ in some direction $u \in H$, and recall that $\rho^{2}\big( \SS^2; \S(\U), u \big) \le 2$ only implies that there exists a direction $w_u \in S_\eta(\SS^2,u)$ in which it is `easy' to grow on $D_u$. How do we plan to infect the sites near the boundary of $D_u$ in direction $-w_u$? Definition~\ref{def:r} appears to give us no help; indeed, for a general two-dimensional update family infecting these sites would not be possible. 

At this point it seems clear that Theorem~\ref{thm:upper} is false, and that we need a new definition of $r(\U)$. However, this turns out not to be the case, for a somewhat subtle reason. Remarkably, it turns out that, while it may indeed be harder to grow in the directions outside the `easiest' hemisphere, it is never \emph{much} harder, and the sense in which this is true is \emph{just} enough to allow for percolation of the entire faces. 

This fact about induced update families lies in contrast to the situation with the original $\U$-bootstrap percolation model, where growth can be genuinely biased (for example, the stable set could be equal to a closed hemisphere of $\SS^{d-1}$). There is thus an important difference between the $\U$-bootstrap process itself, and the corresponding family of induced update families. We will discuss this phenomenon in greater detail in Section~\ref{sec:chat-memphis}, and also in Sections~\ref{sec:memphis} and~\ref{sec:induced}.

The other two problems mentioned above also turn out to create significant difficulties. We overcome both by choosing carefully the set $\QQ$ of directions that we use as the faces of our droplets (see Section~\ref{sec:chat-quasi} and Section~\ref{sec:quasi}); this allows us to infect the sites at the boundary of two (or more) faces, and also (crucially) to infect some (not all) of the sites on the faces that are perpendicular to~$w$. In particular, this implies that the faces of our droplets grow, which allows the droplets to continue growing indefinitely without being likely to meet a face on which percolation fails. However, in order to make this argument work, it turns out that we must choose~$w$, the `easy' direction in which we are growing, more carefully (see Section~\ref{sec:chat-w} and Section~\ref{sec:w}). 

Having laid down these foundations in Sections~\ref{sec:memphis}--\ref{sec:induced}, it will be possible to begin defining the paths of infections that will be used to infect the torus. In order to do so, we need to develop a number of additional technical tools involving the interaction of the lattice with a certain family of polytopes, which we define in Section~\ref{sec:polytopes}. The properties of these polytopes will depend heavily on our choice of the set $\QQ$; in particular, in Section~\ref{sec:deterministic} we prove two deterministic lemmas (Lemmas~\ref{lem:forwards-step-general} and~\ref{lem:sideways-step-general}) that will be used (in Section~\ref{proof:sec}) to define our paths of infections. In order to bound the probability of these paths, we shall use induction on $d$ and $r$, and a useful method that was introduced by Schonmann~\cite{Sch1} in his foundational work on the $r$-neighbour model.

Having completed our brief informal overview, we next discuss in more detail the various issues mentioned above.

\subsection{Induced processes}\label{sec:chat-induced}

Our first important task is to decide what we mean by the `induced process' on a face of a droplet. This may seem straightforward, but there is in fact an important subtlety in our choice. In order to develop some intuition, let us first recall how growth occurs in the $r$-neighbour model; or, more precisely, how it is controlled in Schonmann's proof~\cite{Sch1} that percolation occurs for all fixed $p > 0$. 

In the $r$-neighbour model it suffices to consider droplets of the form $Q_k = [k]^d$, where $[k] = \{1,2,\dots,k\}$. In order to infect the sites of $Q_{k+1} \setminus Q_k$, we begin with the faces of co-dimension $1$, and then work our way downwards. Specifically, observe first that each site of the face $\{k+1\} \times [k]^{d-1}$ has exactly one neighbour in $Q_k$, and so, if all of the sites of $Q_k$ have already been infected, then we may couple growth on this face with the $(r-1)$-neighbour process on a finite subset of $\Z^{d-1}$. Note that, in doing so, we use only those update rules $X \in \cN_r^d$ with $e_1 \notin X$; that is, with $\<x,e_1\> \le 0$ for every $x \in X$. 

Similarly, the face $\{k+1\}^\ell \times [k]^{d-\ell}$ has exactly $\ell$ neighbours in higher-dimensional faces, and so, if all of these faces have already been infected, then we may couple growth on this face with the $(r-\ell)$-neighbour process on a finite subset of $\Z^{d-\ell}$. Once again, note that we use only those $X \in \cN_r^d$ with $\<x,e_i\> \le 0$ for every $x \in X$ and $i \in [\ell]$. 

In the general $\U$-bootstrap setting, the idea is the similar, but we need to work in an arbitrary sub-lattice of $\Z^d$. For a set $W \subset \SS^{d-1}$, let us write
$$\L(W) := W^\perp \cap \Z^d = \big\{ x \in \Z^d : \<x,u\> = 0  \text{ for all } u \in W \big\}$$
for the sub-lattice perpendicular to $W$, and define
\begin{equation}\label{def:HW}
\HH(W) := \big\{ x \in \Z^d \,:\, \< x,u \> \leq 0 \text{ for all } u \in W \big\}.
\end{equation}
We define the update family induced by the set of directions $W$ as follows.

\begin{definition}\label{def:induced}
Let $\U$ be a $d$-dimensional update family. Given $W \subset \SS^{d-1}$, define\footnote{We remark that, in contrast to Definition~\ref{def:MCA}, we may have $\emptyset \in \U[W]$.}
\begin{equation}\label{eq:induced-def}
\U[W] := \big\{ X \cap W^\perp \,:\, X \in \U \,\text{ and }\, X \subset \HH(W) \big\}.
\end{equation}
We call $\,\U[W]$ the \emph{sub-update family of\/ $\U$ induced by $W$}.
\end{definition}

The reader should think of $W$ as being a subset of the directions of the $(d-1)$-dimensional faces of our droplet $D$, and of\, $\U[W]$ as being an approximation of the process that occurs in the 
`$W$-face' of $D$, which is the intersection\footnote{Of course, this intersection may be empty, in which case this analogy fails; however, we shall only be interested in sets for which the intersection is non-empty.} of the faces in the directions of $W$. For example, if $\U = \cN_r^d$ and $W = \{e_1,\ldots,e_\ell\}$, then $\U[W]$ is equivalent\footnote{The two update families are not, strictly speaking, identical, since $\U[W]$ contains all subsets of $\{\pm e_{\ell+1},\ldots,\pm e_d\}$ of size at least $r - \ell$, but the processes that they define are equivalent.} to $\cN_{r-\ell}^{d-\ell}$. We remark that we allow $W = \emptyset$ (note that $\U[\emptyset] = \U$), and write $\U[u]$ for $\U[\{u\}]$. 

We think of $\U[W]$ as acting in the lattice $\L(W)$. This lattice has dimension $\dim(W^\perp)$ whenever $W$ is a set of \emph{rational directions},\footnote{This follows from Lemma~\ref{lem:w-gram-schmidt}, which is proved using the Gram--Schmidt algorithm. Note that we write $\dim(V)$ for the dimension of an affine subspace $V \subset \R^d$, and $\< W \>$ for the span of a set $W$, so, in particular, $\<x\> = \{\lambda x : \lambda \in \R\}$.} meaning that it is a subset of the set
\begin{equation}\label{eq:rational-sphere}
\SS_\Q^{d-1} := \big\{ w \in \SS^{d-1} : \<w\> \cap \Z^d \ne \{\0\} \big\}.
\end{equation}
That is, we say a direction $w$ is rational if the line $\<w\>$ contains a non-zero point of the lattice $\Z^d$. Note that $\SS_\Q^{d-1}$ is a dense\footnote{This is a special case of Lemma~\ref{lem:quasi-rationals-dense}.} and countable subset of $\SS^{d-1}$. In this paper we only need to consider update families induced by sets of rational directions. 

In order to develop some intuition, let us next discuss why the $\U[W]$-process on the lattice $\L(W)$ is a reasonable approximation for the $\U$-process growing on the $W$-face of a droplet $D$. Here, for simplicity, we focus on the case $|W| = 1$; in Section~\ref{sec:chat-quasi} we shall see that, for the special set of directions that we shall consider, general induced process can be constructed from this case by adding elements of $W$ one by one.

Let us therefore assume that $W = \{u\}$, and suppose that we are currently trying to infect a site $x$ of the $W$-face of $D$ (i.e., the face in direction $u$). If $x$ is sufficiently far from the edge of this face\footnote{For a discussion of growth near to the edge of the face, see Section~\ref{sec:chat-quasi}.} then, from the point of view of the process, $D$ may as well be the entire half-space $x + \H_u$. In particular, to infect $x$ it suffices to have infected the set $x + (X \cap \{u\}^\perp)$ for some $X \in \U$ with $X \subset \HH(\{u\})$. 

The reader may find it surprising that we do not need to include those rules $X \in \U$ with $X \not\subset \HH(W)$. The reason is that the crucial property we need -- that the stable set of $\U[u]$ is equivalent\footnote{We remark that here the stable set of $\U[u]$ is assumed to be a subset of $\SS(u) := \SS^{d-1} \cap \{u\}^\perp$.} to the intersection of $\S(\U)$ with a small sphere around $u$ -- is true even without these sets, as we shall see in Lemma~\ref{lem:induced-projection}. In particular, if $u$ belongs to the `easy' open hemisphere of directions (with respect to $\S(\U)$), then the resistance of $\U[u]$ will be at most $r(\U) - 1$, and so it will be possible to use induction to grow under $\U[u]$. 

Why should $\S(\U[u])$ be equivalent to the intersection of $\S(\U)$ with $S_\eta(\SS^{d-1},u)$? The idea is that if $v \in \SS^{d-1}$ is sufficiently close to $u$, then $X \subset \H_v$ if and only if
\begin{equation}\label{eq:induced-motivation}
X \subset \H_u \cup \big( \H_v \cap \{u\}^\perp \big).
\end{equation}
In particular, if $x \in X$ and $x \in \H_v \setminus \H_u$, then~\eqref{eq:induced-motivation} asserts that $x \in \{u\}^\perp$. The reason for this is that lattice points in $y \in \H_v$ with $\<y,u\> > 0$ are all too far away from $\0$ to be in $X$ (which is a finite set). Indeed, since $\eta$ is sufficiently small (as a function of $\U$ and $u$), and $v \in S_\eta(\SS^{d-1},u)$, it follows that we have only `tilted' $\H_u$ very slightly to obtain $\H_v$. The other direction of this implication is proved similarly. 

The equivalence we have just described is of central importance to the proof of Theorem~\ref{thm:upper}, and was in fact our original motivation for the definition of $r(\U)$. It will be proved formally in Lemma~\ref{lem:induced-projection} of Section~\ref{sec:equivalence}, in the more general setting that we need.

\begin{figure}[ht]
  \centering
  \bigskip
  \begin{tikzpicture}[>=latex, scale=0.6]
	\foreach \x in {-9,...,9}
	  \foreach \y in {-2,...,3}
	    \node [open] at (\x,\y) {};
	\foreach \x in {-9,...,9}
	  \foreach \y in {-4,-3}
	    \node [closed] at (\x,\y) {};
	\foreach \x in {-9,...,6}
	  \node [closed] at (\x,-2) {};
	\foreach \x in {-9,...,3}
	  \node [closed] at (\x,-1) {};
 	\foreach \x in {-9,...,-1}
	  \node [closed] at (\x,0) {};
	\node [draw, fill=white, circle, minimum size=0.2cm, inner sep=0pt] at (0,0) {};
	\draw [densely dashed] (-9,0) -- (0,0) -- (9,-3);
	\draw [dotted] (-5.5,0.5) -- (-0.5,0.5) -- (-0.5,-0.5) -- (4.5,-0.5) -- (4.5,-2.5) -- (-5.5,-2.5) -- cycle;
	\draw [->] (-10,0) -- (-10,1) node [above] {$u$};
	\draw [->] (10,-3.33) -- (10.33,-2.33);
	\path (10.33,-2.33) ++(71:0.5) node {$v$};
	\draw (1.5,-2.5) -- (1.5,-4.5) node [below] {$X$};
	\draw (-6.5,-3) -- (-6.5,-4.5) node [below] {$D$};
	\draw (4.2,-0.9) -- (5,-0.5) -- (10,-0.5) node [right] {$x$};
  \end{tikzpicture}
  \caption{The dashed lines mark the boundary of a droplet $D$, with adjacent faces in directions $u$ and $v$. These faces (i.e., the parts of $D$ intersecting the dashed line) have been infected, save for the point at their intersection (shown here as a circle at the centre of the figure; in higher dimensions it would be a face rather than just a point). 
The rule $X$, bounded by the dotted line, cannot be used to infect the final site, because of the uninfected site $x \in X$, which satisfies both $\< x,u \> < 0$ and $\< x,v \> > 0$.}
\label{fig:why-quasi}
\end{figure}
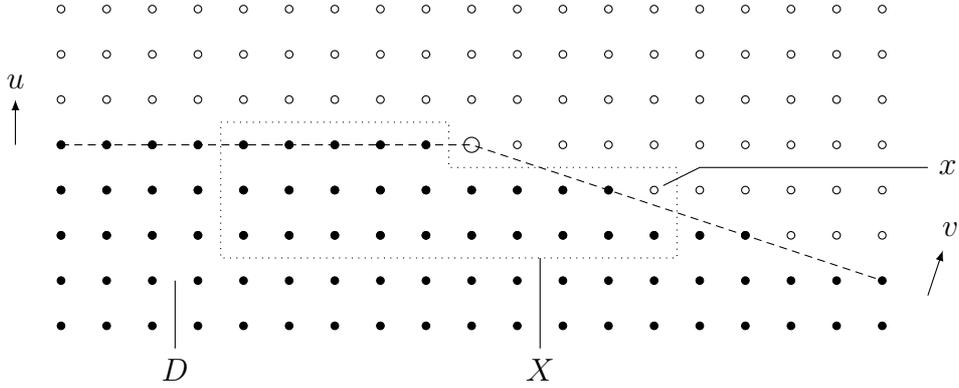

\subsection{Quasistable directions}\label{sec:chat-quasi}

The next choice we need to make is the set of directions to use for the faces of our droplets, which we call the \emph{quasistable set} and denote by $\QQ$. This set needs to have several key properties, and (perhaps surprisingly) showing that a suitable set $\QQ$ exists turns out to be one of the most technically challenging parts of the proof (see Section~\ref{sec:quasi}). In this section we describe (heuristically) the properties we need, and also briefly discuss how we shall go about constructing the set $\QQ$. 

First, the set $\QQ$ should be finite, and each of its elements should be rational. As mentioned above, this allows us to deduce that $\L(W)$ is a $\dim(W^\perp)$-dimensional lattice for any $W \subset \QQ$, and also that successive copies of $\L(W)$ are separated by a minimum distance. It will also allow us to use a union bound over the set $\QQ$. 

We also need $\QQ$ to have two properties related to induced processes. The first of these properties relates to growth near the boundary of a face, and will not be surprising to the reader who is familiar with~\cite{BSU,BDMS}. Intuitively, it says that we are able to infect a site $x$ on the $W$-face of a droplet $D$ using the induced process $\U[W]$, even if $x$ is near the edge of the face.\footnote{For formal statements of what we shall need, see Lemmas~\ref{cor:faces-far-away} and~\ref{lem:surface}.} Note that this is false in general, since (from the point of view of~$x$) the droplet $D$ may no longer be equivalent locally to $(x + \cH(W)) \cap \Z^d$ (see Figure~\ref{fig:why-quasi}). In order to avoid problems caused by `boundary effects'  when infecting the faces of our droplets, we need to choose the elements of $\QQ$ to be `sufficiently close together' (in a sense made precise below) so that whenever the $W$-face of a droplet is non-empty, the situation depicted in the figure does not occur.

The second property of $\QQ$ relating to induced processes is that it should be possible to construct our induced processes inductively. In particular, we shall require
\begin{equation}\label{eq:why-quasi-inductive}
\U\big[ \{u,v\} \big] = \big( \U[u] \big)[v] = \big( \U[v] \big)[u]
\end{equation}
to hold whenever $u$ and $v$ are the directions of adjacent faces of one of our droplets (see Lemma~\ref{lem:commutes} for the more general version of~\eqref{eq:why-quasi-inductive} that we shall actually need). This `abelian' property (which, once again, does not hold for arbitrary directions) is important because it allows us to construct induced processes inductively, via a sequence of induced processes of co-dimension~1. In particular, it will allow us to deduce a general version of the equivalence discussed in Section~\ref{sec:chat-induced} from Lemma~\ref{lem:induced-projection} (see Section~\ref{sec:tech:induced:lemmas}). 

In order to show that $\QQ$ satisfies the two properties described above, it turns out to be sufficient to show that $\QQ$ has the following property (see Figure~\ref{fig:why-quasi}):
\begin{itemize}
\item[(P1)] If $u,v \in \QQ$ are the directions of adjacent faces of a droplet $D$ used in our proof,\footnote{More precisely, we shall show that this is true whenever the Voronoi cells of $u$ and $v$ (with respect to $\QQ$) intersect; see Definition~\ref{def:voronoi} and Lemma~\ref{lem:quasi-innerprod}.}  
then we do not have
\[
\< x,u \> > 0 \qquad \text{and} \qquad \< x,v \> < 0
\]
for any $x \in X \in \U$.
\end{itemize}
For applications of property (P1), see the proofs of Lemmas~\ref{lem:induced-surprise} and~\ref{cor:faces-far-away}.

In two dimensions, Bollob\'as, Smith and Uzzell~\cite{BSU} were able to construct a set $\QQ \subset \SS^1$ satisfying (P1) simply by choosing $\QQ$ to be the set of all unit vectors perpendicular to an element of some update rule (see, e.g.,~\cite[Lemma~3.5]{BDMS}). 
Observe that the naive generalization of this construction to higher dimensions fails immediately, because the set so defined is infinite; moreover, it is not hard to see that an arbitrary sufficiently dense finite subset of this set also fails to satisfy (P1) in general. Indeed, if $u$ and $v$ are adjacent in $\QQ$ and lie on opposite sides of $\{x\}^\perp$ for some $x \in X \in \U$ (which can happen, for example, close to the intersections $\{x,y\}^\perp$), then $\QQ$ fails to satisfy (P1). 

We shall show in Section~\ref{sec:quasi} that a certain (carefully chosen) finite subset of this set has the desired property (P1), and moreover has the following stronger property:
\begin{itemize}
\item[(P2)] For each of the directions $w$ in which we wish to grow our droplets, the direction of every face that intersects the hyperplane 
$\{w\}^\perp$ is perpendicular to~$w$.
\end{itemize}
Roughly speaking, the directions in which we wish to grow will be the centres of the `easy' hemispheres given by Definition~\ref{def:r}. However, the actual choice of these directions is rather more delicate, and is the topic of Section~\ref{sec:chat-w} and Section~\ref{sec:w}.

Recall that the definition of $r(\U)$ only provides us with a single `easy' open hemisphere in which to grow, and observe that growing on faces in this hemisphere might produce a droplet that is much longer in some directions than others. Property (P2) guarantees that in fact the droplet is long only in the direction of $w$, and is roughly `spherical' in all directions in $\{w\}^\perp$; thus, overall, the droplet is roughly `tube-like'. This allows us to avoid having to deal with many different types of polytopes, and thus significantly simplifies the argument in Sections~\ref{sec:polytopes}--\ref{proof:sec}. In this choice we are also motivated by the potential applications to kinetically constrained models, where one often needs to find a path of updates that moves a droplet around without changing its size; see~\cite{Har,HMT1,MMT}.

In Section~\ref{sec:quasi} we construct a set $\QQ$ that satisfies property (P2) via a rather technical induction argument, which (perhaps surprisingly) is significantly harder when $d \ge 4$. In particular, in Lemma~\ref{lem:quasi-new} we show that there exists a set $\QQ$ satisfying  property (P2), and in Lemma~\ref{lem:quasi-innerprod} we show that property (P2) implies property (P1).

\subsection{Choosing a rational direction in which to grow}\label{sec:chat-w}

The alert reader may have noticed the following problem with the approach outlined above. 
Property (P2) of the quasistable set depends on the set of `directions in which we wish to grow our droplets', but this set of directions seems to depend on $\QQ$. In other words, for each $u \in \QQ$ there will be a corresponding face\footnote{Our droplet will also have faces corresponding to many other subsets $W \subset \QQ$, but for simplicity let us ignore this complication for now.} of our droplet, and for that face we only know that there is a single `easy' direction, given by Definition~\ref{def:r}. Since we appear to have no control over these directions, our reasoning appears to circular.

Another inconvenience is that the centre of the `easy' hemisphere of directions, whose existence is guaranteed by Definition~\ref{def:r}, may not be rational (that is, it may not be a member of $\SS_\Q^{d-1}$ as defined in~\eqref{eq:rational-sphere}). This is a problem, because many of our faces will be perpendicular to $w$, and (as discussed in Section~\ref{sec:chat-quasi}) faces must be in rational directions.

The solution to both problems, given in Section~\ref{sec:w} (specifically, in Lemma~\ref{lem:exists-rational-w}), is to show that, for each induced process, we may choose an easy direction $w \in \cL_R \cap \SS^{d-1}$, where 
\begin{equation}\label{def:LR}
\cL_R := \big\{ w \in \R^d : w \in \<x\> \text{ for some } x \in \Z^d \text{ with } \|x\| \le R \big\},
\end{equation}
for some constant $R = R(\U)$. That is, we can choose $w$ so that the line $\<w\>$ intersects the lattice $\Z^d$ in a (non-zero) point that is within a bounded distance of the origin. The key point is that, while the induced processes themselves depend on our choice of $\QQ$, the number of possible easy directions $w$ is bounded by a constant depending only on $R$. We may then choose $\QQ$ so that every such $w$ has the desired property.

The proof of Lemma~\ref{lem:exists-rational-w} is conceptually fairly simple: if the direction $w$ given to us by Definition~\ref{def:r} is not in $\cL_R$, then we `rotate' it slightly until we find a suitable $w'$. The reason this is possible is that the induced resistances $\rho^{d-1}(\SS^{d-1}; \S(\U), u)$ depend only on a bounded number of hyperplanes, and this implies that if $w$ cannot be rotated without changing the resistance of the hemisphere, then it must be `rational' in the required sense. The details, however, are once again somewhat involved.

\subsection{Locally inherited resistance}\label{sec:chat-memphis}

As we mentioned at the start of this section, an important (and somewhat subtle) property of $\U$-bootstrap percolation is as follows: it is never too much `harder' for a droplet to grow (on a face of our main droplet $D$) in any direction than it is in the `easiest' hemisphere. Now that we have formally introduced the notion of an induced process, we can describe this property more precisely.

Recalling Definition~\ref{def:r}, let $H$ be an open hemisphere such that 
\begin{equation}\label{eq:easyH:memphis:chat}
\rho^{d-1}\big(\SS^{d-1},\S(\U),v\big) \leq r - 1
\end{equation}
for every $v \in H$. Let $u \in H$, and recall from Section~\ref{sec:chat-induced} that the set $\S(\U[u])$, considered as a subset of $\SS(u) = \SS^{d-1} \cap \{u\}^\perp$, is equivalent to $\S(\U) \cap S_\eta(\SS^{d-1},u)$. By Definition~\ref{def:r}, it follows that 
\[
r^{d-1}\big( \SS(u); \S(\U[u]) \big) = r^{d-1}\big( S_\eta(\SS^{d-1},u); \S(\U) \big) = \rho^{d-1}\big( \SS^{d-1}; \S(\U), u \big) \leq r - 1,
\]
and hence that there exists an open hemisphere $H'$ of $\SS(u)$ such that 
\begin{equation}\label{eq:uface:memphis:chat}
\rho^{d-2}\big(\SS(u), \S(\U[u]), v \big) \leq r - 2
\end{equation}
for every $v \in H'$. This bound will (eventually) allow us to grow a droplet on the $u$-face of $D$ in the direction of the centre of $H'$ using the induced process $\U[u]$. 

In order to grow in the opposite direction (which will be necessary in order to infect the entire $u$-face of $D$), and also to grow `sideways' (i.e., perpendicular to the centre of $H'$), we need to bound the left-hand side of~\eqref{eq:uface:memphis:chat} for directions $v \not\in H'$. 

To give an idea of how we do so, suppose that the set $\S(\U) \cap S_\eta(\SS^{d-1},u)$ contains an open set $E$ in the usual topology on the sphere $S_\eta(\SS^{d-1},u)$. By Lemma~\ref{lem:inducedwelldef}, varying $\eta$ stretches $E$ into an open set in $\S(\U)$ in the usual topology on $\SS^{d-1}$. Since $u \in H$ and $\eta$ is sufficiently small, this open set is contained in $H$. But this contradicts~\eqref{eq:easyH:memphis:chat}, since $r \le d$, by assumption, 
and any point inside the open set has resistance $d$. It follows that no such open set $E$ exists, and one can then deduce, using~\eqref{eq:rho} and Lemma~\ref{lem:tri}, that
$$\rho^{d-2}\big( S_\eta(\SS^{d-1},u); \S(\U), v \big) \leq d - 2$$
for every $v \in S_\eta(\SS^{d-1},u)$. Again using the equivalence from Section~\ref{sec:chat-induced}, it follows that 
$$\rho^{d-2}\big( \SS(u); \S(\U[u]), v \big) \leq d - 2$$
for every $v \in \SS(u)$. This is our first example of `locally inherited resistance': we have deduced a non-trivial bound on the resistance of \emph{all} directions in $S_\eta(\SS^{d-1},u)$ (not just those in the `easy' hemisphere) from our assumption that $\U$ is critical. In Lemma~\ref{lem:memphis-subcrit} we prove a (straightforward) generalization of this bound, to the setting of an arbitrary sphere $\SS \subset \SS^{d-1}$ and an arbitrary $\SS$-stable set $\T$.

The proof of our second key lemma about locally inherited resistance, Lemma~\ref{lem:memphis}, is significantly more difficult. Roughly speaking, Lemma~\ref{lem:memphis} says that the resistance of a direction outside the easiest hemisphere $H$ is at most one greater than the maximum resistance in $H$. Heuristically, this implies that if a droplet of size $m$ is likely to grow in the direction of the centre of $H$, then a droplet of size roughly $e^m$ is likely to grow in the opposite direction. Since our droplets will typically have logarithmic size (as a function of the size of the corresponding face) when growing in the easiest direction, this will be just enough for the argument to work (see Section~\ref{proof:sec}, and in particular Lemma~\ref{lem:if-tube}).  

\begin{figure}[ht]
  \centering
  \begin{tikzpicture}[>=latex]
    \draw (0,0) circle (3);
    \draw (0:3) arc (0:-180:3 and 0.6);
    \draw [dashed] (0:3) arc (0:180:3 and 0.6);
    \draw (-30:3) arc (0:-180:{3*cos(30)} and {3*cos(30)/5});;
    \draw [dashed] (-30:3) arc (0:180:{3*cos(30)} and {3*cos(30)/5});
    \draw (80:3) arc (0:-180:{3*cos(80)} and {3*cos(80)/5});
    \coordinate (X) at (1.5,{-0.6*sqrt(1-1.5^2/3^2)});
    \draw (X) circle (0.3);
    \draw [dashed] (0,3) -- (0,-3);
    \fill (0,0) circle (0.07);
    \fill (0,-1.5) circle (0.07);
    \fill (0,3) circle (0.07);
    \fill (0,-3) circle (0.07);
    \fill (X) circle (0.07);
    \node at (0,3.5) {$w$};
    \node at (0,-3.5) {$-w$};
    \node at (0.4,0) {$\0$};
    \node at (0.4,-1.5) {$\lambda w$};
    \node at (-3.4,0) {$\SS_0$};
    \node at (-3,-1.5) {$\SS_\lambda$};
    \node at (-2.4,2.4) {$\SS^{k}$};
    \draw (0.3,2.86) to [out=-30,in=-160] (2,2.7) node [right] {$S_\eta(\SS^k,w)$};
    \draw (X) ++(45:0.1) to [out=45,in=180] (3.2,0.5) node [right] {$u$};
    \draw (X) ++(-45:0.3) to [out=-35,in=-180] (3.2,-1) node [right] {$\SS_u$};
  \end{tikzpicture}
  \caption{The setting of Lemma~\ref{lem:axis}: if $C(\SS^{k};\T) \subset \{w\}^\perp$ for some $\SS^{k}$-stable set $\T$ and some $w \in \SS^{k}$, then $r^{k+1}(\SS^{k};\T) \geq r^{k}(\SS_\lambda;\T)$.}
  \label{fig:axis}
\end{figure}
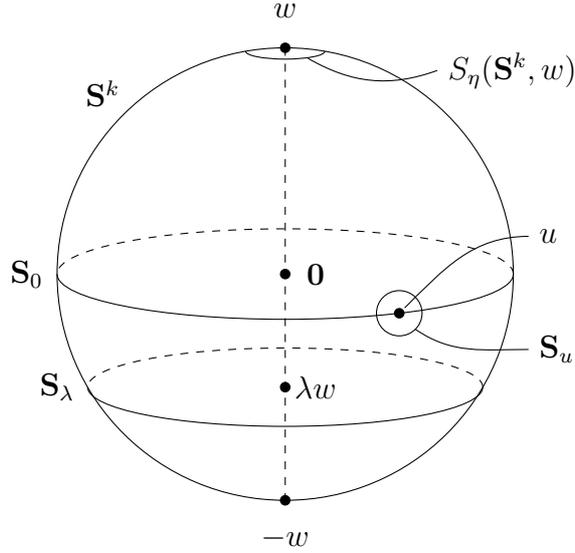

Although the proof of Lemma~\ref{lem:memphis} is technically much more complex than that of Lemma~\ref{lem:memphis-subcrit}, the underlying idea is similar. That is, if $v \in S_\eta(\SS,u)$ and we know the intersection of $\T$ with a neighbourhood of $v$ in $S_\eta(\SS,u)$, then we can deduce the corresponding intersection with a neighbourhood of $v$ in $\SS$ by `stretching' the former in radial directions with respect to $u$ (see Lemma~\ref{lem:memphis-vertical} and Figure~\ref{fig:axis}). To do so, we use the following simple observation: the boundary of every hemisphere as in~\eqref{eq:happy} either passes through both $u$ and $v$, or intersects neither $S_\eta(\SS,u)$ nor $S_\eta(\SS,v)$. In particular, this means that the boundary of any hemisphere of $\T$ that passes through $v$ contains a copy of the vector $v-u$. This restriction on $\T$ forces the intersection of $\T$ with subspheres of $S_\eta(\SS,v)$ obtained by taking slices in directions perpendicular to $v-u$ to be equivalent, and this is exactly what we require for our stretching procedure to work. 

In Section~\ref{sec:induced} we combine Lemmas~\ref{lem:memphis-subcrit} and~\ref{lem:memphis} with the results of Sections~\ref{sec:w} and~\ref{sec:quasi} to deduce the required properties of induced update families (see Lemmas~\ref{lem:good:faces} and~\ref{lem:semigood:faces}). The specific applications of Lemmas~\ref{lem:memphis-subcrit} and~\ref{lem:memphis} can be found in the proofs of Lemmas~\ref{lem:sgood} and~\ref{lem:sgood:extra}. We remark that Lemma~\ref{lem:axis}, which is the main intermediate step in the proof of Lemma~\ref{lem:memphis}, is also used in Section~\ref{sec:w}, and that Lemma~\ref{lem:memphis} also plays an important role in the proof of the lower bounds in Theorem~\ref{thm:universality}, see~\cite{BBMSlower}.

\subsection{Constructing a path of infections}\label{sec:chat-path}

In the second half of the paper, comprising Sections~\ref{sec:polytopes}--\ref{proof:sec} and Appendices~\ref{polytope:app}--\ref{cover:app}, we apply the results described above, which will have been proved in Sections~\ref{sec:memphis}--\ref{sec:induced}, in order to construct a family of `likely' ways in which percolation could occur in the torus $\Z_n^d$. In Section~\ref{proof:sec} we show that, with high probability, a $p$-random initial set $A$ can realise some member of this family of infection paths. We do this by deriving lower bounds on the probability of percolation on the faces of a droplet under the action of induced update families, using induction on the dimension of the face (see Definition~\ref{def:ih} and Lemmas~\ref{lem:ihb} and~\ref{lem:iha}).

Roughly speaking, our construction will be as follows. Let $w \in \cL_R \cap \SS^{d-1}$ be the `easy' direction given by Lemma~\ref{lem:exists-rational-w} (recall that we discussed this lemma in Section~\ref{sec:chat-w}). In order to infect a site $x \in \Z_n^d$, a droplet $D$ will be constructed intersecting the line $\{ x + \lambda w : \lambda \in \R \}$, and then grown in direction $w$ until it infects $x$. Our main task will therefore be to show that such a droplet is likely to exist in $[A]_\U$, and that it is likely to grow sufficiently far in direction $w$. 

In order to show that $D$ is likely to grow in direction $w$, we need to bound from below the probability of percolation on the face of $D$ in direction $u$, for each $u \in \QQ$ such that $\<u,w\> > 0$. To do this, we shall use one of our deterministic results, Lemma~\ref{lem:forwards-step-general}, which says that it suffices to infect the sites in the interior of each lower-dimensional face in the corresponding induced process. This result relies heavily on Lemma~\ref{cor:faces-far-away}, which relates the induced process on a face of $D$ to the original $\U$-process. As discussed above, Lemma~\ref{cor:faces-far-away} is one of the key motivations for and applications of our definition of $\QQ$. 

We can now use the induction hypothesis to bound the probability that the interior of a face is `internally filled' (see Definition~\ref{def:int:filled}) in the corresponding induced process. However, in order to obtain a sufficiently strong bound, we need to adapt an idea used by Schonmann~\cite{Sch1} in his work on the $r$-neighbour model. The idea is to cover a droplet (or a face) with smaller copies of itself (each contained in the larger droplet), and prove that each of these `tiles' is internally filled with probability close to~$1$. By a standard argument, it follows that large components of non-internally filled tiles occur with exponentially small probability, and small components are (deterministically) infected by the internally filled tiles that surround them. We adapt the various parts of this argument to our setting in Lemma~\ref{lem:eating:a:set}, Section~\ref{spherical:IH:sec} and Appendix~\ref{cover:app}. We remark that the proof that each tile is internally filled with high probability uses the results of Section~\ref{sec:induced}, which in turn relies on the results of Sections~\ref{sec:memphis}--\ref{sec:quasi}.

Having shown that a sufficiently large droplet can grow easily in direction $w$, it remains to show that such a droplet is likely to exist in $[A]_\U$. This will be more challenging, since it will require us to control not just `forwards' growth (in the direction of $w$), but also `sideways' growth (perpendicular to $w$). This is necessary because we need a droplet of `size' roughly $\log n$ (in $\{w\}^\perp$) in order to grow $n$ steps in direction $w$, and the largest droplet we can find in $A$ is significantly smaller than this. 

Our main deterministic lemma for dealing with sideways growth is Lemma~\ref{lem:sideways-step-general}, which provides sufficient conditions, stated in terms of the induced processes on the sideways faces of $D$, for the infection of the `forwards end' of a slightly larger droplet. Note that one cannot hope to infect the entire side of $D$, since we only know that there is a single easy direction (fortunately, it is still $w$) on each sideways face of $D$. However, starting from this new seed, we can continue to grow in the easy direction, and eventually be ready for another sideways step. 

The probability that the conditions for sideways growth are satisfied can be controlled using our induction hypothesis (see Lemmas~\ref{lem:forwards-step-probability} and~\ref{lem:sideways-step-probability}), and this lemma therefore allows us (see Lemma~\ref{lem:IHb:likely}) to show that a droplet is likely to form somewhere on (any sufficiently large portion of) the line $\{ x + \lambda w : \lambda \in \R \}$. Applying this argument to an arbitrary face of $D$, we are able to deduce the induction step. 

The proof of Lemma~\ref{lem:sideways-step-general} is the main technical challenge of the second half of the paper, and is built upon a compilation of properties of polytopes, and various operations on polytopes, which are introduced in Sections~\ref{sec:polytopes} and~\ref{buffers:sec}, and proved in Appendices~\ref{polytope:app} and~\ref{buffers:app}. We should mention here that we shall work with a relatively simple family of polytopes, defined by translating, dilating, and stretching a single simple `canonical' family of polytopes; see Section~\ref{our:family:sec}. This will have the advantage of greatly simplifying our induction hypothesis.

\subsection{The equivalence lemma}\label{sec:equivalence}

Having completed our rough outline of the argument, we are ready to begin the formal proof of Theorem~\ref{thm:upper}. The first step, and the subject of this subsection, is to formalize (and generalize) the equivalence outlined in Section~\ref{sec:chat-induced}. In order to do so, we first need to introduce some additional notation. In this subsection (and also in Section~\ref{sec:w}) we work with an arbitrary update family $\F$, which in our applications will be set equal to $\U[W]$ for some $W \subset \QQ$. We work in this more general setting both to simplify the notation, and to emphasize the difference between those results that depend on Definition~\ref{def:induced}, and those that do not. 


Given an arbitrary update family $\F$, we define the \emph{span} of $\F$ to be 
\begin{equation}\label{def:span}
\< \F \> := \Big\< \bigcup_{X \in \F} X \Big\>.
\end{equation}
Observe that if $\F = \U[W]$ then $\<\F\> \subset W^\perp$. Moreover, any update family $\F$ for which $\<\F\> \subset W^\perp$ acts independently on each translate of the lattice $\L(W)$.  

Let us mention that it may be that $\F = \emptyset$ (which corresponds to growth not being possible under $\F$) and it may be that $\emptyset \in \F$ (which corresponds to growth being trivial under $\F$: every uninfected site becomes infected immediately). We say that $\F$ is \emph{trivial} if either $\F = \emptyset$ or $\emptyset \in \F$, and that $\F$ is \emph{non-trivial} otherwise. 

When working with update families satisfying $\<\F\> \subset W^\perp$, it will be convenient to work in the sub-sphere
$$\SS(W) := \SS^{d-1} \cap W^\perp,$$
of $\SS^{d-1}$. We shall also make frequent use of the \emph{projection} $\pi(u,W^\perp)$ of a point $u \in \SS^{d-1} \setminus \< W \>$ onto $\SS(W)$; that is, the unique element of $\SS(W)$ such that 
\begin{equation}\label{def:projection}
u = w + \lambda \cdot \pi(u,W^\perp)
\end{equation}
for some $w \in \< W \>$ and $\lambda > 0$. 
The projection has the following useful property. 

\begin{obs}\label{obs:projection:innerproduct}
Let\/ $W \subset \SS^{d-1}$ and\/ $u \in \SS^{d-1} \setminus \< W \>$. If\/ $x \in W^\perp$, then 
$$\sgn\big( \< x,u \> \big) = \sgn\big( \< x,\pi(u,W^\perp) \> \big).$$
\end{obs}

\begin{proof}
Since $x \in W^\perp$ and $\pi(u,W^\perp) = w + \lambda u$ for some $w \in \< W \>$ and $\lambda > 0$, it follows that $\< x, \pi(u,W^\perp) \> = \lambda \cdot \< x,u \>$, which implies the claimed identity. 
\end{proof}

We are now ready to state and prove our key equivalence lemma. Roughly speaking, the lemma says that the stable set of an update family $\F$ with $\<\F\> \subset W^\perp$ near to a rational direction $u \in \SS(W)$ is equivalent to that of the induced family $\F[u]$.\footnote{The reader who is familiar with~\cite{BSU} may find it useful to think of Lemma~\ref{lem:induced-projection} as a high dimensional analogue of~\cite[Lemma 5.2]{BSU}. In particular, the existence of `$u$-left/right-blocks' in~\cite{BSU} determined the stability or otherwise of each of the two directions in the one-dimensional induced process $\F[u]$.} 
We shall apply the lemma in Section~\ref{sec:induced}, with $\F = \U[W]$ (see Lemma~\ref{lem:induced}).

\begin{lemma}\label{lem:induced-projection}
Let\/ $W \subset \SS_\Q^{d-1}$ and\/ $u \in \SS_\Q^{d-1} \cap W^\perp$, and set\/ $W' := W \cup \{u\}$. If\/~$\F$ is an 
update family such that $\< \F \> \subset W^\perp$, then 
\begin{equation}\label{eq:induced-proj-iff}
v \in \S(\F) \qquad \Leftrightarrow \qquad \pi(v,W'^\perp) \in \S\big( \F[u] \big)
\end{equation}
for every $v \in S_\eta\big( \SS(W), u \big)$. In particular,
\[
\S(\F) \cap S_\eta\big( \SS(W), u \big) \equiv \S\big( \F[u] \big) \cap \SS(W').
\]
\end{lemma}

\begin{proof}
Observe first that 
\begin{equation}\label{eq:SminusSWprimeFu1}
\SS(W') \setminus \S\big( \F[u] \big) = \bigcup_{X \in \F[u]} \bigcap_{x \in X} \big\{ v' \in \SS(W') : \< x,v' \> < 0 \big\},
\end{equation}
since a direction $v' \in \SS(W')$ is unstable for the induced family $\F[u]$ 
if and only if there exists $X \in \F[u]$ such that $\< x,v' \> < 0$ for every $x \in X$. Now, recall that 
\begin{equation}\label{eq:Fu:def:reminder}
\F[u] = \big\{ X \cap \{u\}^\perp : X \in \F \,\text{ and }\, X \subset \HH(u) \big\},
\end{equation}
and that $u \in W^\perp$ and $X \subset W^\perp \cap \Z^d$ for every $X \in \F$. We claim that we can therefore rewrite the right-hand side of~\eqref{eq:SminusSWprimeFu1} as 
\begin{equation}\label{eq:SminusSWprimeFu}
\bigcup_{X \in \F} \bigcap_{x \in X} \Big\{ v' \in \SS(W') : \< x,u \> < 0 \,\text{ or }\, \big( \< x,u \> = 0 \text{ and } \< x,v' \> < 0 \big) \Big\}.
\end{equation}
To see this, note first that those $x \in X \in \F$ with $\< x,u \> < 0$ can be removed from~\eqref{eq:SminusSWprimeFu}, since the corresponding set is $\SS(W')$, and so does not affect the intersection.\footnote{If the reader is concerned that this may remove all elements of some $X \in \F$, note that if $\< x,u \> < 0$ for each $x \in X$, then by~\eqref{eq:Fu:def:reminder} we have $\emptyset \in \F[u]$, and therefore both~\eqref{eq:SminusSWprimeFu1} and~\eqref{eq:SminusSWprimeFu} are equal to $\SS(W')$.} We can also remove those sets $X \in \F$ that are not contained in $\HH(u)$, since the corresponding intersection is empty. We are left with exactly those sets $X \in \F$ such that $X \subset \HH(u)$, and those elements $x \in X \cap \{u\}^\perp$, as claimed.

Next we rewrite~\eqref{eq:SminusSWprimeFu} as
\begin{equation}\label{eq:induced-proj-1}
\bigcup_{X \in \F} \bigcap_{x \in X} \bigcup_{v \in \SS_u} \Big\{ \pi(v,W'^\perp) : \< x,u \> < 0 \,\text{ or }\, \big( \< x,u \> = 0 \text{ and } \< x,v \> < 0 \big) \Big\}.
\end{equation}
where $\SS_u = S_\eta\big( \SS(W), u \big)$. This holds because the homothety $\varphi \colon \SS_u \to \SS(W')$ defined by $\varphi(v) = \pi(v,W'^\perp)$ is a bijection,
and if $\< x,u \> = 0$ for some $x \in X \in \F$, then 
$$\< x,v \> < 0  \qquad \Leftrightarrow \qquad \big\< x,\pi(v,W'^\perp) \big\> < 0$$
for each $v \in \SS_u$, by Observation~\ref{obs:projection:innerproduct}, since $x \in \<\F\> \subset W^\perp$, so if $\< x,u \> = 0$ then $x \in W'^\perp$. 

We next claim that for each $v \in \SS_u$ and $x \in X \in \F$, we have
\begin{equation}\label{eq:induced-proj-tilt}
\< x,u \> < 0 \,\text{ or }\, \big( \< x,u \> = 0 \,\text{ and }\, \< x,v \> < 0 \big) \quad \Leftrightarrow \quad \< x,v \> < 0.
\end{equation}
To see this, note first that if $\< x,u \> < 0$ for some $x \in X \in \F$, then $\< x,v \> < 0$ for every $v \in \SS_u$, by the continuity of $\< x,\cdot \>$, and since $\eta$ is sufficiently small and $\F$ consists of a finite number of finite sets. Similarly, if $\< x,u \> > 0$ then $\< x,v \> > 0$ for every $v \in \SS_u$. 

By~\eqref{eq:induced-proj-tilt}, and recalling that~\eqref{eq:induced-proj-1} is equal to the right-hand side of~\eqref{eq:SminusSWprimeFu1}, we obtain 
\begin{equation}\label{eq:SWminusSWFu:final}
\SS(W') \setminus \S\big( \F[u] \big) = \bigcup_{X \in \F} \bigcap_{x \in X} \Big\{ \pi(v,W'^\perp) : v \in \SS_u \text{ and } \< x,v \> < 0 \Big\}.
\end{equation}
Now,
since a direction $v \in \SS_u$ is unstable for $\F$ if and only if there exists $X \in \F$ such that $\< x,v \> < 0$ for every $x \in X$, we have
\begin{equation}\label{eq:SWminusSWFu:final2}
\SS_u \setminus \S(\F) = \bigcup_{X \in \F} \bigcap_{x \in X} \big\{ v \in \SS_u : \< x,v \> < 0 \big\}.
\end{equation}
Combining~\eqref{eq:SWminusSWFu:final} and~\eqref{eq:SWminusSWFu:final2}, we obtain the equivalence~\eqref{eq:induced-proj-iff}, since $\SS_u = S_\eta\big( \SS(W), u \big)$ and the homothety $\varphi(v) = \pi(v,W'^\perp)$ maps $\SS_u \setminus \S(\F)$ to $\SS(W') \setminus \S\big( \F[u] \big)$.  

Finally, to deduce that the sets $\S(\F) \cap \SS_u$ and $\S\big( \F[u] \big) \cap \SS(W')$  are equivalent, simply note that $\varphi(v) = \pi(v,W'^\perp)$ is a homothety from $\SS_u$ to $\SS(W')$, and apply~\eqref{eq:induced-proj-iff}.
\end{proof}

\subsection{Some constants}

To finish this section, let us introduce some constants that will appear frequently during the proof. First, we define the \emph{radius} of $\U$ to be
\begin{equation}\label{def:radius}
R_0 = R_0(\U) := \max\big\{ \|x\| : x \in X \in \U \big\}.
\end{equation}
We also introduce the following hierarchy of constants:
\begin{equation}\label{def:constants}
0 < \delta \ll \gamma \ll 1 \le R_0 \ll R \ll C.
\end{equation}
where $\ll$ is used informally to indicate the relative sizes of the constants. More precisely, we shall choose $R = R(\U)$ in Lemma~\ref{lem:exists-rational-w}, and then, in Section~\ref{sec:quasi}, define our set $\QQ$ of quasistable directions 
depending on $R$. The constants $\gamma = \gamma(\QQ)$ and $\delta = \delta(\QQ)$ will be defined in Section~\ref{fundamental:sec} (see Definition~\ref{def:gamma} and Lemma~\ref{lem:close-to-many-faces}); in the proof of Lemma~\ref{lem:close-to-many-faces} we shall need $\delta \ll \gamma$. The constant $C$ is chosen last, and will be used to define the family of polytopes (see Definition~\ref{def:cPW}) that we use in Sections~\ref{sec:polytopes}--\ref{proof:sec} to construct our `path of infections'. As a consequence, it will also appear in the functions $t_0(k,s,p)$ and $t_1(k,s,p)$ that we use to define our induction hypothesis (see Definition~\ref{def:ih}). 
The constant $\eps = \eps(\U)$ that we obtain in Theorem~\ref{thm:upper} is of the form $\eps = C^{-O(d)}$.

\section{Locally inherited resistance}\label{sec:memphis}

In this section we prove two lemmas of the following flavour. Let $\SS \subset \SS^{d-1}$ be a sphere  of dimension $k$ (with $1 \leq k \leq d-1$), let $\T$ be an $\SS$-stable set, and let $u \in \SS$. Then, provided an element $v$ of $S_\eta(\SS,u)$ satisfies a certain upper bound condition on its resistance with respect to $\T$ in $\SS$, then it also satisfies a certain upper bound condition on its resistance with respect to $\T$ in $S_\eta(\SS,u)$.

Our main challenge will be to prove Lemma~\ref{lem:memphis}, which states that the resistance of $v$ in $S_\eta(\SS,u)$ is at most its resistance in $\SS$. First, however, we show that a stronger conclusion holds when $\rho^k \big( \SS; \T, v \big) = k$. More precisely, Lemma~\ref{lem:memphis-subcrit} says that if $v$ is not subcritical with respect to $\T$ in $\SS$ (i.e., it has resistance at most $k$), then $v$ is also not subcritical with respect to $\T$ in $S_\eta(\SS,u)$ (that is, its resistance in $S_\eta(\SS,u)$ is at most $k - 1$). This allows us to deduce that none of the induced update families we use to control the growth of a droplet is subcritical.

\begin{lemma}\label{lem:memphis-subcrit}
Let $\SS \subset \SS^{d-1}$ be a $k$-dimensional sphere, where $1 \le k \le d - 1$, and let $\T$ be an $\SS$-stable set. Let $u \in \SS$ and $v \in S_\eta(\SS,u)$, and suppose that
\[
\rho^k \big( \SS; \T, v \big) \le k.
\]
Then 
\[
\rho^{k-1} \big( S_\eta(\SS,u); \T, v \big) \le k - 1.
\]
\end{lemma}

\begin{proof}
Suppose first that $k = 1$, in which case we are required to prove that $v \notin \T$. This~follows from the definition of $\eta$ (see Lemma~\ref{lem:inducedwelldef}), and since $\SS$ is a $1$-dimensional sphere. Indeed, if $v \in \T$ then, since $v \in S_\eta(\SS,u)$, it follows that $w \in \T$ for every $w \in \SS$ that is sufficiently close to $u$ and on the same side of $u$ as $v$. In particular, we have $w \in \T$ for every $w \in \SS$ that is sufficiently close to $v$. But, by~\eqref{eq:rho} and~\eqref{eq:r}, this implies that $\rho^{1}( \SS; \T, v ) = 2$, which is the desired contradiction. 

So let $2 \le k \le d - 1$, and suppose, for a contradiction, that
$$\rho^{k-1} \big( \SS_u; \T, v \big) = k,$$
where $\SS_u := S_\eta(\SS,u)$. By~\eqref{eq:rho}, this implies that $v \in \T$ and 
$$r^{k-1} \big( \SS_v; \T \big) = k,$$
where $\SS_v := S_\theta( \SS_u , v )$ and $\theta > 0$ is sufficiently small (in particular, $\theta \ll \eta$). It follows, by Lemma~\ref{lem:tri} and Definition~\ref{def:tri}, that $\interior_{\SS_v}(H \cap \T) \neq \emptyset$ for every hemisphere $H \subset \SS_v$, and in particular $H \cap \T$ contains a non-empty open set. (Here, and below, `open' refers to the topology of the relevant sphere, in this case the $(k-2)$-sphere $\SS_v$.) 

\begin{figure}[ht]
  \centering
  \begin{tikzpicture}[>=latex]
    \fill (0,0) circle (0.07);
    \fill (8,0) circle (0.07);
    \draw [name path=Seta] (-20:8) arc (-20:20:8);
    \draw [name path=S3] (8,0) circle (1.5);
    \path [name intersections={of=Seta and S3,by={X,Y}}];
    \path (X);
    \pgfgetlastxy{\Xx}{\Xy};
    \draw (X) arc (165:195:{\Xy/sin(15)});
    \path (X) arc (165:170:{\Xy/sin(15)}) coordinate (P1);
    \path (X) arc (165:177:{\Xy/sin(15)}) coordinate (P2);
    \path (X) arc (165:183:{\Xy/sin(15)}) coordinate (P3);
    \path (X) arc (165:190:{\Xy/sin(15)}) coordinate (P4);
    \begin{scope}[on background layer]
      \path [fill,color=gray!70] (P1) arc (170:177:{\Xy/sin(15)}) -- (8,0) -- (P1);
      \path [fill,color=gray!70] (P3) arc (183:190:{\Xy/sin(15)}) -- (8,0) -- (P3);
      \path [fill,color=gray!20] (P1) arc (170:177:{\Xy/sin(15)}) -- (0,0) -- (P1);
      \path [fill,color=gray!35] (P3) arc (183:190:{\Xy/sin(15)}) -- (0,0) -- (P3);
      \path [fill,color=gray!50] (P2) -- (0,0) -- (8,0) -- (P2);
      \path [fill,color=gray!30] (P3) -- (0,0) -- (8,0) -- (P3); 
    \end{scope}
    \draw [dashed] (0,0) -- (8,0);
    \foreach \x in {1,2,3,4}
      \draw [dashed] (0,0) -- (P\x) -- (8,0);
    \draw [line width=3] (P1) arc (170:177:{\Xy/sin(15)});
    \draw [line width=3] (P3) arc (183:190:{\Xy/sin(15)});
    \node at (-0.5,0) {$u$};
    \node at (8.3,0) {$v$};
    \node at (17:8.4) {$\SS_u$};
    \node at ($(8,0)+(-40:2.3)$) {$S_\theta(\SS,v)$};
    \draw [<-] (X) ++(-0.1,-0.3) -- ++(135:1) node [above left] {$\SS_v$};
  \end{tikzpicture}
  \caption{The proof of Lemma~\ref{lem:memphis-subcrit} when $k = 3$. The setting is the 3-sphere $\SS$; also shown are the 2-spheres $\SS_u = S_\eta(\SS,u)$ and $S_\theta(\SS,v)$, and the 1-sphere $\SS_v$. The bold (1-dimensional) arcs on $\SS_v$ are open intervals contained in $\T$; since $\eta$ and $\theta$ are sufficiently small, it follows that the shaded (3-dimensional) wedges are contained in $\T \cap \SS$.}
  \label{fig:memphis-subcrit}
\end{figure}
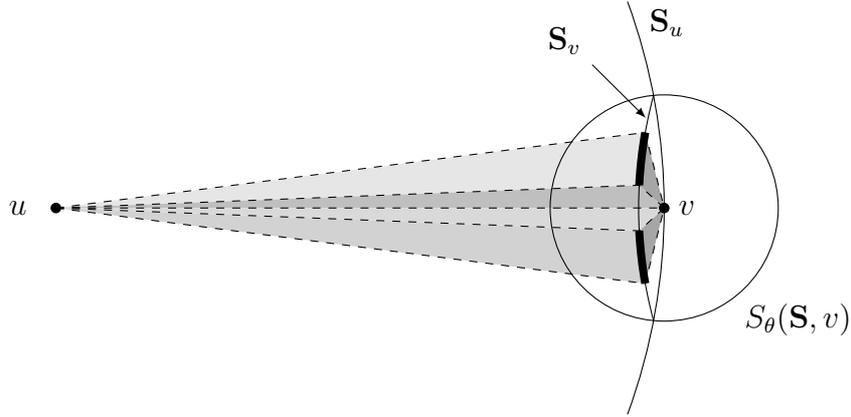

In order to obtain a contradiction, we need to show that $\rho^{k}( \SS; \T, v) = k + 1$, which, by~\eqref{eq:rho}, and since $v \in \T$, is equivalent to 
\begin{equation}\label{eq:memphis-sub-contra}
r^{k} \big( S_{\theta'}(\SS,v); \T \big) = k + 1
\end{equation}
for all sufficiently small $\theta' > 0$. To prove this, by Lemma~\ref{lem:tri} and Definition~\ref{def:tri}, it is enough to show that there exists a non-empty open set in $H \cap \T$ for every hemisphere $H \subset S_{\theta'}(\SS,v)$. The proof of this claim is depicted in Figure~\ref{fig:memphis-subcrit}. Formally, let $H$ be an arbitrary (closed) hemisphere in $S_\theta(\SS,v)$, set $H' := H \cap \SS_u = H \cap \SS_v$, and observe that $H'$ contains a hemisphere of $\SS_v$. It follows from the observations above that $H' \cap \T$ contains a non-empty open set $W$. 

Now, for each $0 < \theta' \le 2\theta$, consider the homothety $\varphi_{\theta'}$ that maps $\SS_v$ to $S_{\theta'}(\SS_u,v)$, and observe that $\varphi_{\theta'}(W)$ is an open subset of $\varphi_{\theta'}(H') \cap \T$ (cf.~the proof of Lemma~\ref{lem:inducedwelldef}). Similarly, for each $0 < \eta' \le 2\eta$, consider the homothety $\chi_{\eta'}$ that maps $\SS_u$ to $S_{\eta'}(\SS,u)$, and observe that $\chi_{\eta'}\big( \varphi_{\theta'}(W) \big)$ is an open subset of $\chi_{\eta'}\big( \varphi_{\theta'}(H') \big) \cap \T$. It follows that the set
$$\bigcup_{0 < \eta' < 2\eta} \bigcup_{0 < \theta' < 2\theta} \chi_{\eta'}\big( \varphi_{\theta'}(W) \big)$$
contains an open subset of $H \cap \T$. 

We have proved that for every hemisphere $H \subset S_\theta(\SS,v)$, there is a non-empty open set in $H \cap \T$. Finally, applying Lemma~\ref{lem:inducedwelldef}, and recalling that $\theta$ was chosen sufficiently small, it follows that the same holds for $S_{\theta'}(\SS,v)$ for every $0 < \theta' \le \theta$, as claimed. As observed above, this proves~\eqref{eq:memphis-sub-contra}, and gives the desired contradiction. 
\end{proof}

The second main lemma of the section is weaker than Lemma~\ref{lem:memphis-subcrit} when $\rho^k( \SS; \T, v ) = k$, but is more general, and is in fact best possible in the level of generality stated. 

\begin{lemma}\label{lem:memphis}
Let\/ $\SS \subset \SS^{d-1}$ be a $k$-dimensional sphere, where $1 \le k \le d - 1$, and let\/ $\T$ be an $\SS$-stable set. If\/ $u \in \SS$ and $v \in S_\eta(\SS,u)$, then
\begin{equation}\label{eq:memphis}
\rho^{k-1} \big( S_\eta(\SS,u); \T, v \big) \le \rho^k \big( \SS; \T, v \big).
\end{equation}
\end{lemma}

In the proof of Lemma~\ref{lem:memphis}, it will be convenient to work in the sphere $\SS^k$ rather than the sphere $\SS$. Since $\T$ is an arbitrary $\SS$-stable set, the following simple observation will allow us to do this.

\begin{obs}\label{lem:moving:to:Sk}
Let\/ $\SS' \subset \SS \subset \SS^{d-1}$ be spheres, with $\ell := \dim(\SS') \geq 0$, and let $\T$ be an\/ $\SS$-stable set. Let $\varphi$ be the homothety that maps the centre of\/ $\SS$ to the origin, and\/ $\SS$ to a subset of\/ $\SS^{d-1}$. Then
$$\rho^\ell \big( \varphi(\SS'); \varphi(\T), \varphi(v) \big) = \rho^\ell \big( \SS'; \T, v \big)$$
for every $v \in \SS'$. 
\end{obs}

For each $1 \le k \le d - 1$, we write $\SS^k$ to denote an arbitrary $k$-dimensional sphere in $\SS^{d-1}$, centred at the origin. Recall from~Definition~\ref{def:happy} that if $\T$ is an $\SS^k$-stable set, then
\begin{equation}\label{eq:happy:reminder}
\T \cap \SS^k = \bigcap_{i=1}^t \bigcup_{u \in Y_i} H_u
\end{equation}
for some finite collection $Y_1,\ldots,Y_t$ of finite families of vectors in $\SS^k$, where $H_u$ is the closed hemisphere of $\SS^k$ centred at $u$, that is, 
$$H_u = \big\{ v \in \SS^k : \< u,v \> \geq 0 \big\}.$$ 
For each $\SS^k$-stable set $\T$, let us choose a minimal representation as in~\eqref{eq:happy:reminder}, and define the \emph{set of centres of $\T$ with respect to $\SS^k$} to be 
\begin{equation}\label{def:CST}
C(\SS^k;\T) := Y_1 \cup \cdots \cup Y_t.
\end{equation}
Note that $C(\SS^k;\T) = \emptyset$ if either $\T \cap \SS^k = \emptyset$ (since we can take $Y_1 = \emptyset$) or $\T \cap \SS^k =\SS^k$ (since we can take $t = 0$). Observe also that if $\T = \S(\F)$ for an update family $\F$, then the set of centres can be taken to be the elements of the update sets. We record this simple fact as an observation, since we shall use it in Section~\ref{sec:w}.

\begin{obs}\label{obs:rules:and:centres}
Let $\F$ be an update family, and let\/ $W \subset \SS_\Q^{d-1}$. Then
$$C(\SS(W);\S(\F)) \subset \<\F\>$$
for some valid choice of $C(\SS(W);\S(\F))$.
\end{obs}

\begin{proof}
This follows from the definition, using the fact that 
$$\S(\F) \cap \SS(W) = \bigcap_{X \in \F} \bigcup_{x \in X} \big\{ u \in \SS(W) : \< x,u \> \geq 0 \big\},$$
by~Lemma~\ref{lem:hemispheres}. 
\end{proof}

Recall from~\eqref{def:projection} the definition of the projection $\pi(u,W^\perp)$ of a vector $u \in \SS^{d-1} \setminus \< W \>$ onto the sphere $\SS(W)$, and note that if $u \in \SS^k$ and $W \subset \SS^k$, then $\pi(u,W^\perp) \in \SS^k \cap W^\perp$. We can now state the following simple but important lemma, which is also used in~\cite{BBMSlower}. It says that resistance is preserved under projection onto subspaces containing $C(\SS^k;\T)$. 

\begin{lemma}\label{lem:memphis-vertical}
Let\/ $1 \le k \le d - 1$, let $\T$ be an\/ $\SS^k$-stable set, and let $W \subset \SS^k$. If\/ $C(\SS^k;\T) \subset W^\perp$, then
\begin{equation}\label{eq:memphis-vertical-rho}
\rho^{k}\big( \SS^k;\T,u \big) = \rho^{k}\big( \SS^{k};\T,v \big),
\end{equation}
for every $u,v \in \SS^{k} \setminus \< W \>$ such that $\pi(u,W^\perp) = \pi(v,W^\perp)$.\medskip
\end{lemma}

\begin{proof}
Let us assume that $v = \pi(u,W^\perp)$; it will suffice to show that~\eqref{eq:memphis-vertical-rho} holds in this case, since it follows that (in general) both sides of~\eqref{eq:memphis-vertical-rho} are equal to $\rho^k\big(\SS^{k};\T,\pi(u,W^\perp) \big)$. If $u = v$ then there is nothing to prove, so let us also assume that $u \not\in W^\perp$. 

Let $\lambda > 0$ and $z \in \<W\>$ be such that $u = \lambda v + z$, and let $w \in C(\SS^k;\T)$. Since $w \in W^\perp$, we have
\begin{equation}\label{eq:memphis-uv-ip}
\< u,w \> = \lambda \< v,w \>.
\end{equation}
Therefore, $u \in \edge H_w$ if and only if $v \in \edge H_w$, where $H_w$ denotes the closed hemisphere of $\SS^{k}$ centred at $w$, and $\edge H_w$ denotes the boundary of $H_w$. It follows that if $\eta > 0$ is sufficiently small, then $\edge H_w$ either contains both $u$ and $v$, or it has empty intersection with both 
$$\SS_u := S_\eta(\SS^{k},u) \qquad \text{and} \qquad  \SS_v := S_\eta(\SS^{k},v).$$

Next, observe that $M_{uv}(\SS_u) = \SS_v$, where $M_{uv}$ is the linear map that rotates $u$ to $v$ in the 2-dimensional plane $P := \< \{u,v\} \>$ and fixes $P^\perp$. We claim that
\begin{equation}\label{eq:memphis-vertical-claim}
x \in \T \cap \SS_u \qquad \Leftrightarrow \qquad M_{uv}(x) \in \T \cap \SS_v.
\end{equation} 
This will suffice to prove~\eqref{eq:memphis-vertical-rho}, since it implies that $\T \cap \SS_u \equiv \T \cap \SS_v$ (cf.~Lemma~\ref{lem:inducedwelldef}), and since $u \in \T$ if and only if $v \in \T$, by~\eqref{eq:memphis-uv-ip} and the definition of $C(\SS^k;\T)$. 

To prove~\eqref{eq:memphis-vertical-claim}, it is enough to show that $\< x,w \> \ge 0$ if and only if $\< M_{uv}(x),w \> \ge 0$, for each $x \in \SS_u$ and $w \in C(\SS^k;\T)$, again by the definition of $C(\SS^k;\T)$. Now, if $\edge H_w$ intersects neither $\SS_u$ nor $\SS_v$ then, by~\eqref{eq:memphis-uv-ip}, we have
$$\< x,w \> \ge 0 \quad \Leftrightarrow \quad \< u,w \> \ge 0  \quad \Leftrightarrow \quad \< v,w \> \ge 0  \quad \Leftrightarrow \quad  \< M_{uv}(x),w \> \ge 0$$
for any $x \in \SS_u$, where the third equivalence holds since $M_{uv}(x) \in \SS_v$. 

We may therefore assume that $\edge H_w$ contains both $u$ and $v$, and hence that $w \in P^\perp$. Since $M_{uv}$ is a linear map that fixes $P^\perp$, it follows that $\< x,w \> = \< M_{uv}(x),w \>$ for every $x \in \SS_u$, and this completes the proof.
\end{proof}

For each $v \in \SS^k$, set $\SS_v := S_\eta(\SS^k,v)$ and define $\varphi_v \colon \SS_v  \to \SS^k \cap \{v\}^\perp$ to be the homothety
\begin{equation}\label{def:phiv:homothety}
\varphi_v(u) := \pi(u,\{v\}^\perp)
\end{equation}
obtained by translating the centre of $\SS_v$ to the origin, and then dilating. 

\begin{lemma}\label{lem:C:perp:inherited}
Let\/ $1 \le k \le d - 1$, let $\T$ be an\/ $\SS^k$-stable set, and let $u \in \SS^{k}$ and $v \in \SS_u$. Then 
$$\T' := \varphi_v(\T \cap \SS_v)$$ 
is a $\varphi_v(\SS_v)$-stable set, and there exists a valid choice of $C\big( \varphi_v(\SS_v); \T' \big)$ with 
\begin{equation}\label{eq:memphis-ind-C}
C\big( \varphi_v(\SS_v); \T' \big) \subset C(\SS^{k};\T) \cap \{u,v\}^\perp.
\end{equation}
\end{lemma}

\begin{proof}
We shall construct a valid choice of $C\big( \varphi_v(\SS_v); \T' \big)$ by removing elements from $C(\SS^k;\T)$. Observe first that, intersecting both sides of~\eqref{eq:happy:reminder} with $\SS_v$, we have 
\begin{equation}\label{eq:happy:reminder:2}
\T \cap \SS_v = \bigcap_{i=1}^t \bigcup_{x \in Y_i} \big( H_x \cap \SS_v \big)
\end{equation}
where $C(\SS^k;\T) = Y_1 \cup \dots \cup Y_t$. Moreover, since $C(\SS^k;\T)$ is finite and $\eta$ is sufficiently small, for each $x \in C(\SS^k;\T)$ we have either $x \in \{u,v\}^\perp$ or $\SS_v \cap \{x\}^\perp = \emptyset$. 

Now, if $x \in \{v\}^\perp$, then $\varphi_v\big( H_x \cap \SS_v \big)$ is the closed hemisphere of $\varphi_v(\SS_v)$ centred at $x$. On the other hand, if $\SS_v \cap \{x\}^\perp = \emptyset$, then $H_x \cap \SS_v \in \{\emptyset,\SS_v\}$, and therefore either the element $x$ can be omitted from $Y_i$ (in the case $H_x \cap \SS_v = \emptyset$) or the set $Y_i$ can be omitted entirely (if $H_x \cap \SS_v = \SS_v$), without changing the right-hand side of~\eqref{eq:happy:reminder:2}.

Thus, by removing elements and sets from $C(\SS^k;\T)$ as described above (including all of those elements not in $\{u,v\}^\perp$), we obtain sets $Y'_1, \ldots, Y'_{t'}$ such that 
$$\T' = \varphi_v\big( \T \cap \SS_v \big) = \bigcap_{i=1}^{t'} \bigcup_{x \in Y'_i} \varphi_v\big( H_x \cap \SS_v \big).$$
This proves that $\T'$ is a $\varphi_v(\SS_v)$-stable set, and moreover that there exists a valid choice of $C\big( \varphi_v(\SS_v); \T' \big)$ such that~\eqref{eq:memphis-ind-C} holds. 
\end{proof}

We shall use Lemma~\ref{lem:C:perp:inherited} to construct a vector $w$ in the span of $\{u,v\}$ such that $C\big( \varphi_v(\SS_v); \T' \big) \subset \{w\}^\perp$, and then apply the following lemma, which is the main technical step in the proof of Lemma~\ref{lem:memphis}. It will be convenient to use the notation 
$$\SS_\lambda := \big\{ u \in \SS^k : \< u, w \> = \lambda \big\}$$ 
for each $\lambda \in (-1,1)$, where $w \in \SS^k$ is fixed in the statement of the lemma.  

\medskip
\pagebreak

\begin{lemma}\label{lem:axis}
Let\/ $1 \le k \le d - 1$, and let\/ $\T$ be an\/ $\SS^{k}$-stable set. If\/ $C(\SS^{k};\T) \subset \{w\}^\perp$ for some $w \in \SS^{k}$, then the following hold.
\begin{enumerate}
\item For each $\lambda \in (-1,1)$, we have
\begin{equation}\label{eq:axis}
r^{k+1} \big( \SS^k; \T \big) =
\begin{cases}
k + 2 & \text{if } \/ r^k \big( \SS_\lambda; \T \big) = k + 1, \\
r^k\big( \SS_\lambda; \T \big) & \text{otherwise.}
\end{cases}
\end{equation}
\item For each $u \in \SS^k \setminus \<w\>$, we have
\begin{equation}\label{eq:axis-rho}
\rho^k\big( \SS^k; \T, u \big) =
\begin{cases}
k + 1 & \text{if } \/\rho^{k-1}\big( \SS_0; \T, v \big) = k, \\
\rho^{k-1}\big( \SS_0; \T, v \big) & \text{otherwise,}
\end{cases}
\end{equation}
where $v := \pi\big( u,\{w\}^\perp \big)$.
\end{enumerate}
\end{lemma}

We remark that the only part of Lemma~\ref{lem:axis} that will be needed in this section is the corollary of part~$(a)$ that says that $r^{k+1} \big( \SS^k; \T \big) \geq r^k\big( \SS_\lambda; \T \big)$; the remaining parts will instead be used in Section~\ref{sec:w}.  

\begin{proof}[Proof of Lemma~\ref{lem:axis}]
The proof is by induction on $k$. To prove the base case $k = 1$, observe that the condition $C(\SS^1;\T) \subset \{w\}^\perp$ implies, by~\eqref{eq:happy:reminder}, that
\begin{equation}\label{eq:axis:basecase:cases}
\T \cap \SS^1 \in \big\{ \emptyset, \{w,-w\}, H_v, H_{-v}, \SS^1 \big\},
\end{equation}
where $v \in \SS^1$ is such that $\< v,w\> = 0$. Now, by Definition~\ref{def:r}, we have $r^1( \SS_\lambda; \SS^1) = 2$ and $r^2( \SS^1; \SS^1 ) = 3$, while 
$$r^1( \SS_\lambda; \T) = r^2(\SS^1; \T) = 1$$ 
in each of the remaining cases of~\eqref{eq:axis:basecase:cases}, since there is an open hemisphere avoiding $\T$, and hence~$(a)$ holds. Part~$(b)$ is also straightforward, since $\rho^0\big( \SS_0; \T, v \big) = 1$ if and only if $v \in \T$, which for each set in~\eqref{eq:axis:basecase:cases} is true if and only if $u \in \T$, since $v = \pi\big( u,\{w\}^\perp \big)$. Thus $\rho^{0}\big( \SS_0; \T, v \big) = 0$ implies that $\rho^{1}\big( \SS; \T, u \big) = 0$, whereas $\rho^{0}\big( \SS_0; \T, v \big) = 1$ implies, by~\eqref{eq:axis:basecase:cases}, that $u$ is contained in an open interval of $\T$, and therefore
$$\rho^{1}\big( \SS; \T, u \big) = r^{1}\big( S_\eta(\SS,u); \T \big) = 2,$$
by~\eqref{eq:rho} and~\eqref{eq:r}, as required. Hence the lemma holds when $k = 1$.

For the induction step, let $k \ge 2$ and assume that the lemma holds for $k - 1$. We may assume that $\T \ne \emptyset$, since otherwise both parts of the lemma hold trivially, and by~\eqref{eq:happy:reminder} it follows that $\{ w, -w \} \subset \T$, since $\{ w, -w \} \subset H_u$ for every $u \in C(\SS^k;\T) \subset \{w\}^\perp$. 

We begin by using Lemma~\ref{lem:memphis-vertical} to prove the following claim.

\begin{claim}\label{claim:reduce:to:S0}
$$r^{k} \big( \SS_\lambda; \T \big) = r^{k} \big( \SS_\mu; \T \big)$$
for every $\lambda,\mu \in (-1,1)$. 
\end{claim}

\begin{clmproof}{claim:reduce:to:S0}
By Lemma~\ref{lem:memphis-vertical}, applied with $W = \{w\}$, we have 
$$\rho^k\big( \SS^k;\T,u \big) = \rho^k\big( \SS^k;\T,u' \big),$$
for every $u \in \SS_\lambda$ and $u' \in \SS_\mu$ with $\pi(u,\{w\}^\perp) = \pi(u',\{w\}^\perp)$. In particular, $u \in \T$ if and only if $u' \in \T$, and therefore the sets $\T \cap \SS_\lambda$ and $\T \cap \SS_\mu$ are equivalent. 
\end{clmproof}

The following claim will also be central to both parts of the lemma, and is where we apply the induction hypothesis.

\begin{claim}\label{claim:axis:final}
For every $u \in \SS_0$,
\[
\rho^k \big( \SS^k; \T, u \big) =
\begin{cases}
k + 1 & \text{if }  \rho^{k-1} \big( \SS_0; \T, u \big) = k, \\
\rho^{k-1} \big( \SS_0; \T, u \big) & \text{otherwise.}
\end{cases}
\]
\end{claim}

\begin{clmproof}{claim:axis:final}
The claim holds trivially if $u \not\in \T$, and if $u \in \T$ then by~\eqref{eq:rho} it is equivalent to 
\begin{equation}\label{eq:axis-rho-ind}
r^k \big( \SS_u; \T \big) =
\begin{cases}
k + 1 & \text{if }  r^{k-1} \big( \SS_u'; \T \big) = k, \\
r^{k-1} \big( \SS_u'; \T \big) & \text{otherwise.}
\end{cases}
\end{equation}
where $\SS_u = S_\eta(\SS^k,u)$ and $\SS_u' = \SS_u \cap \SS_0 = S_\eta(\SS_0,u)$, for some sufficiently small $\eta > 0$.

We shall prove~\eqref{eq:axis-rho-ind} using the induction hypothesis. In order to do so, observe first that, by Lemma~\ref{lem:C:perp:inherited}, 
$$\T' = \varphi_u(\T \cap \SS_u)$$
is a $\varphi_u(\SS_u)$-stable set, where $\varphi_u \colon \SS_u \to \SS^k \cap \{u\}^\perp$ is the homothety defined in~\eqref{def:phiv:homothety}, and there exists a valid choice of $C\big( \varphi_u(\SS_u); \T' \big)$ with 
\[
C\big( \varphi_u(\SS_u); \T' \big) \subset C(\SS^{k};\T) \subset \{w\}^\perp.
\]
Moreover, since $\SS_u' = \SS_u \cap \SS_0$ and $\<u,w\> = 0$, we have
\[
\varphi_u( \SS_u' ) = \big\{ u' \in \varphi_u(\SS_u) : \< u' , w \> = 0 \big\}.
\]
Indeed, recalling that $\varphi_u(x) = \pi(u,\{u\}^\perp)$, we have $x = z + \lambda \cdot \varphi_u(x)$ for some $z \in \< u \>$ and $\lambda > 0$. Since $\<u,w\> = 0$, and $x \in \SS_u$ satisfies $\<x,w\> = 0$ if and only if $x \in \SS_u'$, it follows that $\<\varphi_u(x),w\> = 0$ if and only if $x \in \SS_u'$, as claimed.

Since $\varphi_u(\SS_u)$ is a copy of $\SS^{k-1}$, it follows by the induction hypothesis (with $\lambda = 0$) that
\[
r^k \big( \varphi_u(\SS_u); \T' \big) =
\begin{cases}
k + 1 & \text{if } r^{k-1} \big( \varphi_u(\SS_u'); \T' \big) = k, \\
r^{k-1} \big( \varphi_u(\SS_u'); \T' \big) & \text{otherwise.}
\end{cases}
\]
Hence~\eqref{eq:axis-rho-ind} holds, since $\varphi_u$ is a homothety, and this proves the claim.
\end{clmproof}

Now we can put the components of the proof together. Let us deal with the second part of the lemma first, since it is an immediate consequence of Lemma~\ref{lem:memphis-vertical} and Claim~\ref{claim:axis:final}. Indeed, if $u \in \SS^k \setminus \<w\>$, then by Lemma~\ref{lem:memphis-vertical} we have
\begin{equation}\label{eq:memphis-vertical:axisapp}
\rho^k \big( \SS^k; \T, u \big) = \rho^k \big( \SS^k; \T, v \big),
\end{equation}
where $v := \pi\big( u,\{w\}^\perp \big)$, and by Claim~\ref{claim:axis:final} we have
\[
\rho^k \big( \SS^k; \T, v \big) =
\begin{cases}
k + 1 & \text{if }  \rho^{k-1} \big( \SS_0; \T, v \big) = k, \\
\rho^{k-1} \big( \SS_0; \T, v \big) & \text{otherwise.}
\end{cases}
\]
Hence~\eqref{eq:axis-rho} holds, as required.

Turning to part~$(a)$, observe that, by Claim~\ref{claim:reduce:to:S0}, it will suffice to prove the case $\lambda = 0$. We will first show that for every open hemisphere $H$ of $\SS^k$, there exists $u \in H$ such that
\begin{equation}\label{eq:memphis-ind-finish}
\rho^k \big( \SS^k; \T, u \big) \geq
\begin{cases}
k + 1 & \text{if } r^k \big( \SS_0; \T \big) = k + 1, \\
r^k \big( \SS_0; \T \big) - 1 & \text{otherwise.}
\end{cases}
\end{equation}
By~\eqref{eq:r}, this will prove the claimed lower bound on $r^{k+1} \big( \SS^k; \T \big)$. 

Suppose first that $w \in H$, and recall that $w \in \T$. By~\eqref{eq:rho} and Claim~\ref{claim:reduce:to:S0}, it follows that 
$$\rho^{k} \big( \SS^{k}; \T, w \big) = r^{k} \big( S_\eta(\SS^{k},w); \T \big) = r^{k} \big( \SS_0; \T \big),$$
which suffices to prove both cases of~\eqref{eq:memphis-ind-finish}. We are similarly done if $-w \in H$, so we may assume that $H \cap \{w,-w\} = \emptyset$ and therefore $H \cap \SS_0$ is an open hemisphere of $\SS_0$. 

Choose $u \in H \cap \SS_0$ to maximize $\rho^{k-1}( \SS_0; \T, u)$, and observe that, by~\eqref{eq:r},
\begin{equation}\label{eq:axis-r-final}
\rho^{k-1} \big( \SS_0; \T, u \big) \geq r^k \big( \SS_0; \T \big) - 1.
\end{equation}
In particular, if $r^k \big( \SS_0; \T \big) = k + 1$ then $\rho^{k-1} \big( \SS_0; \T, u \big) = k$. Thus, by~\eqref{eq:axis-r-final} and Claim~\ref{claim:axis:final}, we obtain~\eqref{eq:memphis-ind-finish} as claimed.

Finally, observe that if $r^k \big( \SS_0; \T \big) \leq k$, then by~\eqref{eq:r} there exists an open hemisphere $H \subset \SS_0$ such that
$$\rho^{k-1}(\SS_0;\T,u) \leq r^k(\SS_0;\T) - 1 \leq k - 1$$
for every $u \in H$. By Claim~\ref{claim:axis:final}, it follows that 
\begin{equation}\label{eq:axis-eq-stp}
\rho^k \big( \SS^k; \T, u \big) = \rho^{k-1} \big( \SS_0; \T, u \big) \le r^k(\SS_0;\T) - 1
\end{equation}
for every $u \in H$. Let $H'$ be the open hemisphere of $\SS^k$ such that $H' \cap \SS_0 = H$. By~\eqref{eq:memphis-vertical:axisapp}, it follows from~\eqref{eq:axis-eq-stp} that
$$\rho^k \big( \SS^k; \T, u \big) \leq r^k \big( \SS_0; \T \big) - 1$$
for every $u \in H'$, and therefore $r^{k+1}(\SS^k; \T) \le r^k( \SS_0; \T)$, by~\eqref{eq:r}. 
Since the bound $r^{k+1}( \SS^k; \T) \le k + 2$ holds for every $\T$, this proves the claimed upper bounds, and hence completes the proof of the lemma. 
\end{proof}

We are finally ready to prove the second main lemma of the section. 

\begin{proof}[Proof of Lemma~\ref{lem:memphis}]
We may assume that $\SS = \SS^k$, by Observation~\ref{lem:moving:to:Sk}, and that $v \in \T$ (and thus also $u \in \T$), since otherwise both sides of~\eqref{eq:memphis} are equal to $0$. We are therefore done if $k = 1$, since the right-hand side of~\eqref{eq:memphis} is at most 1 (since the range of $\rho^0$ is $\{0,1\}$), and we have excluded the possibility that the left-hand side is zero. Similarly, if $k = 2$ then we are done if $\rho^{1} \big( S_\eta(\SS,u); \T, v \big) \le 1$, and if $\rho^{1} \big( S_\eta(\SS,u); \T, v \big) = 2$ then we are instead done by Lemma~\ref{lem:memphis-subcrit}. Hence we may assume that $k \ge 3$. 

Let $\SS_u := S_\eta(\SS^k,u)$ and $\SS_v := S_\theta(\SS^k,v)$, where $0 < \theta \ll \eta \ll 1$ are both sufficiently small, and observe that $\SS_u \cap \SS_v = S_\theta(\SS_u, v)$. We shall use Lemma~\ref{lem:axis} to show that
\begin{equation}\label{eq:axis-conseq}
r^k\big( \SS_v; \T \big) \geq r^{k-1}\big( \SS_u \cap \SS_v; \T \big).
\end{equation}
Since this is equivalent to~\eqref{eq:memphis} by~\eqref{eq:rho}, it is enough to prove~\eqref{eq:axis-conseq}. 

Observe first that, by Lemma~\ref{lem:C:perp:inherited}, 
$$\T' = \varphi_v(\T \cap \SS_v)$$ 
is a $\varphi_v(\SS_v)$-stable set, where $\varphi_v \colon \SS_v  \to \SS^k \cap \{v\}^\perp$ is the homothety defined in~\eqref{def:phiv:homothety}, and there exists a valid choice of $C\big( \varphi_v(\SS_v); \T' \big)$ with 
\begin{equation}\label{eq:memphis-ind-C:again}
C\big( \varphi_v(\SS_v); \T' \big) \subset \{u,v\}^\perp.
\end{equation}

Now, set $w' := u - \< u,v \> v$ and $w := w' / \| w' \|$, and note that $w' \in \{v\}^\perp$, so $w \in \varphi_v(\SS_v)$. Moreover, since $w \in \< \{u,v\} \>$ it follows from~\eqref{eq:memphis-ind-C:again} that $C\big( \varphi_v(\SS_v); \T' \big) \subset \{w\}^\perp$. In order to apply Lemma~\ref{lem:axis}, it therefore suffices to prove the following claim.  

\begin{claim}\label{clm:memphis-cond}
$\varphi_v(\SS_u \cap \SS_v) = \big\{ x \in \varphi_v(\SS_v) : \< x, w \> = \lambda \big\}$ for some $\lambda \in (-1,1)$. 
\end{claim}

\begin{clmproof}{clm:memphis-cond}
Recall that $\varphi_v(\SS_v) = \SS^k \cap \{v\}^\perp$. If $x \in \SS_u \cap \SS_v$ then 
$$\< x,w' \> = \< x,u \> - \< u,v \> \< x,v \>$$
is a constant, since $x \in \SS_u$ implies that $\< x,u \>$ is constant (because $\SS_u$ is the intersection of $\SS^k$ with a translate of $\{u\}^\perp$) and $x \in \SS_v$ implies that $\< x,v \>$ is constant (because $\SS_v$ is the intersection of $\SS_u$ with a translate of $\{v\}^\perp$). It follows that $\< \varphi_v(x),w \>$ is also constant, since $\varphi_v(x) = c(x - v')$ for some constant $c$, where $v'$ is the centre of $\SS_v$. 

Now, set $\lambda := \< \varphi_v(x), w\>$, and note that $\lambda \in [-1,1]$, since $w \in \SS^k$ and $\varphi_v(x) \in \SS^k$ for every $x \in \SS_v$. To show that $\lambda \in (-1,1)$, observe that $\SS_u \cap \SS_v$ is a $(k-2)$-sphere, and recall that $k \ge 3$, so $\varphi_v(\SS_u \cap \SS_v)$ cannot be contained in $\{-w,w\}$. 
\end{clmproof}

Applying Lemma~\ref{lem:axis} to the sphere $\varphi_v(\SS_v)$ (which is a copy of $\SS^{k-1}$), we obtain
$$r^k\big( \varphi_v(\SS_v); \T' \big) \geq r^{k-1}\big( \varphi_v(\SS_u \cap \SS_v); \T' \big),$$
which implies~\eqref{eq:axis-conseq}, since $\varphi_v$ is a homothety. This completes the proof of the lemma.
\end{proof}

\section{Finding a rational direction in which to grow}\label{sec:w}

In this section we will prove a key lemma, which provides us with a suitable rational direction in which to grow on each face of our droplet. Our family of quasistable directions will be defined (in Section~\ref{sec:quasi}) in terms of the constant $R$ given by this lemma. 

Recall that $\U$ is a fixed $d$-dimensional update family, and let us write
\begin{equation}\label{def:SW}
\S_W := \S(\U[W]) \cap \SS(W) 
\end{equation}
for the stable set of $\U[W]$ in $\SS(W)$. Recall also from~\eqref{def:LR} that $w \in \cL_R$ if the line $\<w\>$ contains a point of the set $\big\{ x \in \Z^d \setminus \{\0\} : \|x\| \le R \big\}$, and that we say that an update family $\F$ is trivial if $\F = \emptyset$ or $\emptyset \in \F$, and that $\F$ is non-trivial otherwise. 

\begin{lemma}\label{lem:exists-rational-w}
There exists $R = R(\U) > 0$ such that the following holds. If\/ $W \subset \SS_\Q^{d-1}$ is such that\/ $\U[W]$ is non-trivial, then there exists $w \in \cL_R \cap \SS(W)$ such that
\begin{equation}\label{eq:w-rho-condition}
\rho^{k-1}\big( \SS(W); \S_W, u \big) \leq r^k\big( \SS(W); \S_W \big) - 1
\end{equation}
for all $u \in \SS(W)$ such that $\< u,w \> > 0$, where $k = \dim(W^\perp)$.
\end{lemma}

The main point of this lemma is the requirement that $w \in \cL_R$. Indeed, it follows immediately from Definition~\ref{def:r} that there exists $w \in \SS(W)$ such that~\eqref{eq:w-rho-condition} holds for all $u \in \SS(W)$ such that $\< u,w \> > 0$. 

We will deduce Lemma~\ref{lem:exists-rational-w} from the following slightly more general statement, which does not involve the induced update families $\U[W]$. Recall from~\eqref{def:span} that we write $\<\F\>$ for the span of an update family $\F$.

\begin{lemma}\label{lem:exists-rational-w:general}
For each non-trivial update family\/ $\F$, there exists $R > 0$ such that the following holds. If\/ $W \subset \SS_\Q^{d-1}$ is such that $\<\F\> \subset W^\perp$, then there exists $w \in \cL_R \cap \SS(W)$ such that
\begin{equation}\label{eq:w-rho-condition:general}
\rho^{k-1}\big( \SS(W); \S(\F), u \big) \leq r^k\big( \SS(W); \S(\F) \big) - 1
\end{equation}
for all $u \in \SS(W)$ such that $\< u,w \> > 0$, where $k = \dim(W^\perp)$.
\end{lemma}

Applying Lemma~\ref{lem:exists-rational-w:general} with $\F = \U[W]$ for each set $W \subset \SS_\Q^{d-1}$, we obtain Lemma~\ref{lem:exists-rational-w}, since there are only a bounded number of update families with radius at most $R_0$ (see~\eqref{def:radius}), and since the span of $\U[W]$ is contained in $W^\perp$ by Definition~\ref{def:induced}.

The main step in the proof of Lemma~\ref{lem:exists-rational-w:general} will be the following lemma, which deals with the case $\<\F\> = W^\perp$. We will then be able to deduce the general case using Lemma~\ref{lem:axis}.

\begin{lemma}\label{lem:exists-rational-w:span}
For each non-trivial update family\/ $\F$, there exists $R > 0$ such that the following holds. If\/ $W \subset \SS_\Q^{d-1}$ is such that $\<\F\> = W^\perp$, then there exists $w \in \cL_R \cap \SS(W)$ such that
\begin{equation}\label{eq:w-rho-condition:span}
\rho^{k-1}\big( \SS(W); \S(\F), u \big) \leq r^k\big( \SS(W); \S(\F) \big) - 1
\end{equation}
for all $u \in \SS(W)$ such that $\< u,w \> > 0$, where $k = \dim(W^\perp)$.
\end{lemma}

Let us fix, until the end of the proof of Lemma~\ref{lem:exists-rational-w:span}, a non-trivial update family $\F$, and a set $W \subset \SS_\Q^{d-1}$ such that $\<\F\> = W^\perp$. Set $V := \<\F\>$, so $V = W^\perp$.

\subsection{A rational basis for $V^\perp$}

The first step is to prove the following simple (and standard) lemma. Let us say that a vector $x \in \R^d$ is \emph{rational} if it is in $\Q^d$, that is, if each of its coordinates is rational. Similarly, we will say that a subspace of $\R^d$ has a \emph{rational basis} if it has a basis consisting of rational vectors. 

\begin{lemma}\label{lem:Vperp:rational}
$V^\perp$ has a rational basis.
\end{lemma}

We will use the following lemma to prove Lemma~\ref{lem:Vperp:rational}.   

\begin{lemma}\label{lem:w-gram-schmidt}
Let $1 \le j < d$, and let $x_1,\ldots, x_j \in \Q^d$ be linearly independent vectors. There exists a set of rational vectors
$$\{x_{j+1},\dots,x_d\} \subset \{x_1,\dots,x_j\}^\perp$$
such that $x_1,\ldots, x_d$ are linearly independent. 
\end{lemma}

\begin{proof}
By induction, it will suffice to find a single rational vector $x_{j+1} \in \{x_1,\dots,x_j\}^\perp$. To construct $x_{j+1}$, we will use the standard Gram--Schmidt process, which maps a sequence $(u_1,\ldots,u_k)$ of linearly independent vectors in $\R^d$ to a sequence $(v_1,\ldots,v_k)$ of orthogonal vectors in $\R^d$ such that $\< u_1,\ldots,u_i\>  = \< v_1,\ldots,v_i \>$ for each $i \in [k]$. 

To do so, note that at least one of the standard basis vectors $e_1,\dots,e_d$ does not lie in the span of $\{x_1,\ldots,x_j\}$. Suppose that $e_i \not\in \< x_1, \dots, x_j \>$, and apply the Gram--Schmidt process to the sequence $(x_1,\dots,x_j,e_i)$; we obtain a sequence $(v_1,\ldots,v_{j+1})$, where $v_{j+1} \in  \{x_1,\ldots,x_j\}^\perp$. To see that $v_{j+1}$ is rational, we simply need to observe that each step of the Gram--Schmidt process takes a rational vector to a rational vector.\footnote{This follows simply because if $u$ and $v$ are rational vectors, then the orthogonal projection of $v$ onto the line spanned by $u$ (that is, $\lambda u$, where $\lambda = \<u,v\>  / \<u,u\>$) is a rational vector.}
\end{proof}

We can now easily deduce Lemma~\ref{lem:Vperp:rational}, by applying Lemma~\ref{lem:w-gram-schmidt} to a basis of $V$ consisting of elements of the update rules of $\F$. 

\begin{proof}[Proof of Lemma~\ref{lem:Vperp:rational}]
By the definition of $\<\F\>$, and since each update rule is contained in $\Z^d$, there exists  a rational basis for $V$. By Lemma~\ref{lem:w-gram-schmidt}, it follows that there exists a rational basis of $V^\perp$.
\end{proof}

We take the opportunity to prove another simple consequence of the Gram--Schmidt process, which we shall need in Section~\ref{sec:induced}. 

\begin{obs}\label{obs:proj:is:rational}
If\/ $u \in \SS_\Q^{d-1}$ and\/ $W \subset \SS_\Q^{d-1}$, then $\pi(u,W^\perp) \in \SS_\Q^{d-1}$. 
\end{obs}

\begin{proof}
Since $W \subset \SS_\Q^{d-1}$, there exist linearly independent vectors $x_1,\ldots, x_k \in \Q^d$ such that
$$\< x_1,\ldots,x_k \> = \< W\>.$$
By the Gram--Schmidt process, it follows that there exist rational orthogonal vectors $v_1,\ldots,v_k \in \R^d$ such that $\< v_1,\ldots,v_k \> = \<W\>$. Now, applying the Gram--Schmidt step to a rational vector in $\<u\>$, we obtain a rational vector $u' \in \< W \cup \{u\} \>$ such that $u' \in W^\perp$. Noting that $\pi(u,W^\perp) \in \<u'\>$, it follows that $\pi(u,W^\perp) \in \SS_\Q^{d-1}$, as claimed.
\end{proof}

\subsection{The curved cells of $\SS(W)$}

The next step is to partition $\SS(W)$ into `cells' on which the induced resistance with respect to $\F$ is constant. Let $\lambda = \lambda(\F) > 0$ be a sufficiently large constant so that the set
\begin{equation}\label{def:B}
B_\lambda(\0) = \big\{ x \in \Z^d \setminus \{\0\} : \|x\| \le \lambda \big\}
\end{equation}
contains a rational basis of $V^\perp$, and also all elements $x \in X \in \F$, and set $B := B_\lambda(\0)$. In order to avoid clutter, we will omit the dependence on $W$ and $B$ in the notation below. 

\begin{definition}\label{def:w-curved-cell}
Given $\bfdelta = (\delta_x)_{x \in B} \in \{-1,0,1\}^B$, define the \emph{curved $\bfdelta$-cell of\/ $\SS(W)$} to be \[
C(\bfdelta) := \big\{ u \in \SS(W) \,:\, \sgn\big( \< x,u \> \big) = \delta_x \text{ for all } x \in B \big\}.
\]
\end{definition}

Note that $C(\bfdelta)$ will be empty for many values of $\bfdelta$. With that in mind, let
\[
\D = \D(W) := \big\{ \bfdelta \in \{-1,0,1\}^B \,:\, C(\bfdelta) \neq \emptyset \big\}.
\]
The curved cells have the following important property.

\begin{lemma}\label{lem:w-same-cell}
Let $\bfdelta \in \D$ and let $u,v \in C(\bfdelta)$. Then $u \in \S(\F)$ if and only if $v \in \S(\F)$.
\end{lemma}

\begin{proof}
Observe (cf.~Lemma~\ref{lem:hemispheres}) that
$$\S(\F) \cap \SS(W) = \bigcap_{X \in \F} \bigcup_{x \in X} \big\{ u\in \SS(W) : \< x,u \> \geq 0 \big\},$$
since a direction $u \in \SS(W)$ is unstable for $\F$ if and only if there exists $X \in \F$ such that $\< x, u \> < 0$ for every $x \in X$. Now, let $x \in X \in \F$ and note that, since $u$ and $v$ belong to the same cell and $x \in B$, we have $\< x,u \> \geq 0$ if and only if $\< x,v \> \geq 0$. It follows that $u \in \S(\F)$ if and only if $v \in \S(\F)$, as claimed.
\end{proof}

Moreover, using Lemma~\ref{lem:w-same-cell}, we can deduce the following more general fact. 

\begin{lemma}\label{lem:w-same-rhos}
Let $\bfdelta \in \D$ and let $u,v \in C(\bfdelta)$. Then
\begin{equation}\label{eq:w-same-rhos}
\rho^{k-1}\big( \SS(W); \S(\F), u \big) = \rho^{k-1}\big( \SS(W); \S(\F), v \big).
\end{equation}
\end{lemma}

Given $u,v \in \SS(W)$, let $M_{uv}$ denote the canonical rotation of $\SS(W)$ that maps $u$ to $v$ and fixes all directions in $\{u,v\}^\perp$. We will use $M_{uv}$ to compare the neighbourhoods of $u$ and $v$, and hence prove Lemma~\ref{lem:w-same-rhos}. The following lemma provides the property we need. 

\begin{lemma}\label{lem:w-nhood}
Let $\bfdelta \in \D$ and let $u,v \in C(\bfdelta)$. Let $z \in W^\perp$ be such that $u+z \in \SS(W)$. If $\|z\|$ is sufficiently small, then 
$$u+z \in \S(\F) \qquad \Leftrightarrow \qquad v + M_{uv}(z) \in \S(\F).$$
\end{lemma}

\begin{proof}
Let $\bfdelta' \in \D$ be such that $u+z \in C(\bfdelta')$. We claim that $v+M_{uv}(z) \in C(\bfdelta')$ as well. By Lemma~\ref{lem:w-same-cell}, this will suffice to prove the lemma.

Let $Y = \{ x \in B : \delta_x \ne 0 \}$, and observe first that 
$$\sgn\big( \< x,u+z \> \big) = \sgn\big( \< x,u \> \big) = \sgn\big( \< x,v \> \big) = \sgn\big( \< x,v+M_{u,v}(z) \> \big)$$
for every $x \in Y$, since $\|z\|$ is sufficiently small and $u,v \in C(\bfdelta)$. On the other hand, if $x \in B \setminus Y$, then we claim that
$$\<x,u+z\> = \<x,z\> = \<x,M_{uv}(z)\> = \<x,v+M_{uv}(z)\>,$$
since $\<x,u\> = \<x,v\> = \delta_x = 0$, and since $M_{uv}$ fixes $\{u,v\}^\perp$. Indeed, if $z = z_0 + z_1$, where $z_0 \in \< \{u,v\} \>$ and $z_1 \in \{u,v\}^\perp$, then $\< x,z_0 \> = \< x, M_{uv}(z_0) \> = 0$, since $\<x,u\> = \<x,v\> = 0$, and $M_{uv}(z) = M_{uv}(z_0) + z_1$. Hence $\< x, M_{uv}(z) \>= \< x, z_1 \> = \< x, z \>$, as claimed.

We therefore have $\sgn\big( \< x,u+z \> \big) = \sgn\big( \< x,v+M_{u,v}(z) \> \big)$ for every $x \in B$. Since $u+z \in C(\bfdelta')$, this implies that $v+M_{uv}(z) \in C(\bfdelta')$. By Lemma~\ref{lem:w-same-cell}, it follows that $u+z \in \S(\F)$ if and only if $v+M_{uv}(z) \in \S(\F)$, as required.
\end{proof}

Lemma~\ref{lem:w-same-rhos} is now an immediate consequence of Lemma~\ref{lem:w-nhood}. 

\begin{proof}[Proof of Lemma~\ref{lem:w-same-rhos}]
Recall from Definition~\ref{def:r} that $\rho^{k-1}\big( \SS(W); \S(\F), x \big)$ depends only on the intersection of $\S(\F)$ with $S_\eta(\SS(W),x)$. Lemma~\ref{lem:w-nhood} implies that the intersections for $x = u$ and $x = v$ are equivalent, and~\eqref{eq:w-same-rhos} follows.
\end{proof}

We will need two more simple properties of the curved cells $C(\bfdelta)$, which both say (in slightly differences senses) that the cells are `small'. 

\begin{lemma}\label{lem:w-cells-are-small1}
Let $\bfdelta \in \D$. For each $u,v \in C(\bfdelta)$, we have $\< u,v \> > 0$.
\end{lemma}

\begin{proof}
Let $u = \sum_{i=1}^d \lambda_i e_i$ and $v = \sum_{i=1}^d \mu_i e_i$, where $e_1,\ldots,e_d$ are the usual standard basis vectors. Since $\{e_1,\dots,e_d\} \subset B$, we have $\sgn(\lambda_i) = \sgn(\mu_i) = \delta_{e_i}$ for each $i \in [d]$, and therefore $\< u,v \> = \sum_{i=1}^d \lambda_i \mu_i > 0$, where the final inequality holds because $C(\bfdelta) \ne \emptyset$ implies that at least one of the $\delta_{e_i}$ is non-zero.
\end{proof}

Let us write $\Cbar(\bfdelta)$ for the closure of $C(\bfdelta)$.

\begin{lemma}\label{lem:w-cells-are-small}
Let $\bfdelta \in \D$. There exists $u \in \SS(W)$ such that $\< x,u \> > 0$ for all $x \in \Cbar(\bfdelta)$.
\end{lemma}

\begin{proof}
Let $\{v_1,\ldots,v_k\}$ be a basis of $W^\perp$, where each $v_i = \pi(y_i,W^\perp)$ for some $y_i \in B$.\footnote{In fact, since $\<\F\> = W^\perp$, we could just choose a basis of $W^\perp$ using elements of $B$. However, we would like to emphasize that we do not need to use that assumption in this lemma.} Set
$$u := \sum_{i=1}^k \delta_{y_i} v_i,$$
and observe that $u \in W^\perp$, and that $u \ne \0$, since the $v_i$ form a basis of $W^\perp$, and the $\delta_{y_i}$ cannot all be zero (as then $C(\bfdelta) \subset \SS(W)$ would be empty). 

Now, if $x \in \Cbar(\bfdelta)$ then for each $i \in [k]$ we have $\delta_{y_i} \< x,y_i \> \ge 0$, and hence $\delta_{y_i} \< x,v_i \> \ge 0$, by Observation~\ref{obs:projection:innerproduct}, since $x \in W^\perp$. It follows that
\[
\< x,u \> = \sum_{i=1}^k |\< x,v_i \>| > 0,
\]
since $\{v_1,\ldots,v_k\}$ is a basis of $W^\perp$ and $x \in \SS(W)$. Thus, the lemma holds for $u / \|u\|$.
\end{proof}

Let us fix, for each $\bfdelta \in \D$, an element $c_{\bfdelta} \in \SS(W)$ such that $\< x, c_{\bfdelta} \> > 0$ for all $x \in \Cbar(\bfdelta)$. The existence of such a $c_{\bfdelta}$ is guaranteed by Lemma~\ref{lem:w-cells-are-small}.

\subsection{The flat cells of $\SS(W)$}

In what follows, we would like to talk about the (geodesic) convexity of the cells $C(\bfdelta)$ and their closures $\Cbar(\bfdelta)$. In order to avoid a number of technicalities, we shall instead relate the curved cells to a collection of flat analogues, which are gnomonic projections of the curved cells.

\begin{definition}
For each $\bfdelta \in \D$, recall that $\< x, c_{\bfdelta} \> > 0$ for all $x \in \Cbar(\bfdelta)$, and define
\[
H_{\bfdelta} := \big\{ x \in \SS(W) : \< x, c_{\bfdelta} \> > 0 \big\} \quad \text{and} \quad \Pi_{\bfdelta} := \big\{ x \in W^\perp : \< x, c_{\bfdelta} \> = 1 \big\},
\]
and let $\pi_{\bfdelta} \colon H_{\bfdelta} \to \Pi_{\bfdelta}$ denote the gnomonic projection $\pi_{\bfdelta}(x) := x / \< x, c_{\bfdelta} \>$. Now define the \emph{flat cell on $\bfdelta$}, and its closed analogue, by
\[
F(\bfdelta) := \big\{ \pi_{\bfdelta}(x) \,:\, x \in C(\bfdelta) \big\} \qquad \text{and} \qquad \Fbar(\bfdelta) := \big\{ \pi_{\bfdelta}(x) \,:\, x \in \Cbar(\bfdelta) \big\}.
\]
\end{definition}

Note that $\Fbar(\bfdelta)$ is the intersection of the hyperplane $\Pi_{\bfdelta} \subset W^\perp$ with a finite number of closed half-spaces\footnote{To be precise, the set of closed half-spaces with normals in $B$ which define $\Cbar(\bfdelta)$. Since $\pi_{\bfdelta}(x) \in \<x\>$, the boundary of $\Fbar(\bfdelta)$ in $\Pi_{\bfdelta}$ is defined by the same half-spaces.}
 in $W^\perp$, and is therefore a convex 
polytope in $W^\perp$. 

We will use the following classical fact, which is a simple consequence of the Krein--Milman theorem, and was first proved for finite dimensional spaces by Steinitz~\cite{S16}. 

\begin{theorem}[The Krein--Milman theorem]\label{lem:KMthm}
If $K$ is a convex and compact subset of $\R^d$, then $K$ is equal to the convex hull of its extreme points.
\end{theorem}

In order to apply Theorem~\ref{lem:KMthm}, we will need the following simple observation. 

\begin{obs}\label{obs:KM:extreme:points}
If $\bfdelta \in \D$, then $\Fbar(\bfdelta)$ is a convex and compact polytope. 
\end{obs}

\begin{proof}
We observed above that $\Fbar(\bfdelta)$ is a convex polytope in $W^\perp$. It is bounded because
\[
\sup_{x \in \Fbar(\bfdelta)} \| x \| = \sup_{x \in \Cbar(\bfdelta)} \frac{1}{| \< x, c_{\bfdelta} \> |} < \infty,
\]
the final inequality following because $\< x, c_{\bfdelta} \> > 0$ for all $x \in \Cbar(\bfdelta)$, and since $\Cbar(\bfdelta)$ is compact. It follows that $\Fbar(\bfdelta)$ is convex and compact, as claimed. 
\end{proof}

We will need the following simple consequence of Theorem~\ref{lem:KMthm} and Observation~\ref{obs:KM:extreme:points}.

\begin{lemma}\label{lem:KM:extreme:points}
Let $\bfdelta \in \D$ and $w \in \SS^{d-1}$. If\/ $\Fbar(\bfdelta) \cap \{w\}^\perp$ is non-empty, then it is the convex hull of its extreme points. 
\end{lemma}

\begin{proof}
Since $\Fbar(\bfdelta)$ is convex and compact, it follows that  $\Fbar(\bfdelta) \cap \{w\}^\perp$ is also convex and compact. Hence, by Theorem~\ref{lem:KMthm}, it is the convex hull of its extreme points.
\end{proof}

Let us write $E(\bfdelta)$ for the set of extreme points of $\Fbar(\bfdelta)$, and define
\[
\A := \bigg( \bigcup \big\{ A^\perp : A \subset B, \; \dim(A^\perp) = 1 \big\} \bigg) \setminus \{\0\}.
\]
The following lemma will allow us to restrict our attention to the set $\A$. 

\pagebreak

\begin{lemma}\label{lem:w-E-subset-A}
If $\bfdelta \in \D$, then $E(\bfdelta) \subset \A$.
\end{lemma}

\begin{proof}
If $x \in E(\bfdelta)$, then there exists a set $Y \subset B$ such that
$$\Pi_{\bfdelta} \cap \bigcap_{y \in Y} \big\{ x \in \R^d : \<x,y\> = 0 \big\} = \{x\}.$$
Now, by the definition~\eqref{def:B} of $B$, there exists a basis $Z \subset B$ for $V^\perp$. Since $V = W^\perp$ and $\Pi_{\bfdelta} = \{ x \in W^\perp : \< x, c_{\bfdelta} \> = 1 \}$, it follows that the set $A = Y \cup Z \subset B$ satisfies
$$\dim(A^\perp) = 1 \qquad \text{and} \qquad A^\perp \cap \Pi_{\bfdelta} = \{x\},$$ 
as required. 
\end{proof}

Recall from~\eqref{def:LR} that if $R' > 0$ then we write $\cL_{R'}$ for the set of $w \in \R^d$ such that $w \in \<x\>$ for some $x \in \Z^d$ with $\|x\| \le R'$. The set $\A$ has the following crucial property. 

\begin{lemma}\label{lem:A:rational}
There exists $R' > 0$ such that 
$$\A \subset \cL_{R'}.$$
\end{lemma}

\begin{proof}
Let $x \in \A$, and let $A \subset B$ with $\dim(A^\perp) = 1$ be such that $x \in A^\perp$. Let $\{x_1,\ldots,x_{d-1}\} \subset A \subset B$ be a linearly independent set. By Lemma~\ref{lem:w-gram-schmidt}, applied with $j = d - 1$, there exists a rational vector
$$x_d \in \{x_1,\dots,x_{d-1}\}^\perp.$$
Since there are only a bounded number of subsets of $B$, it follows that there exists $R' > 0$ such that $x \in \cL_{R'}$ for every $x \in \A$, as required. 
\end{proof}

\subsection{Rotating from $w$ to $w'$}

We are now ready to state the key technical lemma in the proof of Lemma~\ref{lem:exists-rational-w}. The lemma has two parts: the first will be used (together with Lemma~\ref{lem:w-same-rhos}) to prove~\eqref{eq:w-rho-condition}, while the second will be used (together with Lemma~\ref{lem:A:rational}) to show that $w' \in \cL_R$. For each $w \in \SS^{d-1}$, let us write
$$\bH^+_w = \{x \in \R^d : \<x,w\> > 0\}$$
for the continuous half-space\footnote{Recall that we write $\H_w$ for the discrete half-space in the opposite direction.} in direction $w$. 

\begin{lemma}\label{lem:w-exists-w'}
For each $w \in \SS(W)$ there exists $w' \in \SS(W)$ such that
\begin{equation}\label{eq:w-same-cells-2}
C(\bfdelta) \cap \bH^+_{w'} \neq \emptyset \qquad \Rightarrow \qquad C(\bfdelta) \cap \bH^+_{w} \neq \emptyset,
\end{equation}
for every $\bfdelta \in \D$, and such that
\begin{equation}\label{eq:w-span-w'}
\big\< W^\perp \cap \{w'\}^\perp \cap \A \big\> = W^\perp \cap \{w'\}^\perp.
\end{equation}
\end{lemma}

We will prove Lemma~\ref{lem:w-exists-w'} by induction on the dimension of $\big\< W^\perp \cap \{w\}^\perp \cap \A \big\>$. The induction step is given by the following lemma, which states that when you rotate the hemisphere centred at $w$ just enough to hit a new cell, the boundary (minus the space that is fixed in the rotation) must contain an element of $\A$. 

\begin{lemma}\label{lem:w-rotate-gives-extreme}
Let $U \subset W^\perp$ be a $2$-dimensional subspace and let $w \in U \cap \SS(W)$. Then there exists $w' \in U \cap \SS(W)$ such that 
\begin{equation}\label{eq:w-same-cells-1}
C(\bfdelta) \cap \bH^+_{w'} \neq \emptyset \qquad \Rightarrow \qquad C(\bfdelta) \cap \bH^+_{w} \neq \emptyset,
\end{equation}
for every $\bfdelta \in \D$, and such that
\begin{equation}\label{eq:w-new-extreme}
\big( W^\perp \cap \{w'\}^\perp \cap \A \big) \setminus U^\perp \neq \emptyset.
\end{equation}
\end{lemma}

\begin{proof}
First of all, let us record that~\eqref{eq:w-same-cells-1} is equivalent to
\begin{equation}\label{eq:w-same-cells-F}
F(\bfdelta) \cap \bH^+_{w'} \neq \emptyset \quad \Rightarrow \quad F(\bfdelta) \cap \bH^+_{w} \neq \emptyset,
\end{equation}
since the gnomonic projection $\pi_{\bfdelta}$ preserves the signs of inner products.

Fix an orientation of $U$ and let $M_\theta \colon \R^d \to \R^d$ be the linear map that rotates the oriented $U$ through an angle $\theta$ anticlockwise and fixes $U^\perp$. Define 
$$D(\theta) := \big\{ \bfdelta \in \D : F(\bfdelta) \cap \bH^+_{M_\theta(w)} \neq \emptyset \big\},$$
and let $\theta_0 := \inf\big\{ \theta > 0 : D(\theta) \neq D(0) \big\}$. To see that $\theta_0$ is well-defined, let $\bfdelta'$ be such that $w \in C(\bfdelta')$, and observe that $\< x,w \> > 0$ for all $x \in C(\bfdelta')$, by Lemma~\ref{lem:w-cells-are-small1}, and hence also $\< x,w \> > 0$ for all $x \in F(\bfdelta')$. Since $w \in U$, it follows that $\bfdelta' \in D(0) \setminus D(\pi)$. 

We claim that~\eqref{eq:w-same-cells-F} holds with $w' := M_{\theta_0}(w)$. Note first that if $\theta_0 = 0$ then $w' = w$, and so~\eqref{eq:w-same-cells-F} holds trivially. On the other hand, if $\theta_0 > 0$ then we have
\[
F(\bfdelta) \cap \bH^+_{M_\theta(w)} \neq \emptyset \qquad \Leftrightarrow \qquad F(\bfdelta) \cap \bH^+_{w} \neq \emptyset
\]
for all $0 \le \theta < \theta_0$, by the definitions of $D(\theta)$ and $\theta_0$. Now simply observe that if $F(\bfdelta) \cap \bH^+_{w'} \neq \emptyset$, then $F(\bfdelta) \cap \bH^+_{M_\theta(w)} \neq \emptyset$ for all $\theta < \theta_0$ sufficiently close to $\theta_0$, because the half-space $\bH^+_{w'}$ is open. It follows that~\eqref{eq:w-same-cells-F} holds, as claimed.

It remains to prove~\eqref{eq:w-new-extreme}, which is really the crux of the lemma. Observe first that, by the definition of $D(\theta)$, there exists $\bfdelta_0 \in \D$ such that either $\bfdelta_0 \in D(\theta)$ for all $\theta > \theta_0$ sufficiently close to $\theta_0$, and $\bfdelta_0 \notin D(\theta)$ for all $\theta \in (0,\theta_0)$ (and $\bfdelta_0 \notin D(0)$ if $\theta_0 = 0$), or vice versa. It follows that the set 
$$F := \Fbar(\bfdelta_0) \cap \{w'\}^\perp,$$ 
is non-empty, and therefore, by Lemma~\ref{lem:KM:extreme:points}, $F$ is the convex hull of its extreme points. We claim that there exists an extreme point $e$ of $F$ satisfying
\begin{equation}\label{eq:w-new-extreme-1}
e \not\in U^\perp.
\end{equation} 
To prove~\eqref{eq:w-new-extreme-1}, observe first that there exists $x \in F$ and $\theta_1 \ne \theta_0$ such that 
$$\< x, M_{\theta_1}(w) \> \ne 0.$$
Indeed, by our choice of $\bfdelta_0$, there exists $x \in F(\bfdelta_0)$ that either enters or leaves $\bH^+_{M_\theta(w)}$ at $\theta = \theta_0$, and therefore $x \in \{w'\}^\perp$, and also $\< x, M_{\theta}(w) \>$ varies with $\theta$. It follows that
$$\< M_{-\theta_0}(x), w \> = \< x,w' \> = 0 \ne \< x, M_{\theta_1}(w) \> = \< M_{-\theta_1}(x), w \>,$$ 
where the first and last equalities hold because $M_\theta$ is a rigid transformation, and the second holds since $F \subset \{w'\}^\perp$. Since the rotations $M_\theta$ fix $U^\perp$, it follows that $x \not\in U^\perp$. Now, since $F$ is the convex hull of its extreme points and $x \in F$, it follows that $e \not\in U^\perp$ for some extreme point $e$ of $F$, as claimed. 
 
To complete the proof of~\eqref{eq:w-new-extreme}, we will show that $e$ is an extreme point of $\Fbar(\bfdelta_0)$, which will then imply (by Lemma~\ref{lem:w-E-subset-A}) that $e \in \A$. Suppose, for a contradiction, that there exist distinct $x,y \in \Fbar(\bfdelta_0)$ such that $e$ lies on the straight line between $x$ and $y$. Note that $\{x,y\} \not\subset \{w'\}^\perp$, since $e$ is an extreme point of $F$. Without loss of generality, it follows that $x \in \bH^+_{w'}$ and $y \in \bH^+_{-w'}$. But that is also not possible, because either $F(\bfdelta_0) \cap \bH^+_{w'} = \emptyset$ or $F(\bfdelta_0) \cap \bH^+_{-w'} = \emptyset$, and therefore either $\Fbar(\bfdelta_0) \cap \bH^+_{w'} = \emptyset$ or $\Fbar(\bfdelta_0) \cap \bH^+_{-w'} = \emptyset$. 

It follows that $e$ is an extreme point of $\Fbar(\bfdelta_0)$, and hence that
\[
e \in \Fbar(\bfdelta_0) \cap \A \subset W^\perp \cap \A,
\]
by Lemma~\ref{lem:w-E-subset-A}. Since $F \subset \{w'\}^\perp$, together with~\eqref{eq:w-new-extreme-1} this proves~\eqref{eq:w-new-extreme}.
\end{proof}

Lemma~\ref{lem:w-exists-w'} follows easily from Lemma~\ref{lem:w-rotate-gives-extreme} by induction.

\begin{proof}[Proof of Lemma~\ref{lem:w-exists-w'}]
We will use induction on
\[
d(w) := \dim\Big( \big\< W^\perp \cap \{w\}^\perp \cap \A \big\> \Big).
\]
Note that if $d(w) = k - 1$, then the lemma holds with $w' = w$. On the other hand, if $d(w) \le k - 2$, then there exists a 2-dimensional subspace
\[
U \subset W^\perp \cap \big( W^\perp \cap \{w\}^\perp \cap \A \big)^\perp
\]
such that $w \in U$. By Lemma~\ref{lem:w-rotate-gives-extreme}, it follows that there exists $w' \in U \cap \SS(W)$ such that~\eqref{eq:w-same-cells-2} holds for every $\bfdelta \in \D$, and such that
\begin{equation}\label{eq:w-new-extreme:again}
\big( W^\perp \cap \{w'\}^\perp \cap \A \big) \setminus U^\perp \neq \emptyset.
\end{equation}
Note that $W^\perp \cap \{w\}^\perp \cap \A \subset U^\perp$, by construction, and that 
\[
\big( W^\perp \cap \{w'\}^\perp \cap \A \big) \cap U^\perp = \big( W^\perp \cap \{w\}^\perp \cap \A \big) \cap U^\perp,
\]
since $\{w,w'\} \subset U$. It therefore follows from~\eqref{eq:w-new-extreme:again} that $d(w') > d(w)$. 

Hence, by the induction hypothesis (applied to $w'$), there exists $w'' \in \SS(W)$ such that
$$C(\bfdelta) \cap \bH^+_{w''} \neq \emptyset \quad \Rightarrow \quad C(\bfdelta) \cap \bH^+_{w'} \neq \emptyset \quad \Rightarrow \quad C(\bfdelta) \cap \bH^+_{w} \neq \emptyset$$
for all $\bfdelta \in \D$, and 
$$\big\< W^\perp \cap \{w''\}^\perp \cap \A \big\> = W^\perp \cap \{w''\}^\perp,$$
as required.
\end{proof}

\subsection{The proof of Lemma~\ref{lem:exists-rational-w:span}}

We will next use Lemmas~\ref{lem:w-gram-schmidt},~\ref{lem:w-same-rhos},~\ref{lem:A:rational} and~\ref{lem:w-exists-w'} (and Definition~\ref{def:r}) to find a suitable rational $w' \in \SS(W)$ when $\<\F\> = W^\perp$. 

\begin{proof}[Proof of Lemma~\ref{lem:exists-rational-w:span}]
By Definition~\ref{def:r}, there exists $w \in \SS(W)$ such that 
\begin{equation}\label{eq:w-rho-condition:again}
\rho^{k-1}\big( \SS(W); \S(\F), u \big) \le r^k\big( \SS(W); \S(\F) \big) - 1
\end{equation}
for all $u \in \SS(W) \cap \bH^+_{w}$. By Lemma~\ref{lem:w-exists-w'}, it follows that there exists $w' \in \SS(W)$ such that 
\begin{equation}\label{eq:w-same-cells-2:again}
C(\bfdelta) \cap \bH^+_{w'} \neq \emptyset \quad \Rightarrow \quad C(\bfdelta) \cap \bH^+_{w} \neq \emptyset,
\end{equation}
for all $\bfdelta \in \D$, and such that
\begin{equation}\label{eq:w-span-w':again}
\big\< V \cap \{w'\}^\perp \cap \A \big\> = V \cap \{w'\}^\perp.
\end{equation}
since $V = W^\perp$. Recalling from Lemma~\ref{lem:w-same-rhos} that
\[
\rho^{k-1}\big( \SS(W); \S(\F), u \big) = \rho^{k-1}\big( \SS(W); \S(\F), v \big)
\]
for all $\bfdelta \in \D$ and $u,v \in C(\bfdelta)$, it follows from~\eqref{eq:w-same-cells-2:again} that~\eqref{eq:w-rho-condition:again} holds for all $u \in \SS(W) \cap \bH^+_{w'}$. 

It remains to show that the line $\<w'\>$ intersects $\Z^d$. To see this, recall that $\A \subset \cL_{R'}$ for some $R' > 0$, by Lemma~\ref{lem:A:rational}. By~\eqref{eq:w-span-w':again}, it follows that we may choose a rational basis $\{x_1,\dots,x_{k-1}\} \subset \A$ for the subspace $V \cap \{w'\}^\perp$. Now, by Lemma~\ref{lem:Vperp:rational}, there exists a rational basis $y_{k+1},\dots,y_d$ for $V^\perp$. Hence $\<x_1,\dots,x_{k-1}, y_{k+1},\dots,y_d \> = \{w'\}^\perp$, and it follows, by Lemma~\ref{lem:w-gram-schmidt}, that there exists a rational vector in $\<w'\>$, as claimed. \end{proof}

\subsection{Deducing the general case}

In order to deduce Lemma~\ref{lem:exists-rational-w:general} from Lemma~\ref{lem:exists-rational-w:span}, we will use the following consequence of Lemma~\ref{lem:axis}. 

\begin{lemma}\label{lem:low-dim-to-high-dim}
Let $\F$ be an update family, and let\/ $W \subset W' \subset \SS_\Q^{d-1}$ be such that 
$$\<\F\> \subset W'^\perp \subset W^\perp.$$ 
Set $k := \dim(W^\perp)$ and $k' := \dim(W'^\perp)$, and suppose that $r^k\big( \SS(W); \S(\F) \big) \leq k$. Then 
\begin{equation}\label{eq:low-dim-to-high-dim1}
r^k\big( \SS(W); \S(\F) \big) = r^{k'}\big( \SS(W'); \S(\F) \big) \le k',
\end{equation}
and if $u \in \SS(W) \setminus \<W'\>$, then either
\begin{equation}\label{eq:low-dim-to-high-dim2}
\rho^{k-1}\big( \SS(W); \S(\F), u \big) = \rho^{k'-1}\big( \SS(W'); \S(\F), \pi(u,W'^\perp) \big) \le k' - 1,
\end{equation}
or 
\begin{equation}\label{eq:low-dim-to-high-dim3}
\rho^{k-1}\big( \SS(W); \S(\F), u \big) = k \qquad \text{and} \qquad \rho^{k'-1}\big( \SS(W'); \S(\F), \pi(u,W'^\perp) \big) = k'.
\end{equation}
\end{lemma}

\begin{proof}
Note first that the lemma holds trivially if $k' = k$, since in that case $W'^\perp = W^\perp$ and $\pi(u,W'^\perp) = u$ for every $u \in \SS(W)$. We will first prove the lemma in the case $k' = k - 1$; the general case will then follow easily by induction.

 To do so, we will apply Lemma~\ref{lem:axis} to the $(k-1)$-sphere $\SS(W)$ and the $\SS(W)$-stable set $\S(\F)$.\footnote{Note that $\SS(W)$ is a copy of $\SS^{k-1}$ in the sense used in Section~\ref{sec:memphis}, i.e., it is a $(k-1)$-dimensional sphere in $\SS^{d-1}$, centred at the origin.} Choose $y \in \SS(W) \setminus W'^\perp$ such that $W'^\perp \subset \{y\}^\perp$, and note that if $k' = k - 1$, then $W'^\perp = W^\perp \cap \{y\}^\perp$. 
Now, by Observation~\ref{obs:rules:and:centres}, there is a valid choice of the set of centres of $\S(\F)$ with respect to $\SS(W)$ such that
$$C\big( \SS(W);\S(\F) \big) \subset \<\F\> \subset \{y\}^\perp$$
where the second inclusion holds since $\<\F\> \subset W'^\perp$ and by our choice of $y$. 

By Lemma~\ref{lem:axis}, and since $r^k\big( \SS(W); \S(\F) \big) \leq k$, it follows that
$$r^k \big( \SS(W); \S(\F) \big) = r^{k-1}\big( \SS(W'); \S(\F) \big) \le k - 1,$$
and if $u \in \SS(W) \setminus \<y\>$ then either
$$\rho^{k-1}\big( \SS(W); \S(\F), u \big) = \rho^{k-2}\big( \SS(W'); \S(\F), \pi( u,\{y\}^\perp) \big) \le k - 2,$$
or 
$$\rho^{k-1}\big( \SS(W); \S(\F), u \big) = k \qquad \text{and} \qquad \rho^{k-2}\big( \SS(W'); \S(\F), \pi(u,\{y\}^\perp) \big) = k - 1.$$
Noting that $\pi( u,\{y\}^\perp ) = \pi( u,W'^\perp )$ for $u \in \SS(W) \setminus \<y\>$, since $W'^\perp = W^\perp \cap \{y\}^\perp$, this completes the proof of the lemma in the case $k' = k - 1$. 

We now simply apply induction on $\dim(W^\perp) - \dim(\<\F\>)$. To be precise, let $W'^\perp \subset W''^\perp \subset W^\perp$ with $\dim(W''^\perp) = k - 1$, and observe that $r^{k-1}\big( \SS(W''); \S(\F) \big) \leq k - 1$, by the argument above. By the induction hypothesis, it follows that
$$r^{k-1}\big( \SS(W''); \S(\F) \big) = r^{k'}\big( \SS(W'); \S(\F) \big) \le k'.$$
Moreover, if $u \in \SS(W) \setminus \<W'\>$, then 
$$\pi\big( \pi( u,W''^\perp), W'^\perp \big) = \pi\big( u, W'^\perp \big),$$
and therefore, by the induction hypothesis, either
$$\rho^{k-2}\big( \SS(W''); \S(\F), \pi(u,W''^\perp) \big) = \rho^{k'-1}\big( \SS(W'); \S(\F), \pi(u,W'^\perp) \big) \le k' - 1,$$
or 
$$\rho^{k-2}\big( \SS(W''); \S(\F), \pi(u,W''^\perp) \big) = k - 1 \qquad \text{and} \qquad \rho^{k'-1}\big( \SS(W); \S(\F), \pi(u,W'^\perp) \big) = k'.$$
The result now follows from the case $k' = k - 1$. 
\end{proof}

We are now ready to prove the main lemma of this section, Lemma~\ref{lem:exists-rational-w:general}.

\begin{proof}[Proof of Lemma~\ref{lem:exists-rational-w:general}]
Let $\F$ be a non-trivial update family, and let $W \subset \SS_\Q^{d-1}$ be such that $\<\F\> \subset W^\perp$. Observe first that the left-hand side of~\eqref{eq:w-rho-condition:general} is always at most $k$, by Definition~\ref{def:r}, so if $r^k\big( \SS(W); \S(\F) \big) = k + 1$ then we only need to show that $\cL_R \cap \SS(W)$ is non-empty for some $R = R(\F) > 0$. This follows from our assumption that $\F$ is non-trivial; indeed, there exists a non-empty set $X \in \F$, and $\<\F\> \subset W^\perp$. We may therefore assume that $r^k\big( \SS(W); \S(\F) \big) \le k$.

Fix a set $W \subset W' \subset \SS_\Q^{d-1}$ such that $\<\F\> = W'^\perp$. By Lemma~\ref{lem:exists-rational-w:span}, there exists $R = R(\F) > 0$ and $w \in \cL_R \cap \SS(W')$ such that
\begin{equation}\label{eq:applying:exists-rational-w:span}
\rho^{k'-1}\big( \SS(W'); \S(\F), u' \big) \leq r^{k'}\big( \SS(W'); \S(\F) \big) - 1
\end{equation}
for all $u' \in \SS(W')$ such that $\< u',w \> > 0$, where $k' = \dim(W'^\perp)$. Now, by Lemma~\ref{lem:low-dim-to-high-dim}, 
\begin{equation}\label{eq:applying:low-dim-to-high-dim}
r^k\big( \SS(W); \S(\F) \big) = r^{k'}\big( \SS(W'); \S(\F) \big) \le k',
\end{equation}
and therefore, by~\eqref{eq:applying:exists-rational-w:span},
\begin{equation}\label{eq:exists-rational-w:general:proofstep}
\rho^{k'-1}\big( \SS(W'); \S(\F), u' \big) \le k' - 1
\end{equation}
for all $u' \in \SS(W')$ such that $\< u',w \> > 0$. 

Now, let $u \in \SS(W)$ with $\< u,w \> > 0$, and set $u' := \pi(u,W'^\perp)$. Since $w \in \SS(W')$, we have $\< u',w \> > 0$, by Observation~\ref{obs:projection:innerproduct}. By~\eqref{eq:exists-rational-w:general:proofstep}, it follows that~\eqref{eq:low-dim-to-high-dim3} cannot occur, and therefore, by Lemma~\ref{lem:low-dim-to-high-dim}, we obtain
$$\rho^{k-1}\big( \SS(W); \S(\F), u \big) = \rho^{k'-1}\big( \SS(W'); \S(\F), u' \big).$$ 
Combining this with~\eqref{eq:applying:exists-rational-w:span} and~\eqref{eq:applying:low-dim-to-high-dim}, it follows that
$$\rho^{k-1}\big( \SS(W); \S(\F), u \big) \leq r^k\big( \SS(W); \S(\F) \big) - 1$$
for all $u \in \SS(W)$ such that $\< u,w \> > 0$, as required.
\end{proof}

Finally, let us observe that Lemma~\ref{lem:exists-rational-w:general} implies the existence of the constant $R = R(\U)$ as claimed by Lemma~\ref{lem:exists-rational-w}.

\begin{proof}[Proof of Lemma~\ref{lem:exists-rational-w}]
We simply apply Lemma~\ref{lem:exists-rational-w:general} with $\F = \U[W]$ for each set $W \subset \SS_\Q^{d-1}$. Since each update rule in $\U[W]$ is contained in $W^\perp$,  by~\eqref{eq:induced-def}, it follows that there exists $w \in \cL_R \cap \SS(W)$ such that
$$\rho^{k-1}\big( \SS(W); \S(\F), u \big) \leq r^k\big( \SS(W); \S(\F) \big) - 1$$
for all $u \in \SS(W)$ such that $\< u,w \> > 0$, where $k = \dim(W^\perp)$ and $R$ is a constant depending only on $\F$. Since $\S_W = \S(\U[W]) \cap \SS(W)$, by~\eqref{def:SW}, this implies~\eqref{eq:w-rho-condition}. Moreover, since $\U$ is a finite collection of finite sets, there are only a bounded number of different update families $\U[W]$. It follows that there exists a constant $R = R(\U)$ as claimed. 
\end{proof}

\section{Construction of the quasistable set}\label{sec:quasi}

In this section we construct the set $\QQ$ of quasistable directions, which will be used as the directions of the faces of our droplets. This set will play a crucial role throughout the remainder of the paper; in particular, the droplets that we use (in Section~\ref{proof:sec}) to construct a low energy route to percolation of $\Z_n^d$ will be $\QQ$-droplets. We refer the reader to Section~\ref{sec:chat-quasi} for a discussion of the various properties that we shall require of the set $\QQ$ during the proof of Theorem~\ref{thm:upper}. The definition of $\QQ$ will depend on the constant $R = R(\U)$ given by Lemma~\ref{lem:exists-rational-w}; the reader may, however, prefer to think of the construction in this section as giving an infinite family of sets $\QQ$, one for each $R \in \N$. 

In order to state precisely the key property of $\QQ$, we need to define a graph that encodes which pairs of directions in $\QQ$ are `adjacent'.  

\begin{definition}\label{def:voronoi}
Given a finite set $\X \subset \SS^{d-1}$ and $u \in \X$, the \emph{Voronoi cell} of $u$ with respect to $\X$ is 
$$\Cell_\X(u) := \big\{ w \in \SS^{d-1} : \<u,w\> \ge \<v,w\> \, \text{ for all } \, v \in \X \big\}.$$
The \emph{Voronoi graph} $\Vor(\X)$ has vertex set $\X$ and edge set 
$$E\big(\Vor(\X)\big) := \big\{ uv : \Cell_\X(u) \cap \Cell_\X(v) \ne \emptyset \big\}.$$
\end{definition}

We shall construct droplets from a set of directions by taking a tangents to a suitable sphere for each direction in $\QQ$. Doing this, the faces corresponding to two directions share an edge (or have a lower dimensional intersection) if and only if those directions are adjacent in the Voronoi graph on $\QQ$ (see Section~\ref{sec:polytopes} for the precise statements and proofs of these claims). Recall from~\eqref{def:LR} the definition of $\cL_R$. 

The main aim of this section is to prove the following lemma.

\begin{lemma}\label{lem:quasi-new}
There exists a finite set $\QQ \subset \SS_\Q^{d-1}$, intersecting every open hemisphere of $\SS^{d-1}$, such that if $u \in \QQ$ and $w \in \cL_R$, then
\begin{equation}\label{eq:quasi:property}
\Cell_\QQ(u) \cap \{w\}^\perp \neq \emptyset \qquad \Leftrightarrow \qquad \< u,w \> = 0.
\end{equation}
\end{lemma}

\pagebreak

The condition~\eqref{eq:quasi:property} is the precise version of property (P2) (see Section~\ref{sec:chat-quasi}) that will be used several times during the proof of Theorem~\ref{thm:upper}. If a finite set $\QQ \subset \SS_\Q^{d-1}$ has the properties guaranteed by Lemma~\ref{lem:quasi-new}, then we say that it is \emph{quasistable for range $R$}. 

Before proving Lemma~\ref{lem:quasi-new}, let us deduce from~\eqref{eq:quasi:property} the precise form of property (P1) (see Section~\ref{sec:chat-quasi}) that we shall need, Lemma~\ref{lem:quasi-innerprod}. The following simple property of quasistable sets 
will be useful both in the proof of Lemma~\ref{lem:quasi-innerprod} and also in Sections~\ref{sec:polytopes} and~\ref{buffers:sec}.

\begin{lemma}\label{lem:quasi-cell}
Let $\QQ \subset \SS_\Q^{d-1}$ be quasistable for range $R$, let $u \in \QQ$ and $w \in \cL_R$, and suppose that $\< u,w\> \ne 0$. Then $\< v,w \> \cdot \< u,w \> > 0$ for all $v \in \Cell_\QQ(u)$. 
\end{lemma}

\begin{proof}
Since $\QQ$ is quasistable for range $R$, and since $w \in \cL_R$ and $\< u,w \> \neq 0$, we have from Lemma~\ref{lem:quasi-new} that
\begin{equation}\label{eq:quasi-xperp-empty}
\Cell_\QQ(u) \cap \{w\}^\perp = \emptyset.
\end{equation}
Observe that $\Cell_\QQ(u)$ is path connected (indeed, if $z \in \Cell_\QQ(u)$ then the geodesic from $u$ to $z$ is contained in $\Cell_\QQ(u)$), 
and let $\gamma \colon [0,1] \to \Cell_\QQ(u)$ be a continuous function with $\gamma(0) = u$ and $\gamma(1) = v$. Then since the function $t \mapsto \< \gamma(t), w \>$ is continuous on $[0,1]$, it follows from~\eqref{eq:quasi-xperp-empty} and the intermediate value theorem that $\< u,w \>$ and $\< v,w \>$ are either both strictly positive or both strictly negative.
\end{proof}

We can now easily deduce the following precise version of the property (P1). 

\begin{lemma}\label{lem:quasi-innerprod}
Let $\QQ \subset \SS_\Q^{d-1}$ be quasistable for range $R$. If $uv \in E\big( \Vor(\QQ) \big)$, then there does not exist $w \in \cL_R$ such that
\begin{equation}\label{eq:quasi-innerprod}
\< u,w \> < 0 \qquad \text{and} \qquad \< v,w \> > 0.
\end{equation}
\end{lemma}

\begin{proof}
Let $u,v \in \QQ$ and $w \in \cL_R$ be such that~\eqref{eq:quasi-innerprod} holds. By Definition~\ref{def:voronoi}, it suffices to show that $\Cell_\QQ(u) \cap \Cell_\QQ(v) = \emptyset$, since this implies that $uv \notin E\big( \Vor(\QQ) \big)$. To show this, note that $\< z,w \> < 0$ for all $z \in \Cell_\QQ(u)$ and $\< z,w \> > 0$ for all $z \in \Cell_\QQ(v)$, by Lemma~\ref{lem:quasi-cell} and~\eqref{eq:quasi-innerprod}, and hence $\Cell_\QQ(u) \cap \Cell_\QQ(v) = \emptyset$, as required.
\end{proof}

In two dimensions, it is easy to construct a set satisfying the conclusion of Lemma~\ref{lem:quasi-innerprod}, simply by setting
\begin{equation}\label{eq:quasi-2d}
\QQ \,= \bigcup_{w \,\in\, \Z^d \setminus \{\0\} \,:\, \|w\| \,\le\, R} \big\{ u \in \SS^1 \,:\, \< u,w \> = 0 \big\},
\end{equation}
and this is indeed essentially the construction used in~\cite{BSU,BDMS}. When $d \ge 3$, however, the corresponding set in $\SS^{d-1}$ is infinite, and our construction of a quasistable set 
instead proceeds by choosing a suitable finite subset of this set. The choice of this finite subset is somewhat delicate, and the proof that it satisfies~\eqref{eq:quasi:property} is rather technical; in particular, we shall work in a more general setting in order to facilitate a proof by induction (on the dimension). We therefore first provide some motivation for our approach by describing our construction, and by giving a rough outline of the proof of Lemma~\ref{lem:quasi-new}.

Recall first (from the discussion in Section~\ref{sec:chat-quasi}) that the main difficulty in choosing a set satisfying property (P1) (that is, Lemma~\ref{lem:quasi-innerprod}) occurs when choosing elements of $\QQ$ close to the intersections of two (or more) hyperplanes. Our solution to this difficulty is simple: we define a `buffer' $\cE(H)$ around each such intersection $H = W^\perp$, where $W \subset \cL_R$, and only allow points to be chosen in $\QQ$ if they are either in $H$, or not in $\cE(H)$. The radius of these buffers will be a (rapidly) decreasing function of the dimension of $H$. 

We thus obtain a continuous set $\Qbar \subset \SS^{d-1}$, formed by removing from each sphere $\SS(W)$ the buffers of lower-dimensional spheres $\SS(W')$, and then taking a union over (possibly empty) subsets $W \subset \cL_R$. Our main challenge will be to show that
\begin{equation}\label{eq:quasi:challenge}
d\big( x, \,\Qbar \cap \{w\}^\perp \big) < d\big( x, \,\Qbar \setminus \{w\}^\perp \big)
\end{equation}
for every $w \in \cL_R$ and $x \in \SS^{d-1} \cap \{w\}^\perp$, which is the content of Lemma~\ref{thm:quasi:key}, the main technical lemma of this section. Once we have this fact, it will be straightforward to choose a `dense' finite subset $\QQ \subset \Qbar$, and to show that it satisfies~\eqref{eq:quasi:property}, see Section~\ref{sec:quasi:exists}.

The proof of~\eqref{eq:quasi:challenge} will be by induction on the dimension, and will take up most of the rest of this section. As mentioned before, in order to make the induction work, we shall work in a slightly more general setting, which we define next.

\subsection{Spherical Buffer Systems}

We write $\E^k$ to denote an arbitrary $k$-dimensional Euclidean space.  We remark that the reader should usually think of $\E^k$ as being equal to $\R^k$, but for the induction hypothesis we need the more general setting.
  
\begin{definition}\label{def:SBS}
A \emph{spherical buffer system} is a quadruple $\cB = (\E^{k+1},\SS^k,\cH,\bfdelta)$, where
\begin{itemize}
\item[$(a)$] $\SS^k$ is a $k$-dimensional sphere in $\E^{k+1}$, centred at the origin $O_\cB$ of $\E^{k+1}$;\smallskip
\item[$(b)$] $\cH$ is a finite collection of proper subspaces of $\E^{k+1}$, closed under intersections;\smallskip 
\item[$(c)$] $\bfdelta = (\delta_{-1},\delta_0,\ldots,\delta_{k-1})$ is a sequence of numbers satisfying $$0 < \delta_i < \delta_{i-1}/3$$ 
for each $0 \le i \le k-1$, and with $\delta_{-1}$ equal to the radius of~$\SS^k$.
\end{itemize}
\end{definition}

For the next few definitions, let us fix a spherical buffer system $\cB = (\E^{k+1},\SS^k,\cH,\bfdelta)$. For each $-1 \le i \le k-1$, define 
$$\cH_i = \cH_i(\cB) := \big\{ H \in \cH : \dim(H) = i + 1 \big\}.$$

We need to define two different `buffers' around each $H \in \cH$. 

\begin{definition}\label{def:buffers}
For each $i \in \{0,\dots,k-1\}$ and $H \in \cH_i$, define the \emph{buffer}
\[
 \cE(H) := \big\{ z \in \SS^k : d(z,H) < \delta_i \big\} 
\]
and the \emph{expanded buffer}
\[
 \cE^+(H) := \big\{ z \in \SS^k : d(z,H) < 3\delta_i \big\} 
 \]
of $P$, where $d(.,.)$ denotes the Euclidean distance in $\E^{k+1}$.\end{definition}

\begin{figure}[ht]
  \centering
  \begin{tikzpicture}[>=latex] 
    \draw [name path=A1] (-0.5,3) arc (-110:-80:20);
    \draw [name path=A2] (7,0) arc (0:28:20);
    \draw [name path=A3] (1,0.5) arc (160:130:20);
    \path [name intersections={of=A1 and A2,by=P3}];
    \path [name intersections={of=A2 and A3,by=P1}];
    \path [name intersections={of=A3 and A1,by=P2}];
    \fill (P1) circle (0.1);
    \fill (P2) circle (0.1);
    \fill (P3) circle (0.1);
    \draw [dashed] (-0.45,3.2) arc (-110:-80.2:20);
    \draw [dashed] (-0.55,2.8) arc (-110:-79.8:20);
    \draw [dashed] (6.8,0) arc (0:27.8:20);
    \draw [dashed] (7.2,0) arc (0:28.2:20);
    \draw [dashed] (0.8,0.55) arc (160:129.7:20);
    \draw [dashed] (1.2,0.45) arc (160:130.3:20);
    \draw [dashed] (P1) circle (1.2);
    \draw [dashed] (P2) circle (1.2);
    \draw [dashed] (P3) circle (1.2);
    \draw [<-] (P1) ++(-0.3,0) -- ++(-1.5,0) node [left] {$H_1 \cap H_2$};
    \draw [<-] (P1) ++(0.7,0.1) -- ++(1.5,0) node [right] {$\cE(H_1 \cap H_2)$};
    \draw [<-] (P2) ++(-0.1,0.2) -- ++(-1.2,1.2) node [above left] {$H_1 \cap H_3$};
    \draw [<-] (P2) ++(-0.7,-0.4) -- ++(-1.5,-0.2) node [left] {$\cE(H_1 \cap H_3)$};
    \draw [<-] (P3) ++(0.1,0.2) -- ++(1.2,1.2) node [above right] {$H_2 \cap H_3$};
    \draw [<-] (P3) ++(0.6,-0.5) -- ++(1.5,-0.5) node [right] {$\cE(H_2 \cap H_3)$};
    \coordinate (M1) at ($(-0.5,3)+(70:20)+(-95:20)$);
    \coordinate (M2) at ($(7,0)+(180:20)+(16:20)$);
    \coordinate (M31) at ($(1,0.5)+(-20:20)+(144:20)$);
    \coordinate (M32) at ($(1,0.5)+(-20:20)+(147:20)$);
    \draw [<-] (M1) -- ++(0,1) node [above] {$H_3$};
    \draw [<-] (M1) ++(-0.85,-0.05) -- ++(0,-1) node [below] {$\cE(H_3)$};
    \draw [<-] (M2) -- ++(1,0) node [right] {$H_2$};
    \draw [<-] (M2) ++(0.4,-1) -- ++(1,0) node [right] {$\cE(H_2)$};
    \draw [<-] (M31) -- ++(0.8,-0.6) node [below right] {$H_1$};
    \draw [<-] (M31) ++(0.3,0.6) -- ++(-0.8,0.6) node [left] {$\cE(H_1)$};
    \draw [<->] (P1) ++(0.1,-0.05) -- ++(1,-0.3);
    \draw (P1) ++(0.1,-0.05) ++(0.5,-0.15) arc (-180:-90:1 and 0.6) node [right] {$\delta_0$};
    \draw [-] (M32) -- ++(-0.16,+0.12);
    \draw (M32) ++(-0.2*0.4,0.15*0.4) arc (-35:-90:1 and 0.6) node [left] {$\delta_1$};
  \end{tikzpicture}
  \caption{The figure shows part of the intersection with $\SS^2$ of a spherical buffer system $\cB = (\E^3,\SS^2,\cH,\delta)$. The subspaces $H_1,H_2,H_3 \in \cH$ are 2-dimensional (and so belong to $\cH_1$); their intersections with $\SS^2$ are therefore great circles. The buffers of each of $H_1$, $H_2$ and $H_3$ and their pairwise intersections are also shown.}
  \label{fig:quasi-constr}
\end{figure}
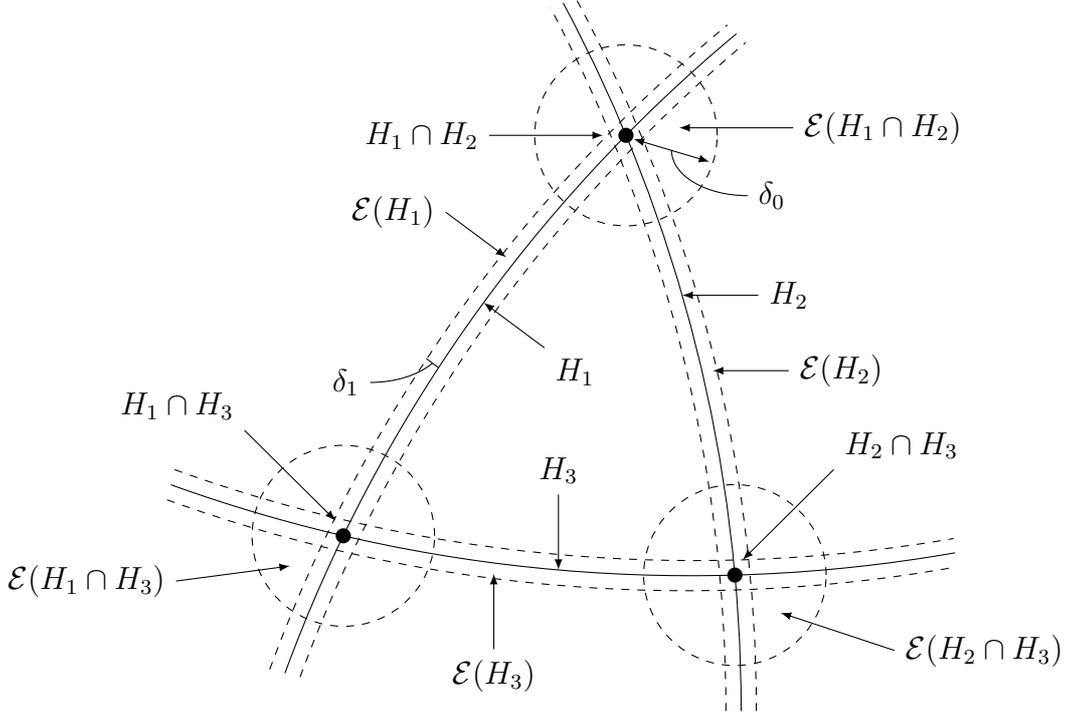

Note that buffers are \emph{open} subsets of $\SS^k$. See Figure~\ref{fig:quasi-constr} for an example of part of a spherical buffer system together with its buffers.

We are now ready to define the set $\Qbar(\cB)$, which can be thought of as a continuous version of the quasistable set. Define $\cH^* := \cH \cup \{ \E^{k+1} \}$ (note that $\E^{k+1} \not\in \cH$), and set $\cE(\{O_\cB\}) := \emptyset$ and moreover $\cE(H) := \emptyset$ for every $H \not\in \cH$.

\begin{definition}\label{def:Qbar}
If $\cB$ is a spherical buffer system, then
\[
 \Qbar(\cB) :=
 \SS^k \cap \bigcup_{H \in \cH^*} \Big( H \setminus \bigcup_{H' \subsetneq H}\cE(H') \Big).
\]
\end{definition}


Our induction hypothesis will be based on the following simple definition. 

\begin{definition}\label{def:SBS:good}
The spherical buffer system $\cB$ is \emph{good\/} if
\[
 \cE^+(H) \cap \cE^+(H') \subset \cE(H \cap H')
\]
for every $H,H' \in \cH$ such that $H \not\subset H'$ and $H' \not\subset H$.
\end{definition}


The following lemma is the main technical result of this section. 

\begin{lemma}\label{thm:quasi:key}
If\/ $\cB = (\E^{k+1},\SS^k,\cH,\bfdelta)$ is a good spherical buffer system and $H \in \cH$, then
 \begin{equation}\label{eq:quasi:key}
  d\big( x, \,\Qbar \cap H \big) < d\big( x, \,\Qbar \setminus H \big)
 \end{equation}
for every\/ $x \in \SS^k \cap H$, where $\Qbar = \Qbar(\cB)$.
\end{lemma}

In Section~\ref{sec:quasi:exists} we shall show that, for any triple $(\E^{k+1},\SS^k,\cH)$ as in Definition~\ref{def:SBS}, we can choose the sequence $\bfdelta$ so as to make $(\E^{k+1},\SS^k,\cH,\bfdelta)$ a good spherical buffer system (see Lemma~\ref{lem:choose:deltas:so:good}), and hence deduce Lemma~\ref{lem:quasi-new} from Lemma~\ref{thm:quasi:key}. 

\subsection{An outline of the proof of Lemma~\ref{thm:quasi:key}}

The proof of Lemma~\ref{thm:quasi:key} is quite technical, and we therefore begin by giving an outline of the argument, to help the reader navigate the details, and to provide some motivation for the more technical definitions introduced during the proof.\footnote{We remark that we do not expect the reader to be able to follow all of the details of this outline on a first reading, but we hope that some readers will nevertheless find it useful as a guide to the proof.} The proof is by induction on $k$. The base case is straightforward\footnote{More precisely, when $k = 0$ the statement is vacuous, and when $k = 1$ it follows easily from Observation~\ref{obs:Qdef3}, below.}, so we shall focus here on the induction step. 

Suppose there exists $y \in \Qbar \setminus H$ with $d(x,y) = d(x,\Qbar)$, let $Q \in \cH$ be minimal such that $y \in Q$, and let $\hat x$ be the closest point of $\SS^k \cap Q$ to $x$. The idea is to try to move $y$ towards $x$ while staying in $\SS^k \cap Q$, by moving along the geodesic from $y$ to $\hat x$. Since we cannot do so without leaving $\Qbar$, it follows that either $y = \hat x$, or there is a buffer in the way, with $y$ on its boundary. This simple idea is made precise in Lemma~\ref{lem:nearest}. 

If $y = \hat x$ then it is not hard to obtain a contradiction, using the observation that if $y \in \Qbar$ and $P \in \cH$, then either $y \in P$ or $y \not\in \cE(P)$ (see Observation~\ref{obs:Qdef3}). To be more precise, we will show that if $P' \in \cH$ is  minimal such that $x \in \cE(P')$ (see Observation~\ref{obs:unique:buffer}), then either $y \in \cE(P')$ or $y \in \cE(P' \cap Q)$. Since $y \notin P'$, because $y \notin H$ and $P' \subset H$ (see~\eqref{eq:showing:PinH}), this will provide the required contradiction. 

To prove that either $y \in \cE(P')$ or $y \in \cE(P' \cap Q)$, we will use the assumption that $\cB$ is good. Indeed, if $P' \subset Q$ then it follows from Observation~\ref{lem:proj} that $y \in \cE(P')$, whereas if $P' \not\subset Q$, then we will show that 
$$y \in \cE^+(P') \cap Q \subset \cE(P' \cap Q).$$
More precisely, we will prove the first inclusion in Lemma~\ref{lem:minbuffer}, and the second holds because $Q \not\subset P'$ (since $y \in Q$ and $y \notin P'$), using our assumption that $\cB$ is good. 

When $y \ne \hat x$, on the other hand, the argument is rather more complicated. Let $P \in \cH$ be minimal such that $y$ is in the boundary of $\cE(P)$. In this case we shall construct (see Definition~\ref{def:BxP} and Lemmas~\ref{lem:Axp:exists} and~\ref{lem:sbs_buffer}) a good spherical buffer system 
$$\cB_{x,P} := \big( A_{x,P}, \SS_{x,P}, \cH_{x,P}, \bfdelta' \big),$$
where $A_{x,P} = \lambda x + P^\perp$ for some $\lambda > 0$, and show that $y \in \SS_{x,P} = \SS^k \cap A_{x,P}$. Our strategy will then depend on whether or not $P \subset H$, and whether or not $x \in P$. 

If $P \not\subset H$ then the proof is similar to the case $y = \hat x$, so assume that $P \subset H$. If $x \in P \subset H$, then we shall find a point of $\Qbar \cap P$ whose distance to $x$ is strictly less than $\delta_j$, where $P \in \cH_j$. Since $y$ lies on the boundary of $\cE(P)$, this will contradict our assumption that $d(x,y) = d(x,\Qbar)$. To find such a point of $\Qbar \cap P$, we escape (one by one) each of the buffers containing $x$ (except $\cE(P)$, since we stay in $P$). For each step, we apply Lemma~\ref{lem:great:escape:step}, which says that we can escape a buffer of dimension $i$ by moving distance at most $2\delta_j^2 / \delta_i$, and in such a way that we do not enter any new buffer. Using our assumption that $\cB$ is good, we can then deduce that we only need to escape one buffer of each dimension strictly less than $j$. For the details, see Section~\ref{sec:GreatEscape}.


Finally, when $x \notin P \subset H$, we apply the induction hypothesis to the good spherical buffer system $\cB_{x,P}$. More precisely, we let $x'$ be one of the closest points of $\SS_{x,P}$ to $x$, show that $x'$ is unique and $x' \in H$, and deduce (by the induction hypothesis) that the closest point of $\Qbar(\cB_{x,P})$ to $x'$ is also in $H$. Since $\Qbar(\cB_{x,P}) = \Qbar(\cB) \cap \SS_{x,P}$ (see Lemma~\ref{lem:sbs:Q}), we will then be able to deduce that $y \in H$, giving us our final contradiction. 

The rest of this section is organised as follows. First, in Section~\ref{sec:obs:warmup}, we make some simple but important observations. Next, in Sections~\ref{sec:BR} and~\ref{sec:BxP}, we construct the `restricted' and `induced' buffer systems $\cB_P$ and $\cB_{x,P}$, show that they are good, and determine the structure of the sets $\Qbar(\cB_P)$ and $\Qbar(\cB_{x,P})$. In Section~\ref{sec:SxP:props} we prove a number of simple but useful properties of $\cB_{x,P}$, in Section~\ref{sec:GreatEscape} we prove the `great escape' lemma, which finds a point of $\Qbar \cap P$ close to $x$, and in Section~\ref{sec:pushing:y} we prove the key fact that either $y = \hat x$, or there is a buffer with $y$ on its boundary. Finally, in Sections~\ref{quasi:technical:proof:sec} and~\ref{sec:quasi:exists}, we complete the proof of Lemma~\ref{thm:quasi:key}, and deduce Lemma~\ref{lem:quasi-new}.

\subsection{Some simple observations to get warmed up}\label{sec:obs:warmup}

We begin with three simple observations, which will be used several times during the proof. 

\begin{obs}\label{obs:unique:buffer}
Let\/ $\cB = (\E^{k+1},\SS^k,\cH,\bfdelta)$ be a good spherical buffer system, and let $x \in \SS^k \setminus \Qbar(\cB)$. Then there exists a unique minimal $P \in \cH$ such that $x \in \cE(P)$. 
\end{obs}

\begin{proof}
By Definition~\ref{def:Qbar}, there exists $P \in \cH$ such that $x \in \cE(P)$. Moreover, since $\cB$ is good, if $x \in \cE(P) \cap \cE(P')$ for some $P,P' \in \cH$ with $P \not\subset P'$ and $P' \not\subset P$, then
\[
x \in \cE^+(P) \cap \cE^+(P') \subset \cE(P \cap P'),
\]
by Definition~\ref{def:SBS:good}. Thus there exists a unique minimal such $P$, as claimed.
\end{proof}

The second observation will provide us with our final contradiction in three of the four cases into which the proof will be divided. 

\begin{obs}\label{obs:Qdef3}
Let $\cB = (\E^{k+1},\SS^k,\cH,\bfdelta)$ be a good spherical buffer system, let $y \in \Qbar(\cB)$, and let $H \in \cH$. Then either $y \in H$ or $y \not\in \cE(H)$. 
\end{obs}

\begin{proof}
By Definition~\ref{def:Qbar}, if $y \in \Qbar(\cB)$ then there exists $P \in \cH^*$ such that 
$$y \in \SS^k \cap P \setminus \bigcup_{H' \subsetneq P} \cE(H').$$ 
Note that if $P \subset H$ then $y \in H$, and that if $H \subsetneq P$ then $y \notin \cE(H)$, so we may assume that $H \not\subset P$ and $P \not\subset H$. Since $\cB$ is good and $H \cap P \subsetneq P$, it follows 
that 
$$y \notin \cE(H \cap P) \supset \cE^+(H) \cap \cE^+(P) \supset \cE(H) \cap P.$$
Since $y \in P$, it follows that $y \notin \cE(H)$, as required.
\end{proof}

The following consequence of Observation~\ref{obs:Qdef3} will play an important role in the proof. 

\medskip
\pagebreak

\begin{obs}\label{obs:compact}
If $\cB = (\E^{k+1},\SS^k,\cH,\bfdelta)$ is a good spherical buffer system and $H \in \cH$, then the sets $\Qbar(\cB)$ and $\Qbar(\cB) \setminus H$ are both compact. 
\end{obs}

\begin{proof}
It follows from Definition~\ref{def:Qbar} that $\Qbar(\cB)$ is compact (even if $\cB$ is not good), since the buffers $\cE(H')$ are open sets. By Observation~\ref{obs:Qdef3}, it follows that $\Qbar(\cB) \setminus H$ is compact (if $\cB$ is good), since $\Qbar(\cB) \setminus H = \Qbar(\cB) \cap \cE(H)^c$. 
\end{proof}

We also need two simple facts about distances to points on a sphere. Given a compact set $U \subset \E^{k+1}$ and a point $x \in \E^{k+1}$, let $M(U,x)$ denote the set of points of $U$ at minimal (Euclidean) distance to $x$. The following observation is a standard geometric fact. 

\begin{obs}\label{obs:spheres:increasing:distance}
Let $\SS \subset \E^{k+1}$ be a sphere (of arbitrary dimension) with centre $z$, and let $x \in \E^{k+1}$. Either 
\begin{itemize}
\item[$(a)$] $M(\SS,x) = \{x'\}$ for some $x' \in \SS$, and the function $y \mapsto d(x,y)$, defined on $\SS$, is increasing in $d(x',y)$, or 
\item[$(b)$] $M(\SS,x) = \SS$ and there exists a subspace $P$ such that $\SS \subset z + P$ and $x \in z + P^\perp$. 
\end{itemize}
\end{obs}

Our final observation is slightly less standard, so we provide a proof. 

\begin{obs}\label{lem:proj}
Let $P \subset Q \subset \E^{k+1}$ be subspaces, let $\SS^k \subset \E^{k+1}$ be the unit sphere centred at the origin, and let $x \in \SS^k$. If $y \in M(\SS^k \cap Q,x)$, then $d(y,P) \le d(x,P)$.
\end{obs}

\begin{proof}
Write $x = (x_1,x_2,x_3) \in P \times (P^\perp \cap Q) \times Q^\perp$. If $\|x_3\| = 1$, then $d(y,P) \le 1 = d(x,P)$, as required, so we may assume that $\|x_3\| < 1$. It follows that $y = (\lambda x_1,\lambda x_2,0)$ where $\lambda \ge 1$ is chosen so that $\|y\| = 1$, and thus $d(x,P)^2 = \|x_2\|^2 + \|x_3\|^2 = 1 - \|x_1\|^2$ and
$d(y,P)^2 = \|\lambda x_2\|^2 = 1 - \|\lambda x_1\|^2$. Since $\lambda \ge 1$, the claimed bound follows.
\end{proof}

\subsection{Restricted buffer systems}\label{sec:BR}

Our first main task is to introduce two different ways in which, given a spherical buffer system $\cB$, another buffer system can be defined on a hyperplane; this will be done in the next two subsections. In order to distinguish these two notions, we shall refer to them as `restricted' and `induced' buffer systems, respectively. The first (and simpler) of the two constructions is as follows.


\begin{definition}\label{def:BR}
Given a spherical buffer system\/ $\cB = (\E^{k+1},\SS^k,\cH,\bfdelta)$ and $P \in \cH^* \setminus \{O_\cB\}$, we define the \emph{restricted} buffer system $\cB_P := (P, \SS^k \cap P, \cH_P, \bfdelta' )$ by setting
\[
 \cH_P := \big\{ H \cap P : P \not\subset H \in \cH \big\}
 \]
and $\bfdelta' := (\delta_{-1},\delta_0,\ldots,\delta_{j-1})$, where $P \in \cH_j$. 
\end{definition}

We remark that the reader may like to think of $\cH_P$ equivalently as
\[
\cH_P := \big\{ H \cap P : H \in \cH \, \text{ and } \, P \cap H \neq P \big\}.
\]
Observe that $\cH_P \subset \cH$, since $P \in \cH$ and $\cH$ is closed under intersections. 

So far we have written $\cE(H)$ for the buffer of $H \in \cH^*$, where $\cB = (\E^{k+1},\SS^k,\cH,\bfdelta)$ is a buffer system. We shall continue to use this notation where the buffer system in question is unambiguously $\cB$, but in other cases (such as when we are considering restricted or induced buffer systems) we shall write $\cE_\cB(H)$ for emphasis. 

The key property of $\cB_P$ that we need is that $\Qbar(\cB_P) = \Qbar(\cB) \cap P$. This is proved in the following lemma, together with the fact that if $\cB$ is good then $\cB_P$ is also good. Both properties follow easily from the definitions.

\begin{lemma}\label{lem:sbs_hyper}
If\/ $\cB = (\E^{k+1},\SS^k,\cH,\bfdelta)$ is a good spherical buffer system and $P \in \cH^* \setminus \{O_\cB\}$, 
then $\cB_P$ is a good spherical buffer system and $\Qbar(\cB_P) = \Qbar(\cB) \cap P$.
\end{lemma}

\begin{proof}
If $P = \E^{k+1}$ then $\cB_P = \cB$, so assume that $P \in \cH_j$ for some $0 \le j \le k - 1$. To show that $\cB_P$ is a spherical buffer system, we need to check that
\begin{itemize}
\item[$(a)$] $\SS^k \cap P$ is a $j$-dimensional sphere, embedded in $P$, with centre $O_\cB$;\smallskip
\item[$(b)$] $\cH_P$ is closed under intersections, and each $H \in \cH_P$ is a proper subspace of $P$;\smallskip
\item[$(c)$] $0 < \delta_i < \delta_{i-1}/3$ for each $0 \le i \le j-1$, and the radius of~$\SS^k \cap P$ is $\delta_{-1}$.
\end{itemize}
Each of these properties follows easily from our assumption that $\cB$ is a spherical buffer system (see Definition~\ref{def:SBS}). 
To show that $\cB_P$ is good, observe that
\begin{equation}\label{eq:BR:buffers}
\cE_{\cB_P}(H \cap P) = \cE_{\cB}(H) \cap P \qquad \text{and} \qquad \cE^+_{\cB_P}(H \cap P)=\cE^+_{\cB}(H) \cap P
\end{equation}
for each $H \in \cH$. Indeed, if $z \in P$ then the Euclidean distance between $z$ and $H \cap P$ is the same in $P$ as it is in $\E^{k+1}$. Note also that if $H \cap P \not\subset H' \cap P$ then $H \not\subset H'$. 
Since $\cB$ is good, and recalling Definition~\ref{def:SBS:good}, it follows that $\cB_P$ is also good. 

It remains to show that $\Qbar(\cB_P) = \Qbar(\cB) \cap P$. 
By Definitions~\ref{def:Qbar} and~\ref{def:BR}, and noting that $\cH_P^* = \{H \cap P : H \in \cH^* \}$, we are required to show that
\begin{equation}\label{eq:using:obs:Qdef2}
\SS^k \cap P \cap \bigcup_{H \in \cH^*} \Big( H \setminus \bigcup_{H' \subsetneq H \cap P} \cE_{\cB_P}(H') \Big) = \SS^k \cap P \cap \bigcup_{H \in \cH^*} \Big( H \setminus \bigcup_{H' \subsetneq H} \cE_{\cB}(H') \Big).
\end{equation}
To see that this holds, recall from~\eqref{eq:BR:buffers} that $\cE_{\cB_P}(H' \cap P) = \cE_{\cB}(H') \cap P$ for every $H' \in \cH$. Since $\cH_P \subset \cH$, as noted above, 
it follows that the right-hand side of~\eqref{eq:using:obs:Qdef2} is contained in the left-hand side. Indeed, if $H' \in \cH_P$ and $H' \subsetneq H \cap P$, then $H' \in \cH$ and $H' \subsetneq H$, so every set $\cE_{\cB_P}(H')$ that is removed from the left is also removed from the right.

Suppose, on the other hand, that $x$ is contained in the left-hand side but not the right-hand side of~\eqref{eq:using:obs:Qdef2}. Let $H \in \cH^*$ be such that $x \in \SS^k \cap H \cap P$, and note that $H \cap P \in \cH$, since $\cH$ is closed under intersections and $P \in \cH$. Thus, since $x$ is not contained in the right-hand side of~\eqref{eq:using:obs:Qdef2}, there must exist $H' \in \cH$ with $H' \subsetneq H \cap P$ such that $x \in \cE_{\cB}(H')$, and hence $x \in \cE_{\cB}(H') \cap P = \cE_{\cB_P}(H')$. Since this holds for every $H \in \cH^*$ containing $x$, this contradicts our assumption that $x$ is contained in the left-hand side of~\eqref{eq:using:obs:Qdef2}, and hence completes the proof of~\eqref{eq:using:obs:Qdef2}. We thus obtain $\Qbar(\cB_P) = \Qbar(\cB) \cap P$, as required. 
\end{proof}

Let us note here the following nice consequence of Lemma~\ref{lem:sbs_hyper}. 

\begin{lemma}\label{lem:Qnonempty}
If\/ $\cB$ is a good spherical buffer system, then $\Qbar(\cB) \ne \emptyset$.
\end{lemma}

\begin{proof}
Set $\cB = (\E^{k+1},\SS^k,\cH,\bfdelta)$, and use induction on the dimension~$k$. Note first that if either $k = 0$ or $\cH^* = \{\E^{k+1}\}$, then $\Qbar(\cB) = \SS^k \ne \emptyset$. We may therefore assume that $k \ge 1$, and that there exists $P \in \cH$. By Lemma~\ref{lem:sbs_hyper}, the spherical buffer system $\cB_P$ satisfies $\Qbar(\cB) \cap P = \Qbar(\cB_P)$ and by the induction hypothesis $\Qbar(\cB_P) \ne \emptyset$.
\end{proof}

To finish this subsection, we make one more simple observation about $\cH_P$. 

\begin{obs}\label{obs:y:not:in:HR}
Let\/ $\cB = (\E^{k+1},\SS^k,\cH,\bfdelta)$ be a spherical buffer system and $y \in \E^{k+1}$. If $P \in \cH^*$ is minimal such that $y \in P$, then $y \not\in H$ for every $H \in \cH_P$. 
\end{obs}

\begin{proof}
Let $H \in \cH_P$, and note that $H \subsetneq P$, by Definition~\ref{def:BR}, and moreover that $H \in \cH$, since $\cH_P \subset \cH$. Since $P \in \cH^*$ was chosen to be minimal such that $y \in P$, it follows that $y \not\in H$, as claimed. 
\end{proof}

\subsection{Induced buffer systems}\label{sec:BxP}

In this subsection we introduce a second method of constructing a buffer system on an affine hyperplane $A$ from a spherical buffer system $\cB$; as noted above, we refer to this construction as the `induced' buffer system on $A$. The hyperplane $A$ that we use will depend on a subspace $P \in \cH$ and a vector $x \in \E^{k+1} \setminus P^\perp$, where $P^\perp$ denotes the subspace perpendicular to $P$ in $\E^{k+1}$. 

Our first task is to define the affine hyperplane $A = A_{x,P}$ that we shall use in the construction; given a set $U \subset \R^d$, let us write $\partial U$ for the boundary of $U$. 

\begin{lemma}\label{lem:Axp:exists}
Let $\cB = (\E^{k+1},\SS^k,\cH,\bfdelta)$ be a spherical buffer system, let $P \in \cH$, and let $x \in \E^{k+1} \setminus P^\perp$. There exists a unique $\lambda > 0$, depending only on $x$ and $P$, such that the affine hyperplane
$$A_{x,P} = \lambda x + P^\perp$$
intersects 
$\partial \cE(P)$. Moreover, $\SS^k \cap A_{x,P} \subset \partial \cE(P)$. 
\end{lemma}

\begin{proof}
Assume, without loss of generality, that $\E^{k+1} = \R^{k+1}$, so in particular $O_\cB = \0$, and that $P = \R^{j+1} \times \{0\}^{k-j}$, i.e., $P$ is the subspace formed by setting the last $k-j$ coordinates equal to zero. Let us write $\SS^j_\delta$ for the $j$-dimensional sphere of radius $\delta$ in $\R^{j+1}$, centred at the origin. Then 
$$\partial \cE(P) = \SS^j_{\rho_j} \times \SS^{k-j-1}_{\delta_j},$$ 
where $\rho_j^2+\delta_j^2=\delta_{-1}^2$. If $x = (x_1,x_2) \in \R^{j+1} \times \R^{k-j}$, then 
\begin{equation}\label{eq:wlog:AxP}
A_{x,P} = \{ \lambda x_1 \} \times \R^{k-j}
\end{equation}
where $\lambda > 0$ is chosen so that $\|\lambda x_1\|=\rho_j$ (note that $x_1\ne 0$, since $x\notin P^\perp$). It follows that $\lambda$ exists, is unique, and only depends on $x$ and $P$. Moreover, we have
\begin{equation}\label{eq:wlog:SxP}
\SS^k \cap A_{x,P} = \{ \lambda x_1\} \times \SS^{k-j-1}_{\delta_j} \subset \partial \cE(P),
\end{equation}
as claimed.
\end{proof}

For each $\cB$, $P$ and $x$ as in Lemma~\ref{lem:Axp:exists}, we write $A_{x,P}$ for the affine hyperplane given by the lemma. We can now define the `induced' buffer systems to which we shall (eventually) apply our induction hypothesis.

\begin{definition}\label{def:BxP}
Let $\cB = (\E^{k+1},\SS^k,\cH,\bfdelta)$ be a spherical buffer system, let $P \in \cH$, and let $x \in \E^{k+1} \setminus P^\perp$. We define the `induced' buffer system $\cB_{x,P} := \big( A_{x,P}, \SS_{x,P}, \cH_{x,P}, \bfdelta' \big)$ by setting
\[
\SS_{x,P} := \SS^k \cap A_{x,P} \qquad \text{and} \qquad \cH_{x,P} := \big\{ H \cap A_{x,P} : P \subset H \in \cH \big\},
\]
and $\bfdelta' = (\delta_j,\ldots,\delta_{k-1})$, where $P \in \cH_j$. 
 \end{definition}  

In order to apply the induction hypothesis to $\cB_{x,P}$, we first need to show that it is a good spherical buffer system. This again follows easily from the definitions. 

\begin{lemma}\label{lem:sbs_buffer}
Let $\cB = (\E^{k+1},\SS^k,\cH,\bfdelta)$ be a good spherical buffer system, let $P \in \cH$, and let $x \in \E^{k+1} \setminus P^\perp$. Then $\cB_{x,P}$ is a good spherical buffer system.
\end{lemma}

\begin{proof}
Note that $P \ne \{O_\cB\}$, since $\E^{k+1} \setminus P^\perp \ne \emptyset$, so we may assume that $P \in \cH_j$ for some $0 \le j \le k - 1$. To show that $\cB_{x,P}$ is a spherical buffer system, we need to check that
\begin{itemize}
\item[$(a)$] $\SS_{x,P}$ is a $(k-j-1)$-dimensional sphere, embedded in $A_{x,P}$;\smallskip
\item[$(b)$] $\cH_{x,P}$ is closed under intersections, $A_{x,P} \not\in \cH_{x,P}$, and each $H \in \cH_{x,P}$ contains the centre of $\SS_{x,P}$;\smallskip
\item[$(c)$] $0 < \delta_i < \delta_{i-1}/3$ for each $j + 1 \le i \le k - 1$, and the radius of~$\SS_{x,P}$ is $\delta_j$.
\end{itemize}
Each of these properties follows easily from our assumption that $\cB$ is a spherical buffer system (see Definition~\ref{def:SBS}) using~\eqref{eq:wlog:SxP} and Definition~\ref{def:BxP}. 

To prove that $\cB_{x,P}$ is good, we shall show that 
\begin{equation}\label{eq:BxP:buffers}
\cE_{\cB_{x,P}}(H \cap A_{x,P}) = \cE_{\cB}(H) \cap A_{x,P}
\end{equation}
for each $P \subset H \in \cH$, and similarly for $\cE^+$. Indeed, if $P \subset H$, then $H \cap A_{x,P}$ is just the orthogonal projection of $H$ onto~$A_{x,P}$, and therefore if $z \in A_{x,P}$, then the Euclidean distance between $z$ and $H$ is the same in $A_{x,P}$ as it is in $\E^{k+1}$. Note also that if $H \cap A_{x,P} \not\subset H' \cap A_{x,P}$ then $H \not\subset H'$. Since $\cB$ is good, and recalling Definition~\ref{def:SBS:good}, it follows immediately that $\cB_{x,P}$ is also good. 
\end{proof}

Our next lemma determines the structure of $\Qbar(\cB_{x,P})$. 

\begin{lemma}\label{lem:sbs:Q}
Let $\cB = (\E^{k+1},\SS^k,\cH,\bfdelta)$ be a good spherical buffer system, let $P \in \cH$, and let $x \in \E^{k+1} \setminus P^\perp$. Then either 
\begin{itemize}
\item[$(a)$] $\Qbar(\cB_{x,P}) = \Qbar(\cB) \cap \SS_{x,P}$\/ and\/ $\SS_{x,P} \cap \cE(H) = \emptyset$\/ for every $P \not\subset H \in \cH$, or\smallskip
\item[$(b)$] $\Qbar(\cB) \cap \SS_{x,P} = \emptyset$\/ and\/ $\SS_{x,P} \subset \cE(H)$\/ for some $H \in \cH$ with $H \subsetneq P$. 
\end{itemize}
\end{lemma}

In the proof of Lemma~\ref{lem:sbs:Q} we will use the following two simple observations. 

\begin{obs}\label{obs:SxP:same:distance}
Let $\cB = (\E^{k+1},\SS^k,\cH,\bfdelta)$ be a spherical buffer system, let $H,P \in \cH$ with $H \subset P$, and let $x \in \E^{k+1} \setminus P^\perp$. If\/ $\SS_{x,P} \cap \cE(H) \ne \emptyset$, then $\SS_{x,P} \subset \cE(H)$.
\end{obs}

\begin{proof}
By~\eqref{eq:wlog:SxP}, the distance of each point of $\SS_{x,P}$ to the set $H$ is the same (it only depends on the distance from $\lambda x_1$ to $H$), and the observation follows immediately. 
\end{proof}

\begin{obs}\label{obs:SxP:sumsets}
Let $\cB = (\E^{k+1},\SS^k,\cH,\bfdelta)$ be a spherical buffer system, let $H,H',P \in \cH$ with $P \subset H' \subsetneq H$, and let $x \in \E^{k+1} \setminus P^\perp$. Then 
$$H' \cap A_{x,P} \subsetneq H \cap A_{x,P}.$$
\end{obs}

\begin{proof}
Set $\tilde{H}' := H' \cap A_{x,P}$ and $\tilde{H} := H \cap A_{x,P}$. Since $P \subset H'$, we have $\tilde{H}' + P \subset H' + P = H'$. Moreover, if 
$y \in H'$, then, writing $y = (y_1,y_2) \in P \times P^\perp$, we have 
$$y = (z,y_2) + (y_1 - z,0) \in \tilde H' + P,$$
where $z$ is the centre of $\SS_{x,P}$. Hence $\tilde H' + P = H'$, and similarly $\tilde H + P = H$. It follows that if $\tilde{H}' = \tilde{H}$, then $H = H'$, contradicting our assumption.
\end{proof}

\begin{proof}[Proof of Lemma~\ref{lem:sbs:Q}]
Suppose first that there exists $P \not\subset H \in \cH$ with $\SS_{x,P} \cap \cE(H) \ne \emptyset$. We claim that 
$$\SS_{x,P} \cap \cE(H \cap P) \ne \emptyset.$$
To prove this, suppose first that $H \not\subset P$. Since $\cB$ is good and $\SS_{x,P} \subset \partial \cE(P)$, by Lemma~\ref{lem:Axp:exists}, it follows that
$$\cE(H \cap P) \supset \cE^+(H) \cap \cE^+(P) \supset \cE(H) \cap \SS_{x,P} \ne \emptyset,$$
as claimed. On the other hand, if $H \subset P$, then $\cE(H \cap P) = \cE(H)$, and so in either case we have $\SS_{x,P} \cap \cE(H \cap P) \ne \emptyset$, as claimed. 

Noting that $H \cap P \in \cH$ and $H \cap P \subsetneq P$, it follows that $\SS_{x,P} \subset \cE(H \cap P)$, by Observation~\ref{obs:SxP:same:distance}, and hence that 
$$\Qbar(\cB) \cap \SS_{x,P} \subset H \cap P,$$
by Observation~\ref{obs:Qdef3}. But $\SS_{x,P} \cap P = \emptyset$, since $\SS_{x,P} \subset \partial \cE(P)$ by Lemma~\ref{lem:Axp:exists}. It therefore follows that in this case we have $\Qbar(\cB) \cap \SS_{x,P} = \emptyset$, as required. 

We may therefore assume that $\SS_{x,P} \cap \cE_\cB(H) = \emptyset$ for every $P \not\subset H \in \cH$, and our task is to show that $\Qbar(\cB_{x,P}) = \Qbar(\cB) \cap \SS_{x,P}$. Let $y \in \Qbar(\cB) \cap \SS_{x,P}$, and observe that, by Definition~\ref{def:Qbar}, there exists $H \in \cH^*$ such that
\begin{equation}\label{eq:sbs:y:property}
y \in \SS_{x,P} \cap H \setminus \bigcup_{H' \subsetneq H} \cE_\cB(H').
\end{equation}
Note also that $P \subset H$, since either $H = \E^{k+1}$, or
$$y \in \SS_{x,P} \cap H \subset \SS_{x,P} \cap \cE_\cB(H),$$
which would contradict our assumption if $P \not\subset H$. Set $\tilde H := H \cap A_{x,P}$, and observe that $\tilde H \in \cH^*_{x,P}$, by Definition~\ref{def:BxP}. We claim that 
\begin{equation}\label{eq:sbs:Q:first:aim}
y \in \SS_{x,P} \cap \tilde H \setminus \bigcup_{H' \subsetneq \tilde H} \cE_{\cB_{x,P}} (H'),
\end{equation}
which, by Definition~\ref{def:Qbar}, would suffice to prove that $y \in \Qbar(\cB_{x,P})$. 

In order to prove~\eqref{eq:sbs:Q:first:aim}, let $H' \in \cH_{x,P}$ with $H' \subsetneq \tilde H$, and recall from Definition~\ref{def:BxP} that there exists $P \subset Q \in \cH$ such that $H' = Q \cap A_{x,P}$. Since $P \subset H \in \cH^*$ and $\cH$ is closed under intersections, it follows that $P \subset H \cap Q \in \cH$, and hence, by~\eqref{eq:BxP:buffers}, that
\begin{equation}\label{eq:BxP:buffers:in:sbs:Q}
\cE_{\cB_{x,P}}(H \cap Q \cap A_{x,P}) = \cE_\cB(H \cap Q) \cap A_{x,P}.
\end{equation}
Moreover, note that $H \cap Q \subsetneq H$, since if $H \subset Q$ then $\tilde H = H \cap A_{x,P} \subset Q \cap A_{x,P} = H'$, contradicting our assumption that $H' \subsetneq \tilde H$. By~\eqref{eq:sbs:y:property}, it follows that 
$$y \not\in \cE_\cB(H \cap Q) \cap A_{x,P} = \cE_{\cB_{x,P}}(H'),$$
by~\eqref{eq:BxP:buffers:in:sbs:Q} and since $H' \subset H$, which implies that $H' = H \cap H' = H \cap Q \cap A_{x,P}$. 

We have proved that $y \not\in \cE_{\cB_{x,P}}(H')$ for every $H' \subsetneq \tilde H$, and hence to deduce~\eqref{eq:sbs:Q:first:aim} it only remains to observe that $y \in \SS_{x,P} \cap H = \SS_{x,P} \cap \tilde H$, since $\tilde H := H \cap A_{x,P}$ and $\SS_{x,P} = \SS^k \cap A_{x,P}$. Since $y$ was arbitrary, we have proved that $\Qbar(\cB)\cap \SS_{x,P} \subset \Qbar(\cB_{x,P})$.

To prove the reverse inclusion, observe first that, by Definition~\ref{def:Qbar}, if $y \in \Qbar(\cB_{x,P})$, then 
\begin{equation}\label{eq:sbs:y:property:reverse}
y \in \SS_{x,P} \cap \tilde H \setminus \bigcup_{H' \subsetneq \tilde H} \cE_{\cB_{x,P}}(H')
\end{equation}
for some $\tilde H \in \cH_{x,P}^*$, and let $P \subset H \in \cH^*$ be such that $\tilde H = H \cap A_{x,P}$. We claim that 
\begin{equation}\label{eq:sbs:Q:second:aim}
y \in \SS_{x,P} \cap H \setminus \bigcup_{H' \subsetneq H} \cE_\cB(H'),
\end{equation}
which, by Definition~\ref{def:Qbar}, will suffice to prove that $y \in \Qbar(\cB) \cap \SS_{x,P}$. 

To prove~\eqref{eq:sbs:Q:second:aim}, note that $y \in \SS_{x,P} \cap H$, and suppose that $y \in \SS_{x,P} \cap \cE_\cB(H')$ for some $H' \in \cH$ with $H' \subsetneq H$. Since $\SS_{x,P} \cap \cE_\cB(H') \ne \emptyset$, it follows (by our assumption) that $P \subset H'$, and hence, by~\eqref{eq:BxP:buffers}, that
$$y \in \cE_\cB(H') \cap A_{x,P} = \cE_{\cB_{x,P}}(H' \cap A_{x,P}).$$
Moreover, we have $H' \cap A_{x,P} \subsetneq \tilde H$, by Observation~\ref{obs:SxP:sumsets}.  This contradicts~\eqref{eq:sbs:y:property:reverse}, and hence completes the proof that $\Qbar(\cB_{x,P}) = \Qbar(\cB) \cap \SS_{x,P}$.
\end{proof}

\subsection{Some simple properties of $\SS_{x,P}$}\label{sec:SxP:props}

We next collate a number of properties of the spheres $\SS_{x,P}$, most of which are fairly straightforward consequences of the definitions. 

Recall that if $U \subset \E^{k+1}$ is compact then $M(U,x)$ denotes the set of points of $U$ at minimal (Euclidean) distance to $x \in \E^{k+1}$. Our first property will be used in the very last step of the proof of Lemma~\ref{thm:quasi:key}, in order to deduce a contradiction from our application of the induction hypothesis to $\cB_{x,P}$. 

\begin{lemma}\label{lem:SxP:closestpoints}
Let $\cB = (\E^{k+1},\SS^k,\cH,\bfdelta)$ be a spherical buffer system, let $P \in \cH$, and let $x \in \E^{k+1} \setminus P^\perp$. If\/ $M(\SS_{x,P},x) = \{y\}$, then 
$$M\big( \Qbar(\cB),x \big) \cap \SS_{x,P} \subset M\big( \Qbar(\cB) \cap \SS_{x,P}, y \big).$$
\end{lemma}

\begin{proof}
By Observation~\ref{obs:spheres:increasing:distance}, and since $M(\SS_{x,P},x) = \{y\}$, it follows that the function $z \mapsto d(x,z)$ is increasing in $d(y,z)$ on $\SS_{x,P}$. This means that if $z,z' \in \Qbar(\cB) \cap \SS_{x,P}$ are such that $d(y,z) < d(y,z')$, then we have $d(x,z) < d(x,z')$, as claimed. 
\end{proof}

Let us take the opportunity to prove another similar consequence of Observation~\ref{obs:spheres:increasing:distance}. 

\begin{lemma}\label{obs:y:closest:to:xhat}
Let $\cB = (\E^{k+1},\SS^k,\cH,\bfdelta)$ be a spherical buffer system, let $x \in \SS^k$, and let $Q \in \cH^*$. If\/ $M(\SS^k \cap Q,x) = \{ \hat x \}$ and $y \in M(\Qbar(\cB) \cap Q,x)$, then $y \in M( \Qbar(\cB) \cap Q,\hat x )$.  
\end{lemma}

\begin{proof}
We may assume that $Q \neq \E^{k+1}$ since otherwise $x = \hat x$. By Observation~\ref{obs:spheres:increasing:distance}, and since $M(\SS^k \cap Q,x) = \{\hat x\}$, it follows that $d(x,z)$ is an increasing function of $d(\hat x,z)$ for all $z \in \SS^k \cap Q$. Hence, if there exists $z \in \Qbar(\cB) \cap Q$ such that $d(\hat x,z) < d(\hat x,y)$, then $d(x,z) < d(x,y)$, contradicting our assumption that $y \in M(\Qbar(\cB) \cap Q,x)$. 
\end{proof}

In order to apply Lemma~\ref{lem:SxP:closestpoints}, we will need the following lemma.

\begin{lemma}\label{lem:SxP:in:H}
Let $\cB = (\E^{k+1},\SS^k,\cH,\bfdelta)$ be a spherical buffer system, let $H,P \in \cH$ with $P \subset H$, and let $x \in H \setminus (P \cup P^\perp)$. Then $M(\SS_{x,P},x) = \{x'\}$ for some $x' \in H$.
\end{lemma}

\begin{proof}
Let $P \in \cH_j$, write $x = (x_1, x_2) \in P \times P^\perp$ and note that $x_1, x_2 \ne 0$. Recall from~\eqref{eq:wlog:SxP} that each element of $\SS_{x,P}$ is of the form $(\lambda x_1,y) \in P \times P^\perp$, where $y \in \SS_{\delta_j}^{k-j-1}$, and observe that
$$M\big( \SS_{\delta_j}^{k-j-1}, x_2 \big) = \big\{ \mu x_2 \big\}$$
for some $\mu > 0$. It follows that if $x' \in M( \SS_{x,P}, x )$, then $x' = (\lambda x_1,\mu x_2)$. But $x \in H$ and $(x_1,0) \in P \subset H$, so $(0,x_2) \in H$ and thus $(\lambda x_1, \mu x_2) \in H$, as claimed. 
\end{proof}

We next prove two simple lemmas which relate $\SS_{x,P}$ to the corresponding sphere for the restricted buffer system $\cB_Q$, and to the sphere $\SS_{y,P}$ when $y \in M(\SS^k \cap Q,x)$. 

\begin{lemma}\label{lem:SxP:BR}
Let $\cB = (\E^{k+1},\SS^k,\cH,\bfdelta)$ be a spherical buffer system, let $Q \in \cH^*$ and $P \in \cH_Q$, and let $x \in Q \setminus P^\perp$. Then 
$$\SS_{x,P}(\cB_Q) = \SS_{x,P}(\cB) \cap Q.$$
\end{lemma}

\begin{proof}
We may assume that $Q \neq \E^{k+1}$, since otherwise the assertion is trivial. Note that, by Definitions~\ref{def:BR} and~\ref{def:BxP}, we have $\SS_{x,P}(\cB_Q) = \SS^k \cap Q \cap A_{x,P}(\cB_Q)$, and also $P \subset Q$. We claim that
$$A_{x,P}(\cB_Q) = A_{x,P}(\cB) \cap Q.$$
To see this, apply Lemma~\ref{lem:Axp:exists} to obtain 
$$A_{x,P}(\cB) = z + P^\perp \qquad \text{and} \qquad A_{x,P}(\cB_Q) = z + \big( P^\perp \cap Q \big)$$
for some $z \in P$, where the fact that we may use $z$ for both sets follows because the value of $\lambda$ (in Lemma~\ref{lem:Axp:exists}) depends only on $x$ and $P$. We thus obtain
$$\SS_{x,P}(\cB_Q) = \SS^k \cap A_{x,P}(\cB) \cap Q  = \SS_{x,P}(\cB) \cap Q,$$
as claimed.
\end{proof}

\begin{lemma}\label{lem:SxP:SyP}
Let $\cB = (\E^{k+1},\SS^k,\cH,\bfdelta)$ be a spherical buffer system, let $P,Q \in \cH$ with $P \subset Q$, and let $x \in \SS^k \setminus P^\perp$. If $y \in M(\SS^k \cap Q,x)$, then $\SS_{x,P} = \SS_{y,P}$.
\end{lemma}

\begin{proof}
Observe that $y = \pi(x,Q)$, and therefore, by~\eqref{def:projection}, we have $y - \mu x \in Q^\perp \subset P^\perp$ for some $\mu > 0$. Noting that $x,y \not\in P^\perp$, let $\lambda, \lambda' > 0$ be such that $A_{x,P} = \lambda x + P^\perp$ and $A_{y,P} = \lambda' y + P^\perp = \lambda' \mu x + P^\perp$, and recall from Lemma~\ref{lem:Axp:exists} that $\lambda$ is unique such that $\lambda x + P^\perp$ intersects $\partial \cE(P)$. It follows that $\lambda' \mu = \lambda$, and $A_{x,P} = A_{y,P}$, as claimed.
\end{proof}

The remaining results in this subsection all rely on the following lemma.

\begin{lemma}\label{lem:xiner}
Let\/ $\cB = (\E^{k+1},\SS^k,\cH,\bfdelta)$ be a spherical buffer system, let\/ $P \in \cH$, and let\/ $x \in \SS^k \setminus P^\perp$ and $y \in \SS_{x,P}$. 
\begin{itemize}
\item[$(a)$] If\/ $H \in \cH$ with $H \subset P$, then either $x \not\in H^\perp$ and
$$x \in \cE(P) \qquad \Leftrightarrow \qquad d(x,H) < d(y,H),$$
or $x \in H^\perp$ and\/ $d(x,H) = d(y,H)$.\smallskip
\item[$(b)$] If\/ $x \in \cE(P)$ and $P \in \cH_j$, then $d(x,y) \le 2\delta_j \le 2\delta_0$.
\end{itemize}
\end{lemma}

\begin{proof}
As usual, we may assume without loss of generality that $\E^{k+1} = \R^{k+1}$ and that $P = \R^{j+1} \times \{0\}^{k-j}$ for some $0 \le j \le k - 1$, since $P \in \cH$ and $\SS^k \setminus P^\perp$ is non-empty. Let $x = (x_1,x_2) \in P \times P^\perp$, so $\SS_{x,P} = \{ \lambda x_1\} \times \SS^{k-j-1}_{\delta_j}$, by~\eqref{eq:wlog:SxP}, where $\lambda > 0$ is chosen so that $\|\lambda x_1\|^2 = \rho_j^2 = \delta_{-1}^2 - \delta_j^2$. Since $x \in \SS^k$, it follows by Definition~\ref{def:buffers} that
\begin{equation}\label{eq:lambda:less:than1}
x \in \cE(P) \qquad \Leftrightarrow \qquad \|x_2\| < \delta_j \quad \text{and} \quad \|x_1\| > \rho_j \qquad \Leftrightarrow \qquad \lambda < 1.
\end{equation}

Let $y = (\lambda x_1,y_2) \in P \times P^\perp$, and observe that if $z \in H \subset P$, then $\langle y,z\rangle = \lambda \langle x,z\rangle$. Therefore, if $z$ is either the closest point in $H$ to $x$ or the closest point in $H$ to $y$, then $\langle x,z \rangle \ge 0$, since if $\langle x,z \rangle < 0$ then $-z$ is closer to $x$ than $z$ (and similarly for $y$). Note also that if $x \in H^\perp$ then $\langle x,z \rangle = \langle y,z \rangle = 0$, and if $x \not\in H^\perp$ then $\langle x,z \rangle > 0$. 

Now, observe that 
$$d(y,z)^2 = \|y\|^2 + \|z\|^2 - 2 \langle y,z\rangle = \delta_{-1}^2 + \|z\|^2 - 2\lambda \langle x,z\rangle,$$
since $y \in \SS_{x,P} \subset \SS^k$, and similarly
$$d(x,z)^2 = \|x\|^2 + \|z\|^2 - 2\langle x,z\rangle = \delta_{-1}^2 + \|z\|^2 - 2\langle x,z\rangle.$$
By~\eqref{eq:lambda:less:than1} and the observations about $\langle x,z \rangle$, it follows that $d(x,H) = d(y,H)$ when $x \in H^\perp$, and that the claimed equivalence holds when $x \not\in H^\perp$, as required.

For part $(b)$, observe that if\/ $x \in \cE(P)$ and $P \in \cH_j$, then both $x$ and $y$ lie within distance $\delta_j$ of the point $w := (\lambda x_1,0)$. Indeed, we have 
$$d(x,w)^2 = \|x_1 - \lambda x_1\|^2 + \|x_2\|^2 = (1 - \lambda)^2 \|x_1\|^2 + \delta_{-1}^2 - \|x_1\|^2,$$ 
and since $\delta_{-1}^2 = \lambda^2 \|x_1\|^2 + \delta_j^2$ and $0 < \lambda < 1$, by~\eqref{eq:lambda:less:than1}, it follows that 
$$d(x,w)^2 = 2\lambda(\lambda - 1) \|x_1\|^2 + \delta_j^2 \le \delta_j^2,$$ 
as claimed. For $d(y,w)$ the same bound is immediate, since $w$ is the centre of $\SS_{x,P}$, which is a sphere of radius $\delta_j$, and $y \in \SS_{x,P}$.
\end{proof}

We next give five applications of Lemma~\ref{lem:xiner}; for the statements to make sense, we need a simple observation. 

\begin{obs}\label{obs:x:notin:Pperp}
Let $\cB = (\E^{k+1},\SS^k,\cH,\bfdelta)$ be a spherical buffer system, let $P \in \cH$, and let $x \in \cE(P)$. Then $x \notin P^\perp$. 
\end{obs}

\begin{proof}
This follows because $x \in \cE(P) \subset \SS^k$ implies\footnote{Recall that $\cE(\{O_\cB\}) = \emptyset$, so $P \ne O_\cB$, and that the radius of~$\SS^k$ is $\delta_{-1}$.} that $d(x,P^\perp)^2 \ge \delta_{-1}^2 - \delta_0^2 > 0$.
\end{proof}

Recall also, from Observation~\ref{obs:unique:buffer}, that if $\cB$ is good and $x \in \SS^k \setminus \Qbar(\cB)$, then there exists a unique minimal $P \in \cH$ such that $x \in \cE(P)$. The following important fact is a simple consequence of Lemmas~\ref{lem:sbs:Q} and~\ref{lem:xiner}. 

\begin{lemma}\label{lem:Qbar:nonempty}
Let\/ $\cB = (\E^{k+1},\SS^k,\cH,\bfdelta)$ be a good spherical buffer system, let\/ $x \in \SS^k \setminus \Qbar(\cB)$, and let\/ $P \in \cH$ be minimal such that\/ $x \in \cE(P)$. Then $\Qbar(\cB) \cap \SS_{x,P} \ne \emptyset$.
\end{lemma}

\begin{proof}
Observe first that $x \notin P^\perp$, by Observation~\ref{obs:x:notin:Pperp}. By Lemma~\ref{lem:sbs:Q}, it follows that if $\Qbar(\cB) \cap \SS_{x,P} = \emptyset$, then there exists $H \in \cH$ such that $H \subsetneq P$ and $\SS_{x,P} \subset \cE(H)$. Now, since $x \in \cE(P)$, by Lemma~\ref{lem:xiner} we have $d(x,H) \le d(y,H)$ for every $y \in \SS_{x,P} \subset \cE(H)$, and hence $x \in \cE(H)$, contradicting the minimality of~$P$.
\end{proof}

We can now use Lemmas~\ref{lem:xiner} and~\ref{lem:Qbar:nonempty} to deduce the following useful fact. 

\begin{lemma}\label{lem:minbuffer}
Let\/ $\cB = (\E^{k+1},\SS^k,\cH,\bfdelta)$ be a good spherical buffer system, let\/ $x \in \SS^k \setminus \Qbar(\cB)$, and let\/ $P \in \cH$ be minimal such that\/ $x \in \cE(P)$. Then 
\begin{equation}\label{eq:closest:in:bigbuffer}
M\big( \Qbar(\cB),x \big) \subset \cE^+(P).
\end{equation}
\end{lemma}

\begin{proof}
By Lemma~\ref{lem:Qbar:nonempty}, we have $\Qbar(\cB) \cap \SS_{x,P} \ne \emptyset$, and therefore, by Lemma~\ref{lem:xiner}, if $P \in \cH_j$ then
$$d\big( x, \Qbar(\cB) \big) \le d\big( x, \Qbar(\cB) \cap \SS_{x,P} \big) \le  2\delta_j.$$
It follows that if $y \in M( \Qbar(\cB),x )$, then 
$$d(y,P) \le d(y,x) + d(x,P) < 2\delta_j + \delta_j,$$ 
and hence $y \in \cE^+(P)$, by Definition~\ref{def:buffers}, as claimed.
\end{proof}

We will use the following two related facts in Sections~\ref{sec:pushing:y} and~\ref{quasi:technical:proof:sec}. 

\begin{lemma}\label{lem:nearest:not:perp}
Let\/ $\cB = (\E^{k+1},\SS^k,\cH,\bfdelta)$ be a good spherical buffer system, and let $P \in \cH$. If $x \in \SS^k \setminus \Qbar(\cB)$ and\/ $y \in M( \Qbar(\cB), x ) \cap \partial\cE(P)$, then $d(x,y) \leq 2\delta_0$ and $x \not\in P^\perp$. 
\end{lemma}

\begin{proof}
By Lemma~\ref{obs:unique:buffer}, there exists a unique minimal $Q \in \cH$ such that $x \in \cE(Q)$. Observe that $\Qbar(\cB) \cap \SS_{x,Q} \neq \emptyset$, by Lemma~\ref{lem:Qbar:nonempty}, and hence if $Q \in \cH_j$, then 
$$d(x,y) = d( x, \Qbar(\cB) ) \le d( x, \Qbar(\cB) \cap \SS_{x,Q} ) \le 2\delta_j \le 2\delta_0,$$
by Lemma~\ref{lem:xiner}. Now if $x \in P^\perp$ then, since $y \in \partial\cE(P)$, we have $d(x,y)^2 \ge \delta_{-1}^2 - \delta_0^2$. Since $\delta_{-1} > 3 \delta_0$, this is a contradiction. 
\end{proof}

\begin{lemma}\label{lem:another:not:perp}
Let\/ $\cB = (\E^{k+1},\SS^k,\cH,\bfdelta)$ be a spherical buffer system, and let $P \in \cH$, $Q \in \cH^*$, $x \in \cE(P) \setminus \Qbar(\cB)$ and $y \in Q \cap M(\Qbar(\cB),x)$. If\/ $\Qbar(\cB) \cap \SS_{x,P} \ne \emptyset$, then $x \not\in Q^\perp$. 
\end{lemma}

\begin{proof}
If $Q = \E^{k+1}$ then the assertion is trivial (because $x \in \SS^k$), so we may assume that $Q \in \cH$. Observe that $x \notin P^\perp$, by Observation~\ref{obs:x:notin:Pperp}, and therefore $d(x,w) \le 2\delta_0$ for every $w \in \SS_{x,P}$, by Lemma~\ref{lem:xiner}. It follows that 
$$d(x,Q) \le d(x,y) \le \max_{w \in \SS_{x,P}} d(x,w) \le 2\delta_0 < \delta_{-1}.$$ 
But if $x \in \SS^k \cap Q^\perp$ then $d(x,Q) = \delta_{-1}$, so $x \notin Q^\perp$, as claimed.
\end{proof}

The following observation will be used in the proof of the next (and final) lemma of this subsection. We omit the straightforward proof.

\begin{obs}\label{obs:plus-intersect-assoc}
Let $P,Q$ be subspaces of\/ $\E^{k+1}$ and let $x \in Q$. Then
\[
(x + P) \cap Q = x + (P \cap Q).
\]
\end{obs}

\medskip
\pagebreak

We need one final property of $\SS_{x,P}$.

\begin{lemma}\label{lem:scover}
Let\/ $\cB$ be a good spherical buffer system, and let $H,P \in \cH$ with $P \not\subset H$. If\/ $x \in H \cap \cE(P)$, then $\SS_{x,P} \subset \cE(H \cap P)$.
\end{lemma}

\begin{proof}
Note that $x \not\in P^\perp$, by Observation~\ref{obs:x:notin:Pperp}. Suppose first that $\dim(H) \le \dim(P)$. Since $x \in H \cap \cE(P)$ and the $\delta_j$ are decreasing, it follows that if $P \in \cH_j$, then 
\begin{equation}\label{eq:scover:using:xiner}
\SS_{x,P} \subset \big\{ z \in \SS^k : d(x,z) \le 2\delta_j \big\} \subset \cE^+(H),
\end{equation}
by Lemma~\ref{lem:xiner}. Since $x \in \cE(P)$, it also follows from~\eqref{eq:scover:using:xiner} that $\SS_{x,P} \subset \cE^+(P)$. If $H \not\subset P$ then, since $\cB$ is good, it now follows that
$$\SS_{x,P} \subset \cE^+(H) \cap \cE^+(P) \subset \cE(H \cap P),$$ 
as required. On the other hand, if $H \subset P$, then we have $\dim(H) < \dim(P)$ (since $P \not\subset H$), and therefore it follows from the first inclusion in~\eqref{eq:scover:using:xiner} that $\SS_{x,P} \subset \cE(H)$.

We may therefore assume that $\dim(H) > \dim(P)$, and hence that $\dim(P^\perp \cap H) > 0$, and moreover that $H \not\subset P$. Recalling that $H \subset \cE^+(H)$, and that $\SS_{x,P} \subset \cE^+(P)$, by Lemma~\ref{lem:Axp:exists}, it follows that
\begin{equation}\label{eq:scover:1}
\SS_{x,P} \cap H \subset \cE^+(P) \cap \cE^+(H) \subset \cE(H \cap P),
\end{equation}
since $\cB$ is good. We claim that moreover 
\begin{equation}\label{eq:scover:2}
\SS_{x,P} \cap H = \SS^k \cap A_{x,P} \cap H \ne \emptyset.
\end{equation}
To see this, recall that $A_{x,P} = \lambda x + P^\perp$ for some $\lambda > 0$, and therefore, since $x \in H$, we have
$$A_{x,P} \cap H = (\lambda x + P^\perp) \cap H = \lambda x + (P^\perp \cap H),$$ 
by Observation~\ref{obs:plus-intersect-assoc}. Hence $A_{x,P} \cap H$ is an affine space of dimension at least $1$. Since $A_{x,P} \cap H$ contains the centre of $\SS_{x,P}$, 
it also intersects $\SS^k$, and so~\eqref{eq:scover:2} holds. 

Finally, it follows from~\eqref{eq:scover:1} and~\eqref{eq:scover:2} that $\SS_{x,P} \cap \cE(H \cap P) \ne \emptyset$, and therefore, by Observation~\ref{obs:SxP:same:distance}, we have $\SS_{x,P} \subset \cE(H \cap P)$, as required.
\end{proof}

\subsection{The great escape}\label{sec:GreatEscape}

When $x \in P$, all points of $\SS_{x,P}$ are equidistant from~$x$. The next lemma (the `great escape') shows that in this case, if $\Qbar(\cB) \cap \SS_{x,P} \ne \emptyset$, then there are points of $\Qbar(\cB) \cap P$ that are closer to $x$ than any point of $\SS_{x,P} \subset \partial \cE(P)$. 

\begin{lemma}\label{lem:great:escape}
Let\/ $\cB=(\E^{k+1},\SS^k,\cH,\bfdelta)$ be a good spherical buffer system, let $P \in \cH_j$, and let\/ $x \in \SS^k \cap P$. If\/ $\Qbar(\cB) \cap \SS_{x,P} \ne \emptyset$, then
\[
d(x,\Qbar(\cB) \cap P) < \delta_j.
\]
In particular, $M(\Qbar(\cB),x) \subset P$. 
\end{lemma}

We will construct an escape route inductively, applying the following lemma at each step to find a nearby point on the boundary of one of the buffers containing $x$. 

\begin{lemma}\label{lem:great:escape:step}
Let\/ $\cB=(\E^{k+1},\SS^k,\cH,\bfdelta)$ be a good spherical buffer system, let $H \in \cH_i$ and $P \in \cH_j$ with $H \subsetneq P$, and let\/ $x,y \in \SS^k \cap P$ with $d(x,H) \le d(y,H)$. If\/ $\Qbar(\cB) \cap \SS_{x,P} \ne \emptyset$, then there exists $z \in \SS^k \cap P \setminus \cE(H)$ such that
$$d(y, z) \le \frac{2\delta_j^2}{\delta_i},$$
and $d(y,H') \le d(z,H')$ for every $H' \in \cH$ such that either $H' \subset H$ or $H \subset H' \subset P$.
\end{lemma}

\begin{proof}
If $y \not\in \cE(H)$ then the conditions are satisfied with $z = y$, so we may assume that $y \in \cE(H)$, and therefore $0 \le i < j$. The key observation is that 
\begin{equation}\label{eq:SxP:EH:disjoint}
\SS_{x,P} \cap \cE(H) = \emptyset.
\end{equation}
To see this, recall that $\Qbar(\cB) \cap \SS_{x,P} \neq \emptyset$, so we can fix $w \in \Qbar(\cB) \cap \SS_{x,P}$. Note that $w \not\in H$, since $H \subset P$ and $\SS_{x,P} \cap P = \emptyset$. By Observation~\ref{obs:Qdef3}, it follows that $w \not\in \cE(H)$, and therefore $\SS_{x,P} \not\subset \cE(H)$. This then implies~\eqref{eq:SxP:EH:disjoint}, by Observation~\ref{obs:SxP:same:distance}. 

In order to use~\eqref{eq:SxP:EH:disjoint} to construct $z$, let $\lambda > 0$ be such that $A_{x,P} = \lambda x + P^\perp$, and write $x = (x_1,x_2,0) \in H \times (P \cap H^\perp) \times P^\perp$. By~\eqref{eq:SxP:EH:disjoint}, we have
$$\delta_i^2 \le d(H,\SS_{x,P})^2 = \|\lambda x_2\|^2 + \delta_j^2,$$
and since $x \in \SS^k \cap P$ we have $\lambda < 1$, by~\eqref{eq:lambda:less:than1}. It follows that 
\begin{equation}\label{eq:great-escape-y2}
\|y_2\|^2 = d(y,H)^2 \ge d(x,H)^2 = \|x_2\|^2 \ge \delta_i^2 - \delta_j^2,
\end{equation}
where $y = (y_1,y_2,0) \in H \times (P \cap H^\perp) \times P^\perp$, and the second step is by assumption. 

Now, define $z := (\mu_1 y_1,\mu_2 y_2,0) \in H \times (P \cap H^\perp) \times P^\perp$, where $0 < \mu_1 < 1 < \mu_2$ are chosen so that $z \in \SS^k \cap \partial \cE(H)$. This is possible because the conditions are equivalent to $\|\mu_1 y_1\|^2 + \|\mu_2 y_2\|^2 = \|y_1\|^2 + \|y_2\|^2$ and $\|\mu_2 y_2\| = \delta_i$, and since $\|y_2\| < \delta_i$, because we assumed at the start of the proof that $y \in \cE(H)$.  Observe that
$$(1 - \mu_1^2) \|y_1\|^2 = (\mu_2^2 - 1) \| y_2 \|^2 \le \delta_i^2 - (\delta_i^2 - \delta_j^2) = \delta_j^2$$
since $\|y_1\|^2 + \|y_2\|^2 = \|y\|^2 = \|\mu_1 y_1\|^2 + \| \mu_2 y_2\|^2$, and by~\eqref{eq:great-escape-y2}, and hence that 
$$1 - \mu_1 \le \frac{\delta_j^2}{\|y_1\|^2} \qquad \text{and} \qquad \mu_2 - 1 \le \frac{\delta_j^2}{\|y_2\|^2}.$$
Since $\|y_1\|^2 \ge \delta_{-1}^2 - \delta_i^2$ (because $\|y_2\| < \delta_i$) and $\|y_2\|^2 \ge \delta_i^2 - \delta_j^2$, it follows that 
$$d(y,z)^2 = (1 - \mu_1)^2 \|y_1\|^2 + (\mu_2 - 1)^2 \|y_2\|^2 \le \frac{\delta_j^4}{\|y_1\|^2} + \frac{\delta_j^4}{\|y_2\|^2} \le \frac{2\delta_j^4}{\delta_i^2},$$
as claimed, since $0 \le i < j$, so $\delta_{-1} > 3\delta_i > 9\delta_j > 0$. 

For the final part of the lemma, observe first that if $H \subset H' \subset P$, then
$$d(z,H') = \mu_2 d(y,H') \ge d(y,H'),$$
since $\mu_2 > 1$. On the other hand, if $H' \subset H$, then note that
$$d(z,H')^2 = d(\mu_1 y_1,H')^2 + \|\mu_2 y_2 \|^2 = \mu_1^2 d(y_1,H')^2 + \mu_2^2\| y_2 \|^2,$$
and therefore 
$$d(z,H')^2 - d(y,H')^2 = (\mu_1^2 - 1) d(y_1,H')^2 + (\mu_2^2 - 1) \| y_2 \|^2.$$
Since $(1 - \mu_1^2) \|y_1\|^2 = (\mu_2^2 - 1) \| y_2 \|^2$, it follows that
$$d(z,H')^2 - d(y,H')^2 = (\mu_1^2 - 1) \big( d(y_1,H')^2 - \|y_1\|^2 \big) \geq 0,$$
completing the proof.
\end{proof}

To deduce Lemma~\ref{lem:great:escape}, we apply Lemma~\ref{lem:great:escape:step} once for each buffer containing~$x$. 

\begin{proof}[Proof of Lemma~\ref{lem:great:escape}]
First, note that if $x \in \Qbar(\cB)$ then $d(x,\Qbar(\cB) \cap P) = 0$, so we may assume that $x \not\in \Qbar(\cB)$. By Observation~\ref{obs:unique:buffer}, it follows that there exists a unique minimal $H_0 \in \cH$ such that $x \in \cE(H_0)$. 

We claim that $H_0 \subsetneq P$. Indeed, $H_0 \subset P$ follows from the minimality of $H_0$, since $x \in P \subset \cE(P)$, and if $H_0 = P$ then by Definition~\ref{def:Qbar} we would have $x \in \Qbar(\cB)$, since $x \in \SS^k \cap P \setminus \bigcup_{P' \subsetneq P} \cE(P')$. Let $0 \le i_0 < j$ be such that $H_0 \in \cH_{i_0}$. 

By Lemma~\ref{lem:great:escape:step} (applied with $y = x$), there exists $z_1 \in \SS^k \cap P \setminus \cE(H_0)$, with 
$$d(x,z_1) \le 2\delta_j^2 / \delta_{i_0},$$ 
such that $d(x,H') \le d(z_1,H')$ for every $H' \in \cH$ with $H' \subset H_0$ or $H_0 \subset H' \subset P$. 

We now iterate the above argument until we find a $z \in \Qbar(\cB) \cap P$ with $d(x,z) < \delta_j$. To be precise, suppose we have found, for some $\ell \ge 1$, sequences 
$$H_0 \subsetneq H_1 \subsetneq \cdots \subsetneq H_{\ell-1} \subsetneq P \qquad \text{and} \qquad z_1,\ldots,z_\ell \in \SS^k \cap P$$ 
such that, setting $z_0 := x$, the following hold for each $0 \le t \le \ell - 1$:
\begin{itemize}
\item[$(a)$] $H_t \in \cH$ is minimal such that $z_t \in \cE(H_t)$, and $z_{t+1} \not\in \cE(H_t)$;\smallskip
\item[$(b)$] $d(z_t,z_{t+1}) \le 2\delta_j^2 / \delta_{i_t}$, where $H_t \in \cH_{i_t}$;\smallskip \item[$(c)$] $d(z_t,H') \le d(z_{t+1},H')$ for every $H' \in \cH$ with $H' \subset H_t$ or $H_t \subset H' \subset P$. 
\end{itemize}
If $z_\ell \in \Qbar(\cB)$ then set $z := z_\ell$, and observe that 
\begin{equation}\label{eq:escape:sum}
d(x,z) \le \sum_{i = i_0}^{j-1} \frac{2\delta_j^2}{\delta_i} < \delta_j,
\end{equation}
as required, so assume that $z_\ell \not\in \Qbar(\cB)$ and (by Observation~\ref{obs:unique:buffer}) let $H_\ell \in \cH$ be minimal such that $z_\ell \in \cE(H_\ell)$. We plan to apply Lemma~\ref{lem:great:escape:step} with $y = z_\ell$ and $H = H_\ell$, so we need to check that the conditions of the lemma hold. 

\begin{claim}\label{clm:escape}
$H_{\ell-1} \subsetneq H_\ell \subsetneq P$ and $d(x,H_\ell) \le d(z_\ell,H_\ell)$. 
\end{claim}

\begin{clmproof}{clm:escape}
Note first that $H_\ell \subsetneq P$, since $z_\ell \in P \subset \cE(P)$ and $z_\ell \not\in \Qbar(\cB)$. Indeed, $H_\ell \in \cH$ was chosen be minimal such that $z_\ell \in \cE(H_\ell)$, and if $H_\ell = P$ then $z_\ell \in \SS^k \cap P \setminus \bigcup_{P' \subsetneq P} \cE(P') \subset \Qbar(\cB)$. We also have $H_\ell \ne H_{\ell-1}$, since $z_\ell \in \cE(H_\ell) \setminus \cE(H_{\ell-1})$. 

Next, observe that if $H_\ell \subset H_{\ell-1}$, then $d(z_{\ell-1},H_\ell) \le d(z_\ell,H_\ell)$, by property~$(c)$ applied with $H' = H_\ell$. Since $z_\ell \in \cE(H_\ell)$, it follows that $z_{\ell-1} \in \cE(H_\ell)$, contradicting the minimality of $H_{\ell-1}$. It follows that $H_\ell \not\subset H_{\ell-1}$. 

Now, recall that $\cB$ is good, and observe that $z_\ell \in \cE^+(H_{\ell-1})$, since $z_{\ell-1} \in \cE(H_{\ell-1})$ and $d(z_{\ell-1}, z_\ell) < \delta_j < \delta_{i_{\ell-1}}$, by property~$(b)$. It follows that if $H_{\ell-1} \not\subset H_\ell$, then
$$z_\ell \in \cE^+(H_{\ell-1}) \cap \cE(H_\ell) \subset \cE(H_{\ell-1} \cap H_\ell),$$
contradicting the minimality of $H_\ell$. Hence $H_{\ell-1} \subsetneq H_\ell$, as claimed. 

Finally, observe that, since $H_0 \subsetneq H_1 \subsetneq \cdots \subsetneq H_{\ell} \subsetneq P$, we have
$$d(x,H_\ell) \le d(z_1,H_\ell) \le \cdots \le d(z_\ell,H_\ell),$$
by property~$(c)$ with $H' = H_\ell$. 
\end{clmproof}

By Claim~\ref{clm:escape}, the conditions of Lemma~\ref{lem:great:escape:step} are satisfied with $y = z_\ell$ and $H = H_\ell$, and therefore there exists $z_{\ell+1} \in \SS^k \cap P \setminus \cE(H_\ell)$, with 
$$d(z_\ell,z_{\ell+1}) \le 2\delta_j^2 / \delta_{i_\ell},$$ 
such that $d(z_\ell,H') \le d(z_{\ell+1},H')$ for every $H' \in \cH$ with $H' \subset H_\ell$ or $H_\ell \subset H' \subset P$. This completes the inductive step, and since the dimension of the subspaces $H_t$ is strictly increasing in $t$, we must eventually find a $z \in \Qbar(\cB) \cap P$ with 
$$d(x,\Qbar(\cB) \cap P) \le d(x,z) < \delta_j,$$ 
as required. Finally, to deduce that $M(\Qbar(\cB),x) \subset P$, simply recall from Observation~\ref{obs:Qdef3} that all points of $\Qbar(\cB) \setminus P$ lie outside $\cE(P)$. 
\end{proof}

\subsection{Pushing $y$ towards $x$}\label{sec:pushing:y}

In this subsection we prove the following key lemma, which we shall apply to the restricted buffer system $\cB_Q$, where $Q \in \cH^*$ is minimal such that $y \in Q$. Recall from Observation~\ref{obs:y:not:in:HR} that $y \not\in H$ for every $H \in \cH_Q$. 

\begin{lemma}\label{lem:nearest}
Let\/ $\cB = (\E^{k+1},\SS^k,\cH,\bfdelta)$ be a good spherical buffer system with $k \geq 1$. Let $x \in \SS^k \setminus \Qbar(\cB)$ and\/ $y \in M( \Qbar(\cB), x )$, and suppose that $y \not\in H$ for every $H \in \cH$. Then there exists $P \in \cH$ such that\/ $y \in \SS_{x,P} \subset \partial\cE(P)$. 
\end{lemma}

\begin{proof}
Observe first that $x \ne y$, since $y \in \Qbar(\cB)$ and $x \not\in \Qbar(\cB)$. Since $y$ is one of the closest points of $\Qbar(\cB)$ to $x$, it follows that if we move along a geodesic in $\SS^k$ from $y$ to $x$ then we must immediately leave $\Qbar(\cB)$. Note that geodesics exist in $\SS^k$ because we have assumed that $k \geq 1$. Recall that, by Definition~\ref{def:Qbar}, and since $\E^{k+1} \in \cH^*$, we have
\begin{equation}\label{eq:not:Q:in:buffer}
\SS^k \setminus \Qbar(\cB) \subset \bigcup_{H \in \cH} \cE(H).
\end{equation}
Moreover, since $y \in \Qbar(\cB)$ and $y \not\in H$ for every $H \in \cH$, it again follows by Definition~\ref{def:Qbar} that
\begin{equation}\label{eq:y:not:in:a:buffer}
y \in \SS^k \setminus \bigcup_{H \in \cH} \cE(H).
\end{equation}
Together with our observation above about moving along a geodesic from $y$ to $x$, it follows from~\eqref{eq:not:Q:in:buffer} and~\eqref{eq:y:not:in:a:buffer} that $y$ lies in the boundary of at least one buffer. Choose a minimal $P \in \cH$ with $y \in \partial\cE(P)$. Note that $x \not\in P^\perp$, by Lemma~\ref{lem:nearest:not:perp}. We claim that $y \in \SS_{x,P}$. Since $\SS_{x,P} \subset \partial\cE(P)$, by Lemma~\ref{lem:Axp:exists}, this claim will prove the lemma.  

Our proof that $y \in \SS_{x,P}$ comes in two steps: in the first, we show that 
\begin{equation}\label{eq:M:in:SxP}
M\big( \partial\cE(P) \cap (y+P), x \big) = \big\{ x' \big\}
\end{equation} 
for some $x' \in \SS_{x,P}$; in the second, we show that $x' = y$. In both steps, it will be convenient to assume 
that $\E^{k+1} = \R^{k+1}$ and $P = \R^{j+1} \times \{0\}^{k-j}$, so
$$\partial \cE(P) = \SS^j_{\rho_j} \times \SS^{k-j-1}_{\delta_j},$$ 
where $\rho_j^2+\delta_j^2=\delta_{-1}^2$ (cf.~the proof of Lemma~\ref{lem:Axp:exists}). Let $y = (y_1,y_2) \in P \times P^\perp$, and note that $\|y_2\| = \delta_j$ (since $y \in \partial\cE(P)$), and that
\begin{equation}\label{eq:partial:cap:y+P}
\partial\cE(P) \cap (y+P) = \SS^j_{\rho_j} \times \{y_2\}.
\end{equation}
Recalling~\eqref{eq:wlog:SxP}, and noting that $x_1 \ne 0$, since $x \not\in P^\perp$, it follows that the (unique) closest point to $x = (x_1,x_2) \in P \times P^\perp$ in this set is 
$$x' = (\lambda x_1,y_2) \in \{ \lambda x_1\} \times \SS^{k-j-1}_{\delta_j} = \SS_{x,P},$$
where $\|\lambda x_1\| = \rho_j$ and $\lambda > 0$. Thus $x' \in \SS_{x,P}$, and hence we have~\eqref{eq:M:in:SxP}, as claimed. 

Now suppose that $y \ne x'$. Let us deal with the case $j = 0$ separately, since in general we shall want to take geodesics in $\SS_{\rho_j}^j \times \{y_2\}$, and these do not exist if $j = 0$. Fortunately this case is straightforward: indeed, now $\SS_{\rho_j}^j \times \{y_2\} = \{-\lambda x_1, \lambda x_1\} \times \{y_2\}$, and we are worried about the case $y_1 = -\lambda x_1$. But this would imply that
\[
d(x,y) \geq (1+\lambda)\|x_1\| > \rho_j.
\]
However, by Lemma~\ref{lem:nearest:not:perp}, we also have $d(x,y) \leq 2 \delta_0$. Since $\rho_j^2 = \delta_{-1}^2 - \delta_j^2 > 4\delta_0^2$, this is a contradiction, and it follows that if $j = 0$, then $y = x'$, as claimed.

We may therefore assume that $j \geq 1$. By the same reasoning as in the case $j = 0$, we may moreover assume that $x'$ and $y$ are not antipodal points in the sphere $\SS^j_{\rho_j} \times \{y_2\}$. Let $G$ be the (unique) geodesic in $\SS^j_{\rho_j} \times \{y_2\}$ from $y$ to $x'$, and let $Q$ be the affine span of $G$. Note that $Q$ is a 2-dimensional affine space, and is the same as the affine span of the non-collinear points $x'$, $y$ and $(0,y_2)$, because $(0,y_2)$ is the centre of the sphere $\SS^j_{\rho_j} \times \{y_2\}$. 

Now, since $x'$ is the (unique) closest point of $\SS^j_{\rho_j} \times \{y_2\}$ to $x$, and therefore, by Observation~\ref{obs:spheres:increasing:distance}, the function $z \mapsto d(x,z)$ is increasing in $d(x',z)$ for $z \in \SS^j_{\rho_j} \times \{y_2\}$. It follows that every point $z$ of the geodesic $G$ (other than $y$ itself), being strictly closer to $x'$ than $y$, is also strictly closer to $x$ than $y$. Moreover, since $\SS^j_{\rho_j} \times \{y_2\} = \partial\cE(P) \cap (y+P)$, by~\eqref{eq:partial:cap:y+P}, we have $G \subset Q \subset y + P$, and therefore the geodesic $G$ is parallel to $P'$ for every $P \subset P' \in \cH$. In particular, this implies that $G$ does not cross the boundary of $\cE(P')$. 

Since every point of $G \setminus \{y\}$ is strictly closer to $x$ than $y$, and recalling~\eqref{eq:not:Q:in:buffer} and~\eqref{eq:y:not:in:a:buffer}, and that $y \in M( \Qbar(\cB), x )$, it follows that $y \in \partial \cE(P')$ for some $P \not\subset P' \in \cH$. Since $P$ is minimal such that $y \in \partial \cE(P)$, we also have $P' \not\subset P$, and since $\cB$ is good, it follows that
$$y \in \partial \cE(P) \cap \partial\cE(P') \subset \cE^+(P) \cap \cE^+(P') \subset \cE(P\cap P'),$$ 
contradicting~\eqref{eq:y:not:in:a:buffer}. This implies that $y = x'$, and hence $y \in \SS_{x,P}$, as required.
\end{proof}

\subsection{The proof of Lemma~\ref{thm:quasi:key}}\label{quasi:technical:proof:sec}

We are finally ready to prove the key technical lemma of the section, Lemma~\ref{thm:quasi:key}. 

\begin{proof}[Proof of Lemma~\ref{thm:quasi:key}]
The proof is by induction on $k$. Note that the statement is vacuous when $k = 0$, since each member of $\cH_{-1}$ is equal to $\{O_\cB\}$, so let $k \ge 1$ and assume that the lemma holds for all smaller values of $k$. Let $\cB = (\E^{k+1},\SS^k,\cH,\bfdelta)$ be a good spherical buffer system, let $H \in \cH$, and let $x \in \SS^k \cap H$. 

Observe first that if $x \in \Qbar = \Qbar(\cB)$ then we are done by Observation~\ref{obs:Qdef3}, since if $y \in \Qbar \setminus H$ then $y \not\in \cE(H)$, and hence 
$$d\big( x, \,\Qbar \cap H \big) = 0 < \delta_{k-1} \le d\big( x, \,\Qbar \setminus H \big).$$
We may therefore assume that $x \notin \Qbar$, and hence, by Observation~\ref{obs:unique:buffer} and since $\cB$ is good, that there exists a unique minimal $P' \in \cH$ such that $x \in \cE(P')$. Observe that $x \in H \subset \cE(H)$, and therefore either $P' = H$ or $H \not\subset P'$, by the minimality of $P'$. Since $\cB$ is good, it follows that if $P' \not\subset H$, then
\begin{equation}\label{eq:showing:PinH}
x \in \cE(P') \cap H \subset \cE^+(P') \cap \cE^+(H) \subset \cE(P' \cap H),
\end{equation}
contradicting the minimality of $P'$. We therefore have $P' \subset H$. 

Suppose, for a contradiction, that $d( x, \Qbar \setminus H ) = d( x, \Qbar )$, and recall that $\Qbar \setminus H$ is compact, by Observation~\ref{obs:compact},  since $\cB$ is good. It follows that there exists $y \in M(\Qbar,x)$ with $y \not\in H$. Note that $y \in \cE^+(P')$, by Lemma~\ref{lem:minbuffer}.

Now let $Q \in \cH^*$ be minimal such that $y \in Q$. Note that we may have $Q = \E^{k+1}$, but we cannot have $Q = \{O_\cB\}$. Observe that $y \notin P'$, since $P' \subset H$ and $y \notin H$, and therefore $Q \not\subset P'$, since $y \in Q$. It follows that if $Q \in \cH_0$ (that is, if $\dim(Q) = 1$), then $Q \not\subset P'$ and $P' \not\subset Q$ (since $\cE(P')$ is non-empty, so $P' \ne \{O_\cB\}$), and also $P' \cap Q = \{O_\cB\}$, so 
$$y \in \cE^+(P') \cap Q \subset \cE(P' \cap Q) = \emptyset,$$
a contradiction. We may therefore assume that $\dim(Q) \ge 2$. Let $\cB_Q$ be the spherical buffer system restricted to $Q$, as defined in Definition~\ref{def:BR}. Recall from Lemma~\ref{lem:sbs_hyper} that $\cB_Q$ is a good spherical buffer system (since $Q \neq \{O_\cB\}$), and $\Qbar(\cB_Q) = \Qbar(\cB) \cap Q$.  

Let $\hat x \in M(\SS^k \cap Q,x)$. By Lemma~\ref{lem:Qbar:nonempty} we have $\Qbar(\cB) \cap \S_{x,P'} \ne \emptyset$, and therefore $x \notin Q^\perp$, by Lemma~\ref{lem:another:not:perp}. It follows that
$$M\big( \SS^k \cap Q,x \big) = \big\{ \hat x \big\},$$
by Observation~\ref{obs:spheres:increasing:distance}. Let us next eliminate the case $y = \hat x$.

\begin{claim}\label{clm:quasi:case1}
We have $y \ne \hat x$, and in particular $\hat x \not\in \Qbar(\cB_Q)$. 
\end{claim}

\begin{clmproof}{clm:quasi:case1}
Recall that $y \notin P'$ and $Q \not\subset P'$. If also $P'\not\subset Q$ then, since $\cB$ is good, we have 
$$y \in \cE^+(P') \cap Q \subset \cE(P' \cap Q).$$
But $y \notin P' \cap Q$, so by Observation~\ref{obs:Qdef3} this contradicts the fact that $y \in \Qbar(\cB)$. 

We may therefore assume that $P' \subsetneq Q$. If $y = \hat x \in M(\SS^k \cap Q, x)$, then it follows by Observation~\ref{lem:proj} that $d(y,P') \le d(x,P')$, and hence $y \in \cE(P')$. But $y \notin P'$, so by Observation~\ref{obs:Qdef3} we again obtain a contradiction to the fact that $y \in \Qbar(\cB)$. 

Finally, since $M( \SS^k \cap Q,x ) = \{ \hat x \}$ and $\Qbar(\cB_Q) = \Qbar(\cB) \cap Q$, if $\hat x \in \Qbar(\cB_Q)$ then it would follow that $M( \Qbar(\cB) \cap Q,x ) = \{ \hat x \}$, and therefore $y = \hat x$.  
\end{clmproof}

We would like to apply Lemma~\ref{lem:nearest} to the good spherical buffer system $\cB_Q$, and the vectors $\hat x \in (\SS^k \cap Q) \setminus \Qbar(\cB_Q)$ (which holds by Claim~\ref{clm:quasi:case1}) and $y$. To verify the conditions of the lemma, recall that $\dim(Q) \ge 2$, and observe that 
$$y \in M\big( \Qbar(\cB_Q),\hat x \big),$$  
by Lemma~\ref{obs:y:closest:to:xhat}, since $\Qbar(\cB_Q) = \Qbar(\cB) \cap Q$. Moreover, since $Q \in \cH^*$ is minimal such that $y \in Q$, it follows by Observation~\ref{obs:y:not:in:HR} that $y \not\in H'$ for every $H' \in \cH_Q$.

We may therefore apply Lemma~\ref{lem:nearest}, and deduce that there exists $P \in \cH_Q$ (so, in particular, $P \in \cH$ and $P \subset Q$, and moreover $\hat{x} \not\in P^\perp$) with 
\begin{equation}\label{eq:y:in:SxP:final}
y \in \SS_{\hat x,P}(\cB_Q) = \SS_{\hat x,P}(\cB) \cap Q \subset \partial\cE(P),
\end{equation}
where the equality holds by Lemma~\ref{lem:SxP:BR}. It follows that $x \not\in P^\perp$, by Lemma~\ref{lem:nearest:not:perp}, and hence that $\SS_{\hat x,P}(\cB) = \SS_{x,P}(\cB)$, by Lemma~\ref{lem:SxP:SyP} if $Q \in \cH$, and since $\hat x = x$ if $Q = \E^{k+1}$. Thus, by~\eqref{eq:y:in:SxP:final}, we have $y \in \SS_{x,P}$.

To obtain the desired contradiction, we split into three cases.

\medskip
\noindent Case 1: $P \not\subset H$.
\medskip

We claim first that $x \notin \cE(P)$. To see this, observe that if $x \in H \cap \cE(P)$, then since $P \not\subset H$, it follows by Lemma~\ref{lem:scover} that $y \in \SS_{x,P} \subset \cE(H \cap P)$. But $y \notin H \cap P$, since $y \in \partial \cE(P)$, and by Observation~\ref{obs:Qdef3} this implies that $y \not\in \Qbar(\cB)$, which is a contradiction. Hence, since $x \in H$, it must be that $x \notin \cE(P)$, as claimed.

Next, we claim that  $P' \subsetneq P$. To see this, note that $P \not\subset P'$, since $P' \subset H$ but $P \not\subset H$, and that $y \in \SS_{x,P} \subset \partial \cE(P)$, by Lemma~\ref{lem:Axp:exists}. Thus, if $P' \not\subset P$, then 
$$y \in \cE^+(P) \cap \cE^+(P') \subset \cE(P \cap P'),$$ 
since $\cB$ is good. But $y \notin P \cap P'$, since $y \in \partial \cE(P)$ (and thus $y \notin P$), and hence, by Observation~\ref{obs:Qdef3}, it again follows that $y \not\in \Qbar(\cB)$.  

Now, since $x \notin \cE(P) \cup P^\perp$ and $P' \subset P$, it follows by Lemma~\ref{lem:xiner} that $d(w,P') \le d(x,P')$ for every $w \in\SS_{x,P}$. Recalling that $x \in \cE(P')$, it follows that $\SS_{x,P} \subset \cE(P')$, and hence that $y \in \cE(P')$. But $y \notin P'$, since $P' \subset H$ and $y \notin H$, and hence, by Observation~\ref{obs:Qdef3}, we again deduce that $y \not\in \Qbar(\cB)$, which gives us our desired contradiction.

\medskip
\noindent Case 2: $x \in P \subset H$.
\medskip

In this case we apply Lemma~\ref{lem:great:escape} to $\cB$. To do so, note that $x \in \SS^k \cap P$ and $y \in \Qbar(\cB) \cap \SS_{x,P}$. It follows that 
$$M(\Qbar(\cB),x) \subset P \subset H,$$ 
and hence $y \in H$, which is a contradiction. 

\medskip
\noindent Case 3: $x \notin P \subset H$.
\medskip

In this final case we apply the induction hypothesis to the spherical buffer system $\cB_{x,P}$, which is good by Lemma~\ref{lem:sbs_buffer}. Since $P \subset H$ and $x \in H \setminus (P \cup P^\perp)$, it follows by Lemma~\ref{lem:SxP:in:H} that $M(\SS_{x,P},x) = \{x'\}$ for some $x' \in \SS_{x,P} \cap H$. Moreover, $H \cap A_{x,P} \in \cH_{x,P}$, since $P \subset H \in \cH$. By the induction hypothesis, it follows that
\begin{equation}\label{eq:IH:for:BxP}
d\big( x', \,\Qbar(\cB_{x,P}) \cap H \big) < d\big( x', \,\Qbar(\cB_{x,P}) \setminus H \big).
\end{equation}
Now, since $y \in \Qbar(\cB) \cap \SS_{x,P}$, we have $\Qbar(\cB_{x,P}) = \Qbar(\cB)\cap \SS_{x,P}$, by Lemma~\ref{lem:sbs:Q}, and thus
$$y \in M\big( \Qbar(\cB), x \big) \cap \SS_{x,P} \subset M\big( \Qbar(\cB) \cap \SS_{x,P}, x' \big) = M\big( \Qbar(\cB_{x,P}), x' \big),$$
the inclusion following from Lemma~\ref{lem:SxP:closestpoints}, since $M(\SS_{x,P},x) = \{x'\}$. Hence, by~\eqref{eq:IH:for:BxP}, we deduce that $y \in H$, which is our final contradiction. This completes the proof of Lemma~\ref{thm:quasi:key}.
\end{proof}

\subsection{Deducing that $\QQ$ exists}\label{sec:quasi:exists}

In order to apply Lemma~\ref{thm:quasi:key}, we first need to observe that we can choose $\bfdelta$ so as to make our spherical buffer system good. From now on, $\SS^k$ will (as earlier) denote a $k$-dimensional unit sphere centred at the origin of $\R^d$.  

\begin{lemma}\label{lem:choose:deltas:so:good}
Let $\cH$ be a finite collection of proper subspaces of $\R^{k+1}$, closed under intersections. Then there exists a choice of\/ $\bfdelta \in \R^{k+1}$ such that\/ $\cB = (\R^{k+1},\SS^k,\cH,\bfdelta)$ is a good spherical buffer system.
\end{lemma}

\begin{proof}
We set $\delta_{-1} := 1$ (i.e., the radius of $\SS^k$), and define the remaining $\delta_i$ inductively. Recall that for $\cB$ to be a spherical buffer system we need
\begin{equation}\label{eq:choose:deltas:decreasing}
0 < \delta_i < \delta_{i-1}/3
\end{equation}
for each $0 \le i \le k-1$, and for $\cB$ to be good we need
$$\cE^+(H) \cap \cE^+(H') \subset \cE(H \cap H')$$
for every $H,H' \in \cH$ such that $H \not\subset H'$ and $H' \not\subset H$.

To define $\delta_0$, note first that if $H,H'\in \cH$ with $H \cap H' = \{\0\}$, then for any $z \in \SS^k$, we have $d(z,H) + d(z,H') > 0$. By compactness, it follows that there exists $\delta \in (0,1/3)$ such that $d(z,H) + d(z,H') > 6\delta$ for all $z \in \SS^k$. We define $\delta_0$ to be the minimum of these
values of $\delta$ over the (finite) set of choices of such pairs $H,H' \in \cH$, noting that this satisfies~\eqref{eq:choose:deltas:decreasing}. It follows that
$$\cE^+(H) \cap \cE^+(H') = \emptyset$$
for every pair $H,H' \in \cH$ with $H \cap H' = \{\0\}$.

Assume now that we have defined $\delta_{-1},\dots,\delta_i$ satisfying~\eqref{eq:choose:deltas:decreasing} and suppose that $H,H' \in \cH$ are such that $H \not\subset H'$, $H' \not\subset H$ and $H \cap H' \in \cH_i$. Note that
$$\bigcap_{\delta > 0} \, \overline{B}_{\delta}(H) \cap \overline{B}_{\delta}(H') \cap \SS^k = H \cap H' \cap \SS^k,$$
where $\overline{B}_{\delta}(P) := \big\{ z \in \R^d : d(z,P) \le \delta \big\}$ for each $P \in \cH$. Since $\cE(H \cap H')$ is an open set in $\SS^k$ containing $H \cap H' \cap \SS^k$, it follows by compactness that there exists $\delta \in (0,\delta_i/3)$ such that
$$\overline{B}_{3\delta}(H) \cap \overline{B}_{3\delta}(H') \cap \SS^k \subset \cE(H \cap H').$$
Take $\delta_{i+1}$ to be the minimum $\delta$ over all such pairs $H,H' \in \cH$. Now, observe that
$$\dim(H \cap H') = i + 1 < \min\big\{ \dim(H), \dim(H') \big\},$$
and therefore, since $\delta_{-1},\dots,\delta_{i+1}$ satisfy~\eqref{eq:choose:deltas:decreasing}, we have
$$\cE^+(H)\cap \cE^+(H') \subset \overline{B}_{3\delta_{i+1}}(H) \cap \overline{B}_{3\delta_{i+1}}(H') \cap \SS^k,$$
and hence $\cE^+(H) \cap \cE^+(H') \subseteq \cE(H \cap H')$ for all pairs $H,H' \in \cH$ with $H \cap H' \in \cH_i$. Continuing in this way until we have defined the entire vector $\bfdelta$ completes the proof.
\end{proof}

\pagebreak

In order to complete the proof of Lemma~\ref{lem:quasi-new}, we need to choose a suitable finite subset $\QQ \subset \Qbar(\cB) \cap \SS_\Q^{d-1}$, where $\cB = (\R^d,\SS^{d-1},\cH,\bfdelta)$ is the following spherical buffer system. Recall that the constant $R = R(\U)$ was fixed in Lemma~\ref{lem:exists-rational-w}, and define $\cH = \cH(R)$ to be the family of subspaces of $\R^d$ obtained via intersections of subspaces of the form $\{w\}^\perp$ with $w \in \cL_R$. By Lemma~\ref{lem:choose:deltas:so:good}, there exists $\bfdelta(R)$ such that the spherical buffer system
\begin{equation}\label{def:our:buffer:system}
\cB(R) = (\R^d,\SS^{d-1},\cH(R),\bfdelta(R))
\end{equation}
is good. Set $\Qbar := \Qbar(\cB(R))$. In order to guarantee that we can choose a sufficiently `dense' set of rational directions in $\Qbar$, we require the following lemma.

\begin{lemma}\label{lem:quasi-rationals-dense}
For each $H \in \cH(R)^*$, the set $H \cap \SS_\Q^{d-1}$ is dense in $H \cap \SS^{d-1}$.
\end{lemma}

\begin{proof}
By the definition of $\cH(R)^*$, there exists a finite (possibly empty) set $W \subset \cL_R$ such that $H = W^\perp$. By Lemma~\ref{lem:Vperp:rational}, it follows that $W^\perp$ has a rational basis. Therefore, in order to approximate $x \in H \cap \SS^{d-1} = \SS(W)$ by a point of $H \cap \SS_\Q^{d-1}$, we simply approximate each coordinate of $x$ (with respect to the rational basis of $W^\perp$) by a rational number, and then project onto the sphere $\SS(W)$, recalling from Observation~\ref{obs:proj:is:rational} that the point we obtain by doing so is rational. 
\end{proof}

We are finally ready to prove the main result of this section.

\begin{proof}[Proof of Lemma~\ref{lem:quasi-new}]
Let $\cB(R) = (\R^d,\SS^{d-1},\cH(R),\bfdelta(R))$ be the spherical buffer system defined in~\eqref{def:our:buffer:system}, and recall that $\cB(R)$ is good, since we chose $\bfdelta(R)$ using Lemma~\ref{lem:choose:deltas:so:good}. Applying Lemma~\ref{thm:quasi:key} to $\cB(R)$, it follows that there exists $\theta > 0$ such that
\begin{equation}\label{eq:quasi:theta}
d\big( x, \,\Qbar \cap H \big) < d\big( x, \,\Qbar \setminus H \big) - \theta
\end{equation}
for all $H \in \cH(R)$ and $x \in \SS^{d-1} \cap H$, since $\cH(R)$ is finite and $\SS^{d-1} \cap H$ is compact.

We construct $\QQ$ from $\Qbar$ in two steps: first, for each $H \in \cH(R)^*$, we greedily choose a finite subset $\X_H$ of the interior of $\Qbar \cap H$ with respect to $\SS^{d-1} \cap H$, that is, of the set
$$\Y_H := \SS^{d-1} \cap H \setminus \bigcup_{H' \subsetneq H} \overline{\cE}(H'),$$
where $\overline{\cE}(H')$ denotes the closure of the buffer $\cE(H')$ with respect to $\SS^{d-1} \cap H$, with the following property:
\begin{itemize}
\item[$(a)$] For all $w \in \Y_H$, there exists $u \in \X_H$ such that $d(u,w) \leq \theta/2$.
\end{itemize}
Thus, if there exists $w \in \Y_H$ such that $d(u,w) > \theta/2$ for all $u \in \X_H$, then add $w$ to $\X_H$, and repeat. This procedure terminates because there exists a finite cover of $\Y_H$ with closed balls of radius $\theta/6$, and each such ball can contain only one element of $\X_H$.

We now simply adjust each element of $\X_H$ slightly, using Lemma~\ref{lem:quasi-rationals-dense}, to obtain a set $\X'_H \subset \Y_H \cap \SS_\Q^{d-1}$ of rational directions such that:
\begin{itemize}
\item[$(b)$] For all $w \in \Y_H$, there exists $u \in \X'_H$ such that $d(u,w) < \theta$.
\end{itemize}
Note that for this step it was important that every point of $\X_H$ is an interior point of $\Qbar \cap H$ with respect to $\SS^{d-1} \cap H$, so that we could take an open ball around an element of $\Qbar \cap H$ (again, open with respect to $\SS^{d-1} \cap H$) and apply Lemma~\ref{lem:quasi-rationals-dense}.

\pagebreak

Finally, define
$$\QQ := \bigcup_{H \in \cH(R)^*} \X'_H.$$
We claim that $\QQ$ is quasistable for range $R$; that is, it is a finite set of rational directions that intersects every open hemisphere of $\SS^{d-1}$, and such that~\eqref{eq:quasi:property} holds for all if $u \in \QQ$ and $w \in \cL_R$. This will follow easily from the following claim.

\begin{claim}\label{claim:quasi:property:proof}
$$d\big( x, \,\QQ \cap H \big) < d\big( x, \,\QQ \setminus H \big)$$
for every $H \in \cH(R)$ and $x \in \SS^{d-1} \cap H$. 
\end{claim}

\begin{clmproof}{claim:quasi:property:proof}
Recalling that $\Qbar \cap H$ is compact, let $y \in M(\Qbar \cap H,x)$ and let $P \in \cH(R)$ be minimal such that $y \in P$. Note that $P \subset H$, since $\cH(R)$ is closed under taking intersections, and that 
\begin{equation}\label{eq:y:in:bigger:YP}
y \in \SS^{d-1} \cap P \setminus \bigcup_{H' \subsetneq P} \cE(H'),
\end{equation}
by Observation~\ref{obs:Qdef3} and the minimality of $P$, and since $y \in \Qbar$. Now, if $d(y,\Y_P) = 0$, then by property $(b)$ and~\eqref{eq:quasi:theta}, and since $\QQ \subset \Qbar$, it follows that 
$$d\big( x, \QQ \cap H \big) \le d(x,u) < d(x,y) + \theta < d( x, \Qbar \setminus H) \le d( x, \QQ \setminus H)$$
for some $u \in \X'_P$, as claimed. It will therefore suffice to show that $d(y,\Y_P) = 0$. 

To do so, it will be convenient to set $\delta_{H'} := \delta_i$ for each $i \in \{0,\dots,d-1\}$ and $H' \in \cH_i$. Note that $d(y,H') \ge \delta_{H'}$ for all $H' \subsetneq P$, by~\eqref{eq:y:in:bigger:YP}, and that if $y \notin \Y_P$, then 
\begin{equation}\label{eq:y:is:on:the:edge}
d(y,H') = \delta_{H'}
\end{equation}
for some $H' \subsetneq P$. Let $H'$ be minimal such that~\eqref{eq:y:is:on:the:edge} holds, and define
$$y' := \frac{y + \eps z}{\|y + \eps z\|},$$
where $y = (w,z) \in H' \times (P \cap (H')^\perp)$, and $\eps > 0$ is sufficiently small. Note that $y' \in \SS^{d-1} \cap P$; we claim that moreover $y' \notin \overline{\cE}(H'')$ for all $H'' \subsetneq P$, and hence that $y' \in \Y_P$. Since $\eps$ can be chosen arbitrarily small, this will suffice to prove that $d(y,\Y_P) = 0$. 

Suppose, for a contradiction, that $y' \in \overline{\cE}(H'')$ for some $H'' \subsetneq P$. We claim first that $H' \subset H''$. To see this, recall first that $d(y,H'') \ge \delta_{H''}$, and note that therefore $d(y,H'') = \delta_{H''}$, since we chose $\eps$ sufficiently small. Thus, by the minimality of $H'$, we cannot have $H'' \subset H'$. Now, since $\B(R)$ is good, if $H' \not\subset H''$ then it would follow that $y \in \cE^+(H') \cap \cE^+(H'') \subset \cE(H' \cap H'')$, contradicting~\eqref{eq:y:in:bigger:YP}. We therefore have $H' \subset H''$. 

Now, to compare $d(y,H'')$ and $d(y',H'')$, recall that $y = (w,z) \in H' \times (P \cap (H')^\perp)$, let $z = (z_1,z_2) \in (H'' \cap (H')^\perp) \times (P \cap (H'')^\perp)$, and observe that
$$d(y',H'') = \frac{(1+\eps)\|z_2\|}{\|y + \eps z\|}.$$
Next, note that, since $y \in \SS^{d-1}$ and $\<y,z\> = \<w + z,z\> = \|z\|^2$, we have
\[
\|y + \eps z\|^2 = 1 + 2\eps\<y,z\> + \eps^2\|z\|^2 = 1 + 2\eps\|z\|^2 + \eps^2\|z\|^2 < (1 + \eps\|z\|)^2,
\]
where the last inequality holds since $\|z\| = \delta_{H'} < 1$. It follows that
$$d(y',H'') > \bigg( \frac{1+\eps}{1 + \eps\|z\|} \bigg) \|z_2\| > \|z_2\| = d(y,H'') \ge \delta_{H''},$$
since $y \in \SS^{d-1} \setminus \cE(H'')$, by~\eqref{eq:y:in:bigger:YP}. But this contradicts our assumption that $y' \in \overline{\cE}(H'')$, and hence completes the proof of the claim.
\end{clmproof}

To complete the proof of Lemma~\ref{lem:quasi-new}, we need to show that $\QQ$ intersects every open hemisphere of $\SS^{d-1}$, and that~\eqref{eq:quasi:property} holds for every $u \in \QQ$ and $w \in \cL_R$. The first of these two properties follows from the fact that $\QQ$ contains all of the $2d$ nearest neighbours $\{\pm e_1,\ldots,\pm e_d\}$ of the origin in $\Z^d$. Indeed, if $W = \{e_1,\ldots,e_d\} \setminus \{e_i\}$, then $W^\perp \cap \SS^{d-1} = \{e_i,-e_i\}$. It follows that $e_i$ and $-e_i$ are isolated points of $\Qbar$, and hence are chosen in $\QQ$.

To prove that~\eqref{eq:quasi:property} holds, note first that if $\<u,w\> = 0$ then $u \in \Cell_\QQ(u) \cap \{w\}^\perp$, and therefore one direction is trivial. To prove the other direction, let $x \in \Cell_\QQ(u) \cap \{w\}^\perp$ and suppose that $\<u,w\> \ne 0$. By Claim~\ref{claim:quasi:property:proof}, applied with $H = \{w\}^\perp$, it follows that
$$d(x,v) = d\big( x, \,\QQ \cap \{w\}^\perp \big) < d\big( x, \,\QQ \setminus \{w\}^\perp \big) \le d(x,u)$$
for some $v \in \QQ \cap \{ w\}^\perp$. However, by Definition~\ref{def:voronoi}, we have $\< x,u \> \ge \< x,v \>$ for all $v \in \QQ$, and hence $d(x,u) \le d(x,v)$. This contradiction proves~\eqref{eq:quasi:property}, and hence completes the proof of Lemma~\ref{lem:quasi-new}.
\end{proof}

\section{The resistance of induced update families}\label{sec:induced}

In this section we apply the results of the previous three sections to prove three key lemmas about the resistance of induced update families. The definition of resistance (Definition~\ref{def:r}) is specifically designed for the torus; indeed, $r^k$ is always defined in terms of the `easiest' open hemisphere. However, the faces of droplets are non-toral, so we actually need to control the difficulty of growing in an arbitrary direction on a face. The two main lemmas of this section (Lemmas~\ref{lem:good:faces} and~\ref{lem:semigood:faces}) address this problem. 

Let us fix, for the rest of the paper, a quasistable set $\QQ \subset \SS_\Q^{d-1}$ for\/ range $R = R(\U)$, that is, a finite set of directions, intersecting every open hemisphere of $\SS^{d-1}$, such that 
\begin{equation}\label{eq:quasi:property:repeat}
\Cell_\QQ(u) \cap \{w\}^\perp \neq \emptyset \qquad \Leftrightarrow \qquad \< u,w \> = 0
\end{equation}
for each $u \in \QQ$ and $w \in \cL_R$. Recall that such a set is guaranteed to exist by Lemma~\ref{lem:quasi-new}, and that the constant $R = R(\U)$ was chosen in Lemma~\ref{lem:exists-rational-w}.

Several of the lemmas proved in this section will hold for cliques\footnote{We include the empty set in the family of cliques, and also the set $\{w\}$ for each $w \in \QQ$.} in the Voronoi graph $\Vor(\QQ)$ (see Definition~\ref{def:voronoi}), rather than for arbitrary sets $W \subset \QQ$. This is because the update families induced by cliques behave (in a certain sense) `as expected' (see Lemma~\ref{lem:induced-surprise}) due to the properties of $\QQ$ (in particular, Lemma~\ref{lem:quasi-innerprod}). In Section~\ref{sec:polytopes} we will introduce a family of polytopes that we will use to control the growth of a droplet, and show that each face of  each polytope in our family corresponds to a clique in $\Vor(\QQ)$. 

\pagebreak

We will begin with a key definition, which captures the information we need about growing in all directions (not just easy directions) in induced processes. Recall from~\eqref{def:SW} that we write $\S_W$ for the stable set in $\SS(W)$ of the induced process $\U[W]$.

\begin{definition}\label{def:sgood}
Let $W \subset \QQ$ and let $k := \dim(W^\perp)$.
\begin{itemize}
\item We say $W$ is \emph{$(s,w)$-semi-good}, where $1 \leq s \leq k$ and $w \in \SS(W)$, if
\begin{equation}\label{eq:ssemigood}
\rho^{k-1} \big( \SS(W); \S_W, v \big) \leq s - 1
\end{equation}
for all $v \in \SS(W)$ such that $\< v, w \> > 0$. \label{ftn:sgood} \smallskip
\item We say $W$ is \emph{$s$-good} if it is $(s,w)$-semi-good for some $w \in \SS(W)$, and also
\begin{equation}\label{eq:sgood}
\rho^{k-1}\big( \SS(W); \S_W, v \big) \leq \min\big\{ s, k - 1 \big\}
\end{equation}
for all $v \in \SS(W)$.
\end{itemize}
\end{definition}

The first of our three key lemmas about induced processes provides us with a rational direction $w$ such that $W$ is $(s,w)$-semi-good, and is a simple consequence of Lemma~\ref{lem:exists-rational-w}. This lemma will play a crucial role in Section~\ref{proof:sec}. 

\begin{lemma}\label{lem:Q-exists-rational-w}
Let $W \subset \QQ$, and let $1 \le s \le k := \dim(W^\perp)$. If\/ $\U[W]$ is non-trivial, and
$$r^k\big( \SS(W); \S_W \big) \le s,$$
then there exists $w \in \cL_R \cap \SS(W)$ such that $W$ is $(s,w)$-semi-good.
\end{lemma}

\begin{proof}
By Lemma~\ref{lem:exists-rational-w} and our choice of $R = R(\U)$, there exists $w \in \cL_R \cap \SS(W)$ such that
$$\rho^{k-1}\big( \SS(W); \S_W, u \big) \leq r^k\big( \SS(W); \S_W \big) - 1 \le s - 1$$
for all $u \in \SS(W)$ such that $\< u,w \> > 0$. This implies that $W$ is $(s,w)$-semi-good, as claimed.
\end{proof}

We will also need the following easy consequence of Lemma~\ref{lem:Q-exists-rational-w}. 

\begin{lemma}\label{lem:sgood-rational-w}
Let $W \subset \QQ$, let $1 \le s \le \dim(W^\perp)$, and suppose that\/ $\U[W]$ is non-trivial. If\/ $W$ is $s$-good, then there exists $w \in \cL_R \cap \SS(W)$ such that $W$ is $(s,w)$-semi-good.
\end{lemma}

\begin{proof}
Since $W$ is $s$-good, it follows, by Definition~\ref{def:sgood}, that $W$ is $(s,w')$-semi-good for some $w' \in \SS(W)$, and therefore, setting $k := \dim(W^\perp)$, that
$$\rho^{k-1} \big( \SS(W); \S_W, v \big) \leq s - 1$$
for all $v \in \SS(W)$ such that $\< v, w' \> > 0$. By~\eqref{eq:r}, it follows that $r^k\big( \SS(W); \S_W \big) \le s$, and hence, by Lemma~\ref{lem:Q-exists-rational-w}, there exists $w \in \cL_R \cap \SS(W)$ as required.
\end{proof}

We will next state the two main lemmas of this section, whose proofs will require significantly more work. They allow us to deduce properties of the update family induced by a clique $W'$ from similar properties of the family induced by a sub-clique $W \subset W'$. The first applies when $W$ is $s$-good for some $1 \le s \le \dim(W^\perp)$. 

\medskip
\pagebreak

\begin{lemma}\label{lem:good:faces}
Let $1 \le s \le k \le d$, and let $W \subset \QQ$ be an $s$-good clique with $\dim(W^\perp) = k$. Let $W \subset W' \subset \QQ$ be a clique, and set $k' := \dim( W'^\perp )$ and $s' := \min\{s,k'\}$.
\begin{enumerate}
\item If\/ $k' = 0$ then\/ $\emptyset \in \U[W']$.\vskip0.1cm
\item If\/ $k' > 0$ then\/ $W'$ is $s'$-good.
\end{enumerate}
\end{lemma}

The second lemma applies when $W$ is $(s,w)$-semi-good for some $1 \le s \le \dim(W^\perp)$ and some $w \in \SS(W)$. 

\begin{lemma}\label{lem:semigood:faces}
Let $1 \le s \le k \le d$, and let $W \subset \QQ$ be a clique with $\dim(W^\perp) = k$. Let\/ $w \in \SS(W)$, let $W' \subset \QQ$ be a clique with $W \subsetneq W' \subset \{w\}^\perp$, and set $k' = \dim( W'^\perp )$ and $s' = \min\{s,k'\}$. If\/ $W$ is $(s,w)$-semi-good, then $W'$ is $(s',w)$-semi-good.
\end{lemma}

The section is organized as follows: first, in Section~\ref{sec:induced:lemmas}, we will prove a number of fundamental properties of induced families; next, in Section~\ref{sec:tech:induced:lemmas}, we will prove a key technical lemma; in Section~\ref{sec:induced:induction:steps}, we will apply these results in the case $|W' \setminus W| = 1$, and in Section~\ref{sec:induced:proofs} we will complete the proofs of Lemmas~\ref{lem:good:faces} and~\ref{lem:semigood:faces}.

\subsection{Some simple lemmas about induced families}\label{sec:induced:lemmas}

We begin with a simple but key consequence of Lemma~\ref{lem:quasi-innerprod}. This lemma is one of the main motivations for the construction of the family $\QQ$ of quasistable directions. We remark that the lemma would be false without the assumption that $W \cup \{u\}$ is a clique. 

\begin{lemma}\label{lem:induced-surprise}
Let $W \subset \QQ$ and $u \in \QQ$ be such that $W' = W \cup \{u\}$ is a clique. Then
\begin{equation}\label{eq:induced-surprise}
X \subset \HH(W') \qquad \Leftrightarrow \qquad X \subset \HH(W) \quad \text{and} \quad X \cap W^\perp \subset \HH(u)
\end{equation} 
for each $X \in \U$.
\end{lemma}

\begin{proof}
Recall from~\eqref{def:HW} that $\HH(W)$ is the set of $x \in \Z^d$ such that $\< x,u \> \le 0$ for all $u \in W$. One direction of~\eqref{eq:induced-surprise} is therefore immediate, since if $X \subset \HH(W') = \HH(W) \cap \HH(u)$ then clearly $X \subset \HH(W)$ and $X \cap W^\perp \subset X \subset \HH(u)$. To prove the reverse implication, note that $uv \in E\big(\Vor(\QQ)\big)$ for all $v \in W \setminus \{u\}$, since $W \cup \{u\}$ is a clique. Note also that $X \subset \cL_R$. By Lemma~\ref{lem:quasi-innerprod}, it follows that  
\begin{equation}\label{eq:quasi-innerprod:surprise:app}
\< x,v \> < 0  \qquad \Rightarrow \qquad \< x,u \> \leq 0
\end{equation}
for every $x \in X$ and $v \in W$. Now, if $x \in X \subset \HH(W)$ and $x \notin W^\perp$, then there exists $v \in W$ such that $\< x,v \> < 0$. By~\eqref{eq:quasi-innerprod:surprise:app}, it follows that $x \in \HH(u)$, and hence if $X \cap W^\perp \subset \HH(u)$ then $X \subset \HH(u)$, as required.
\end{proof}

We will need the following simple consequence of Lemma~\ref{lem:induced-surprise} in Section~\ref{proof:sec}. 

\begin{lemma}\label{lem:emptyset:for:subfaces}
Let\/ $W \subset W' \subset \QQ$ be cliques. If\/ $\emptyset \in \U[W]$, then\/ $\emptyset \in \U[W']$.
\end{lemma}

\begin{proof}
It is clearly enough, by induction, to prove the case $| W' \setminus W | = 1$, and this case follows from applying Lemma~\ref{lem:induced-surprise}. Indeed, by Definition~\ref{def:induced}, the condition $\emptyset \in \U[W]$ implies that $X \cap W^\perp = \emptyset$ for some $X \in \U$ with $X \subset \HH(W)$, and therefore $X \cap W'^\perp = \emptyset$ and $X \subset \HH(W')$, the former since $W'^\perp \subset W^\perp$ and the latter by Lemma~\ref{lem:induced-surprise}.
\end{proof}

The following commutativity lemma is another immediate consequence of Lemma~\ref{lem:induced-surprise}. It allows us to form the update family induced by a clique $W$ by successively forming update families induced by the elements of $W$ in turn.

\begin{lemma}\label{lem:commutes}
Let $W \subset \QQ$ and $u \in \QQ$. If\/ $W \cup \{u\}$ is a clique, then
\begin{equation}\label{eq:inducedplusone}
\big( \U[W] \big) [u] = \U\big[ W \cup \{u\} \big].
\end{equation}
\end{lemma}

\begin{proof}
Set $W' = W \cup \{u\}$, and recall that 
$$\U[ W' ] = \big\{ X \cap W'^\perp : X \in \U \,\text{ and }\, X \subset \HH(W') \big\}.$$
Similarly, since $(W')^\perp = W^\perp \cap \{u\}^\perp$, we have 
$$\big( \U[W] \big) [u] = \big\{ X \cap W'^\perp : X \in \U, \,\, X \subset \HH(W) \,\text{ and }\, X \cap W^\perp \subset \HH(u) \big\}.$$
The lemma now follows immediately from Lemma~\ref{lem:induced-surprise}.
\end{proof}

Observe that the induced update family $\U[W]$ depends on the set $W$, and not just the subspace $W^\perp$. For example, $\U[W]$ will generally be different from $\U[-W]$, where $-W = \{ -w : w \in W\}$. Nevertheless, our next lemma shows that adding a new direction $u \in \<W\>$ to a clique $W \subset \QQ$ does not change the induced update family, as long as $u \in \QQ$ and $W \cup \{u\}$ is also a clique. This is a straightforward consequence of Lemmas~\ref{lem:quasi-innerprod} and~\ref{lem:commutes}, and will enable us to deal with this `degenerate' case. 

\begin{lemma}\label{lem:new-elt-in-span}
Let $W \subset \QQ$ and $u \in \QQ \cap \< W \>$. If\/ $W' = W \cup \{u\}$ is a clique, then 
$$\U[W] = \U[W'].$$
\end{lemma}

\begin{proof}
Since $W' = W \cup \{u\}$ is a clique, it follows from Definition~\ref{def:induced} and Lemma~\ref{lem:commutes} that 
$$\U[W'] = \big(\U[W]\big)[u] = \big\{ X \cap \{u\}^\perp : X \in \U[W] \,\text{ and }\, X \subset \HH(u) \big\}.$$
Now, since $u \in \< W \>$, it follows that $X \subset W^\perp \subset \{u\}^\perp$ for each $X \in \U[W]$, and therefore $X \cap \{u\}^\perp = X$ and $X \subset \HH(u)$, as required.
\end{proof}

For the next lemma, imagine two adjacent faces of our droplet, one corresponding to a set $W$, and the other corresponding to a direction $u$. The lemma says that if the projection of $u$ onto $W^\perp$ is unstable in $\U[W]$, then growth occurs trivially in the induced process corresponding to the boundary of the two faces. 


\begin{lemma}\label{lem:unstable:projection}
Let\/ $W \subset \QQ$ and\/ $u \in \QQ \setminus \<W\>$. If\/ $W' = W \cup \{u\}$ is a clique, then 
$$\pi(u,W^\perp) \notin \S_W \qquad \Rightarrow \qquad \emptyset \in \U[W'].$$
\end{lemma}

\begin{proof}
Since $\pi(u,W^\perp) \notin \S_W$, there exists a set $Y \in \U[W]$ such that $\<x, \pi(u,W^\perp) \> < 0$ for every $x \in Y$. Since $Y \subset W^\perp$, it follows that $\<x,u \> < 0$ for every $x \in Y$, by Observation~\ref{obs:projection:innerproduct}, and in particular that $Y \subset \HH(u)$. Let $X \in \U$ be such that $X \subset \HH(W)$ and $Y = X \cap W^\perp$. Observe that $X \subset \HH(W')$, by Lemma~\ref{lem:induced-surprise}, and that $X \cap W'^\perp = \emptyset$, since $\<x,u \> < 0$ for every $x \in Y$. It follows that $\emptyset \in \U[W']$, as claimed.
\end{proof}

\pagebreak

Finally, let us note one more simple fact, which we will use in Section~\ref{sec:tech:induced:lemmas}.

\begin{lemma}\label{lem:proj:induces:same}
Let $W \subset \QQ$, and set $\F = \U[W]$. If $u \in \QQ \setminus \< W \>$, then
$$\F[u] = \F[u'],$$
where $u' = \pi(u,W^\perp)$.
\end{lemma}

\begin{proof}
We need to show that
$$\big\{ X \cap \{u\}^\perp : X \in \F \,\text{ and }\, X \subset \HH(u) \big\} = \big\{ X \cap \{u'\}^\perp : X \in \F \,\text{ and }\, X \subset \HH(u') \big\}.$$
To prove this, recall that $X \subset W^\perp$ for every $X \in \F$ (since $\F = \U[W]$), and therefore 
$$\sgn\big( \< x,u \> \big) = \sgn\big( \< x,u' \> \big)$$
for every $x \in X \in \F$, by Observation~\ref{obs:projection:innerproduct}. It follows that $X \cap \{u\}^\perp = X \cap \{u'\}^\perp$, and $X \subset \HH(u)$ if and only if $X \subset \HH(u')$, as required.
\end{proof}

\subsection{A technical lemma about induced families}\label{sec:tech:induced:lemmas}

We are now ready to prove the following lemma, which will be our main tool in the proofs of Lemmas~\ref{lem:good:faces} and~\ref{lem:semigood:faces}. We will prove it by applying Lemma~\ref{lem:induced-projection} with $\F = \U[W]$, and using Lemma~\ref{lem:commutes}.

\begin{lemma}\label{lem:induced}
Let $W \subset \QQ$ and $u \in \QQ \setminus \< W \>$. If\/ $W' = W \cup \{u\}$ is a clique, then
\begin{equation}\label{eq:induced-families-iff}
v \in \S_W \qquad \Leftrightarrow \qquad \pi(v,W'^\perp) \in \S_{W'}
\end{equation}
for every $v \in S_\eta\big( \SS(W), \pi(u,W^\perp) \big)$. 
\end{lemma}

\begin{proof} 
Set $u' = \pi(u,W^\perp)$, and recall that $W \cup \{u\} \subset \QQ \subset \SS_\Q^{d-1}$. By Observation~\ref{obs:proj:is:rational}, it follows that $u' \in \SS_\Q^{d-1}$. Setting $W'' = W \cup \{u'\}$, and applying Lemma~\ref{lem:induced-projection} with $\F = \U[W]$, we obtain
\begin{equation}\label{eq:projection:app}
v \in \S_W  \qquad \Leftrightarrow \qquad v \in \S(\F)  \qquad \Leftrightarrow \qquad \pi(v,W''^\perp) \in \S\big( \F[u'] \big)
\end{equation}
for every $v \in S_\eta\big( \SS(W), u' \big)$, since $\S_W = \S(\F) \cap \SS(W)$ and $\< \F \> \subset W^\perp$, by~\eqref{eq:induced-def}. 

To prove~\eqref{eq:induced-families-iff}, we therefore need to show that
\begin{equation}\label{eq:induced-stp}
\pi(v,W''^\perp) \in \S\big( \F[u'] \big) \qquad \Leftrightarrow \qquad \pi(v,W'^\perp) \in \S_{W'}.
\end{equation}
To prove~\eqref{eq:induced-stp}, observe first that 
$$\big(W \cup \{u\}\big)^\perp = \big(W \cup \{u'\}\big)^\perp,$$
since if $x \in W^\perp$ then $\< x,u \> = 0$ if and only if $\< x,u' \> = 0$, by Observation~\ref{obs:projection:innerproduct}. It follows that $\pi(v,W'^\perp) = \pi(v,W''^\perp)$, and it will therefore suffice to show that
$$\S\big( \F[u'] \big) \cap \SS(W') = \S_{W'},$$
which, by the definition of $\S_{W'}$, will follow from the identity
\begin{equation}\label{eq:families-equal}
\F[u'] = \F[u] = \U[W']. 
\end{equation}
The first equality in~\eqref{eq:families-equal} holds by Lemma~\ref{lem:proj:induces:same}, and the second by Lemma~\ref{lem:commutes}, since $W' = W \cup \{u\}$ is a clique. This proves~\eqref{eq:induced-stp}, which in turn implies~\eqref{eq:induced-families-iff}, as required.
\end{proof}

\pagebreak

We will use the following immediate consequences of Lemma~\ref{lem:induced}. 

\begin{lemma}\label{cor:induced}
Let $W \subset \QQ$ and $u \in \QQ \setminus \< W \>$. If\/ $W' = W \cup \{u\}$ is a clique, then
\begin{equation}\label{induced:cor:r}
r^{k-1} \big( S_\eta\big( \SS(W), \pi(u,W^\perp) \big); \S_W \big) \, = \, r^{k-1} \big( \SS(W'); \S_{W'} \big),
\end{equation}
where $k = \dim(W^\perp)$, and
\begin{equation}\label{induced:cor:rho}
\rho^{k-2} \big( S_\eta\big(\SS(W),\pi(u,W^\perp)\big); \S_W, v \big) = \rho^{k-2}\big( \SS(W'); \S_{W'}, v' \big),
\end{equation}
for all $v \in S_\eta\big(\SS(W),\pi(u,W^\perp)\big)$, where $v' = \pi(v,W'^\perp)$.
\end{lemma}

\begin{proof}
To prove~\eqref{induced:cor:r}, it will suffice to show that
\[
\S_W \cap S_\eta\big( \SS(W), \pi(u,W^\perp) \big) \equiv \S_{W'}.
\]
To see this, simply consider the homothety $\varphi \colon S_\eta\big( \SS(W), \pi(u,W^\perp) \big) \to \SS(W')$ defined by $\varphi(v) = \pi(v,W'^\perp)$, and apply~\eqref{eq:induced-families-iff}. 
Since the homothety $\varphi$ that was used to prove the equivalence maps $v$ to $v'$, we also obtain~\eqref{induced:cor:rho}. 
\end{proof}

\subsection{The induction steps}\label{sec:induced:induction:steps}

We will prove Lemmas~\ref{lem:good:faces} and~\ref{lem:semigood:faces} by induction, adding vertices of $W' \setminus W$ one by one. In order to do so, we first prove three lemmas about the case $|W' \setminus W| = 1$; these lemmas will provide the induction steps. 

The first of these three lemmas deals with the case $s = 1$. It follows easily from Lemmas~\ref{lem:unstable:projection} and~\ref{lem:induced}, and will be used in the proofs of both Lemma~\ref{lem:good:faces} and Lemma~\ref{lem:semigood:faces}. 

\begin{lemma}\label{lem:1good}
Let\/ $W \subset \QQ$ and\/ $u \in \QQ \setminus \<W\>$ be such that\/ $W' = W \cup \{u\}$ is a clique.
Let~$w \in \SS(W)$ and suppose that\/ $W$ is $(1,w)$-semi-good. Then the following hold:
\begin{enumerate}
\item If\/ $\< u,w \> > 0$ then\/ $\emptyset \in \U[W']$.\vskip0.1cm
\item If\/ $\< u,w \> = 0$ then\/ $W'$ is $(1,w)$-semi-good.
\end{enumerate}
\end{lemma}

\begin{proof}
Recall from Definition~\ref{def:sgood} that, since $W$ is $(1,w)$-semi-good, we have
\begin{equation}\label{eq:ssemigood-prf}
\rho^{k-1} \big( \SS(W); \S_W, v \big) = 0
\end{equation}
and hence $v \notin \S_W$, for all $v \in \SS(W)$ such that $\< v, w \> > 0$, where $k = \dim(W^\perp)$. 

Suppose first that $\< u,w \> > 0$, and note that $\< \pi(u,W^\perp), w \> > 0$, by Observation~\ref{obs:projection:innerproduct}, since $w \in W^\perp$. By~\eqref{eq:ssemigood-prf} and since $\pi(u,W^\perp) \in \SS(W)$, it follows that $\pi(u,W^\perp) \notin \S_W$, and hence that $\emptyset \in \U[W']$, by Lemma~\ref{lem:unstable:projection}, as required.

For part~$(b)$,  observe first that if $w \in \SS(W)$ and $\< u,w \> = 0$, then $w \in \SS(W')$, so $\dim(W'^\perp) > 0$. We need to show that $W'$ is $(1,w)$-semi-good, which means that $v' \notin \S_{W'}$ for all $v' \in \SS(W')$ such that $\< v', w \> > 0$. To do so, observe first that, by Lemma~\ref{lem:induced}, we have
$$v \in \S_W \qquad \Leftrightarrow \qquad v' \in \S_{W'},$$
where $v \in S_\eta\big( \SS(W), \pi(u,W^\perp) \big)$ is such that $v' = \pi(v,W'^\perp)$. 

Now, by Observation~\ref{obs:projection:innerproduct}, we have $\< v,w \> > 0$ if and only if $\< v', w \> > 0$, since $w \in W'^\perp$. It follows that $v \in \SS(W)$ and $\< v,w \> > 0$ for any such $v' \in \SS(W')$. Hence, by~\eqref{eq:ssemigood-prf}, we have $v \notin \S_W$, and therefore $v' \notin \S_{W'}$, as required.
\end{proof}

We will next prove the following analogue of Lemma~\ref{lem:1good} for the case $s \ge 2$. The proof uses similar ideas, but is more complicated -- in particular, we will need to use the results of Section~\ref{sec:memphis}. We will use part $(b)$ to prove Lemma~\ref{lem:semigood:faces}, and part $(a)$ in Section~\ref{proof:sec}. 

\begin{lemma}\label{lem:sgood}
Let\/ $W \subset \QQ$ and\/ $u \in \QQ \setminus \<W\>$ be such that\/ $W' := W \cup \{u\}$ is a clique. Let~$w \in \SS(W)$ and $2 \le s \le k := \dim( W^\perp )$, and set $s^* := \min\{s,k-1\}$. If\/ $W$ is $(s,w)$-semi-good, then the following hold:
\begin{enumerate}
\item If\/ $\< u,w \> > 0$ then\/ $W'$ is $(s-1)$-good.\vskip0.1cm
\item If\/ $\< u,w \> = 0$ then\/ $W'$ is $(s^*,w)$-semi-good.
\end{enumerate}
\end{lemma}

\begin{proof}
Recall from Definition~\ref{def:sgood} that $W$ being $(s,w)$-semi-good means that
\begin{equation}\label{eq:semigood:app}
\rho^{k-1} \big( \SS(W); \S_W, v \big) \le s - 1
\end{equation}
for all $v \in \SS(W)$ such that $\< v, w \> > 0$. Since $\SS(W)$ is a $(k-1)$-dimensional sphere and $2 \le s \le k \le d$, it follows from Lemma~\ref{lem:memphis-subcrit} that
\begin{equation}\label{eq:sgood-memphis1}
\rho^{k-2} \Big( S_\eta\big( \SS(W), \pi(u,W^\perp) \big); \S_W, v \Big) \le k - 2
\end{equation}
and from Lemma~\ref{lem:memphis}, that 
\begin{equation}\label{eq:sgood-memphis2}
\rho^{k-2} \Big( S_\eta\big( \SS(W), \pi(u,W^\perp) \big); \S_W, v \Big) \le s - 1 
\end{equation}
for every $v \in S_\eta\big( \SS(W),\pi(u,W^\perp) \big)$ such that $\< v, w \> > 0$. 

For part $(a)$, observe that if $\pi(u,W^\perp) \in \S_W$, then by Lemma~\ref{cor:induced} and~\eqref{eq:rho} we have
\begin{align}
r^{k-1} \big( \SS(W'); \S_{W'} \big) & \, = \, r^{k-1} \big( S_\eta\big( \SS(W), \pi(u,W^\perp) \big); \S_W \big) \nonumber \\
& \, = \, \rho^{k-1} \big( \SS(W); \S_W, \pi(u,W^\perp) \big), \label{eq:sgood-res}
\end{align}
and similarly $r^{k-1} \big( \SS(W'); \S_{W'} \big) = 1$ if $\pi(u,W^\perp) \not\in \S_W$. Note that if $\< u,w \> > 0$, then $\< \pi(u,W^\perp),w \> > 0$, by Observation~\ref{obs:projection:innerproduct} and since $w \in W^\perp$, and therefore
$$r^{k-1} \big( \SS(W'); \S_{W'} \big) \le s - 1,$$
by~\eqref{eq:semigood:app} and~\eqref{eq:sgood-res}, and since $s \ge 2$. Hence, recalling~\eqref{eq:r} and Definition~\ref{def:sgood}, there exists $w' \in \SS(W')$ such that $W'$ is $(s-1,w')$-semi-good. 

In order to prove that $W'$ is $(s-1)$-good, it therefore suffices to show that 
\begin{equation}\label{eq:sgood:reminder}
\rho^{k-2}\big( \SS(W'); \S_{W'}, v' \big) \le \min\big\{ s - 1, k - 2 \big\}
\end{equation}
for all $v' \in \SS(W')$. To do so, observe first that, by Lemma~\ref{cor:induced}, 
\begin{equation}\label{eq:sgood-two-rhos}
\rho^{k-2}\big( \SS(W'); \S_{W'}, v' \big) = \rho^{k-2} \big( S_\eta\big(\SS(W),\pi(u,W^\perp)\big); \S_W, v \big),
\end{equation}
where $v \in S_\eta\big(\SS(W),\pi(u,W^\perp)\big)$ is such that $v'= \pi(v,W'^\perp)$. Note also that for every $v' \in \SS(W')$, there exists a unique $v \in S_\eta\big(\SS(W),\pi(u,W^\perp)\big)$ with $v' = \pi(v,W'^\perp)$.

To deduce~\eqref{eq:sgood:reminder} from~\eqref{eq:sgood-two-rhos}, recall from~\eqref{eq:sgood-memphis1} and~\eqref{eq:sgood-memphis2} that 
\begin{equation}\label{eq:sgood:memphis:synthesis}
\rho^{k-2} \Big( S_\eta\big( \SS(W), \pi(u,W^\perp) \big); \S_W, v \Big) \le \min\big\{ s - 1, k - 2 \big\}
\end{equation}
for every $v \in S_\eta\big( \SS(W), \pi(u,W^\perp) \big)$ such that $\< v, w \> > 0$. Now, since $\< u,w \> > 0$ and $w \in W^\perp$, we have $\< \pi(u,W^\perp),w \> > 0$, by Observation~\ref{obs:projection:innerproduct}, and hence $\< v, w \> > 0$ for every $v \in S_\eta\big( \SS(W), \pi(u,W^\perp) \big)$, because $\eta$ is sufficiently small. Combining~\eqref{eq:sgood-two-rhos} and~\eqref{eq:sgood:memphis:synthesis}, we obtain~\eqref{eq:sgood:reminder}. As noted above, this implies that $W'$ is $(s-1)$-good, as required. 

For part~$(b)$, we need to show that if $\< u,w \> = 0$, then $W'$ is $(s^*,w)$-semi-good, where $s^* = \min\{s,k-1\}$. That is, we need to show that
\begin{equation}\label{eq:ssemigood:yetagain}
\rho^{k-2} \big( \SS(W'); \S_{W'}, v' \big) \le s^* - 1 = \min\big\{ s - 1, k - 2 \big\}
\end{equation}
for every $v' \in \SS(W')$ such that $\< v', w \> > 0$. Observe first that~\eqref{eq:sgood-two-rhos} and~\eqref{eq:sgood:memphis:synthesis} also hold in this case, that is, by Lemma~\ref{cor:induced}, we have
\begin{equation}\label{eq:sgood-two-rhos:again}
\rho^{k-2}\big( \SS(W'); \S_{W'}, v' \big) = \rho^{k-2} \big( S_\eta\big(\SS(W),\pi(u,W^\perp)\big); \S_W, v \big),
\end{equation}
where $v \in S_\eta\big(\SS(W),\pi(u,W^\perp)\big)$ is such that $v'= \pi(v,W'^\perp)$, and by~\eqref{eq:sgood-memphis1} and~\eqref{eq:sgood-memphis2},  
\begin{equation}\label{eq:sgood:memphis:synthesis:again}
\rho^{k-2} \Big( S_\eta\big( \SS(W), \pi(u,W^\perp) \big); \S_W, v \Big) \le \min\big\{ s - 1, k - 2 \big\}
\end{equation}
for every $v \in S_\eta\big( \SS(W), \pi(u,W^\perp) \big)$ such that $\< v, w \> > 0$. In particular, if $\< v', w \> > 0$ then~\eqref{eq:sgood:memphis:synthesis:again} holds for the (unique) $v \in S_\eta\big(\SS(W),\pi(u,W^\perp)\big)$ such that $v'= \pi(v,W'^\perp)$, by Observation~\ref{obs:projection:innerproduct}, since $w \in W^\perp$ and $\<u,w \> = 0$, so $w \in W'^\perp$. Thus, combining~\eqref{eq:sgood-two-rhos:again} and~\eqref{eq:sgood:memphis:synthesis:again}, we obtain~\eqref{eq:ssemigood:yetagain}, as required. This completes the proof of the lemma.
\end{proof}

In order to prove Lemma~\ref{lem:good:faces}, we will need the following consequence of Lemma~\ref{lem:sgood}.

\begin{lemma}\label{lem:sgood:extra}
Let\/ $W \subset \QQ$ and\/ $u \in \QQ \setminus \< W \>$ be such that\/ $W' = W \cup \{u\}$ is a clique. Let $s \le k = \dim( W^\perp )$ with $s^* = \min\{s,k-1\} \ge 1$. If\/ $W$ is $s$-good, then\/ $W'$ is $s^*$-good.
\end{lemma}

\begin{proof}
Note that if $W$ is $s$-good, then it is $(s^*+1,w)$-semi-good for every $w \in \SS(W)$. Since $u \not\in \< W \>$, it follows that $W$ is $(s^*+1,w)$-semi-good for some $w \in \SS(W)$ with $\<u,w\> > 0$. Noting that $2 \le s^* + 1 \le k$, it follows by part~$(a)$ of Lemma~\ref{lem:sgood} that $W'$ is $s^*$-good, as required.
\end{proof}

\subsection{The proofs of Lemmas~\ref{lem:good:faces} and~\ref{lem:semigood:faces}}\label{sec:induced:proofs}

We are finally ready to prove the two main lemmas of the section. The two proofs are similar to one another; the first uses Lemmas~\ref{lem:new-elt-in-span},~\ref{lem:1good} and~\ref{lem:sgood:extra}, while the second uses Lemmas~\ref{lem:new-elt-in-span},~\ref{lem:1good} and~\ref{lem:sgood}.

\begin{proof}[Proof of Lemma~\ref{lem:good:faces}]
Let $W' \setminus W = \{u_1,\dots,u_\ell\}$, and note that if $\ell = 0$ then there is nothing to prove.  For each $0 \le i \le \ell$, set 
$$W_i := W \cup \{u_1,\dots,u_i\}, \qquad k_i := \dim(W_i^\perp) \qquad \text{and} \qquad s_i := \min\{s,k_i\}.$$ 
We claim that either $k_i > 0$ and $W_i$ is $s_i$-good, or $k_i = 0$ and $\emptyset \in \U[W_i]$.

The proof is by induction on $i$; note that the case $i = 0$ follows from our assumptions, since $W_0 = W$ is $s$-good and $k_0 = k \ge 1$. So let $1 \le i \le \ell$, and assume that the claim holds for $i - 1$. Suppose first that $u_i \in \< W_{i-1} \>$. In this case $W_i^\perp = W_{i-1}^\perp$, and thus $k_i = k_{i-1}$, and moreover $\U[W_i] = \U[W_{i-1}]$ by Lemma~\ref{lem:new-elt-in-span}, so the claim follows from the induction hypothesis. 

\pagebreak

Let us therefore assume that $u_i \notin \< W_{i-1} \>$, and note that, in particular, this implies that $k_i = k_{i-1} - 1$, and therefore $s_i = \min\{s_{i-1},k_{i-1}-1\}$. If $s_i \ge 1$, then it follows from Lemma~\ref{lem:sgood:extra} that $W_i$ is $s_i$-good, which completes the induction step if $k_{i-1} \ge 2$. 

If $k_{i-1} = 1$, on the other hand, then by the induction hypothesis we have $W_{i-1}^\perp = \<w\>$ for some $w \in \SS^{d-1}$, and $W_{i-1}$ is $1$-good. By Definition~\ref{def:sgood}, it follows that $W_{i-1}$ is both $(1,w)$-semi-good and $(1,-w)$-semi-good. Moreover, since $u_i \notin \< W_{i-1} \>$, we have $\< u_i,w \> \ne 0$. Hence, applying Lemma~\ref{lem:1good} (for either $w$ or $-w$, depending on the sign of $\< u_i,w \>$), it follows that $\emptyset \in \U[W_i]$, as claimed.

This proves the induction step, and since $W_\ell = W'$, the lemma follows.
\end{proof}

A slight modification of the argument above gives our second main lemma. 

\begin{proof}[Proof of Lemma~\ref{lem:semigood:faces}]
Let $W' \setminus W = \{u_1,\dots,u_\ell\}$, and note that $\ell > 0$, by assumption.  For each $0 \le i \le \ell$, set 
$$W_i := W \cup \{u_1,\dots,u_i\}, \qquad k_i := \dim(W_i^\perp) \qquad \text{and} \qquad s_i := \min\{s,k_i\}.$$ 
Note that $k_i \ge k' \ge 1$, since $W' \subset \{w\}^\perp$. We claim that $W_i$ is $(s_i,w)$-semi-good.

The proof is by induction on $i$; note that the case $i = 0$ follows from our assumptions, since $W_0 = W$ is $(s,w)$-semi-good. So let $1 \le i \le \ell$, and assume that the claim holds for $i - 1$. Suppose first that $u_i \in \< W_{i-1} \>$, so  $W_i^\perp = W_{i-1}^\perp$ and thus $k_i = k_{i-1}$. Since $\U[W_i] = \U[W_{i-1}]$, by Lemma~\ref{lem:new-elt-in-span}, the claim follows from the induction hypothesis. 

Let us therefore assume that $u_i \notin \< W_{i-1} \>$, and note that, in particular, this implies that $k_i = k_{i-1} - 1$, and therefore $s_i = \min\{s_{i-1},k_{i-1}-1\}$. Recall that $W' \subset \{w\}^\perp$, and note that therefore $\<u_i,w\> = 0$. By Lemma~\ref{lem:sgood}, it follows that if $s_{i-1} \ge 2$, then $W_i$ is $(s_i,w)$-semi-good, as required. If $s_{i-1} = 1$, on the other hand, then $W_{i-1}$ is $(1,w)$-semi-good, so by Lemma~\ref{lem:1good} (and since  $\<u_i,w\> = 0$) $W_i$ is $(1,w)$-semi-good, as claimed.

This proves the induction step, and since $W_\ell = W'$, the lemma follows.
\end{proof}

This concludes the first half of the paper: we have now proved all of the crucial properties of the stable sets of induced processes, and are ready to begin the process of constructing paths of infections. The central objects in our construction will be a certain family of polytopes, which will be defined in the next section.

\section{Polytopes}\label{sec:polytopes}

In this section we will introduce the family of polytopes that we shall use to prove Theorem~\ref{thm:upper}. The first step is to define two simpler families of `canonical' polytopes: `spherical' polytopes (see Section~\ref{spherical:sec}) and `tubular' polytopes (see Section~\ref{tubular:sec}). We will also state a few simple but fundamental properties of these polytopes; since the proofs are standard, but somewhat technical, we postpone most of them to Appendix~\ref{polytope:app}.

Recall that $\QQ$ is a fixed (finite) set of quasistable directions for range $R = R(\U)$, so $\QQ$~intersects every open hemisphere of $\SS^{d-1}$, and~\eqref{eq:quasi:property:repeat} holds for each $u \in \QQ$ and $w \in \cL_R$. 

\subsection{Spherical polytopes}\label{spherical:sec}

The most basic polytope that we shall study is
\begin{equation}\label{def:P}
P(\emptyset) := \bigcap_{u \in \QQ} \big\{ x \in \R^d : \< x,u \> \leq 1 \big\},
\end{equation}
Observe that $P(\emptyset)$ is a bounded (and hence compact) $d$-dimensional polytope, because $\QQ$ intersects every open hemisphere of $\SS^{d-1}$. Now, for each set $W \subset \QQ$, define
\begin{equation}\label{eq:P(W)}
P(W) := P(\emptyset) \cap \bigcap_{u \in W} \big\{ x \in \R^d : \< x,u \> = 1 \big\}.
\end{equation}
Note that if $P(W)$ is non-empty, then it is a face of $P$. We will show (see Lemma~\ref{lem:faces-to-cliques}) that if $P(W) \neq \emptyset$, then $W$ is a clique\footnote{The converse is (unfortunately) false, and $\QQ$ may contain cliques corresponding to empty faces. However, we will show in Section~\ref{maximal:clique:sec} that there is a natural family of cliques, every member of which corresponds to a non-empty face, see Definition~\ref{def:curlyW} and Lemma~\ref{lem:clique-of-correct-dim}.} in the Voronoi graph $\Vor(\QQ)$ (see Definition~\ref{def:voronoi}). For each (possibly empty) clique $W \subset \QQ$, let us write
\begin{equation}\label{eq:QQW}
N_\QQ(W) := \big\{ u \in \QQ \setminus W \,:\, W \cup \{u\} \text{ is a clique} \big\}
\end{equation}
for the set of common neighbours of $W$ in $\Vor(\QQ)$. 

Our first aim is to prove the following lemma.

\begin{lemma}\label{lem:P-nbrs}
Let $W \subset \QQ$, and suppose that $P(W) \neq \emptyset$. Then $W$ is a clique, and 
\begin{equation}\label{eq:P-nbrs}
P(W) = \bigcap_{u \in W} \big\{ x \in \R^d : \< x,u \> = 1 \big\} \cap \bigcap_{u \in N_\QQ(W)} \big\{ x \in \R^d : \< x,u \> \leq 1 \big\}.
\end{equation}
\end{lemma}

We begin with the following simple observation.

\begin{lemma}\label{lem:facets-to-cells}
Let $u \in \QQ$ and $x \in \R^d$ be such that $\< x,u \> = 1$. Then 
$$x \in P(\emptyset) \qquad \Leftrightarrow \qquad x / \| x \| \in \Cell_\QQ(u).$$
\end{lemma}

\begin{proof}
Recall from Definition~\ref{def:voronoi} that if $w \in \SS^{d-1}$ then $w \in \Cell_\QQ(u)$ if and only if $\< w,u \> \ge \< w,v \>$ for every $v \in \QQ$. It follows that $x / \| x \| \in \Cell_\QQ(u)$ if and only if $\< x,v \> \le \< x,u \> = 1$ for every $v \in \QQ$, and by~\eqref{def:P} this holds if and only if $x \in P(\emptyset)$, as required.
\end{proof}

We next use Lemma~\ref{lem:facets-to-cells} to show that faces of $P(\emptyset)$ correspond to cliques in $\QQ$.

\begin{lemma}\label{lem:faces-to-cliques}
If\/ $W \subset \QQ$ and $P(W) \neq \emptyset$, then $W$ is a clique. 
\end{lemma}

\begin{proof}
Recall from~\eqref{eq:P(W)} that if $x \in P(W)$, then $x \in P(\emptyset)$ and $\< x,u \> = 1$ for every $u \in W$, and therefore $x / \| x \| \in \bigcap_{u \in W} \Cell_\QQ(u)$, by Lemma~\ref{lem:facets-to-cells}. Now, by Definition~\ref{def:voronoi}, the set $\Cell_\QQ(u) \cap \Cell_\QQ(v)$ is non-empty, for distinct $u,v \in \QQ$, if and only if $uv \in E\big(\Vor(\QQ)\big)$. It follows that $\bigcap_{u \in W} \Cell_\QQ(u)$ being non-empty implies that $W$ is a clique.
\end{proof}

We can now easily deduce Lemma~\ref{lem:P-nbrs}. 

\begin{proof}[Proof of Lemma~\ref{lem:P-nbrs}]
Let us write $P'(W)$ for the right-hand side of~\eqref{eq:P-nbrs}. Note first that $P(W) \subset P'(W)$, by~\eqref{eq:P(W)} and since $N_\QQ(W) \subset \QQ$. We therefore need to show that if $x \in P'(W)$ and $P(W) \neq \emptyset$ then $x \in P(W)$. 

Let $x \in P'(W)$ and $y \in P(W)$, and (recalling from~\eqref{eq:P(W)} that $P(W)$ is compact) let $\lambda \ge 0$ be maximal such that $z := y + \lambda(x-y) \in P(W)$. If $\lambda \ge 1$ then, since $P(W)$ is convex, it follows that $x \in P(W)$, as required.  On the other hand, if $\lambda < 1$ then there exists $v \in \QQ$ such that $\< z,v \> = 1$ and $\< x,v \> > 1$. It follows that $z \in P(W \cup \{v\})$, and therefore $W \cup \{v\}$ is a clique, by Lemma~\ref{lem:faces-to-cliques}. However, since $x \in P'(W)$ and $\< x,v \> > 1$, we must have $v \not\in W \cup N_\QQ(W)$, so this is a contradiction. 
\end{proof}

\subsection{Tubular polytopes}\label{tubular:sec}

When working in a lattice $\L(W)$, for some $W \subset \QQ$, typically we shall only know that our droplets are likely to grow in a certain direction $w \in \SS(W)$ (cf.~Definition~\ref{def:sgood}). In this situation, we will not be able to control the growth of the infected set using polytopes of the form $a + t \cdot P(W)$; instead, we shall need to use a different family of polytopes, that are formed by `stretching' $P(W)$ in direction $w$. 

To define these polytopes, set $\delta(u,w) := \1[ \< u,w \> \geq 0 ]$. Now, given $w \in \R^d$, let\footnote{We emphasize that in this definition $w$ is not assumed to be a unit vector.}
\begin{equation}\label{eq:P(emptyset,w)}
P(\emptyset,w) := \bigcap_{u \in \QQ} \big\{ x \in \R^d : \big\< x - \delta(u,w) w, \, u \big\> \leq 1 \big\},
\end{equation}
and for each set $W \subset \QQ$, let
\begin{equation}\label{eq:P(W,w)}
P(W, w) := P(\emptyset,w) \cap \bigcap_{u \in W} \big\{ x \in \R^d : \big\< x - \delta(u,w) w, \, u \big\> = 1 \big\}.
\end{equation}
Recall from~\eqref{def:LR} that
$$\cL_R := \big\{ w \in \R^d : w \in \<x\> \text{ for some } x \in \Z^d \text{ with } \|x\| \le R \big\},$$
The following is the tubular analogue of Lemma~\ref{lem:P-nbrs}.

\begin{lemma}\label{lem:tube-nbrs}
Let $W \subset \QQ$ and $w \in \cL_R$. If\/ $P(W,w) \neq \emptyset$, then $W$ is a clique, and
\begin{align*}
& P(W,w) = \bigcap_{u \in W} \Big\{ x \in \R^d : \big\< x - \delta(u,w) w, \, u\big\> = 1 \Big\} \nonumber \\ 
& \hspace{5cm} \cap \bigcap_{u \in N_\QQ(W)} \Big\{ x \in \R^d : \big\< x - \delta(u,w) w, \, u \big\> \le 1 \Big\}.
\end{align*}
\end{lemma}

The proof of Lemma~\ref{lem:tube-nbrs} will be given in the appendix 
(see Lemma~\ref{lem:tube-nbrs:app}). However, we shall state here two of the lemmas used in the proof, since they will also be needed later on. The first of these two lemmas is less obvious than it looks; in particular, it requires the full power of the definition of $\QQ$. To highlight this, we give a sketch of the proof here; the full details are given in the appendix (see Lemma~\ref{lem:tube-is-tube:app}). 

\begin{lemma}\label{lem:tube-is-tube}
Let $W \subset \QQ$ and let $w \in \cL_R$. If $w \in W^\perp$, then
\begin{equation}\label{eq:tube-is-tube-emptyset}
P(W,w) = \bigcup_{\lambda \in [0,1]} \big( P(W) + \lambda w \big).
\end{equation}
\end{lemma}

\begin{proof}[Sketch proof]
It suffices to prove the lemma in the case $W = \emptyset$, since the general statement follows by intersecting with the set $\bigcap_{u \in W} \big\{ x \in \R^d : \< x,u \> = 1 \big\}$ (using the assumption that $w \in W^\perp$). It is straightforward to verify that $P(\emptyset) + \lambda w \subset P(\emptyset,w)$ for each $\lambda \in [0,1]$, so we shall concentrate on the other inclusion; that is, showing that if $x \in P(\emptyset,w)$, then $x \in P(\emptyset) + \lambda w$ for some $\lambda \in [0,1]$. We consider here only the (most interesting) case $0 \le \<x,w\> \le \|w\|^2$; the proof in the other cases is similar, but simpler.

Let $z \in \{w\}^\perp$ be such that $x = \lambda w + z$, and note that $\lambda = \<x,w\> / \|w\|^2 \in [0,1]$, and that we may assume that $z \ne \0$, because $\lambda w \in P(\emptyset) + \lambda w$. Let $u \in \QQ$ be such that $z / \|z\| \in \Cell_\QQ(u)$. Since $w \in \cL_R$, it follows by~\eqref{eq:quasi:property:repeat} that $\< w,u \> = 0$, and therefore $\< x,u \> = \< z,u \>$. Since $x \in P(\emptyset,w)$, and recalling~\eqref{eq:P(emptyset,w)}, it follows that $\< z,u \> \le 1$. 

Now, observe that $\< z,u \> > 0$, since $\QQ$ intersects every open hemisphere of $\SS^{d-1}$. By Lemma~\ref{lem:facets-to-cells}, it follows that $z / \< z,u \> \in P(\emptyset)$, and hence $z \in P(\emptyset)$, since $P(\emptyset)$ is convex and $\0 \in P(\emptyset)$. Since $x = \lambda w + z$, this implies that $x \in P(\emptyset) + \lambda w$, as required.
\end{proof}

When $w \not\in W^\perp$, on the other hand, $P(W,w)$ is rather less interesting: it is just a translate of $P(W)$. The proof of the following lemma is not especially enlightening, so we defer the details to the appendix (see Lemma~\ref{lem:tube-faces-to-P-faces:app}).

\begin{lemma}\label{lem:tube-faces-to-P-faces}
Let $W \subset \QQ$ and $w \in \cL_R$, and suppose that $P(W) \neq \emptyset$ and $w \notin W^\perp$. If $\< u,w \> > 0$  for some $u \in W$, then 
$$P(W,w) = P(W) + w,$$ 
and otherwise $P(W,w) = P(W)$.
\end{lemma}

\subsection{Maximal cliques, and the dimension of a face}\label{maximal:clique:sec}

When working with the polytopes $P(W)$ and $P(W,w)$, it will often be important to know not only that they are non-empty, but that they have dimension $\dim(W^\perp)$. We next define a family of cliques for which we shall be able to prove that this is indeed the case; this family will play an important role in Sections~\ref{fundamental:sec}--\ref{proof:sec}. 

\begin{definition}\label{def:curlyW}
We define the set of \emph{maximal cliques} in $\Vor(\QQ)$ to be
$$\W := \big\{ W \subset \QQ : \text{$W$ is a clique and $P(W') \neq P(W)$ for every $W \subsetneq W' \subset \QQ$} \big\}.$$
We shall also write $\W_k$ for the set of $W \in \W$ such that $\dim(W^\perp) = k$. 
\end{definition}
 
Note that if $P(W)$ is non-empty, then there exists $W \subset W' \in \W$ with $P(W') = P(W)$. In the appendix we shall prove the following stronger statement (see Lemma~\ref{lem:exists:WinW:trivial:app}).

\begin{lemma}\label{lem:exists:WinW:trivial}
Let $W \subset \QQ$ with $P(W) \ne \emptyset$. There exists $W' \in \W$ with 
$$W \subset W' \subset W \cup N_\QQ(W) \qquad \text{and} \qquad P(W,w) = P(W',w)$$
for every $w \in \cL_R$. 
\end{lemma}

As mentioned above, the crucial property of cliques $W \in \W$ is that the polytopes $P(W,w)$ have the same dimension as the space $W^\perp$. Let us write $\aff(X)$ for the affine span of a set $X \subset \R^d$. The following lemma is proved in the appendix (see Lemma~\ref{lem:clique-of-correct-dim:app}).

\begin{lemma}\label{lem:clique-of-correct-dim}
Let $W \in \W$ and $w \in \cL_R$. Then $P(W,w) \neq \emptyset$ and 
$$\dim\big(\aff\big( P(W,w) \big) \big) = \dim(W^\perp).$$
\end{lemma}

Lemma~\ref{lem:clique-of-correct-dim} has the following useful consequence.

\begin{lemma}\label{lem:faces:lower:dim}
Let $W \in \W$. If\/ $W \subsetneq W' \in \W$, then 
$$\dim\big( W'^\perp \big) < \dim(W^\perp).$$
\end{lemma}

\begin{proof}
Since $W \in \W$, we have $P(W') \ne P(W)$. But $P(W')$ is formed by intersecting $P(W)$ with some hyperplanes, so $P(W') \ne P(W)$ implies that the dimension of the affine span of $P(W')$ is strictly less than that of the affine span of $P(W)$. By Lemma~\ref{lem:clique-of-correct-dim} (applied with $w = \0$), the claim follows. 
\end{proof}

\subsection{Our family of polytopes}\label{our:family:sec}

We are now ready to introduce the family of polytopes that we shall work with throughout the rest of the proof of Theorem~\ref{thm:upper}. Many of the lemmas that we will prove about the polytopes in this family actually hold in much greater generality, but it will be convenient (in particular, to simplify the notation, and our induction hypothesis) to restrict our attention to this family. 

\begin{definition}\label{def:cPW}
For each $W \subset \QQ$, define\footnote{Recall from~\eqref{def:constants} that $C = C(\QQ) > 0$ is a sufficiently large constant. We require $t > C$ (rather than $t \ge C$) so that the closed interior of a polytope $P \in \cP(W)$ is also in $\cP(W)$, see Definition~\ref{lem:interiors-are-canonical}. We require $a \in \Z^d + W^\perp$ so that the lattice $a + \L(W)$ is non-empty (and therefore has dimension $\dim(W^\perp)$).} 
$$\cP(W) := \Big\{ (W,w,a,t,\tau) : w \in \cL_R \cap \SS^{d-1}, \; a \in \Z^d + W^\perp, \; t > C \,\text{ and } \, \tau \ge 0 \Big\},$$
and for each quintuple $(W,w,a,t,\tau) \in \cP(W)$, define a polytope
\begin{equation}\label{def:PWattau}
P(W,w;a,t,\tau) := a + t \cdot P\big( W, (\tau/t)w \big).
\end{equation}
\end{definition}

Abusing notation slightly, we write $P \in \cP(W)$ to mean that $P = P(W,w;a,t,\tau)$ for some $(W,w,a,t,\tau) \in \cP(W)$. We moreover write $w(P)$, $a(P)$, $t(P)$ and $\tau(P)$ for the corresponding elements of the quintuple associated with $P$, and define
\begin{equation}\label{def:cPWw}
\cP(W,w) := \big\{ P \in \cP(W) : w(P) = w \big\},
\end{equation}
\begin{equation}\label{def:cPWwt}
\cP(W,w;t) := \big\{ P \in \cP(W,w) : t(P) = t \big\},
\end{equation}
and
\begin{equation}\label{def:cPWwtt}
\cP(W,w;t,\tau) := \big\{ P \in \cP(W,w;t) : \tau(P) = \tau \big\}.
\end{equation}
If two polytopes $P,Q \in \cP(W)$ are equal as subsets of $\R^d$, then (abusing notation further) we shall sometimes write $P = Q$, and we trust that this will not cause confusion. For example, observe that if $P = P(W,w;a,t,\tau)$ with $w \in W^\perp$, and
\begin{equation}\label{def:Pminus}
P^- := P\big( W,-w;a+\tau w, t, \tau \big),
\end{equation}
then $P = P^-$, by Lemma~\ref{lem:tube-is-tube} and~\eqref{def:PWattau}. The following consequence of Lemmas~\ref{lem:tube-nbrs} and~\ref{lem:clique-of-correct-dim} will be useful, and is proved in Lemma~\ref{lem:P(W,w,a,t,tau):def2:app}. 

\begin{lemma}\label{lem:P(W,w,a,t,tau):def2}
If\/ $W \in \W$ and $P = P(W,w;a,t,\tau) \in \cP(W)$, then 
\begin{align}
& P = \bigcap_{u \in W} \Big\{ x \in \R^d : \big\< x - a - \delta(u,w) \tau w, \, u \big\> = t \Big\} \nonumber \\ 
& \hspace{4cm} \cap \bigcap_{u \in N_\QQ(W)} \Big\{ x \in \R^d : \big\< x - a - \delta(u,w) \tau w, \, u \big\> \le t \Big\}.\label{eq:P(W,w,a,t,tau):def2}
\end{align}
\end{lemma}

Finally, if $P = P(W,w;a,t,\tau) \in \cP(W)$ and $W \subset W' \subset \QQ$, then let us write
\begin{equation}\label{def:Delta:face}
\Delta(P,W') := P(W',w;a,t,\tau)
\end{equation}
for the $W'$-face of $P$. Note, in particular, that $\Delta(P,W') \subset P$.

\section{The bootstrap process in a polytope}\label{fundamental:sec}

In this section we prove two fundamental lemmas about the bootstrap process in a polytope. The first provides a connection between different induced processes, and will allow us, in Section~\ref{sec:deterministic}, to prove one of our key deterministic results, Lemma~\ref{lem:surface}. Recall from~\eqref{def:constants} that $C = C(\QQ) > 0$ is a sufficiently large constant.

\begin{lemma}\label{cor:faces-far-away}
Let\/ $W \subset W' \in \W$, and let\/ $X' \in \U[W']$. There exists $X \in \U[W]$ such that if $P \in \cP(W)$ and $x \in \Delta(P,W')$, then
$$x + X' \subset \Delta(P,W') \quad \Rightarrow \quad x + X \subset P.$$
\end{lemma}

Our second fundamental lemma will allow us to complete the infection of a polytope once it is `almost' entirely infected. In order to state this lemma, we will need to define what it means for a polytope $P$ to be `internally filled' in the $\U[W]$-process; this notion will also play a key role in Sections~\ref{sec:deterministic} and~\ref{proof:sec}, including in the precise statement of our induction hypothesis (see Definition~\ref{def:ih}). We remark that the definition we introduce here is slightly different from the definition of `$\U$-internally filled' in~\cite{BSU,BDMS}. 

Roughly speaking, $P$ is internally filled in the $\U[W]$-process if every point of $P \cap \Z^d$ is infected by the $\U[W]$-process with initial set $P \cap A$. However, there is an important additional condition, which is that the growth under $\U[W]$ of the initial infection $P \cap A$ is constrained to take place \emph{inside $P$} (in particular, routes to the full infection of $P$ that pass via the infection of sites outside $P$ and then back inside $P$ do not count). The reason for imposing this constraint is that $P$ itself will typically be a face of a higher dimensional polytope, and so the use of the induced update family $\U[W]$ is only valid inside $P$.

\begin{definition}\label{def:int:filled}
Let $W \in \W$ and $P \subset \R^d$. Given $B \subset P$, define $B_0 := B \cap \Z^d$ and
\[
B_{t+1} := B_t \cup \big\{ x \in P \cap \Z^d \,:\, x + X \subset B_t \,\text{ for some }\, X \in \U[W] \big\}
\]
for each $t \ge 0$. We write $[B]^P_{\U[W]} := \bigcup_{t\geq 0} B_t$ for the set of eventually infected sites in this restricted process. We say that $P$ is \emph{internally filled} by $A$ in the $\U[W]$-process if
\[
\big[ P \cap A \big]_{\U[W]}^P = P \cap \Z^d,
\]
and write $I_W^\bullet(P)$ for the event that $P$ is internally filled by $A$ in the $\U[W]$-process. 
\end{definition}

We also need the following definition from~\cite{BDMS,BBMSlower}: a \emph{strongly connected component} of a finite set $K \subset \R^d$ is a connected component of the graph $G$ with vertex set $K$, and edge set 
\begin{equation}\label{def:strongly:conn}
E(G) = \big\{  xy : \|x-y\| \le 2R_0 \big\},
\end{equation}
where $R_0 = R_0(\U)$ is the radius of $\U$, see~\eqref{def:radius}. Recall Definition~\ref{def:sgood}, and also from~\eqref{def:constants} that $\delta = \delta(\QQ) > 0$ is a sufficiently small constant (chosen so that Lemma~\ref{lem:close-to-many-faces} holds).

We can now state the second main lemma of this section.

\begin{lemma}\label{lem:eating:a:set}
Let $W \in \W$, set $k := \dim(W^\perp)$, and suppose that $W$ is $k$-good. Let $P \in \cP(W)$ and let $K \subset P \cap \Z^d$. If every strongly connected component of $K$ has diameter at most $\delta \cdot t(P)$, then $\big[ P \setminus K \big]^P_{\U[W]} = P \cap \Z^d$.
\end{lemma}

\pagebreak

Lemma~\ref{lem:eating:a:set} will play an important role in Section~\ref{proof:sec}, where we shall use it (as part of an adaptation of the `Schonmann trick' from~\cite{Sch1}) in order to obtain an exponential failure probability. As in~\cite{Sch1}, this trick will play a key role in the induction argument.

\subsection{The distance between faces of a polytope}

In the proofs of Lemmas~\ref{cor:faces-far-away} and~\ref{lem:eating:a:set} we shall need to control the distance between non-adjacent faces of our polytopes. We will next define a constant $\gamma = \gamma(\QQ)$ that allows us to do so, and state some facts involving this constant. The proofs of these facts are deferred to Appendix~\ref{gamma:app}. 

\begin{definition}\label{def:gamma}
Define
$$\gamma = \gamma(\QQ) := \min \Big\{ D(W,u) \,:\, W \subset \QQ, \; P(W) \neq \emptyset, \; u \in \QQ \; \text{ and }\, P(W \cup \{u\}) = \emptyset \Big\},$$
where
$$D(W,u) := \min \big\{ \| x - y \| \,:\, x \in P(W) \,\text{ and }\, \< y,u \> = 1 \big\},$$
which is well-defined because $P(W)$ is compact and $\big\{ y \in \R^d \,:\, \< y,u \> = 1 \big\}$ is closed.
\end{definition}

\begin{remark}\label{rmk:gamma}
Observe that $\gamma > 0$, since if $D(W,u) = 0$ for some $u \in \QQ$ and $W \subset \QQ$, then there exists $x \in P(W)$ with $\< x,u \> = 1$, which implies that $x \in P(W \cup \{u\})$. \end{remark}

In the proof of Lemma~\ref{cor:faces-far-away} we shall use the following lemma. When reading the statement of the lemma, one should imagine that $x$ is a vertex (of some rescaled lattice) that we wish to infect on the face $P(W \cup \{u\},w)$ of the polytope $P(W,w)$, and that $y$ is a (rescaled) element of some rule $X \in \U[W]$.

\begin{lemma}\label{lem:faces-far-away}
Let $W \subset \QQ$ and $w \in \cL_R$, and let $u \in N_\QQ(W)$. Let $x \in P(W \cup \{u\},w)$, and suppose that $y \in W^\perp$ is such that $\| y \| \le \gamma$ and 
\begin{equation}\label{eq:faces-far-away:condition}
\< y,v \> \le 0 \qquad \text{for every} \qquad v \in \{u\} \cup N_\QQ(W \cup \{u\}).
\end{equation}
Then $x + y \in P(W,w)$. 
\end{lemma}

The proof of Lemma~\ref{lem:faces-far-away} uses Lemmas~\ref{lem:faces-to-cliques},~\ref{lem:tube-nbrs},~\ref{lem:tube-is-tube} and~\ref{lem:tube-faces-to-P-faces}, see Appendix~\ref{gamma:app}. 

For the proof of Lemma~\ref{lem:eating:a:set}, we shall use the following two lemmas. The first is an immediate consequence of Definition~\ref{def:gamma}. 

\begin{lemma}\label{lem:far:from:nonfaces}
Let\/ $W \subset \QQ$ and\/ $u \in \QQ$. If\/ $P(W \cup \{u\}) = \emptyset$, then 
$$\<x,u\> \le 1 - \gamma$$ 
for every $x \in P(W)$.
\end{lemma}

\begin{proof}
Given $x \in P(W)$, let $y$ be the orthogonal projection of $x$ onto the hyperplane $\{z \in \R^d : \<z,u\> = 1\}$. By Definition~\ref{def:gamma}, we have
\[
\gamma \leq D(W,u) \le \|x - y\| = 1 - \<x,u\>,
\]
where the equality holds since $\<x,u\> \leq 1$ for every $x \in P(W)$, by~\eqref{eq:P(W)}
\end{proof}

Recall that $\delta = \delta(\QQ) > 0$ is a sufficiently small constant. The following lemma may be taken to be the definition of $\delta$.  

\begin{lemma}\label{lem:close-to-many-faces}
Let $W \subset \QQ$ and $T \subset N_\QQ(W)$ be such that $P(W\cup\{u\}) \neq \emptyset$ for all $u \in T$. If there exists $x \in P(W)$ such that 
$$\< x,u \> \ge 1 - 2\delta$$
for every $u \in T$, then $W \cup T$ is a clique and $P(W \cup T) \neq \emptyset$.
\end{lemma}

Lemma~\ref{lem:close-to-many-faces} is proved in Appendix~\ref{gamma:app} (see Lemma~\ref{lem:close-to-many-faces:app}). 

\subsection{The proof of Lemma~\ref{cor:faces-far-away}}

Lemma~\ref{cor:faces-far-away} is a fairly straightforward consequence of Lemmas~\ref{lem:quasi-innerprod} and~\ref{lem:faces-far-away}. Let's first make a simple observation, which provides us with the set $X$ that we will use to prove Lemma~\ref{cor:faces-far-away}. 

\begin{obs}\label{obs:faces-far-away}
Let\/ $W \subset W' \subset \QQ$, and let\/ $X' \in \U[W']$. There exists $X \in \U[W]$ such that $X \subset \HH(W')$ and $X' = X \cap W'^\perp$.
\end{obs}

\begin{proof}
Recall from Definition~\ref{def:induced} 
that if $X' \in \U[W']$, then there exists $Y \in \U$, with $Y \subset \HH(W')$, such that $X' = Y \cap W'^\perp$. We claim that the set $X = Y \cap W^\perp$ satisfies $X \in \U[W]$ and $X' = X \cap W'^\perp$. Indeed, $X' = X \cap W'^\perp$ holds because $W^\perp \cap W'^\perp = W'^\perp$, and $X \in \U[W]$ holds because $Y \subset \HH(W') \subset \HH(W)$, in both cases because $W \subset W'$. 
\end{proof}

\begin{proof}[Proof of Lemma~\ref{cor:faces-far-away}]
Let $X \in \U[W]$ be the set given by Observation~\ref{obs:faces-far-away}, so $X \in \U[W]$, $X \subset \HH(W')$ and $X' = X \cap W'^\perp$. Let $P \in \cP(W)$, let $x \in \Delta(P,W')$ satisfy $x + X' \subset \Delta(P,W')$, and let $y \in X$. Our aim is to show that $x + y \in P$.

We divide the proof into two cases, according to whether or not $y \in (W' \setminus W)^\perp$. 

Suppose first that we do have $y \in (W' \setminus W)^\perp$. Then 
$$y \in X \cap (W' \setminus W)^\perp \subset W^\perp \cap (W' \setminus W)^\perp = W'^\perp,$$ 
and hence $y \in X' = X \cap W'^\perp$. Since $x + X' \subset \Delta(P,W')$ by assumption, it follows that $x + y \in \Delta(P,W') \subset P$, as required.

So suppose instead that $y \notin (W' \setminus W)^\perp$, and observe that, since $\< y,u \> \le 0$ for every $u \in W'$, because $y \in X \subset \HH(W')$, we must have $\< y,u \> < 0$ for some $u \in W' \setminus W$. In this case we shall prove that $x + y \in P$ using Lemma~\ref{lem:faces-far-away}. The key fact is that, since $\<y,u\> < 0$ and $y \in \cL_R$, it follows from Lemma~\ref{lem:quasi-innerprod} that 
\begin{equation}\label{eq:faces-far-away-cond}
\< y,v \> \le 0 \qquad \text{for every} \qquad v \in \{u\} \cup N_\QQ(W \cup \{u\}),
\end{equation}
since $uv \in E\big(\Vor(\QQ)\big)$ for every $v \in N_\QQ(W \cup \{u\})$.

To complete the proof, we need to formalize the scaling that we shall use to apply Lemma~\ref{lem:faces-far-away}, and then verify the remaining conditions of that lemma. Thus, let $P = P(W,w;a,t,\tau)$, and recall from Definition~\ref{def:cPW} that $t > C \ge R_0/\gamma$, where $\gamma = \gamma(\QQ)$ is the constant defined in Definition~\ref{def:gamma}. Note that $x + y \in P$ if and only if
\begin{equation}\label{eq:newsurface-ip-cond}
t^{-1} \big( x + y - a \big) \in P(W,\hat{w}),
\end{equation}
by~\eqref{def:PWattau}, where $\hat{w} = (\tau / t) w \in \cL_R$. It therefore suffices to prove that~\eqref{eq:newsurface-ip-cond} holds. We shall apply Lemma~\ref{lem:faces-far-away} to the points $x' := t^{-1}(x - a)$ and $y' := t^{-1} y$.

To check that the conditions of the lemma hold, note first that $u \in N_\QQ(W)$, since $u \in W' \setminus W$ and $W'$ is a clique. Next, observe that
$$x \in \Delta(P,W') = a + t \cdot P(W',\hat{w}) \subset a + t \cdot P(W \cup \{u\},\hat{w}),$$
and thus $x' \in P(W \cup \{u\},\hat{w})$. Observe also that $y \in X \subset W^\perp$, and that $\| y' \| \le \gamma$, since $\| y \| \le R_0$ and $t \ge R_0/\gamma$. Since we have already verified~\eqref{eq:faces-far-away:condition}, it therefore follows by Lemma~\ref{lem:faces-far-away} that $x' + y' \in P(W,\hat{w})$, and hence that~\eqref{eq:newsurface-ip-cond} holds, as required.
\end{proof}

\subsection{The proof of Lemma~\ref{lem:eating:a:set}}

We shall deduce Lemma~\ref{lem:eating:a:set} from Lemmas~\ref{lem:far:from:nonfaces} and~\ref{lem:close-to-many-faces}, together with the following consequence of Lemma~\ref{lem:good:faces}. 

\begin{lemma}\label{lem:exists:emptyset:corner}
Let\/ $W \in \W$, set\/ $k := \dim(W^\perp)$, and suppose that\/ $W$ is $k$-good. Then for every $W \subset W'' \in \W$, there exists $W'' \subset W' \in \W$ such that $\emptyset \in \U[W']$.
\end{lemma}

\begin{proof} 
Set $k' := \dim( W''^\perp )$, and apply Lemma~\ref{lem:good:faces} with $s = k$. If $k' = 0$, then $\emptyset \in \U[W'']$, as required. We may therefore assume that $k' > 0$, in which case $W''$ is $k'$-good.

We claim that there exists $W' \in \W$ with $W'' \subset W'$ and $\dim(W'^\perp) = 0$. We will then apply Lemma~\ref{lem:good:faces} again to deduce that $\emptyset \in \U[W']$. To define $W'$, choose an arbitrary vertex (i.e., an extreme point) $x$ of the convex polytope $P(W'')$, and set 
$$W' := \big\{ u \in \QQ : \< x,u \> = 1 \big\}.$$ 
Observe that $W'' \subset W'$ (by~\eqref{eq:P(W)}), and that $W'$ is a clique, by Lemma~\ref{lem:faces-to-cliques}, since $x \in P(W')$. Moreover, since $x$ is an extreme point of $P(W'')$, it follows that $W'$ is maximal such that $P(W') = \{x\}$, and hence $W' \in \W$, by Definition~\ref{def:curlyW}.

By Lemma~\ref{lem:clique-of-correct-dim}, it follows that $\dim(W'^\perp) = 0$. Hence, applying Lemma~\ref{lem:good:faces} again, we deduce that $\emptyset \in \U[W']$, as required.
\end{proof}

We are now ready to prove Lemma~\ref{lem:eating:a:set}. 

\begin{proof}[Proof of Lemma~\ref{lem:eating:a:set}]
Observe first that it will suffice to prove the lemma in the case $\tau(P) = 0$. To see why this is the case, recall from Lemmas~\ref{lem:tube-is-tube} and~\ref{lem:tube-faces-to-P-faces} that  the tubular polytope $P(W,w;a,t,\tau)$ is a union of copies of the spherical polytope $P(W,w;a,t,0)$. Moreover, if $P = \bigcup_{i \in I} P_i$ and $\big[ P_i \cap A \big]_{\U[W]}^{P_i} = P_i$ for each $i \in I$, then $\big[ P \cap A \big]_{\U[W]}^P = P$. 

Let us therefore assume that $P = P(W,w;a,t,0)$, and let $K'$ be a strongly connected component of $K$ with diameter at most $\delta t$. Recalling Lemma~\ref{lem:P(W,w,a,t,tau):def2} and~\eqref{def:radius}, set
\begin{equation}\label{eq:eating-T}
T := \Big\{ u \in N_\QQ(W) : \< x - a, \, u \> \ge t - 2R_0 \; \text{ for some } x \in K' \Big\}.
\end{equation}
We think of the elements of $T$ as corresponding to the faces of $P$ that are `close' to $K'$. 

We shall first use Lemmas~\ref{lem:far:from:nonfaces},~\ref{lem:close-to-many-faces} and~\ref{lem:exists:emptyset:corner} to prove the following claim. We will use the claim to deduce the existence of an update rule $X \in \U[W]$ that will enable us to infect the sites of $K'$ one-by-one, even if they are close to the corners of $P$.

\begin{claim}\label{clm:eating-a-set-1}
There exists a clique $W' \supset W \cup T$ such that\/ $\emptyset \in \U[W']$.
\end{claim}

\begin{clmproof}{clm:eating-a-set-1}
In order to apply Lemma~\ref{lem:exists:emptyset:corner}, we first need to use Lemma~\ref{lem:close-to-many-faces} to show that $W \cup T$ is a clique. Observe that if $P(W \cup \{u\}) = \emptyset$ for some $u \in T$, then 
$$\< x - a, \, u \> \le (1 - \gamma)t < t - 2R_0$$ 
for every $x \in K' \subset P = a + t \cdot P(W)$, by Lemma~\ref{lem:far:from:nonfaces}, and since $t > C \ge 2R_0 / \gamma$. We therefore have $P(W \cup \{u\}) \neq \emptyset$ for all $u \in T$. Now, in order to apply Lemma~\ref{lem:close-to-many-faces}, note that since $K'$ has diameter at most $\delta t$, we have
$$\<x-a,u\> \ge t - 2R_0 - \delta t > (1 - 2\delta)t$$
for every $u \in T$ and $x \in K'$, since $t > C \ge 2R_0/\delta$. Noting that $K'$ is non-empty, it follows by Lemma~\ref{lem:close-to-many-faces} that $W \cup T$ is a clique and $P(W \cup T) \neq \emptyset$. 

Now, applying Lemma~\ref{lem:exists:WinW:trivial} to the set $W \cup T$, we obtain a set $W'' \in \W$ with 
$$W \cup T \subset W'' \qquad \text{and} \qquad P(W'') = P(W \cup T).$$
Finally, recalling that $W$ is $k$-good, by Lemma~\ref{lem:exists:emptyset:corner} we obtain a clique $W'' \subset W' \in \W$ such that $\emptyset \in \U[W']$. Since $W \cup T \subset W'' \subset W'$, this proves the claim.
\end{clmproof}

Consider the $\U[W]$-process in $P$. In the next claim, we shall show that we can infect the elements of $K'$ one-by-one in increasing order of their inner product with 
$$v := \sum_{u \in W'} u,$$ 
with ties broken arbitrarily. Before stating the claim formally, recall from~\eqref{eq:induced-def} that since $\emptyset \in \U[W']$, there must exist a set $Y \in \U$ such that $Y \subset \HH(W')$ and $Y \cap W'^\perp = \emptyset$. Set $X := Y \cap W^\perp$ and observe that $X \subset Y \subset \HH(W') \setminus W'^\perp \subset \HH(W)$, and therefore
\begin{equation}\label{eq:eating:a:set:X:properties}
X \in \U[W] \qquad \text{and} \qquad X \subset \HH(W') \setminus W'^\perp.
\end{equation}
Let us fix a set $X$ satisfying~\eqref{eq:eating:a:set:X:properties}. The following claim shows that each element $y \in K'$ can be infected (in the $\U[W]$-process, and in fact only using the set $X \in \U[W]$) by the set $P \setminus K'$, together with those elements of $K'$ that have smaller inner product with $v$. 

\begin{claim}\label{clm:eating:from:the:corner}
If $y \in K'$, then 
$$y + X \subset (P \setminus K) \cup \big\{ z \in K' : \<z,v\> < \<y,v\> \big\}.$$
\end{claim}

\begin{clmproof}{clm:eating:from:the:corner}
We will first show that $y + X \subset P$. By Lemma~\ref{lem:P(W,w,a,t,tau):def2}, and recalling that $y \in K' \subset P$ and $X \subset W^\perp$, since $X \in \U[W]$, to do so it will suffice to show that 
\begin{equation}\label{eq:need:to:be:in:P}
\big\< x + y - a, \, u \big\> \le t
\end{equation}
for every $x \in X$ and $u \in N_\QQ(W)$. If $u \in T$, then this holds because $\< x,u \> \leq 0$ and $\<y-a,u\> \leq t$, the first since $T \subset W'$ and $X \subset \HH(W')$, and the second because $y \in P$. If~$u \not\in T$, on the other hand, then it follows from~\eqref{eq:eating-T} that
$$\big\< y - a, \, u \big\> < t - 2R_0,$$
since $y \in K'$. Recalling from~\eqref{def:radius} that $\|x\| \le R_0$ for every $x \in X$, we obtain~\eqref{eq:need:to:be:in:P}, and hence $y + X \subset P$, as claimed. 

Now, let $x \in X$, and suppose that $x+y \in K$. By the definition~\eqref{def:strongly:conn} of a strongly connected component, and since $\|x\| \le R_0$ and $y \in K'$, it follows that $x+y \in K'$. We~claim that moreover $\<x+y,v\> < \<y,v\>$, i.e., that $\<x,v\> < 0$. To see this, recall from~\eqref{eq:eating:a:set:X:properties} that $X \subset \HH(W') \setminus W'^\perp$, and therefore $\< x,u \> \le 0$ for every $u \in W'$, and moreover $\< x,u \> < 0$ for some $u \in W'$. Thus $\< x,v \> = \sum_{u \in W'} \< x,u \> < 0$, 
as claimed. 
\end{clmproof}

Now, let $y \in K'$ and suppose that we have already infected all elements $z \in K'$ with $\<z,v\> < \<y,v\>$. Then, by Claim~\ref{clm:eating:from:the:corner}, the set $y + X$ is entirely infected. Recalling that $X \in \U[W]$, it follows that $y$ is also infected in the $\U[W]$-process in $P$, and hence the entire set $K'$ is contained in the closure $\big[ P \setminus K \big]^P_{\U[W]}$. Since $K'$ was an arbitrary strongly connected component of $K$, the lemma follows.
\end{proof}

\section{Interiors, extensions, buffers, and growth sequences}\label{buffers:sec}

In this section we define several notions of the `interior' and `extension' of the polytopes introduced in Section~\ref{sec:polytopes}. Various basic properties of these notions are stated, with the (relatively straightforward) proofs being given in Appendix~\ref{buffers:app}. 

\subsection{The interior of a polytope}

To begin, given $W \subset \QQ$ and a polytope $P \in \cP(W)$, let us define the \emph{interior} of $P$ to be\footnote{Recall from Definition~\ref{def:cPW} and~\eqref{def:Delta:face} the definitions of the family of polytopes $\cP(W)$, and of the $W'$-face $\Delta(P,W')$ of a polytope $P$. We define the interior of our canonical polytopes $P(W,w)$ similarly, by setting $\interior\big( P(W,w) \big) := t^{-1} \big( \interior\big( t \cdot P(W,w) \big) \big)$ for any $t > C$.}
\begin{equation}\label{def:intP}
\interior(P) := P  \setminus \bigcup_{u \in N_\QQ(W)} \Delta\big( P,W \cup \{u\} \big).
\end{equation}
Note that if $W \in \W$ then this coincides with the usual definition of the interior in the Euclidean space $\aff(P)$, and $\interior(P)$ is non-empty, by Lemmas~\ref{lem:clique-of-correct-dim} and~\ref{lem:faces:lower:dim}. We remark that $\interior(P) \not\in \cP(W)$, since it is not closed (unless $\dim(W^\perp) = 0$, in which case $\interior(P) = P$). It will therefore frequently be necessary to work instead with the following polytope, which is in $\cP(W)$ and contains the same lattice points as $\interior(P)$.

\begin{definition}\label{def:closed:interior}
For each $W \in \W$ and each polytope $P = P(W,w;a,t,\tau) \in \cP(W)$, choose an arbitrary $y \in \interior\big( P(W) \big)$ and $\eps > 0$ sufficiently small, and define 
\begin{equation}\label{def:clint}
\clint(P) := P(W,w;a+\eps y,t-\eps,\tau).
\end{equation}
We call $\clint(P)$ the \emph{closed interior} of $P$. 
\end{definition}

We remark that the reader should not worry that this definition (and also those below) depends on the choice of $y \in \interior\big( P(W) \big)$ and $\eps > 0$: any such $y$ and (sufficiently small) $\eps$ will do. The following lemma (see Lemma~\ref{lem:interiors-are-canonical:app}) motivates the definition; in fact, one can think of $\clint(P)$ as being an arbitrary polytope satisfying the conclusion of the lemma.

\begin{lemma}\label{lem:interiors-are-canonical}
Let $W \in \W$ and $P \in \cP(W)$. Then\/ $\clint(P) \in \cP(W)$, 
\begin{equation}\label{eq:interior:droplets:equal}
\clint(P) \subset \interior(P) \qquad \text{and} \qquad \clint(P) \cap \Z^d = \interior(P) \cap \Z^d.
\end{equation}
\end{lemma}
 
The closed interior will play an important role in Sections~\ref{sec:deterministic} and~\ref{proof:sec}. In order to prove the deterministic lemmas in Section~\ref{sec:deterministic}, we shall also need the following `shifted' version, which is translated so that it intersects a face of $P$. 

\pagebreak

\begin{definition}\label{def:shifted:interior}
For each $W \in \W$, each polytope $P = P(W,w;a,t,\tau) \in \cP(W)$, and each $W \subset W' \in \W$, choose $y \in \interior\big( P(W') \big)$ and $\eps > 0$ sufficiently small, and define
\begin{equation}\label{def:clint2}
\clint(P \to W') := P(W,w;a+\eps y,t-\eps,\tau).
\end{equation}
We call $\clint(P \to W')$ the \emph{$W'$-shifted closed interior} of $P$. 
\end{definition}

This polytope has the useful property that its $W'$-face contains the same lattice points as the interior of $\Delta(P,W')$ (see Lemmas~\ref{lem:clint:Delta:app1} and~\ref{lem:clint:Delta:app2}). 

\begin{lemma}\label{lem:clint:Delta}
Let $W \in \W$ and $P \in \cP(W)$, and let $W \subset W' \in \W$. Then
$$P^\circ \in \cP(W), \qquad P^\circ \subset P \qquad \text{and} \qquad \Delta\big( P^\circ,W' \big) \subset \interior\big( \Delta(P,W') \big),$$ 
where $P^\circ := \clint(P \to W')$. Moreover,
$$\Delta\big( P^\circ,W' \big) \cap \Z^d = \interior \big( \Delta(P,W') \big) \cap \Z^d.$$
\end{lemma}

The $W'$-shifted closed interior has another important property (see Lemma~\ref{lem:x-notin-Delta:Wv:app}): it only intersects faces of $P$ corresponding to subsets of $W'$. 

\begin{lemma}\label{lem:x-notin-Delta:Wv}
Let $W \in \W$ and $P \in \cP(W)$, and let $W \subset W' \in \W$. If $x \in \clint(P \to W')$, then
$$x \notin \Delta\big( P, W \cup \{v\} \big)$$
for every $v \in \QQ \setminus W'$.
\end{lemma}

\subsection{Forwards and sideways faces}

In Section~\ref{sec:deterministic}, we shall divide growth into two types: `forwards growth' (in direction $w$) 
and `sideways growth' (on the faces of $P$ that are perpendicular to $w$). Next we define the families of faces corresponding to these two types of growth. First, for each $W \in \W$ and $w \in \cL_R \cap \SS^{d-1}$, set
$$N_\QQ(W,w) := \big\{ u \in N_\QQ(W) : \< u,w \> > 0 \big\}.$$
It will be notationally convenient to define the forwards and sideways `faces' of a polytope in $\cP(W,w)$ to be the  corresponding sets of (maximal) cliques. 

\begin{definition}\label{def:forwards:sideways:faces}
Let $W \in \W$ and $w \in \cL_R \cap \SS^{d-1}$.  
\begin{itemize}
\item[$(a)$] The \emph{forwards faces} of a polytope $P \in \cP(W,w)$ are
\begin{equation}\label{def:WtoP}
\W^\to(P) := \big\{ W' \in \W \,:\, W \subsetneq W' \, \text{ and } \, W' \cap N_\QQ(W,w) \ne \emptyset \big\}.
\end{equation}
\item[$(b)$] The \emph{sideways faces} of a polytope $P \in \cP(W,w)$ are
\begin{equation}\label{def:Wperp}
\W^\perp(P) := \big\{ W' \in \W \,:\, W \subsetneq W' \subset \{w\}^\perp \big\}.
\end{equation}
\end{itemize}
\end{definition}

In particular, note that the sets $\W^\to(P)$ and $\W^\perp(P)$ only depend on $W$ and $w$, and that $\W^\perp(P)$ is empty unless $w \in W^\perp$. When growing forwards, our task will be to infect the `forwards boundary' of a polytope, which is defined as follows.

\pagebreak

\begin{definition}\label{def:forwards:boundary}
Let $W \in \W$ and $P \in \cP(W)$. The \emph{forwards boundary} of $P$ is the set
\begin{equation}\label{def:Delta}
\Delta(P) := \bigcup_{W' \in \W^\to(P)} \interior\big( \Delta(P,W') \big).
\end{equation}
\end{definition}

Let us note here (see Lemma~\ref{lem:Delta:forwards:app} for the proof) that this set is equal to the union of all of the `co-dimension 1' forwards faces of $P$. 

\begin{lemma}\label{lem:Delta:forwards}
Let $W \in \W$ and $P \in \cP(W)$. Then
$$\Delta(P) = \bigcup_{u \in N_\QQ(W,w)} \Delta( P, W \cup \{u\}).$$
\end{lemma}

Lemma~\ref{lem:Delta:forwards} is a straightforward consequence of the definitions and the following easy lemma (see Lemma~\ref{lem:max:face:interior:app}), which will also be used in Section~\ref{sec:deterministic}. 

\begin{lemma}\label{lem:max:face:interior} 
Let $W \subset \QQ$ with $P(W) \ne \emptyset$, let $P \in \cP(W)$, and let $x \in P$. If\/ $W' \subset \QQ$ is maximal such that 
$$W \subset W' \in \W \qquad \text{and} \qquad x \in \Delta(P,W'),$$ 
then $W' \setminus W \subset N_\QQ(W)$ and $x \in \interior\big( \Delta(P,W') \big)$.
\end{lemma}

Let us also note here the following property of the forwards faces, which will be needed in Section~\ref{proof:sec}. 

\begin{lemma}\label{lem:forwardsfaces:tau}
Let $W \in \W$ and $P \in \cP(W)$. If\/ $W' \in \W^\to(P)$, then there exists $Q \in \cP(W')$ with $\tau(Q) = 0$ such that $Q = \Delta(P,W')$ (as subsets of $\R^d$).
\end{lemma}

Lemma~\ref{lem:forwardsfaces:tau} follows easily from Lemma~\ref{lem:tube-faces-to-P-faces} and the definitions, see Lemma~\ref{lem:forwardsfaces:tau:app}.

\subsection{Extending and retracting a polytope}

In order to define sequences of growing droplets, we shall use two concepts of the `extension' of a polytope: one for growing only in direction $w$, and one for growing in all directions. We will only need these notions when $\dim(W^\perp) \ne 0$, and usually only when moreover $w \in W^\perp$. 

\begin{definition}\label{def:ext:P}
Given $W \in \W \setminus \W_0$ and a polytope $P = P(W,w;a,t,\tau) \in \cP(W)$, the \emph{extension} of $P$ is 
\begin{equation}\label{def:ext}
\ext(P) := P\big( W,w; a - \eps y, t + \eps, \tau \big)
\end{equation}
where $y \in \interior\big( P(W) \big)$, and $\eps > 0$ is minimal such that $\ext(P) \cap \Z^d \ne P \cap \Z^d$.
\end{definition} 

We remark that, as in Definitions~\ref{def:closed:interior} and~\ref{def:shifted:interior}, it does not matter that the definition of $\ext(P)$ depends on the choice of $y$: we may choose any element $y \in \interior\big( P(W) \big)$. We prove in Lemma~\ref{lem:ext:exists} that there exists $\eps > 0$ such that $\ext(P) \cap \Z^d \ne P \cap \Z^d$.

Growing from $P$ to $\ext(P)$ will be one of the two basic steps we use to grow a droplet. In order to use extensions the following basic properties will be needed: $P$ is contained in its extension, and every lattice point in $\ext(P) \setminus P$ lies on one of the faces. 

\pagebreak

\begin{lemma}\label{lem:exterior}
Let $W \in \W$, let $P \in \cP(W)$, and set $P' := \ext(P)$. Then
$$P \subset P' \qquad \text{and} \qquad \interior(P') \cap \Z^d \subset P.$$
\end{lemma}

The proof of Lemma~\ref{lem:exterior} is straightforward (see Lemmas~\ref{lem:exterior:app} and~\ref{lem:exterior:faces:app}). When growing in direction $w$, we shall instead use the following (simpler) notions. 

\begin{definition}\label{def:fext:P}
Let $W \in \W$, and let $P = P(W,w;a,t,\tau) \in \cP(W)$ with $w \in W^\perp$.
\begin{itemize}
\item[$(a)$] The \emph{forwards extension} of $P$ is
\begin{equation}\label{def:fext}
\fext(P) := P\big( W,w;a,t,\tau' \big)
\end{equation}
where $\tau' > \tau$ is minimal such that $\fext(P) \cap \Z^d \ne P \cap \Z^d$.\smallskip
\item[$(b)$] The \emph{forwards retraction} of $P$ is
$$\fret(P) := P\big( W,w;a,t,\tau' \big)$$
where $0 < \tau' < \tau$ and $\tau - \tau'$ is sufficiently small. If $\tau = 0$, then $\fret(P) := P$.
\end{itemize}
\end{definition}

We shall use the following simple facts about the forwards extension and retraction. Note that $\fret(P) \subset P \subset \fext(P)$. The first property is that all of the lattice points in $\fext(P) \setminus P$ lie on the forwards faces of $\fext(P)$, and similarly all of the lattice points in $P \setminus \fret(P)$ lie on the forwards faces of $P$.

\begin{lemma}\label{lem:fext}
Let $W \in \W$ and $P \in \cP(W)$, with $w(P) \in W^\perp$. Then
$$\fext(P) \cap \Z^d \subset P \cup \Delta\big( \fext(P) \big) \qquad \text{and} \qquad P \cap \Z^d \subset \fret(P) \cup \Delta(P).$$
\end{lemma}

The second property is also straightforward; it is moreover not hard to see that the additional conditions on $P$ in this lemma are both necessary. 

\begin{lemma}\label{lem:fret}
Let $W \in \W$ and $P \in \cP(W)$, with $w(P) \in W^\perp$. If\/ $\tau(P) > 0$ and $\Delta(P) \cap \Z^d \ne \emptyset$, then  
$$\fext\big( \fret(P) \big) = P.$$
\end{lemma}

Lemmas~\ref{lem:fext} and~\ref{lem:fret} both follow easily from the definitions; see Lemmas~\ref{lem:fext:app} and~\ref{lem:fret:app} for the details. 

\subsection{Growth sequences}\label{sec:growth:seq}

We shall use $\ext$ and $\fext$ in Sections~\ref{sec:deterministic} and~\ref{proof:sec} to define sequences of growing droplets; let us next introduce the two basic constructions that will be used to do this. The first, which is very simple, only allows us to grow `forwards', whereas the second will be used when we also need to grow `sideways'. 

\begin{definition}\label{def:forwards:growth:sequence}
The \emph{forwards growth sequence} $\G$ with seed $Q$, where $Q \in \cP(W,w)$ for some $W \in \W$ and $w \in \cL_R \cap \SS(W)$, is defined by setting $Q_0 := Q$ and 
\begin{equation}\label{def:sequence:of:growing:polytopes:forward}
Q_j := \fext( Q_{j-1} ) 
\end{equation}
for each $j \ge 1$.
\end{definition}

The first important property of forwards growth sequences is as follows; see Lemma~\ref{lem:poly-seq-eats-end:app}. This property allows us to infect the `forwards end' of a polytope step by step. 

\begin{lemma}\label{lem:poly-seq-eats-end}
Let $W \in \W$, and let $P,Q \in \cP(W)$ be such that
$$Q \subset P, \qquad a(P) = a(Q), \qquad t(P) = t(Q) \qquad \text{and} \qquad w(P) = w(Q) \in W^\perp.$$
Let $\G = (Q_j)_{j \ge 0}$ be the forwards growth sequence with seed $Q$, and let $m$ be maximal such that $Q_m \subset P$. Then
\[
P \cap \Z^d \subset Q_m.
\]
\end{lemma}

We remark that it may not be true that $P \subset Q_m$, for example if $\Delta(P) \cap \Z^d$ is empty. This construction will be used in the proofs of Lemmas~\ref{lem:sideways-step-general} and~\ref{lem:extensions:cover:tube}. In Section~\ref{proof:sec} a union bound will be used to control the probability that at some step we fail to grow from $Q_{j-1}$ to $Q_j$. To bound the number of steps, we shall use the following lemma.

\begin{lemma}\label{lem:growing:linearly:forwards}
There exists a constant $\xi = \xi(\QQ) > 0$ such that the following holds. Let $W \in \W$ and $w \in \cL_R \cap \SS(W)$, and let $\G$ be the forwards growth sequence with seed $Q \in \cP(W,w)$. Then
$$\tau(Q_j) \ge \tau(Q_0) + \xi \cdot j$$
for all $j \ge 1/\xi$. 
\end{lemma}

To see why Lemma~\ref{lem:growing:linearly:forwards} should be true, observe that in each step one (or more) of the forwards faces $\Delta(Q_j,W')$, where $W' \in \W^\to(Q_j)$, intersects a new translate of the lattice $\L(W')$. Since there are only a bounded number of faces, and each passes through a bounded number of translates of $\L(W')$ when $\tau$ increases by $1$ (say), the claimed bound follows (see Lemma~\ref{lem:growing:linearly:forwards:app} for the details). We remark that some lower bound on $j$ is needed for the lemma to be true, since individual steps of a growth sequence can be arbitrarily small. 

In Section~\ref{proof:sec} (see the proof of Lemma~\ref{lem:IHb:likely}) we shall sometimes also need to grow `sideways'. The following construction will allow us to do so.  

\begin{definition}\label{def:growth:sequence}
A \emph{growth sequence} $\G$ with seed $Q$, where $Q \in \cP(W,w)$ for some $W \in \W$ and $w \in \cL_R \cap \SS(W)$, is a collection of polytopes $Q^{(i)}_j$ such that $Q^{(1)}_0 := Q$,  
\begin{equation}\label{def:sequence:of:growing:polytopes}
Q^{(i)}_j := \fext\big( Q^{(i)}_{j-1} \big) \qquad \text{and} \qquad Q^{(i+1)}_0 := \ext\big( Q^{(i)}_* \big)
\end{equation}
for each $i,j \ge 1$, for some $Q^{(i)}_* \in \cP(W,w)$ with $t(Q^{(i)}_*) = t(Q^{(i)}_0)$ and $Q^{(i)}_* \subset Q^{(i)}_{\ell(i)} \setminus Q^{(i)}_0$, where $\ell(i)$ will depend on the particular application.\footnote{We remark that $Q^{(i)}_*$ is not determined by the polytopes $Q^{(i)}_j$; what matters is that there \emph{exists} a polytope satisfying the stated conditions. In particular, the definition implies that each $\ell(i)$ is sufficiently large so that the set $Q^{(i)}_{\ell(i)} \setminus Q^{(i)}_0$ contains a suitable polytope $Q^{(i)}_*$. It will be important that $Q^{(i)}_*$ is disjoint from $Q^{(i)}_0$ because we will usually not know that the whole of $Q^{(i)}_0$ is infected, see Section~\ref{subsec:buffers}.}  
\end{definition}

Note that for each fixed $i \ge 1$, the sequence $( Q^{(i)}_j )_{j \ge 0}$ is a forwards growth sequence with seed $Q^{(i)}_0$. We shall need analogues of Lemmas~\ref{lem:poly-seq-eats-end} and~\ref{lem:growing:linearly:forwards} for sideways growth; the first of these is again quite straightforward (see Lemma~\ref{lem:growth-seq-eats-sides:app}).

\begin{lemma}\label{lem:growth-seq-eats-sides}
Let $W \in \W$, $w \in \cL_R \cap \SS(W)$ and $t > C$, and let $Q \in \cP(W,w)$ with $t(Q) \le t$.  Let $\G$ be a growth sequence with seed $Q$, and let $m$ be maximal such that $t\big( Q^{(m)}_0 \big) \le t$. Then
$$P \cap \Z^d \subset Q^{(m)}_*$$
for some $P \in \cP(W,w;t)$.
\end{lemma}

Unfortunately, the analogue of Lemma~\ref{lem:growing:linearly:forwards} is slightly more complicated, since there exist growth sequences for which $t(Q^{(m)}_0) / m$ is arbitrarily small. We therefore need to choose the polytopes $Q^{(i)}_*$ with a little care. To do this, the following definition will help; recall from Lemma~\ref{lem:growing:linearly:forwards} that the average increase in $\tau$ in a forwards step is at least $\xi$. 

\begin{definition} 
We say that a polytope $P$ is a \emph{grower} if either 
$$t\big( \ext(P) \big) \ge t(P) + \xi,$$
or there exists $W' \in \W^\perp(P)$ such that 
$$\Delta\big( \ext(P), W' \big) \cap \Z^d \ne \emptyset.$$
We say that a growth sequence $\G$ is \emph{happy} if $Q^{(i)}_*$ is a grower for every $i \in \N$.
\end{definition}

The following lemma will be used in the proof of Lemma~\ref{lem:IHb:likely}.  

\begin{lemma}\label{lem:growing:linearly:sideways}
There exists a constant $\xi' = \xi'(\QQ) > 0$ such that the following holds. Let $W \in \W$ and $w \in \cL_R \cap \SS(W)$, and let $\G$ be a happy growth sequence with seed $Q \in \cP(W,w)$. Then
$$t(Q^{(i)}_0) \ge t(Q) + \xi' \cdot i$$
for every $i \ge 1/\xi'$. 
\end{lemma}

The proof of Lemma~\ref{lem:growing:linearly:sideways} is similar to that of Lemma~\ref{lem:growing:linearly:forwards}; see Lemma~\ref{lem:growing:linearly:sideways:app} for the details. We remark that it is straightforward to construct a happy growth sequence, using the following lemma (see Lemma~\ref{lem:finding:a:grower:app}). 

\begin{lemma}\label{lem:finding:a:grower}
Let $W \in \W$ and $w \in \cL_R \cap \SS(W)$, and let $Q \in \cP(W,w)$. Then there exists $Q' \in \cP(W,w)$, with
$$t(Q) = t(Q'), \qquad |\tau(Q) - \tau(Q')| \le 1 \qquad \text{and} \qquad a(Q) - a(Q') = \mu w$$
for some $0 \le \mu \le 1$, such that $Q'$ is a grower. 
\end{lemma}

\subsection{Buffers}\label{subsec:buffers}

When growing on the sideways faces of a polytope, it will only be known that the corresponding cliques are $(s,w)$-semi-good for some $s \ge 1$ and $w \in \cL_R \cap \SS^{d-1}$, and it will therefore not be possible to infect the entire extension. Fortunately, in order to grow forwards, i.e., from a polytope $P$ to $\fext(P)$, we do not need all of $P$, but only the points `close to' the forwards faces. To be precise, we shall infect the following subset. 

\pagebreak

\begin{definition}\label{def:buffer}
Let $W \in \W$ and $P = P(W,w;a,t,\tau) \in \cP(W)$, with $w \in W^\perp$. For each $x \in P$, define 
$$\tau_P(x) := \inf\big\{ \tau^* \ge 0 : x \in P\big( W,w;a,t,\tau^* \big) \big\}.$$
The \emph{forwards buffer} of $P$ is the set
$$B(P) := \big\{ x \in P : \tau_P(x) > \tau - C \big\}.$$
\end{definition}

Note that if $\tau(P) \le C$, then $B(P) = P$. Let us first observe that $B(P)$ contains all points of $P$ within distance $R_0$ of $\Delta(P)$. This follows easily from the definition, using the fact that $C = C(\QQ)$ is a sufficiently large constant (see Lemma~\ref{lem:newbuffer:contains:oldbuffer:app}).

\begin{lemma}\label{lem:newbuffer:contains:oldbuffer}
Let $W \in \W$, $w \in \cL_R \cap \SS(W)$ and $P \in \cP(W,w)$. Then
$$\bigcup_{x \in \Delta(P)} \big\{ y \in P : \|x - y\| \le R_0 \big\} \subset B(P).$$
\end{lemma}

We shall need the following three properties of buffers. They will be used in Section~\ref{sec:deterministic} to prove our two main deterministic lemmas. First, for forwards growth, we shall use the following lemma. 

\begin{lemma}\label{lem:forwards:event}
Let $W \in \W$, $w \in \cL_R \cap \SS(W)$ and $P \in \cP(W,w)$, and set $P' := \fext(P)$. Then
\begin{equation}\label{eq:forwards:event}
B(P') \cap \Z^d \subset B(P) \cup \Delta(P').
\end{equation}
\end{lemma}

For sideways growth (that is, perpendicular to $w$) we shall instead use the following two lemmas. The first is similar to Lemma~\ref{lem:forwards:event}. 

\begin{lemma}\label{lem:fext:fret:buffers}
Let $W \in \W$, $w \in \cL_R \cap \SS(W)$ and $P \in \cP(W,w)$, and set $P' := \fret(P)$. Then 
$$B(P) \cap \Z^d \subset B(P') \cup \Delta(P).$$
\end{lemma}

Our second lemma for sideways growth is a little more technical, and requires some additional notation. Given $W \in \W$, $w \in \cL_R \cap \SS(W)$ and a polytope $P \in \cP(W,w)$, we define the \emph{forward half} of $P$ to be the polytope 
\begin{equation}\label{def:forward:half}
F(P) := P \cap \big( P + \tau(P) w / 2 \big).
\end{equation}
We remark that $F(P) \in \cP(W,w)$, and moreover if $P = P(W,w;a,t,\tau)$ then 
$$F(P) = P\big( W,w; a + \tau w/2,t,\tau/2 \big),$$
see Lemma~\ref{lem:Pzero:properties1} for the details. We can now state the final lemma of the section. 

\begin{lemma}\label{lem:sideways:event}
Let $W \in \W$, $w \in \cL_R \cap \SS(W)$ and $P \in \cP(W,w)$, and set $P' := \ext(P)$ and $P'' = \fret(P')$. If $\tau(P) \ge 5C$, then
\begin{equation}\label{eq:sideways:event}
B(P'') \cap \Z^d \subset P \cup \bigcup_{W' \in \W^\perp(P)} \interior\big( \Delta(F,W') \big),
\end{equation}
where $F := F(P')$. 
\end{lemma}

Lemmas~\ref{lem:forwards:event},~\ref{lem:fext:fret:buffers} and~\ref{lem:sideways:event} are proved in Appendix~\ref{app:buffers:sec}, see Lemmas~\ref{lem:forwards:event:app},~\ref{lem:fext:fret:buffers:app} and~\ref{lem:sideways:event:app}.

\section{Deterministic growth of droplets}\label{sec:deterministic}

The aim of this section, the last before the proof of Theorem~\ref{thm:upper}, is to prove two key lemmas about the deterministic growth of droplets. The first of the two deterministic growth lemmas is for forwards growth, and is the subject of Section~\ref{sec:growth-forwards}. The second of the lemmas is for sideways growth, and is given in Section~\ref{sec:growth-sideways}. We begin, however, in Section~\ref{fun:sec}, by proving a lemma that will be important in both of their proofs.

\subsection{A fundamental lemma}\label{fun:sec}

Given $W \subset W' \subset \QQ$ and $P \in \cP(W)$, define the \emph{$W'$-buffer} of $P$ to be the set
\begin{equation}\label{def:buffer:face}
B(P,W') := \bigcup_{x \in \Delta(P,W')} \big\{ y \in P : \|x - y\| \le R_0 \big\}.
\end{equation}
Observe that $B(P,W') \subset B(P)$ for every $W' \in \W^\to(P)$, by 
Lemmas~\ref{lem:Delta:forwards} and~\ref{lem:newbuffer:contains:oldbuffer}, since each such clique $W'$ contains an element $u \in N_\QQ(W,w)$. 

The following lemma is a consequence of Lemma~\ref{cor:faces-far-away}, and its proof is in fact the only time that we apply that lemma. Roughly speaking, the lemma says that if we wish to infect the face $\Delta(P,W')$ of a polytope in the $\U[W]$-process, and we know that the rest of the $W'$-buffer of $P$ is already infected, then it suffices to consider the $\U[W']$-process restricted to the face $\Delta(P,W')$. 

\begin{lemma}\label{lem:surface}
Let $W \subset W' \in \W$, and let $P \in \cP(W)$. Then
\begin{equation}\label{eq:face}
\big[ \Delta \cap A \big]^{\Delta}_{\U[W']} \subset \big[ (B \setminus \Delta) \cup (\Delta \cap A) \big]^B_{\U[W]},
\end{equation}
where 
$\Delta := \Delta(P,W')$ and $B := B(P,W')$.
\end{lemma}

\begin{proof}
By Lemma~\ref{cor:faces-far-away}, for each $X' \in \U[W']$ there exists $X \in \U[W]$ such that if $x \in \Delta$ then 
$$x + X' \subset \Delta \quad \Rightarrow \quad x + X \subset P.$$
By~\eqref{def:buffer:face}, and recalling that $\|y\| \le R_0$ for every $y \in X$, it follows that if $x \in \Delta$, then 
$$x + X' \subset \Delta \quad \Rightarrow \quad x + X \subset B,$$
and hence every site that is infected in the $\U[W']$-process on $\Delta$ with initial set $\Delta \cap A$ is also infected in the $\U[W]$-process with initial set $(B \setminus \Delta) \cup (\Delta \cap A)$, as required.
\end{proof}

We shall use Lemma~\ref{lem:surface} to prove the two main results of this section: Lemma~\ref{lem:forwards-step-general}, which deals with `forwards' growth (in direction $w$), and Lemma~\ref{lem:sideways-step-general}, which deals with `sideways' growth (perpendicular to $w$). The proof of Lemma~\ref{lem:sideways-step-general}, in particular, will be rather technical. These two lemmas will be our main deterministic tools in Section~\ref{proof:sec}. 

\subsection{Forwards deterministic growth}\label{sec:growth-forwards}

In order to state our key deterministic lemma for forwards growth, we need an additional definition, which is chosen (cf.~Definition~\ref{def:forwards:boundary} and Lemmas~\ref{lem:interiors-are-canonical} and~\ref{lem:fext}) to encode the deterministic property of $A$ that will be used to grow from $P$ to $\fext(P)$. Recall from Definition~\ref{def:int:filled} that we write $I_W^\bullet(P)$ for the event that a set $P \subset \R^d$ is internally filled in the $\U[W]$-process, i.e., that $\big[ P \cap A \big]_{\U[W]}^P = P \cap \Z^d$.

\begin{definition}\label{def:forwards-edge-filled}
Given $W \in \W$ and $P \in \cP(W)$, we say that $P$ is \emph{forwards edge-filled} by $A$ if the event 
$$I_{W'}^\bullet\Big( \clint\big( \Delta(P,W') \big) \Big)$$ 
holds for every $W' \in \W^\to(P)$. 
\end{definition}

Recall from Definition~\ref{def:buffer} the definition of the buffer $B(P)$. The following lemma is the key deterministic property of forwards growth.  

\begin{lemma}\label{lem:forwards-step-general}
Let $W \in \W$, $w \in \cL_R \cap \SS(W)$ and $P \in \cP(W,w)$, and set $P' := \fext(P)$. If $P'$ is forwards edge-filled by $A$, then
\begin{equation}\label{eq:forwards-event}
B(P') \cap \Z^d \subset \big[ B(P) \cup \big( P' \cap A \big) \big]^{P'}_{\U[W]}.
\end{equation}
\end{lemma}

Recall from Lemma~\ref{lem:forwards:event} that
$$B(P') \cap \Z^d \subset B(P) \cup \Delta(P'),$$
and from~\eqref{def:Delta} that $\Delta(P') \subset P'$.  Lemma~\ref{lem:forwards-step-general} will therefore follow easily from the following lemma, which will also be used later in the proof of Lemma~\ref{lem:sideways-step-general}.

\begin{lemma}\label{lem:multi-surface}
Let $W \in \W$, and let $P \in \cP(W)$. If $P$ is forwards edge-filled by $A$, then
$$\Delta(P) \cap \Z^d \subset \big[ \big( B(P) \setminus \Delta(P) \big) \cup \big( \Delta(P) \cap A \big) \big]^{B(P)}_{\U[W]}.$$
\end{lemma}

We will infect $\Delta(P)$ in stages, starting with the interiors of the faces of $P$ of highest dimension, and working our way down the dimensions. For each $j \ge 0$, define
\begin{equation}\label{def:Delta:j}
\Delta_j(P) := \bigcup_{\substack{W' \in \W^\to(P) \\ \dim(W'^\perp) \ge j}} \interior\big( \Delta(P,W') \big),
\end{equation}
and note that $\Delta_0(P) = \Delta(P)$, by Definition~\ref{def:forwards:boundary}, and that if $\dim(W^\perp) = k$, then $\Delta_k(P) = \emptyset$, since $\dim(W'^\perp) \le k - 1$ for every $W' \in \W^\to(P)$, by~\eqref{def:WtoP} and Lemma~\ref{lem:faces:lower:dim}.

\begin{proof}[Proof of Lemma~\ref{lem:multi-surface}]
Since $\Delta_0(P) = \Delta(P)$ and $\Delta_k(P) = \emptyset$, in order to prove the lemma it is enough to show that
\begin{equation}\label{eq:multi-surface-ind}
\Delta_j(P) \cap \Z^d \subset \big[ \Delta_{j+1}(P) \cup \big( B(P) \setminus \Delta(P) \big) \cup \big( \Delta(P) \cap A \big) \big]^{B(P)}_{\U[W]}
\end{equation}
for each $0 \le j < k$.  In particular, by~\eqref{def:Delta:j}, it suffices to prove the following claim.

\begin{claim}\label{claim:multi-surface}
If\/ $0 \le j < k$ and $W' \in \W^\to(P)$, with $\dim(W'^\perp) = j$, then
\begin{equation}\label{eq:multi-surface-stp}
\interior\big( \Delta(P,W') \big) \cap \Z^d \subset \big[ \Delta_{j+1}(P) \cup \big( B(P) \setminus \Delta(P) \big) \cup \big( \Delta(P) \cap A \big) \big]^{B(P)}_{\U[W]}.
\end{equation}
\end{claim}

\begin{clmproof}{claim:multi-surface}
Recalling~\eqref{def:clint2}, set $P^\circ := \clint(P \to W') \in \cP(W)$ and observe that
\begin{equation}\label{eq:multi-surface-canonical2}
\interior\big( \Delta(P,W') \big) \cap \Z^d = \clint \big( \Delta(P,W') \big) \cap \Z^d = \Delta(P^\circ,W') \cap \Z^d,
\end{equation}
by Lemma~\ref{lem:interiors-are-canonical} (applied to the polytope $\Delta(P,W') \in \cP(W')$) and Lemma~\ref{lem:clint:Delta}. Moreover, since $P$ is forwards edge-filled by $A$, we have
$$\clint\big( \Delta(P,W') \big) \cap \Z^d = \big[ \clint\big( \Delta(P,W') \big) \cap A \big]^{\clint(\Delta(P,W'))}_{\U[W']},$$
by Definitions~\ref{def:int:filled} and~\ref{def:forwards-edge-filled}. By~\eqref{eq:multi-surface-canonical2}, it follows that
\begin{equation}\label{eq:multi-surface:step1}
\interior\big( \Delta(P,W') \big) \cap \Z^d = \big[ \Delta(P^\circ,W') \cap A \big]^{\Delta(P^\circ,W')}_{\U[W']}.
\end{equation}
Now, by Lemma~\ref{lem:surface} (applied to the polytope $P^\circ \in \cP(W)$), we have 
$$\big[ \Delta(P^\circ,W') \cap A \big]^{\Delta(P^\circ,W')}_{\U[W']} \subset \big[ \big( B \setminus \Delta(P^\circ,W') \big) \cup \big( \Delta(P^\circ,W') \cap A \big) \big]^B_{\U[W]},$$
where $B := B(P^\circ,W')$, and therefore, by~\eqref{eq:multi-surface:step1}, 
$$\interior\big( \Delta(P,W') \big) \cap \Z^d \subset \big[ \big( B \setminus \Delta(P^\circ,W') \big) \cup \big( \Delta(P^\circ,W') \cap A \big) \big]^B_{\U[W]}.$$
Now, observe that, by~\eqref{def:Delta} and~\eqref{def:buffer:face} and Lemmas~\ref{lem:clint:Delta} and~\ref{lem:newbuffer:contains:oldbuffer}, 
\begin{equation}\label{eq:multi-surface-Delta-B}
\Delta(P^\circ,W') \subset \Delta(P)  \qquad \text{and} \qquad B = B(P^\circ,W') \subset B(P).
\end{equation}
Thus, in order to complete the proof of~\eqref{eq:multi-surface-stp}, it remains to show that
\begin{equation}\label{eq:multi-surface:step3}
\big( B \setminus \Delta(P^\circ,W') \big) \cap \Z^d \subset \Delta_{j+1}(P) \cup \big( B(P) \setminus \Delta(P) \big). 
\end{equation}

To prove~\eqref{eq:multi-surface:step3}, let $x \in B \cap \Z^d$ (and hence $x \in B(P)$, by~\eqref{eq:multi-surface-Delta-B}), and suppose that $x \in \Delta(P) \setminus \Delta_{j+1}(P)$. We are required to show that $x \in \Delta(P^\circ,W')$. By~\eqref{def:Delta} and~\eqref{def:Delta:j}, we have
$$x \in \interior\big( \Delta(P,W'') \big) \cap \Z^d$$
for some $W'' \in \W^\to(P)$ with $\dim(W''^\perp) \le j$. Moreover, since $x \in B \subset P^\circ$, we have
$$x \notin \Delta\big( P, W \cup \{v\} \big)$$
for all $v \in \QQ \setminus W'$, by Lemma~\ref{lem:x-notin-Delta:Wv}, so $W'' \subset W'$. Since $\dim(W'^\perp) = j$, it follows that $\dim(W''^\perp) = \dim(W'^\perp)$, and so, by Lemma~\ref{lem:faces:lower:dim}, we have $W'' = W'$. Therefore, by~\eqref{eq:multi-surface-canonical2},
$$x \in \interior\big( \Delta(P,W'') \big) \cap \Z^d = \interior\big( \Delta(P,W') \big) \cap \Z^d = \Delta(P^\circ,W') \cap \Z^d,$$ 
as required. This proves~\eqref{eq:multi-surface:step3}, and hence completes the proof of the claim.
\end{clmproof}

As observed above, Claim~\ref{claim:multi-surface} implies that~\eqref{eq:multi-surface-ind} holds for every $0 \le j < k$, and hence completes the proof of the lemma. 
\end{proof}
 
We can now easily deduce our main deterministic lemma for growth in direction $w$.

\begin{proof}[Proof of Lemma~\ref{lem:forwards-step-general}]
By Lemma~\ref{lem:forwards:event}, we have
\begin{equation}\label{eq:forwards-step-B}
B(P') \cap \Z^d \subset B(P) \cup \Delta(P').
\end{equation}
Moreover, by Lemma~\ref{lem:multi-surface}, since $P'$ is forwards edge-filled by $A$, we have
$$\Delta(P') \cap \Z^d \subset \big[ \big( B(P') \setminus \Delta(P') \big) \cup \big( \Delta(P') \cap A \big) \big]^{B(P')}_{\U[W]}.$$

Since $\Delta(P') \subset B(P') \subset P'$, 
it follows, by~\eqref{eq:forwards-step-B}, that
\[
\Delta(P') \cap \Z^d \subset \big[ B(P) \cup ( P' \cap A ) \big]^{P'}_{\U[W]}.
\]
Applying~\eqref{eq:forwards-step-B} once again, we obtain~\eqref{eq:forwards-event}, as required.
\end{proof}

\subsection{Sideways deterministic growth}\label{sec:growth-sideways}

For sideways growth, we need two further definitions relating to `internal filling', which will be analogues of Definitions~\ref{def:int:filled} and~\ref{def:forwards-edge-filled} for the setting of sideways growth. The first of these definitions will form part of the induction hypothesis in the next section (see Definition~\ref{def:ih}). Recall from~\eqref{def:cPWwtt} the definitions of the families of polytopes $\cP(W,w;t)$ and $\cP(W,w;t,\tau)$. 

\begin{definition}\label{def:int-half-filled}
Given $W \in \W$ and $P = P(W,w;a,t,\tau) \in \cP(W)$, we say that $P$ is \emph{internally half-filled} by $A$, written $I_W^\circ(P)$, if
$$P' \cap \Z^d \subset \big[ P \cap A \big]_{\U[W]}^P$$
for some $P' \in \cP(W,w;t)$ such that $P' \subset P$.
\end{definition}

The second definition is rather technical, which is an unfortunate consequence of the `directed' nature of growth on faces perpendicular to $w$. Set $\alpha := 1/8d$, and recall from~\eqref{def:Wperp} the definition of the family $\W^\perp(P)$ of sideways faces of $P$.

\begin{definition}\label{def:sideways-edge-filled}
Given $W \in \W$ and $P = P(W,w;a,t,\tau) \in \cP(W)$, we say that $P$ is \emph{sideways edge-filled} by $A$ if for every $W' \in \W^\perp(P)$, the following holds with $Q^* := \clint\big( \Delta(P,W') \big)$ and $t^* = t(Q^*)$: 
\begin{enumerate}
\item every polytope $Q \in \cP\big( W',w;t^*,\alpha\tau \big)$ with $Q \subset Q^*$ is internally half-filled;\vskip0.1cm
\item every polytope $Q \in \cP\big( W',w;t^* \big)$ with $Q \subset Q^*$ is forwards edge-filled.
\end{enumerate}
\end{definition}

This following lemma is our main deterministic lemma for sideways growth. Recall from~\eqref{def:constants} that $C > 0$ is a large constant. We say that $P$ is \emph{long} if $\tau(P) \ge C \cdot t(P)$. 

\begin{lemma}\label{lem:sideways-step-general}
Let $W \in \W$, $w \in \cL_R \cap \SS(W)$ and $P \in \cP(W,w)$, and set $P' := \ext(P)$. Suppose that $P'$ is long, and is forwards edge-filled and sideways edge-filled by $A$. Then
\begin{equation}\label{eq:sideways-event}
B(P') \cap \Z^d \subset \big[ P \cup \big( P' \cap A \big) \big]_{\U[W]}^{P'}.
\end{equation}
\end{lemma}

Our proof of Lemma~\ref{lem:sideways-step-general} will proceed by showing first that $P$ can grow sideways (perpendicular to $w$) and then that it can grow forwards (in the direction of $w$). To grow forwards, we shall use Lemma~\ref{lem:forwards-step-general}; to grow sideways, we need the following lemma.

\begin{lemma}\label{lem:sideways-step-part2}
Let $W \in \W$, $w \in \cL_R \cap \SS(W)$ and $P \in \cP(W,w)$, and set $P' := \ext(P)$ and~$P'' := \fret(P')$. If $P'$ is long and sideways edge-filled by $A$, then
\begin{equation}\label{eq:sideways-step-clm1}
B(P'') \cap \Z^d \subset \big[ P \cup \big( P' \cap A \big) \big]_{\U[W]}^{P'}.
\end{equation}
\end{lemma}

To prove Lemma~\ref{lem:sideways-step-part2}, we shall (partially) fill in the sides of $P''$, starting with those of highest dimension, moving in the direction of $w$ as the dimension decreases. Let us fix, until the end of the proof of Lemma~\ref{lem:sideways-step-part2}, $W$, $w$, $P$, $P'$ and~$P''$ as in the statement of the lemma, and set $k := \dim(W^\perp)$ and $\tau' := \tau(P')$. Recall that, by Lemma~\ref{lem:sideways:event}, 
$$B(P'') \cap \Z^d \subset P \cup \bigcup_{W' \in \W^\perp(P)} \interior\big( \Delta(F,W') \big),$$ 
where $F = F(P')$ is the forward half of $P'$, see~\eqref{def:forward:half}.  
It will therefore suffice to show that 
\begin{equation}\label{eq:sideways:step:aim}
\interior\big( \Delta(F,W') \big) \cap \Z^d \subset \big[ P \cup \big( P' \cap A \big) \big]_{\U[W]}^{P'}
\end{equation}
for every $W' \in \W^\perp(P)$. For each $0 \le j \le k$, define\footnote{It follows from Lemma~\ref{lem:tube-is-tube} that $P_j \in \cP(W,w)$; see Lemma~\ref{lem:Pzero:properties1}.}\begin{equation}\label{def:Pjs}
P_j := P' \cap \big( P' + 4(k - j) \alpha \tau' \cdot w \big),
\end{equation}
and observe that
\begin{equation}\label{eq:Pjs:nested}
F \subset P_0 \subset P_1 \subset \cdots \subset P_k = P'.
\end{equation}
We shall prove by induction on $k - j$ that
\begin{equation}\label{eq:sideways:IH:Deltas}
\Delta^\perp_j(P_j) \cap \Z^d \subset \big[ P \cup (P' \cap A) \big]_{\U[W]}^{P'}
\end{equation}
for each $0 \leq j \leq k$, where 
\begin{equation}\label{def:Delta:perp}
\Delta^\perp_j(P_j) := \bigcup_{\substack{W' \in \W^\perp(P) \\ \dim(W'^\perp) \ge j}} \interior\big( \Delta(P_j,W') \big).
\end{equation}
Note that the case $j = k$ of~\eqref{eq:sideways:IH:Deltas} holds because $\Delta^\perp_k(P_k) = \emptyset$, since $\dim(W'^\perp) \le k - 1$ for every $W' \in \W^\perp(P)$, by~\eqref{def:Wperp} and Lemma~\ref{lem:faces:lower:dim}, and that the case $j = 0$ implies~\eqref{eq:sideways:step:aim}. Our main challenge will therefore be to prove the following lemma, which provides the induction step. 

\begin{lemma}\label{lemma:inductionstep:sidewaysgrowth:det}
Let $0 \le j \le k - 1$, and let $W' \in \W^\perp(P)$ with $\dim(W'^\perp) = j$. Then
\begin{equation}\label{eq:sideways-step-part2-stp1}
\interior\big( \Delta(P_j,W') \big) \cap \Z^d \subset \big[ \Delta^\perp_{j+1}(P_{j+1}) \cup P \cup (P' \cap A) \big]_{\U[W]}^{P'}.
\end{equation}
\end{lemma}

The proof of Lemma~\ref{lemma:inductionstep:sidewaysgrowth:det} is unfortunately rather technical, involving the introduction of several further polytopes (and their interiors), the most important of which will be 
$$Q := \interior\big( \Delta(P_{j+1},W') \big).$$
To help the reader negotiate this proliferation of polytopes, we have used variants of `$P$' for $k$-dimensional polytopes (that are contained in a translation of $W^\perp$) and variants of `$Q$' for $j$-dimensional polytopes (that are contained in a translation of $W'^\perp$).

Lemma~\ref{lemma:inductionstep:sidewaysgrowth:det} will be proved in two steps; the first is an application of Lemma~\ref{lem:surface}. 

\begin{lemma}\label{lem:sideways:surface}
Let $0 \le j \le k - 1$, and let $W' \in \W^\perp(P)$ with $\dim(W'^\perp) = j$. Then
\begin{equation}\label{eq:sideways:surface}
\big[ Q \cap A \big]_{\U[W']}^{Q} \subset \big[ \Delta^\perp_{j+1}(P_{j+1}) \cup P \cup (P' \cap A) \big]_{\U[W]}^{P'},
\end{equation}
where $Q := \interior\big( \Delta(P_{j+1},W') \big)$. 
\end{lemma}

\begin{proof}
Set $P^\circ := \clint(P_{j+1} \to W')$, and recall that $P^\circ \subset P_{j+1}$ and  
$$Q \cap \Z^d = \Delta( P^\circ,W') \cap \Z^d,$$
by Lemma~\ref{lem:clint:Delta}. Thus, applying Lemma~\ref{lem:surface} to $P^\circ \in \cP(W)$, we obtain
$$\big[ Q \cap A \big]^{Q}_{\U[W']} \subset \big[ \big( B(P^\circ,W') \setminus Q \big) \cup ( Q \cap A ) \big]^{B(P^\circ,W')}_{\U[W]}.$$
Since $Q \cup P^\circ \subset P_{j+1} \subset P'$, in order to prove~\eqref{eq:sideways:surface} it therefore suffices to show that
\begin{equation}\label{eq:sideways:step3}
B(P^\circ,W') \cap \Z^d \subset \Delta^\perp_{j+1}(P_{j+1}) \cup P \cup Q. 
\end{equation}

To prove~\eqref{eq:sideways:step3}, observe first that, since $P_{j+1} \subset P'$ and by Lemma~\ref{lem:exterior}, 
\begin{equation}\label{eq:interior:in:interior:in:P}
\interior(P_{j+1}) \cap \Z^d \subset \interior(P') \cap \Z^d \subset P.
\end{equation}
Now let $x \in B(P^\circ,W') \cap \Z^d$, and note that $x \in P_{j+1}$. If $x \not\in P$, then it follows from~\eqref{eq:interior:in:interior:in:P} that $x \in P_{j+1} \setminus \interior(P_{j+1})$, and thus $x \in \Delta( P_{j+1}, W \cup \{u\} )$ for some $u \in N_\QQ(W)$.

Let $W \cup \{u\} \subset W'' \in \W$ be maximal such that $x \in \Delta(P_{j+1},W'')$, and observe that 
$$x \in \interior \big( \Delta(P_{j+1},W'') \big),$$
by Lemma~\ref{lem:max:face:interior}. We also have $W'' \subset W'$, by Lemma~\ref{lem:x-notin-Delta:Wv}, since $x \in P^\circ$, so
$$x \notin \Delta\big( P_{j+1}, W \cup \{v\} \big)$$
for every $v \in \QQ \setminus W'$. Therefore $W \subsetneq W'' \subset W' \in \W^\perp(P)$, and hence $W'' \in \W^\perp(P)$.

We are now done, since if $W'' = W'$, then 
$$x \in \interior \big( \Delta(P_{j+1},W') \big) = Q,$$
and if $W'' \ne W'$, then $\dim(W''^\perp) > \dim(W'^\perp) = j$, by Lemma~\ref{lem:faces:lower:dim}, and hence 
$$x \in \interior\big( \Delta(P_{j+1},W'') \big) \subset \Delta^\perp_{j+1}(P_{j+1}),$$
as required. This proves~\eqref{eq:sideways:step3}, and hence completes the proof of the lemma.
\end{proof}

The second step uses Lemma~\ref{lem:forwards-step-general}, as well as our assumptions (in Lemma~\ref{lem:sideways-step-part2}) that $P'$ is long and sideways edge-filled by $A$. 

\begin{lemma}\label{lem:sideways-step-part1}
Let $0 \le j \le k - 1$, and let $W' \in \W^\perp(P)$ with $\dim(W'^\perp) = j$. Then
\begin{equation}\label{eq:sideways-stage-event}
\interior\big( \Delta(P_j,W') \big) \cap \Z^d \subset \big[ Q \cap A \big]_{\U[W']}^{Q},
\end{equation}
where $Q = \interior\big( \Delta(P_{j+1},W') \big)$. 
\end{lemma}

In outline, the proof of this lemma is straightforward: we shall use property~$(a)$ of Definition~\ref{def:sideways-edge-filled} to find a polytope $Q_0 \in \cP(W',w,t^*)$ such that 
$$Q_0 \cap \Z^d \subset \big[ Q \cap A \big]_{\U[W']}^{Q} \setminus P_j,$$ 
and then grow this polytope using property~$(b)$ of the definition and Lemma~\ref{lem:forwards-step-general}. However, since checking the details carefully requires some (tedious) technical calculations, we postpone a few of the details to Appendix~\ref{sideways:app}. 

\begin{proof}[Proof of Lemma~\ref{lem:sideways-step-part1}]
Recall that, since $P$ and $P'$ satisfy the conditions of Lemma~\ref{lem:sideways-step-part2}, the polytope $P' = \ext(P)$ is long and sideways edge-filled by $A$. Fix $W' \in \W^\perp(P)$, and set $Q^* := \clint\big( \Delta(P',W') \big)$ and $t^* := t(Q^*)$. We begin with the following simple claim. 

\begin{claim}\label{claim:Qplus:minusQ}
There exists a polytope $Q' \in \cP\big( W',w; t^*, \alpha \tau' \big)$ with 
$$Q' \subset \big( Q^* \cap Q \big) \setminus P_j.$$ 
\end{claim}

\begin{clmproof}{claim:Qplus:minusQ}
To see that the claim is plausible, recall that $P'$ is long, and that
$$\tau(P_{j+1}) = \tau(P_j) + 4\alpha \tau'.$$
Checking the details is straightforward, but requires a slightly tedious calculation. For completeness, we provide the details in Appendix~\ref{sideways:app}. 
\end{clmproof}

Since $P'$ is sideways edge-filled by $A$, it follows from Definition~\ref{def:sideways-edge-filled} that $Q'$ is internally half-filled. By Definition~\ref{def:int-half-filled}, it follows that
\begin{equation}\label{eq:covertube:startpoint}
Q_0 \cap \Z^d \subset \big[ Q' \cap A \big]_{\U[W']}^{Q'} \subset \big[ Q \cap A \big]_{\U[W']}^{Q}
\end{equation}
for some polytope $Q_0 \in \cP\big( W',w;t^* \big)$ with $Q_0 \subset Q' \subset Q$.

To complete the proof, we use Lemma~\ref{lem:forwards-step-general} and Definition~\ref{def:sideways-edge-filled} to extend $Q_0$ in direction $w$ all the way to the end of $Q^*$, and hence infect the entire interior of $\Delta(P_j,W')$. To be precise, recalling Definition~\ref{def:forwards:growth:sequence}, let $\G$ be the forwards growth sequence with seed $Q_0$, so for each $i \ge 1$ we have
$$Q_i = \fext(Q_{i-1}) \in \cP( W',w;t^*).$$
Now, since $P'$ is sideways edge-filled by $A$, it follows from Definition~\ref{def:sideways-edge-filled} that if $Q_i \subset Q^*$, then $Q_i$ is forwards edge-filled. By Lemma~\ref{lem:forwards-step-general}, it follows that
\begin{equation}\label{eq:forwards-step:covertube}
B(Q_i) \cap \Z^d \subset \big[ B(Q_{i-1}) \cup \big( Q_i \cap A \big) \big]^{Q_i}_{\U[W']}
\end{equation}
for every $i \ge 1$ such that $Q_i \subset Q^*$. Let $m \in \N$ be maximal such that $Q_m \subset Q^*$. Recall that $B(Q_0) \subset Q_0 \subset Q \cap Q^*$, and therefore\footnote{To see that $Q_m \subset Q$, recall that $P_{j+1}$ contains the forward half of $P'$.} $Q_i \subset Q_m \subset Q \cap Q^*$ for every $0 \le i \le m$. It follows, by~\eqref{eq:forwards-step:covertube} and induction, that
\begin{equation}\label{eq:sideways-final-1}
B(Q_i) \cap \Z^d \subset \big[ Q_0 \cup \big( Q \cap A \big) \big]^{Q}_{\U[W']}
\end{equation}
for every $0 \le i \le m$. Note also that
\begin{equation}\label{eq:sideways-final-2}
Q_m \cap \Z^d \subset \bigg( Q_0 \cup \bigcup_{i=1}^m B(Q_i) \bigg) \cap \Z^d,
\end{equation}
since $Q_i = \fext(Q_{i-1})$, so by Lemma~\ref{lem:fext} we have 
\[
\big( Q_i \setminus Q_{i-1} \big) \cap \Z^d \subset \Delta(Q_i) \subset B(Q_i).
\]

To complete the proof of~\eqref{eq:sideways-stage-event}, we therefore only need to prove the following claim.

\begin{claim}\label{claim:extensions:cover:tube}
$$\interior\big( \Delta(P_j,W') \big) \cap \Z^d \subset Q_m.$$
\end{claim}

Indeed, by~\eqref{eq:covertube:startpoint},~\eqref{eq:sideways-final-1} and~\eqref{eq:sideways-final-2}, it will then follow that
\begin{equation}\label{eq:sideways-stage-event:repeated}
\interior\big( \Delta(P_j,W') \big) \cap \Z^d \subset \big[ Q_0 \cup (Q \cap A) \big]_{\U[W']}^Q  \subset \big[ Q \cap A \big]_{\U[W']}^Q,
\end{equation}
as required. We provide a sketch of the proof of Claim~\ref{claim:extensions:cover:tube} below; the full (slightly tedious) details can be found in Appendix~\ref{sideways:app}. 

\begin{clmproof}{claim:extensions:cover:tube}
Applying Lemma~\ref{lem:poly-seq-eats-end} to the polytopes $Q_0$ and 
$$Q^* \cap \big( Q^* + c w \big) = P\big( W',w,a(Q_0),t^*,\tau_1 \big),$$ 
where $c$ and $\tau_1$ are chosen so that this is the case (see Lemma~\ref{lem:Pzero:properties1}), we obtain 
$$Q^* \cap \big( Q^* + c w \big) \cap \Z^d \subset Q_m.$$
Now, since $Q_0 \subset Q \setminus P_j$, and by Lemma~\ref{lem:interiors-are-canonical}, we have 
$$\interior\big( \Delta(P_j,W') \big) \cap \Z^d \subset Q^* \cap \big( Q^* + c w \big),$$
which completes the proof of the claim.
\end{clmproof}

As noted above, combining Claim~\ref{claim:extensions:cover:tube} with~\eqref{eq:covertube:startpoint},~\eqref{eq:sideways-final-1} and~\eqref{eq:sideways-final-2}, we obtain~\eqref{eq:sideways-stage-event:repeated}, as required.  
\end{proof}

Lemma~\ref{lemma:inductionstep:sidewaysgrowth:det} follows immediately from Lemmas~\ref{lem:sideways:surface} and~\ref{lem:sideways-step-part1}. 

\begin{proof}[Proof of Lemma~\ref{lemma:inductionstep:sidewaysgrowth:det}]
Set $Q := \interior\big( \Delta(P_{j+1},W') \big)$. By Lemma~\ref{lem:sideways-step-part1} we have 
$$\interior\big( \Delta(P_j,W') \big) \cap \Z^d \subset \big[ Q \cap A \big]_{\U[W']}^{Q},$$
and by Lemma~\ref{lem:sideways:surface} we have
$$\big[ Q \cap A \big]_{\U[W']}^{Q} \subset \big[ \Delta^\perp_{j+1}(P_{j+1}) \cup P \cup (P' \cap A) \big]_{\U[W]}^{P'},$$
so~\eqref{eq:sideways-step-part2-stp1} follows. 
\end{proof}

The deduction of Lemma~\ref{lem:sideways-step-part2} was already sketched earlier, but for the reader's convenience let us repeat the details. 

\begin{proof}[Proof of Lemma~\ref{lem:sideways-step-part2}]
Recall that $P' = \ext(P)$ and~$P'' = \fret(P')$, and suppose that $P'$ is long and sideways edge-filled by $A$. By Lemma~\ref{lemma:inductionstep:sidewaysgrowth:det}, and recalling~\eqref{eq:Pjs:nested} and~\eqref{def:Delta:perp}, we have
$$\Delta^\perp_j(P_j) \cap \Z^d \subset \big[ \Delta^\perp_{j+1}(P_{j+1}) \cup P \cup (P' \cap A) \big]_{\U[W]}^{P'}$$
for each $0 \le j \le k - 1$. Since $\Delta^\perp_k(P_k) = \emptyset$, by~\eqref{def:Wperp} and Lemma~\ref{lem:faces:lower:dim}, and by~\eqref{eq:Pjs:nested}, it follows that 
$$\interior\big( \Delta(F,W') \big) \big) \cap \Z^d \subset \big[ P \cup (P' \cap A) \big]_{\U[W]}^{P'}$$
for every $W' \in \W^\perp(P)$. 
Since
$$B(P'') \cap \Z^d \subset P \cup \bigcup_{W' \in \W^\perp(P)} \interior\big( \Delta(F,W') \big) \big),$$ 
by Lemma~\ref{lem:sideways:event}, this completes the proof of the lemma. 
\end{proof}

We can now easily deduce Lemma~\ref{lem:sideways-step-general} from Lemmas~\ref{lem:forwards-step-general} and~\ref{lem:sideways-step-part2}. 

\begin{proof}[Proof of Lemma~\ref{lem:sideways-step-general}]
Recall that $W \in \W$ and $w \in \cL_R \cap \SS(W)$, and that $P \in \cP(W,w)$ and $P' = \ext(P)$, and set $P'' := \fret(P')$. Since $P'$ is long and sideways edge-filled, it follows from Lemma~\ref{lem:sideways-step-part2} that
\begin{equation}\label{sideways-step-general:proof2}
B(P'') \cap \Z^d \subset \big[ P \cup \big( P' \cap A \big) \big]_{\U[W]}^{P'}.
\end{equation}
Now, by Lemmas~\ref{lem:fret} and~\ref{lem:fext:fret:buffers}, either
$$B(P') \cap \Z^d \subset B(P'') \cap \Z^d \qquad \text{or} \qquad \fext(P'') = P'.$$
If $B(P') \cap \Z^d \subset B(P'') \cap \Z^d$ then we are done. On the other hand, if $\fext(P'') = P'$, then since $P'$ is forwards edge-filled, we have 
$$B(P') \cap \Z^d \subset \big[ B(P'') \cup \big( P' \cap A \big) \big]^{P'}_{\U[W]}.$$
by Lemma~\ref{lem:forwards-step-general}. Combining this with~\eqref{sideways-step-general:proof2} gives~\eqref{eq:sideways-event}, as required. 
\end{proof}

\subsection{Growing both forwards and backwards}

To finish this section, we prove one further lemma, which will be used in Section~\ref{proof:sec} in order to grow both forwards and `backwards' in a polytope, and hence deduce that an internally half-filled polytope is in fact internally filled. It is a straightforward consequence of Lemma~\ref{lem:forwards-step-general}.

\begin{lemma}\label{lem:extensions:cover:tube}
Let $W \in \W$, let $w \in \cL_R \cap \SS(W)$, let $P \in \cP(W,w)$, and let $A \subset \Z^d$. Suppose that $I_W^\circ(P)$ holds, and that every $Q \in \cP(W,w) \cup \cP(W,-w)$ with $t(Q) = t(P)$ and $Q \subset P$ is forwards edge-filled by $A$. Then the event $I_W^\bullet(P)$ holds. \end{lemma}

\begin{proof}
Let $P = P(W,w;a,t,\tau)$, and recall that, by Definition~\ref{def:int-half-filled}, if $I_W^\circ(P)$ holds then there exists a polytope $Q_0 = P(W,w;a_0,t,\tau_0) \subset P$ such that $Q_0 \cap \Z^d \subset [P \cap A]^P_{\U[W]}$. Now, let $Q_0' := Q_0^- \in \cP(W,-w)$ as in~\eqref{def:Pminus}, so that we have $Q_0' = Q_0$ (as subsets of $\R^d$). Let $\G$ and $\G'$ be the forwards growth sequences with seeds $Q_0$ and $Q_0'$ respectively, so
$$Q_i = \fext(Q_{i-1}) \qquad \text{and} \qquad Q'_i = \fext(Q'_{i-1})$$
for each $i \ge 1$. Let $m \in \N$ be maximal such that $Q_m \subset P$, and observe that, since $Q_i \in \cP(W,w)$ and $t(Q_i) = t$, it follows from our assumptions that, for each $1 \le i \le m$, $Q_i$ is forwards edge-filled by $A$. Hence, by Lemma~\ref{lem:forwards-step-general}, we have
$$B(Q_i) \cap \Z^d \subset \big[ B(Q_{i-1}) \cup \big( P \cap A \big) \big]^{P}_{\U[W]}$$
for every $i \in [m]$. Moreover, by Lemma~\ref{lem:fext}, we have 
$$\big( Q_i \setminus Q_{i-1} \big) \cap \Z^d \subset \Delta( Q_i )  \subset B( Q_i )$$
for every $i \in \N$, and therefore
$$Q_m \cap \Z^d \subset \big[ Q_0 \cup ( P \cap A ) \big]^{P}_{\U[W]}.$$
Similarly $Q_i' \in \cP(W,-w)$ and $t(Q_i') = t$, so by the same argument it follows that 
$$Q'_{m'} \cap \Z^d \subset \big[ Q_0 \cup ( P \cap A ) \big]^{P}_{\U[W]},$$
where $m' \in \N$ is maximal such that $Q'_{m'} \subset P$. By Lemma~\ref{lem:poly-seq-eats-end}, we have
$$P \cap \Z^d = Q_{m} \cup Q'_{m'},$$
and, recalling that $Q_0 \subset \big[ P \cap A \big]^P_{\U[W]}$, we therefore obtain
$$P \cap \Z^d = \big[ P \cap A \big]_{\U[W]}^P.$$
as required.
\end{proof}

\pagebreak

\section{The proof of Theorem~\ref{thm:upper}}\label{proof:sec}

In this section we complete the proof of our main theorem by constructing, for each $W \in \W$, a `low-energy' route to infecting the sites of $\Z^d$ inside each (sufficiently large) polytope in $\cP(W)$. We construct these routes inductively (our induction hypothesis is given in Definition~\ref{def:ih}, below), using the results of Sections~\ref{sec:induced}--\ref{sec:deterministic}. 

In order to state the induction hypothesis, which will be the focus of most of this section, we need to define two functions. Recall from~\eqref{def:constants} that $C$ is a sufficiently large constant depending on $\U$ and, for each $1 \le k \le d$, set\footnote{Since we cannot remove the dependence of our bounds on $C$, we have made no attempt to optimize the constants in this section, opting instead to simplify the presentation.} 
$$\lambda(k) := (8d)^k \cdot C^{3d}.$$
Now, for each $0 \leq s \leq k \leq d$, define 
\begin{equation}\label{def:t0}
t_0(k,s,p) := 
\begin{cases} 
\exp_{(s-1)}\big( p^{-\lambda(k)} \big) & \text{if } s \geq 1 \\ 
C & \text{if } s = 0, 
\end{cases}
\end{equation}
and for each $1 \leq s \leq k \leq d$, set
\begin{equation}\label{def:t1}
t_1(k,s,p) := p^{-t_0(k-1,s-1,p)^{3d}}.
\end{equation}
Recall Definitions~\ref{def:sgood},~\ref{def:curlyW},~\ref{def:int:filled} and~\ref{def:int-half-filled}. Our induction hypothesis is as follows.

\begin{definition}\label{def:ih}
For each $1 \leq s \leq k \leq d$, let $\ih(k,s) = \iha(k,s) \wedge \ihb(k,s)$ be the statement that for every $W \in \W_k$, and all sufficiently small $p > 0$, the following two properties hold: 
\begin{itemize}
\item $\iha(k,s)$: If $W$ is $s$-good, then
\begin{equation}\label{eq:upper-case-a}
\P_{2^k p} \big( I_W^\bullet(P) \big) \ge 1 - \exp\bigg( - \frac{t(P)}{t_0(k,s,p)} \bigg)
\end{equation}
for every $P \in \cP(W)$ such that $t(P) \ge t_0(k,s,p)$ and $\tau(P) = 0$.\medskip
\item $\ihb(k,s)$: If $W$ is $(s,w)$-semi-good for some $w \in \cL_R \cap \SS(W)$, then 
\begin{equation}\label{eq:upper-case-b}
\P_{2^k p} \big( I_W^\circ(P) \big) \ge 1 - e^{- t(P)}
\end{equation}
for every $P \in \cP(W,w)$ such that $\tau(P) \ge t(P)^{2k} \cdot t_1(k,s,p)$.
\end{itemize}
\end{definition}

Before beginning the proof of the induction hypothesis, let us note some simple properties of the functions $t_0(k,s,p)$ and $t_1(k,s,p)$. Each inequality follows easily from~\eqref{def:t0} and~\eqref{def:t1}, so we postpone the details to Appendix~\ref{sideways:app}. 

\begin{obs}\label{obs:t:properties}
Let\/ $1 \le s \le k \le d$, and let\/ $p > 0$ be sufficiently small. Then \begin{equation}\label{eq:t:properties:a}
t_0(k,s,p) \ge t_1(k,s,p)^{8d}.
\end{equation}
If $s^* = \min\{s, k-1\}$, then
\begin{equation}\label{eq:t:properties:c}
t_0(k,s,p) \ge t_0(k-1,s^*,p)^{8d},
\end{equation}
and if $s \ge 2$, then
\begin{equation}\label{eq:t:properties:b}
t_0(k-1,s-1,p) \ge 2 \cdot \log t_1(k-1,s^*,p).
\end{equation}
Moreover, $t_0(1,1,p) \ge t_1(1,1,p) \ge p^{-2R_0}$.
\end{obs}

We remark that we shall also use the bound $t_1(k,s,p) \cdot p^{O(t_0(k-1,s-1,p)^{2d})} \gg 1$, which follows immediately\footnote{Throughout this section, all constants that are implicit in our uses of $O$-notation are allowed to depend on $\U$ (and on $R$ and $\QQ$), but \emph{not} on the constant $C$.} from the definitions, in the proof of Lemma~\ref{lem:IHb:likely}. 

\subsection{The base cases: supercritical growth}\label{basecase:IH:sec}

The base cases of our induction will consist of the statements $\ih(1,1)$ and $\ihb(k,1)$, both of which correspond to `supercritical' growth on a face. We begin with the case $d = 1$, which follows easily from the definitions.

\begin{lemma}\label{lem:ihbase}
$\ih(1,1)$ holds.
\end{lemma}

\begin{proof}
Let $W \in \W$ with $\dim(W^\perp) = 1$, and note that $\SS(W) = \{w,-w\}$ for some $w \in \SS^{d-1}$. To prove $\ihb(1,1)$, suppose that $w \in \cL_R$ and that $W$ is $(1,w)$-semi-good, and observe that, by Definition~\ref{def:sgood}, this means that $\rho^{0} \big( \SS(W); \S_W, w \big) = 0$. By Definition~\ref{def:r}, it follows that $w \not\in \S_W$, and (recalling~\eqref{def:SW}) this implies that there exists a rule $X \in \U[W]$ such that $\< x,w \> < 0$ for all $x \in X$.

Now, let $P \in \cP(W,w)$ with $\tau(P) \ge t(P)^2 \cdot t_1(1,1,p)$, and observe that $P$ is a line segment (by Lemma~\ref{lem:clique-of-correct-dim}) of length $\Theta(\tau(P))$, and moreover that $|P \cap \Z^d| = \Theta(\tau(P))$, where the implicit constants depend only on $R$.\footnote{Indeed, the constants only depend on $w$, and for each $R$ there are a finite number of $w \in \cL_R \cap \SS^{d-1}$.} Since $\|x\| \le R_0$ for every $x \in X$, and recalling that $X \subset W^\perp$ and $\dim(W^\perp) = 1$, it follows that if $A$ contains $R_0$ consecutive sites of $\Z^d$ in $P \setminus F(P)$ (recall from~\eqref{def:forward:half} that $F(P)$ is the forward half of $P$), then $P$ is internally half-filled. Since
$$\tau(P) \ge t(P)^2 \cdot t_1(1,1,p) \ge t(P)^2 \cdot p^{-2R_0},$$ 
by Observation~\ref{obs:t:properties}, it follows that
\begin{equation}\label{eq:ih11}
1 - \P_{2p}\big( I_W^\circ(P) \big) \le \big( 1 - (2p)^{R_0} \big)^{\Omega(\tau(P))} \le \exp\big( - \Omega\big( p^{R_0} \tau(P) \big) \big) \le e^{-t(P)}
\end{equation}
for all sufficiently small $p > 0$, as required. This completes the proof of $\ihb(1,1)$.

The proof of $\iha(1,1)$ is very similar. Indeed, if $W$ is $1$-good and $\dim(W^\perp) = 1$ then, by Definition~\ref{def:sgood}, we have $\rho^{0}\big( \SS(W); \S_W, v \big) = 0$ for every $v \in \SS(W) = \{w,-w\}$. By~\eqref{eq:rho}, 
it follows that $w,-w \not\in \S_W$, and (as above) this implies that there exist rules $X,X' \in \U[W]$ such that $\< x,w \> < 0$ for all $x \in X$, and $\< x,w \> > 0$ for all $x \in X'$. 

Now, let $P \in \cP(W)$ with $t(P) \ge t_0(1,1,p)$ and $\tau(P) = 0$, and observe (as above) that $P$ is a line of length $\Theta(t(P))$. Since $W \subset \QQ$, it follows that $|P \cap \Z^d| = \Theta(t(P))$, where the implicit constants depend only on $\QQ$. Moreover, if $A$ contains $R_0$ consecutive sites of $P \cap \Z^d$, then $P$ is internally filled, using the rules $X$ and $X'$. Since $t_0(1,1,p) \ge p^{-2R_0}$, it follows that
$$1 - \P_{2p}\big( I_W^\bullet(P) \big) \le \big( 1 - p^{R_0} \big)^{\Omega(t(P))} \le \exp\bigg( - \frac{t(P)}{t_0(1,1,p)} \bigg),$$
completing the proof of $\iha(1,1)$. 
\end{proof}

Next we prove $\ihb(k,1)$ for $k \ge 2$, which requires a slightly different (and simpler) argument than the one that will be used for the case $s \ge 2$. The reason for this is that when $s = 1$ we do not need to (and, in fact, cannot) grow our droplets sideways. To grow forwards, we need the following simple consequence of Lemma~\ref{lem:1good}.

\begin{lemma}\label{lem:forwards-step:supercritical}
Let $W \in \W_k$ and $w \in \cL_R \cap \SS(W)$ be such that $W$ is $(1,w)$-semi-good. If~$P \in \cP(W,w)$, then $P$ is forwards edge-filled. 
\end{lemma}

\begin{proof}
By Lemma~\ref{lem:1good}, we have $\emptyset \in \U[W \cup \{u\}]$ for every $u \in N_\QQ(W,w)$, since $w \in W^\perp$, so $\<u,w\> > 0$ implies that $u \notin \< W \>$. By~\eqref{def:WtoP} and Lemma~\ref{lem:emptyset:for:subfaces}, it follows that $\emptyset \in \U[W']$ for every $W' \in \W^\to(P)$, and hence that 
$$I_{W'}^\bullet\Big( \clint\big( \Delta(P,W') \big) \Big)$$ 
holds (deterministically) for every $W' \in \W^\to(P)$, for every set $A$. By Definition~\ref{def:forwards-edge-filled}, it follows that $P$ is forwards edge-filled, as required.
\end{proof}

We can now prove $\ihb(k,1)$. To do so, we simply find a constant-size `seed' on each translate of the line $\<w\>$, and use Lemmas~\ref{lem:forwards-step-general} and~\ref{lem:forwards-step:supercritical} to show that these seeds grow forwards to infect a suitable polytope. In the proof we do not need to approximate the constants too precisely, since our bounds will hold with room to spare. 

\begin{lemma}\label{lem:ihbase:supercritical}
$\ihb(k,1)$ holds for every $1 \le k \le d$.
\end{lemma}

\begin{proof}
Note that, since $\ihb(1,1)$ holds by Lemma~\ref{lem:ihbase}, we may assume that $k \ge 2$. Let $W \in \W$ and $w \in \cL_R \cap \SS(W)$ be such that $\dim(W^\perp) = k$ and 
$W$ is $(1,w)$-semi-good. Let $P = P(W,w;a,t,\tau)$, with $\tau \ge t^{2k} \cdot t_1(k,1,p)$, and set\footnote{We think of $P'$ as being the `middle third' of $P$, cf.~the definition~\eqref{def:forward:half} of the forwards half of $P$.}
\begin{equation}\label{def:middle:third}
P' = P\big(W,w;a+\tau w / 3,t,\tau/3 \big).
\end{equation}
We shall show that $P' \cap \Z^d \not\subset \big[ P \cap A \big]_{\U[W]}^P$ with probability (in $\Pr_{2^kp}$) at most $e^{-t}$. 

Let $E$ denote the event that for each $x \in P' \cap \Z^d$, there exists $a(x) \in \Z^d + W^\perp$ such that 
$$x \in P\big(W,w;a(x),t',\tau\big) \qquad \text{and} \qquad P\big(W,w;a(x),t',0\big) \cap \Z^d \subset P \cap A,$$
where $t' = \min\{t,2C\}$. To bound the probability of $E$, let $a_0(x) \in \Z^d + W^\perp$ be such that $x \in P(W,w;a_0(x),t',0) \subset P$, and observe that, by Lemma~\ref{lem:tube-is-tube} (and since the diameter of $t' \cdot P(W,w)$ is bounded), there exists a set of $\Omega(\tau)$ values of $a(x)$ on the line $a_0(x) + \< w \>$  such that $x \in P(W,w;a(x),t',\tau)$ and such that the corresponding polytopes $P(W,w;a(x),t',0)$ are disjoint and contained in $P$. Now, observe that $P(W,w;a(x),t',0)$ contains fewer than $C^{2d}$ vertices of the lattice $\Z^d$, and therefore these are all contained in the $(2^kp)$-random set $A$ with probability at least $p^{C^{2d}}$. It follows that the probability that there exists $x \in P' \cap \Z^d$ for which no such $a(x)$ exists is at most
$$O\big( t^k \cdot \tau \big) \cdot \big( 1 - p^{C^{2d}} \big)^{\Omega(\tau)} \le e^{-t}$$
if $p$ is sufficiently small, since $|P' \cap \Z^d| = O\big( t^k \cdot \tau \big)$ and
$$\tau \ge t^{2k} \cdot t_1(k,1,p) = t^{2k} \cdot  p^{-C^{3d}}.$$

We next claim that the event $E$ implies that $P' \cap \Z^d \subset \big[ P \cap A \big]_{\U[W]}^P$. To see this, let $x \in P'$, and set $Q_0 := P(W,w;a(x),t',0)$, so $Q_0 \subset P \cap A$ and $x$ is contained in some member of the forwards growth sequence with seed $Q_0$. Since $W$ is $(1,w)$-semi-good, it follows by Lemma~\ref{lem:forwards-step:supercritical} that $Q_i = \fext(Q_{i-1})$ is forwards edge-filled (by the empty set) for every $i \in \N$. By Lemma~\ref{lem:forwards-step-general}, it follows that if $Q_i \subset P$, then
$$B(Q_i) \cap \Z^d \subset \big[ B(Q_{i-1}) \big]^P_{\U[W]}.$$
Hence, if $m \in \N$ is minimal such that $x \in Q_m$, then 
$$x \in Q_m \subset \big[ P \cap A \big]^P_{\U[W]}.$$
Since $x \in P' \cap \Z^d$ was arbitrary, it follows that $P' \cap \Z^d \subset \big[ P \cap A \big]_{\U[W]}^P$, as required.
\end{proof}

The remainder of this section is divided into three parts. In Section~\ref{tubular:IH:sec} we prove the induction step for $\ihb(k,s)$; in Section~\ref{spherical:IH:sec}, we deduce the induction step for $\iha(k,s)$; and in Section~\ref{final:proof:sec} we use the induction hypothesis $\ihb(d,r)$ to prove Theorem~\ref{thm:upper}.

\subsection{Tubular droplets}\label{tubular:IH:sec}

The aim of this subsection is to prove the following lemma.

\begin{lemma}\label{lem:ihb}
Let $2 \leq s \leq k \leq d$. Suppose that $\ih(k',s')$ holds for all $1 \leq s' \leq k' \leq d$ such that $s' \leq s$, $k' \leq k$ and $(k',s') \neq (k,s)$. Then $\ihb(k,s)$ holds.
\end{lemma}

In order to avoid repetition, let us fix $2 \le s \le k \le d$ until the end of the proof of Lemma~\ref{lem:ihb}, and assume that $\ih(k',s')$ holds for all $1 \leq s' \leq k' \leq d$ such that $s' \leq s$, $k' \leq k$ and $(k',s') \neq (k,s)$. 

Recall the statements of Lemmas~\ref{lem:forwards-step-general} and~\ref{lem:sideways-step-general}. We begin the proof by using the induction hypothesis to bound the probability that a polytope is forwards edge-filled.

\begin{lemma}\label{lem:forwards-step-probability}
Let $W \in \W_k$ and $w \in \cL_R \cap \SS(W)$ be such that $W$ is $(s,w)$-semi-good. If~$P \in \cP(W,w)$, then
\begin{equation}\label{eq:forwards:probability}
\Pr_{2^{k-1} p}\big( P \text{ is forwards edge-filled by } A \big) \ge 1 - O(1) \cdot \exp\bigg( - \frac{t(P)}{t_0(k-1,s-1,p)} \bigg).
\end{equation}
\end{lemma}

\begin{proof}
Recall from Definition~\ref{def:forwards-edge-filled} that $P$ is forwards edge-filled by $A$ if the event 
$$I_{W'}^\bullet\Big( \clint\big( \Delta(P,W') \big) \Big)$$ 
holds for every $W' \in \W^\to(P)$. Fix $W' \in \W^\to(P)$, and set $k' := \dim(W'^\perp)$. We shall apply the induction hypothesis to a polytope $Q \in \cP(W')$ such that $Q = \clint\big( \Delta(P,W') \big)$ (as subsets of $\R^d$) and $\tau(Q) = 0$. To prove that such a polytope exists, we first apply Lemma~\ref{lem:forwardsfaces:tau}, which implies that, since $W' \in \W^\to(P)$, there exists $Q' \in \cP(W')$ with $\tau(Q') = 0$ such that $Q' = \Delta(P,W')$ (as subsets of $\R^d$). We may now choose\footnote{To be precise, we simply use the same $y \in \interior\big( P(W') \big)$ and $\eps > 0$ to define both $\clint(Q')$ and $\clint\big( \Delta(P,W') \big)$.} the closed interior of $Q'$ so that $\clint(Q') = \clint\big( \Delta(P,W') \big)$, and by Definition~\ref{def:closed:interior} and Lemma~\ref{lem:interiors-are-canonical}, it follows that $Q := \clint(Q')$ satisfies $Q \in \cP(W')$ and $\tau(Q) = 0$, as claimed. 

Now, recalling~\eqref{def:WtoP}, let $u \in W' \cap N_\QQ(W,w)$, and observe that $u \notin \< W \>$ and $k' < k$, by Lemma~\ref{lem:faces:lower:dim}. Since $s \ge 2$ and $u \in N_\QQ(W,w)$, it follows from Lemma~\ref{lem:sgood} that the set $W \cup \{u\}$ is $(s-1)$-good. If $k'= 0$, then it follows by Lemma~\ref{lem:good:faces} that $\emptyset \in \U[W']$, and in this case we are done, as in the proof of Lemma~\ref{lem:forwards-step:supercritical}. We may therefore assume that $k' > 0$, and hence, by Lemma~\ref{lem:good:faces}, that $W'$ is~$s'$-good, where $s' := \min\{s-1,k'\}$. 

We may assume that $t(P) > t_0(k-1,s-1,p) \ge t_0(k',s',p)$, since otherwise the bound~\eqref{eq:forwards:probability} holds trivially. By~\eqref{def:clint}, it follows\footnote{Recall that we chose $\eps > 0$ sufficiently small in~\eqref{def:clint}.} that $t(Q) = t(P) - \eps \ge t_0(k',s',p)$. Hence, by $\ih_a(k',s')$, we obtain
\begin{equation}\label{eq:forwards:oneface}
\P_{2^{k'} p}\big( I_{W'}^\bullet(Q) \big) \ge 1 - \exp\bigg( - \frac{t(Q)}{t_0(k',s',p)} \bigg).
\end{equation}
Taking a union bound over the $O(1)$ sets $W' \in \W^\to(P)$, the lemma follows.
\end{proof}

A similar argument bounds the probability that a polytope is sideways edge-filled. 

\begin{lemma}\label{lem:sideways-step-probability}
Let $W \in \W_k$ and $w \in \cL_R \cap \SS(W)$ be such that $W$ is $(s,w)$-semi-good. If~$P \in \cP(W,w)$ and 
\begin{equation}\label{eq:sideways-step-conditions}
\tau(P) \ge 8d \cdot t(P)^{2k-2} \cdot t_1(k-1,s^*,p),
\end{equation}
where $s^* = \min\{s, k-1\}$, then
$$\Pr_{2^{k-1} p}\big( P \text{ is sideways edge-filled by } A \big) \ge 1 - O\big( \tau(P)^2 \big) \cdot \exp\bigg( - \frac{t(P)}{t_0(k-2,s^*-1,p)} \bigg).$$
\end{lemma}

The following observation will be used in the proof of Lemma~\ref{lem:sideways-step-probability}. 

\begin{lemma}\label{lem:counting:polytopes}
Let $W \in \W$ and $P = P(W,w;a,t,\tau) \in \cP(W)$. There are $O( \tau^2 + 1 )$ sets $Q \cap \Z^d$ such that $Q \in \cP(W,w;t)$ and $Q \subset P$.
\end{lemma}

\begin{proof}
To count the sets $Q \cap \Z^d$ as described, observe that each polytope $Q \in \cP(W,w;t)$ with $Q \subset P$ can be obtained from $P$ by changing $\tau$ and translating in direction $w$. By Lemma~\ref{lem:growing:linearly:forwards}, adjacent values of $\tau$ where the set changes differ, on average, by at least $\xi$, and we therefore have at most $O( \tau + 1 )$ choices for the set $Q \cap \Z^d$, given $a(Q)$. Similarly (e.g., by applying the same argument to the polytope $P^- = P(W,-w;a+\tau w, t, \tau)$), adjacent values of $a$ where the set changes also differ, on average, by at least $\xi$, and we therefore have at most $O( \tau + 1 )$ choices for $Q \cap \Z^d$, given $\tau(Q)$. It follows that we have at most $O( \tau^2 + 1 )$ choices for the set $Q \cap \Z^d$, as claimed.
\end{proof}

\begin{proof}[Proof of Lemma~\ref{lem:sideways-step-probability}]
Recall from Definition~\ref{def:sideways-edge-filled} that $P$ is sideways edge-filled by $A$ if the following events hold for every $W' \in \W^\perp(P)$:
\begin{enumerate}
\item every polytope $Q \in \cP\big( W',w;t^*,\alpha\tau(P) \big)$ with $Q \subset Q^*$ is internally half-filled;\smallskip
\item every polytope $Q \in \cP\big( W',w;t^* \big)$ with $Q \subset Q^*$ is forwards edge-filled;
\end{enumerate}
where $Q^* = \clint\big( \Delta(P,W') \big)$ and $t^* = t(Q^*)$. We bound the probability of property $(a)$ using  the induction hypothesis, and the probability of property $(b)$ using Lemmas~\ref{lem:forwards-step:supercritical} and~\ref{lem:forwards-step-probability}.\footnote{Note that although 
we stated Lemma~\ref{lem:forwards-step-probability} for the pair ($k,s)$, exactly the same proof works for all pairs $(k',s')$ with $k' \le k$, $s' \le s$ and $2 \le s' \le k'$.} To do this, let us fix $W' \in \W^\perp(P)$, and set $k' := \dim(W'^\perp)$. Recall from~\eqref{def:Wperp} that $W \subsetneq W' \subset \{w\}^\perp$ and $W' \in \W$, and therefore $1 \le k' < k$, by Lemma~\ref{lem:faces:lower:dim}.

Since $W$ is $(s,w)$-semi-good and $W' \in \W^\perp(P)$, it follows by Lemma~\ref{lem:semigood:faces} that $W'$ is $(s',w)$-semi-good, where $s' = \min\{s,k'\}$. Now, since $k' < k$ and $s' \le s^* = \min\{s, k-1\}$, and recalling that $\alpha = 1/8d$, it follows from~\eqref{eq:sideways-step-conditions} that
$$\alpha \cdot \tau(P) \ge (t^*)^{2k'} \cdot t_1(k',s',p).$$ 
Since $w \in W'^\perp$, it therefore follows by $\ih_b(k',s')$ that
$$\P_{2^{k'} p}\big( I_{W'}^\circ(Q) \big) \ge 1 - e^{-t^*}$$
for each $Q \in \cP\big( W',w;t^*,\alpha\tau(P) \big)$ with $Q \subset Q^*$. Moreover, by Lemmas~\ref{lem:forwards-step:supercritical} and~\ref{lem:forwards-step-probability}, 
$$\Pr_{2^{k'-1} p}\big( Q \text{ is forwards edge-filled by } A \big) \ge 1 - O(1) \cdot \exp\bigg( - \frac{t^*}{t_0(k'-1,s'-1,p)} \bigg)$$
for each $Q \in \cP( W',w;t^* )$ with $Q \subset Q^*$. 

Finally, by Lemma~\ref{lem:counting:polytopes}, there are $O\big( \tau(P)^2 + 1 \big)$ choices for the set $Q \cap \Z^d$ such that $Q \in \cP( W',w;t^* )$ with $Q \subset Q^*$. Taking a union bound over these sets, and the $O(1)$ choices of $W' \in \W^\perp(P)$, and noting that $t^* \ge t(P) - 1$, the lemma follows.
\end{proof}

We can now show that the event $I_W^\circ(P)$ occurs with high probability. It will be straightforward to `bootstrap' this result to give the exponential bound we require. 

\begin{lemma}\label{lem:IHb:likely}
Let $W \in \W_k$ and $w \in \cL_R \cap \SS(W)$, and let $P \in \cP(W,w)$ be such that
\begin{equation}\label{eq:IHb:likely:condition}
\tau(P) \ge t(P)^{2k-1} \cdot t_1(k,s,p).
\end{equation}
If\/ $W$ is $(s,w)$-semi-good, then\/ $\P_{2^k p} \big( I_W^\circ(P) \big) \to 1$ as $p \to 0$. 
\end{lemma}

The deduction of Lemma~\ref{lem:IHb:likely} from the lemmas above is not difficult, but the details are a little technical, so let us first give a brief sketch of the construction. Roughly speaking, we find a `seed' polytope $Q$ contained in $A$, and then grow it forwards and sideways using Lemmas~\ref{lem:forwards-step-general} and~\ref{lem:sideways-step-general}. We shall bound, for each possible choice of $Q$, the probability that we fail to grow in a given step with Lemmas~\ref{lem:forwards-step-probability} and~\ref{lem:sideways-step-probability}, using sprinkling to maintain independence. 

Recall from Definition~\ref{def:growth:sequence} that a growth sequence $\G$ with seed $Q \in \cP(W,w)$ is obtained by setting $Q^{(1)}_0 = Q$ and defining 
$$Q^{(i)}_j := \fext\big( Q^{(i)}_{j-1} \big) \qquad \text{and} \qquad Q^{(i+1)}_0 := \ext\big( Q^{(i)}_* \big)$$
for each $i,j \ge 1$, where $Q^{(i)}_* \in \cP(W,w)$ satisfies $t(Q^{(i)}_*) = t(Q^{(i)}_0)$ and $Q^{(i)}_* \subset Q^{(i)}_{\ell(i)} \setminus Q^{(i)}_0$. Here $\ell(i)$ and $Q^{(i)}_*$ will be chosen so that we can apply Lemma~\ref{lem:sideways-step-probability}. 

In the proof of Lemma~\ref{lem:IHb:likely} we will show that, with high probability, 
$$B( Q^{(i)}_j ) \cap \Z^d \subset \big[ B(Q^{(1)}_0) \cup (P \cap A) \big]_{\U[W]}^P$$
for every $i$ and $j$ such that $Q^{(i)}_j \subset P$. To deduce the lemma, we shall prove that there exists $m \in \N$, with $Q^{(m)}_{\ell(m)} \subset P$, such that 
$$P' \cap \Z^d \subset \bigcup_{j = 1}^{\ell(m)} B( Q^{(m)}_j )$$
for some polytope $P' \in \cP(W,w)$ with $P' \subset P$ and $t(P') = t(P)$. 

\begin{proof}[Proof of Lemma~\ref{lem:IHb:likely}]
As noted above, sprinkling will be used to maintain independence between the two stages of our construction. Let us therefore set $q := 2^{k-1}p$, and let $A_1$ and $A_2$ be independent $q$-random subsets of $\Z^d$, noting that therefore $A := A_1 \cup A_2$ is a $q'$-random subset of $\Z^d$, where $q' := 2q-q^2 < 2^k p$.
 
Let $P = P(W,w;a,t,\tau)$, and let $E$ denote the event that there exists a polytope 
$Q \in \cP(W,w)$, with $a(P) - a(Q) \in \<w\>$, such that 
\begin{equation}\label{eq:Q0:properties}
t(Q) = t(1) \qquad \text{and} \qquad Q \cap \Z^d \subset A_1 \cap P\big( W,w;a+\tau w/3,t,\tau/3 \big),
\end{equation}
(cf.~\eqref{def:middle:third}), where
\begin{equation}\label{eq:Q0:properties2}
t(1) := \min\big\{ t_0(k-1,s-1,p)^2, \, t(P) \big\}.
\end{equation}
Note that there exists a collection of $\Omega\big( \tau(P) / t(1) \big) \geq \Omega\big( \tau(P) / t(P) \big)$ disjoint polytopes $Q \in \cP(W,w;t(1))$ in the middle third of $P$, and that for each such polytope $Q$, the event $Q \subset A_1$ occurs with probability $p^{O(t(1)^d)}$. Since
\begin{equation}\label{eq:finding:Q0}
\frac{\tau(P)}{t(P)} \cdot p^{O(t(1)^d)} \ge t_1(k,s,p) \cdot p^{O(t_0(k-1,s-1,p)^{2d})} \gg 1,
\end{equation}
by~\eqref{eq:IHb:likely:condition} and~\eqref{def:t1}, it follows by Chernoff's inequality that $\Pr(E) \to 1$ as $p \to 0$. 

Fix such a polytope $Q$, and let us assume that $t(Q) = t_0(k-1,s-1,p)^2$, since if $t(Q) = t(P)$ then we are already done. We claim that there exists a happy growth sequence $\G$ with seed $Q$ such that
\begin{equation}\label{eq:Qstar:length}
\tau\big( Q^{(i)}_* \big) = 9d \cdot t( Q^{(i)}_0 )^{2k-2} \cdot t_1(k-1,s^*,p) + O(1)
\end{equation}
for every $i \in \N$, where $s^* := \min\{s, k-1\}$. Indeed, to construct $\G$ we simply let $\ell(i) \in \N$ be minimal such that there exists a grower $Q^{(i)}_* \subset Q^{(i)}_{\ell(i)} \setminus Q^{(i)}_0$ with $t(Q^{(i)}_*) = t( Q^{(i)}_0 )$ and such that~\eqref{eq:Qstar:length} holds, noting that, by Lemma~\ref{lem:finding:a:grower}, the condition that $Q^{(i)}_*$ is a grower is not difficult to satisfy. Since $\G$ is happy, it follows by Lemma~\ref{lem:growing:linearly:sideways} that
\begin{equation}\label{eq:speed:of:growing:sideways}
t(Q^{(i)}_0) \ge t(Q) + \xi' \cdot i
\end{equation}
for every $i \ge 1/\xi'$. Now, if $m \in \N$ is maximal such that $t(Q^{(m)}_0) \le t = t(P)$, then\footnote{Note that $t(Q^{(m)}_0)$ is not necessarily equal to $t(P)$, since we did not assume that $P$ is the minimal polytope containing the set $P \cap \Z^d$.} 
\begin{equation}\label{eq:Qm:contains:poly}
P' \cap \Z^d \subset Q^{(m)}_*
\end{equation}
for some $P' \in \cP(W,w;t)$, by Lemma~\ref{lem:growth-seq-eats-sides}, and
\begin{equation}\label{eq:m-at-most-t(P)}
m = O\big( t(P) \big),
\end{equation}
by~\eqref{eq:speed:of:growing:sideways}. The following claim will suffice to complete the proof of the lemma. 

\begin{claim}\label{claim:filling:buffers:Qsequence}
With high probability as $p \to 0$, we have
\begin{equation}\label{eq:filling:buffers:Qsequence}
Q^{(m)}_* \cap \Z^d \subset \big[ Q \cup (P \cap A_2) \big]_{\U[W]}^P.
\end{equation}
\end{claim}

\begin{clmproof}{claim:filling:buffers:Qsequence}
We shall first show that $Q^{(m)}_* \subset P$. To do so, recall that $Q \cap \Z^d$ is contained in the middle third of $P$, and that $t(Q^{(m)}_0) \le t$ and $a(P) - a(Q) \in \<w\>$. Now, by~\eqref{eq:IHb:likely:condition},~\eqref{eq:Qstar:length} and~\eqref{eq:m-at-most-t(P)}, and since $t(Q_0^{(i)}) \leq t(P)$ for all $i$, we have
\begin{align*}
\sum_{i = 1}^m \tau\big( Q^{(i)}_* \big) &\le m \cdot O\big( t(P)^{2k-2} \big) \cdot t_1(k-1,s^*,p) \\
&\leq O\big( t(P)^{2k-1} \big) \cdot t_1(k-1,s,p) \\[+0.3cm]
&\ll \tau(P),
\end{align*}
so $Q^{(m)}_* \subset P$, and moreover $Q^{(i)}_j \subset P$ for every $1 \le i \le m$ and $0 \le j \le \ell(i)$.

Now, by Lemma~\ref{lem:forwards-step-general}, if $Q^{(i)}_j \subset P$ and $Q^{(i)}_j$ is forwards edge-filled by $A_2$ for some $i,j \ge 1$, then 
\begin{equation}\label{eq:forwards:Qs}
B( Q^{(i)}_j ) \cap \Z^d \subset \big[ B( Q^{(i)}_{j-1} ) \cup ( P \cap A_2 ) \big]^{P}_{\U[W]}.
\end{equation}
For each $i \in \N$, set $t(i) := t( Q^{(i)}_0 )$. By Lemma~\ref{lem:forwards-step-probability}, it follows that the probability that~\eqref{eq:forwards:Qs} fails to hold for some $i \in \N$ and $j \in [\ell(i)]$ with $Q^{(i)}_j \subset P$ is at most 
\begin{equation}\label{eq:filling:buffers:Qsequence:bound1}
O(1) \cdot \sum_{i \ge 1} \ell(i) \cdot \exp\bigg( - \frac{t(i)}{t_0(k-1,s-1,p)} \bigg).
\end{equation}
To bound this sum, observe first that, by Lemma~\ref{lem:growing:linearly:forwards} and~\eqref{eq:Qstar:length}, we have
$$\ell(i) = O\big( t(i)^{2k-2} \cdot t_1(k-1,s^*,p) \big)$$
for every $i \ge 1$. Since 
$$t(i) \ge t(1) + \xi' \cdot i = t_0(k-1,s-1,p)^2 + \xi' \cdot i$$
for every $i \ge 1/\xi'$, by~\eqref{eq:speed:of:growing:sideways}, it follows that~\eqref{eq:filling:buffers:Qsequence:bound1} is at most\footnote{Here we use the fact that $\int_{c^2}^\infty x^{O(1)} e^{-x/c} \,dx = c^{O(1)} e^{-c}$.}
\begin{equation}\label{eq:filling:buffers:Qsequence:bound2}
t_0(k-1,s-1,p)^{O(1)} \cdot t_1(k-1,s^*,p) \cdot \exp\big( - t_0(k-1,s-1,p) \big). 
\end{equation}
Now, by Observation~\ref{obs:t:properties} (in particular, by~\eqref{eq:t:properties:b}),~\eqref{eq:filling:buffers:Qsequence:bound2} tends to zero as $p \to 0$. Thus, with high probability,~\eqref{eq:forwards:Qs} holds for all $i \in \N$ and $j \in [\ell(i)]$ such that $Q^{(i)}_j \subset P$.

Now, by Lemma~\ref{lem:fext} and the definition of $\G$, we have 
$$Q^{(i)}_* \cap \Z^d \subset \bigcup_{j = 1}^{\ell(i)} \big( Q^{(i)}_{j} \setminus Q^{(i)}_{j-1} \big) \cap \Z^d \subset \bigcup_{j = 1}^{\ell(i)} \Delta( Q^{(i)}_j ).$$
In particular, recalling that $\Delta(P) \subset B(P)$, it follows from~\eqref{eq:forwards:Qs} that
\begin{equation}\label{eq:using:forwards:Qs}
Q^{(i)}_* \cap \Z^d \subset \big[ B( Q^{(i)}_{0} ) \cup ( P \cap A_2 ) \big]^{P}_{\U[W]}
\end{equation}
for every $1 \le i \le m$. It therefore remains to show that, with high probability,  
\begin{equation}\label{eq:sideways:Qs}
B( Q^{(i)}_0 ) \cap \Z^d \subset \big[ Q^{(i-1)}_* \cup ( P \cap A_2 ) \big]^{P}_{\U[W]}
\end{equation}
for every $2 \le i \le m$. Since $Q^{(i)}_0$ is long, by~\eqref{eq:Qstar:length}, and $Q^{(i)}_0 \subset P$, this will follow from Lemma~\ref{lem:sideways-step-general} if each polytope $Q^{(i)}_0$ is forwards and sideways edge-filled by $A_2$. 

We have already bounded the probability that $Q^{(i)}_0$ is not forwards edge-filled by $A_2$ above. Moreover, by Lemma~\ref{lem:sideways-step-probability}, and using the bound~\eqref{eq:Qstar:length}, the probability that $Q^{(i)}_0$ is not sideways edge-filled by $A_2$ is at most 
$$O\big( \ell(i-1)^2 \big)  \cdot \exp\bigg( - \frac{t(i)}{t_0(k-2,s^*-1,p)} \bigg),$$
since $t( Q^{(i)}_0 ) = t(i)$ and $\tau( Q^{(i)}_0 ) = O( \ell(i-1) )$. Summing over $i$, and noting that $s^* \le s$, it follows that the probability that $Q^{(i)}_0$ is not sideways edge-filled by $A_2$ for some $2 \le i \le m$ is bounded above by~\eqref{eq:filling:buffers:Qsequence:bound2}. Therefore, by the calculation above, with high probability,~\eqref{eq:sideways:Qs} holds for every $2 \le i \le m$. 

Finally, combining~\eqref{eq:using:forwards:Qs} and~\eqref{eq:sideways:Qs}, we obtain~\eqref{eq:filling:buffers:Qsequence}, as required.
\end{clmproof}

By Claim~\ref{claim:filling:buffers:Qsequence} and~\eqref{eq:Qm:contains:poly}, and recalling that $Q \cap \Z^d \subset A_1$, it follows  that, with high probability,
\begin{equation}\label{eq:Qm:contains:P}
P' \cap \Z^d \subset Q^{(m)}_* \cap \Z^d \subset \big[ Q \cup (P \cap A_2) \big]^P_{\U[W]} \subset \big[ P \cap A \big]^P_{\U[W]} \subset P
\end{equation}
for some $P' \in \cP(W,w;t)$. Since $t(P) = t$, we can moreover choose $P' \subset P$, and hence the event $I_W^\circ(P)$ holds with high probability, as required.
\end{proof}

It is now straightforward to deduce $\ihb(k,s)$ from Lemma~\ref{lem:IHb:likely}. Indeed, to obtain an  exponentially small failure probability we simply make multiple (independent) attempts to find a polytope $P' \in \cP\big(W,w;t(P) \big)$ with $P' \cap \Z^d \subset [P \cap A]_{\U[W]}^P$. 

\begin{proof}[Proof of Lemma~\ref{lem:ihb}]
Let $W \in \W_k$ and $w \in \cL_R \cap \SS(W)$, and let $P \in \cP(W,w)$, with
$$\tau(P) \ge t(P)^{2k} \cdot t_1(k,s,p).$$ 
Suppose that $W$ is $(s,w)$-semi-good, and (recalling that $\delta > 0$ is a small constant) note that $P$ contains at least $\delta \cdot t(P)$ disjoint copies of any polytope $P' \in \cP(W,w)$ with $t(P') = t(P)$ and 
$$\tau(P') = t(P)^{2k-1} \cdot t_1(k,s,p) \le \tau(P) / t(P).$$ 
By Lemma~\ref{lem:IHb:likely}, each of these is internally half-filled with probability $1 - o(1)$, and if any of them is internally half-filled then $P$ is internally half-filled. Hence
$$\P_{2^k p} \big( I_W^\circ(P) \big) \ge 1 - e^{- t(P)},$$
as required.
\end{proof}

\subsection{Spherical droplets}\label{spherical:IH:sec}

The aim of this subsection is to prove the following lemma.

\begin{lemma}\label{lem:iha}
Let $1 \le s \le k \le d$ with $k \ge 2$. Suppose that $\ih(k',s')$ holds for all $1 \leq s' \leq k' \le k$ such that $s' \leq s$ and $(k',s') \neq (k,s)$. Then $\iha(k,s)$ holds.
\end{lemma}

As in the previous subsection, let us fix $1 \le s \le k \le d$ with $k \ge 2$ until the end of the proof of Lemma~\ref{lem:iha}, and assume that $\ih(k',s')$ holds for all $1 \le s' \le k' \le k$ with $s' \leq s$ and $(k',s') \neq (k,s)$. By Lemmas~\ref{lem:ihbase:supercritical} and~\ref{lem:ihb}, it follows that $\ihb(k,s)$ holds.

In order to prove the exponential bound that we need in~\eqref{eq:upper-case-a}, we shall use a variant of the renormalization trick of Schonmann~\cite{Sch1}. Roughly speaking, the idea is to cover our large polytope $P \in \cP(W)$ with copies of a smaller polytope $Q \in \cP(W,w)$ that is internally filled with high probability (see Lemmas~\ref{lem:if-tube} and~\ref{lem:cover:exists}). We then apply a standard percolation argument to show that (with very high probability) the `connected components' of non-internally-filled copies of $Q$ are all small, and use Lemma~\ref{lem:eating:a:set} to fill in the gaps created by these small components (see Lemmas~\ref{cor:eating:polytopes} and~\ref{lem:schonmann}). 

To begin, we shall use the induction hypothesis and $\ihb(k,s)$ to prove the following lemma, which provides us with a suitable `small' polytope $Q$. 

\begin{lemma}\label{lem:if-tube}
Let $W \in \W_k$, let $w \in \cL_R \cap \SS(W)$, and let $Q \in \cP(W,w)$ satisfy
\begin{equation}\label{eq:droplet:size}
t(Q) \ge t_0(k,s,p)^{1/4d}  \qquad \text{and} \qquad \tau(Q) = t(Q)^{2k} \cdot t_1(k,s,p).
\end{equation}
If $W$ is $s$-good and $(s,w)$-semi-good, then $\P_{2^k p} \big( I_W^\bullet(Q) \big) \to 1$ as $p \to 0$.
\end{lemma}

We shall prove Lemma~\ref{lem:if-tube} using Lemma~\ref{lem:extensions:cover:tube}, together with the following simple consequence of Lemmas~\ref{lem:forwards-step:supercritical} and~\ref{lem:forwards-step-probability}, which allows us to bound the probability that a polytope is forwards edge-filled in an arbitrary (rather than just the easiest) direction.

\begin{lemma}\label{lem:forwards-step-probability2}
Let $W \in \W_k$ be $s$-good, and let $Q \in \cP(W)$. Then
\begin{equation}\label{eq:forwards:probability:alldirections}
\Pr_{2^{k-1} p}\big( Q \text{ is forwards edge-filled by } A \big) \ge 1 - O(1) \cdot \exp\bigg( - \frac{t(Q)}{t_0(k-1,s^*,p)} \bigg),
\end{equation}
where $s^* = \min\{s, k-1\}$. 
\end{lemma}

\begin{proof}
Recall from Definition~\ref{def:sgood} that since $W$ is $s$-good, it is $(s^*+1,w)$-semi-good for every $w \in \SS(W)$. The claimed bound now follows immediately from Lemma~\ref{lem:forwards-step:supercritical} (when $s = 1$) and Lemma~\ref{lem:forwards-step-probability} (in the case $s \ge 2$).
\end{proof}

We are now ready to prove that $Q$ is internally filled with high probability.

\begin{proof}[Proof of Lemma~\ref{lem:if-tube}]
Let $E$ denote the event that every $Q' \in \cP(W,w) \cup \cP(W,-w)$ with $t(Q') = t(Q)$ and $Q' \subset Q$ is forwards edge-filled by $A$. By Lemma~\ref{lem:extensions:cover:tube}, the event $I_W^\circ(Q) \cap E$ implies (deterministically) that $Q$ is internally filled by $A$. Moreover, since $W$ is $(s,w)$-semi-good, by $\ihb(k,s)$ and~\eqref{eq:droplet:size} we have
\begin{equation}\label{eq:using:IHB:tube}
\P_{2^k p} \big( I_W^\circ(Q) \big) \to 1 
\end{equation}
as $p \to 0$. It therefore suffices to show that $E$ occurs with high probability as $p \to 0$.

To bound the probability of $E$, let $Q' \in \cP(W,w) \cup \cP(W,-w)$ with $t(Q') = t(Q)$ and $Q' \subset Q$, and observe that, since $W$ is $s$-good, we have
$$\Pr_{2^{k-1} p}\big( Q' \text{ is forwards edge-filled by } A \big) \ge 1 - O(1) \cdot \exp\bigg( - \frac{t(Q)}{t_0(k-1,s^*,p)} \bigg)$$
by Lemma~\ref{lem:forwards-step-probability2}, where $s^* := \min\{s, k-1\}$. Since, by Lemma~\ref{lem:counting:polytopes}, there are at most $O\big( \tau(Q)^2 \big)$ choices for the set $Q' \cap \Z^d$, it follows that
\begin{equation}\label{eq:tube:bounding:E}
\Pr_{2^{k-1} p}(E) \ge 1 - O\big( \tau(Q)^2 \big) \cdot \exp\bigg( - \frac{t(Q)}{t_0(k-1,s^*,p)} \bigg).
\end{equation}
To bound the right-hand side of~\eqref{eq:tube:bounding:E}, observe first that 
\begin{equation}\label{eq:tube-bounding-t}
t(Q) \ge t_0(k,s,p)^{1/4d} \ge t_0(k-1,s^*,p)^{2},
\end{equation}
where the first inequality holds by~\eqref{eq:droplet:size}, and the second by~\eqref{eq:t:properties:c}. It follows that
$$\tau(Q)^2 \cdot \exp\bigg( - \frac{t(Q)}{t_0(k-1,s^*,p)} \bigg) \le t(Q)^{4k} \cdot t_1(k,s,p)^2 \cdot \exp\big( - t(Q)^{1/2} \big) \to 0$$
as $p \to 0$, where in the first step we used~\eqref{eq:droplet:size} and~\eqref{eq:tube-bounding-t} , and in the second we used~\eqref{eq:droplet:size} and~\eqref{eq:t:properties:a}, which together imply that $t(Q) \ge t_1(k,s,p)^2 \gg 1$. 

Combining this with~\eqref{eq:using:IHB:tube} and~\eqref{eq:tube:bounding:E}, it follows that the event $I_W^\circ(Q) \cap E$ occurs with high probability, and therefore that $\P_{2^k p} \big( I_W^\bullet(Q) \big) \to 1$ as $p \to 0$, as required.
\end{proof}

As mentioned above, in order to deduce from Lemma~\ref{lem:if-tube} the exponential failure bound that we need in~\eqref{eq:upper-case-a}, we shall use a variant of the renormalization trick of Schonmann~\cite{Sch1}. In the next few lemmas we develop the tools that we will require for this method. Firstly, we define the concept of a (perfect) cover of a set $P$ with copies of a set $Q$; we remark that in this section $P$ and $Q$ will always be polytopes, but in the proof in Appendix~\ref{cover:app} we will need to cover slightly more general sets.

\begin{definition}\label{def:cover}
Given sets $P,Q \subset \R^d$, we say that $\C$ is a \emph{cover of $P$ by copies of $Q$} if
\begin{enumerate}
\item each member of $\C$ is a translate of $Q$;\smallskip
\item every element of $P$ is contained in some member of $\C$.
\end{enumerate}
Moreover, we say that $\C$ is a \emph{perfect cover of $P$ by copies of $Q$} if in addition
\begin{enumerate}
\item[$(c)$] each member of $\C$ is contained in $P$.
\end{enumerate}
Given a cover $\C$ of $P$ by copies of $Q$, define
a graph $G_\C$ on vertex set $\C$ as follows: 
\begin{equation}\label{def:GC:edges}
E(G_\C) = \big\{ Q_1Q_2 : d(Q_1,Q_2) \le 2R_0 \big\}.
\end{equation}
We will use the (standard) notation $\Delta(G_\C)$ to denote the maximum degree of $G_\C$. 
\end{definition}

In particular, recalling~\eqref{def:strongly:conn}, observe that the diameter of a strongly connected component of a finite set $K \subset P$ can be bounded in terms of the maximum length of a path in $G_\C[\B]$, where $\B$ is the set of vertices of $\C$ that intersect $K$. 

Recall from Definition~\ref{def:ih} that when proving $\iha(k,s)$ our polytope $P \in \cP(W)$ will be such that $t(P) \ge t_0(k,s,p)$. We shall construct a perfect cover of $P$ with copies of the minimal polytope $Q$ satisfying the conditions of Lemma~\ref{lem:if-tube}, whose diameter satisfies 
\begin{equation}\label{eq:diamQ:small:enough}
\diam(Q) = O\big( t(Q) + \tau(Q) \big) = O\big( t_0(k,s,p)^{k/2d} \cdot t_1(k,s,p) \big) \ll t_0(k,s,p)^{2/3},
\end{equation}
where the final step follows from~\eqref{eq:t:properties:a}. In particular, note that $\diam(Q) \ll t(P)$. 

In order to carry out our strategy, we need to know that bounded degree perfect covers exist. The following lemma is proved in Appendix~\ref{cover:app}.

\begin{lemma}\label{lem:cover:exists}
There exists $\Delta > 0$ depending only on $\QQ$ such that the following holds. Let $W \in \W$ and $w \in \cL_R \cap \SS^{d-1}$, and let $P,Q \in \cP(W,w)$ satisfy
\begin{equation}\label{eq:PQ:condition}
\diam(Q) \le \frac{t(P)}{\Delta}.
\end{equation}
Then there exists a perfect cover $\C$ of $P$ by copies of $Q$ such that $\Delta( G_\C ) \le \Delta$. 
\end{lemma}

\pagebreak

We remark that we have stated (and proved) Lemma~\ref{lem:cover:exists} in greater generality than we need for our application, since we expect it to be useful in future work. 

Given a set $A \subset \Z^d$ and a perfect cover $\C$ of $P$ by copies of $Q$, for some polytopes $P$ and $Q$, let us define 
\begin{equation}\label{def:B:not:spanned}
\B = \B(A,\C) := \big\{ Q' \in \C : I_W^\bullet(Q') \text{ does not hold} \, \big\}. 
\end{equation}
The next step is to apply Lemma~\ref{lem:eating:a:set} to show that if $P$ is not internally filled by $A$, then $G_\C[\B]$ must contain a long path. Since the polytope $Q$ to be used for our perfect cover is internally filled with high probability, by Lemma~\ref{lem:if-tube}, we shall be able to show that the existence of such a path is extremely unlikely (see Lemma~\ref{lem:schonmann}). Recall from~\eqref{def:constants} and Lemma~\ref{lem:eating:a:set} that $\delta = \delta(\QQ) > 0$ is a sufficiently small constant. 

\begin{lemma}\label{cor:eating:polytopes}
Let $W \in \W_k$ be $k$-good, let $P,Q \in \cP(W)$, and let $\B = \B(A,\C)$ for some set $A \subset \Z^d$ and some perfect cover $\C$ of $P$ by copies of $Q$. If every path in $G_\C[\B]$ has at most 
$$\frac{\delta \cdot t(P)}{2 \cdot \diam(Q)}$$ 
vertices, then $I_W^\bullet(P)$ holds. 
\end{lemma}

\begin{proof}
Set $K := \bigcup_{Q' \in \B} Q' \cap \Z^d$, and observe that $K \subset P \cap \Z^d$, and that 
$$(P \setminus K) \cap \Z^d \subset \big[ P \cap A \big]^P_{\U[W]},$$
since $\C$ is a perfect cover and $I_W^\bullet(Q')$ holds for every $Q' \in \C \setminus \B$. It follows that if $I_W^\bullet(P)$ does not hold, then $\big[ P \setminus K \big]^P_{\U[W]} \ne P \cap \Z^d$, and hence, by Lemma~\ref{lem:eating:a:set}, that there exists a strongly connected component $K'$ of $K$ with diameter at least $\delta \cdot t(P)$. 

Let $x,y \in K'$ with $d(x,y) \ge \delta \cdot t(P)$ and, recalling~\eqref{def:strongly:conn} and~\eqref{def:GC:edges}, observe that there must exist a path $(Q_1,\ldots,Q_\ell)$ in $\G[\B]$ such that $x \in Q_1$ and $y \in Q_\ell$. It follows that
$$\delta \cdot t(P) \le d(x,y) \le \ell \cdot \big( \diam(Q) + 2R_0 \big) \le 2\ell \cdot \diam(Q),$$
since $\diam(Q) > t(Q) > C$, by Definition~\ref{def:cPW}. It follows that $G_\C[\B]$ contains a path with at least $(\delta/2) \cdot t(P) / \diam(Q)$ vertices, as required.
\end{proof}

We are finally ready to prove the following Schonmann-type lemma. Let $\beta > 0$ be a sufficiently small constant (in particular, we need it to satisfy $\beta < ( e^{2/\delta} \Delta )^{-(\Delta + 1)}$), and recall that the (minimal) polytope $Q$ from Lemma~\ref{lem:if-tube} satisfies the conditions below. 

\begin{lemma}\label{lem:schonmann}
Let $W \in \W_k$ be $k$-good. If\/\footnote{Strictly speaking, we mean here that the inequality~\eqref{eq:schonmann:condition} holds for all translates of $Q$.}
\begin{equation}\label{eq:schonmann:condition}
\P_{2^k p}\big(I_W^\bullet(Q) \big) \ge 1 - \beta
\end{equation}
where $Q \in \cP(W)$ is such that $\diam(Q) \le t_0(k,s,p)^{2/3}$, then
$$\P_{2^k p}\big(I_W^\bullet(P) \big) \ge 1 - \exp\bigg( - \frac{t(P)}{t_0(k,s,p)} \bigg)$$
for every $P \in \cP(W)$ with $t(P) \ge t_0(k,s,p)$ and $\tau(P) = 0$. 
\end{lemma}

\pagebreak

\begin{proof}
Observe first that, by Lemma~\ref{lem:cover:exists}, and since 
$$\diam(Q) \le t_0(k,s,p)^{2/3} \ll t_0(k,s,p) \le t(P),$$
by~\eqref{eq:diamQ:small:enough}, there exists a perfect cover $\C$ of $P$ by copies of $Q$ such that $\Delta(G_\C) \le \Delta$.\footnote{To be precise, if $w(P) \ne w(Q)$ then we apply Lemma~\ref{lem:cover:exists} to the polytope $P' = P(W,w';a,t,0)$, where $P = P(W,w;a,t,0)$ and $w' = w(Q)$, which is equal to $P$ as a subset of $\R^d$.} Let $A \subset \Z^d$ be $(2^k p)$-random set, and set $\B = \B(A,\C)$, as in~\eqref{def:B:not:spanned}. Note that each of the copies of $Q$ in $\C$ is internally filled with probability at least $1 - \beta$, by~\eqref{eq:schonmann:condition}, and that these events are independent for non-intersecting copies of $Q$.  

By Lemma~\ref{cor:eating:polytopes}, if the polytope $P$ is not internally filled by $A$, then there exists a path in $G_\C[\B]$ with 
$$m \ge \frac{\delta \cdot t(P)}{2 \cdot \diam(Q)}$$ 
vertices. Since $\Delta(G_\C) \le \Delta$, there are at most $|\C| \cdot \Delta^m$ possible such paths, and each contains an independent set in $\G_\C$ of size at least $m / (\Delta+1)$. Since the vertices of an independent set in $\G_\C$ are disjoint copies of $Q$, it follows that the expected number of such paths in $G_\C[\B]$ is at most
$$|\C| \cdot \Delta^m \cdot \beta^{m / (\Delta+1)} \le t(P)^{O(1)} \exp\bigg( - \frac{t(P)}{\diam(Q)} \bigg)$$
since $\beta < ( e^{2/\delta} \Delta )^{-(\Delta+1)}$ and $\tau(P) = 0$, which implies that $|\C| = t(P)^{O(1)}$. 

Since $\diam(Q) \le t_0(k,s,p)^{2/3}$ and $t(P) \ge t_0(k,s,p)$, it follows by Markov's inequality that
$$\P_p\big(I_W^\bullet(P)^c \big) \le \exp\bigg( - \frac{t(P)}{t_0(k,s,p)} \bigg),$$
as required.
\end{proof}

We are now ready to complete the proof of the induction step for $\iha(k,s)$.

\begin{proof}[Proof of Lemma~\ref{lem:iha}]
Let $W \in \W_k$ be $s$-good, and let $P \in \cP(W)$ with 
$$t(P) \ge t_0(k,s,p) \qquad \text{and} \qquad \tau(P) = 0.$$ 
Note that $\U[W] \ne \emptyset$, since $W$ is $k$-good, and if $\emptyset \in \U[W]$ then~\eqref{eq:upper-case-a} holds trivially; we may therefore assume that $\U[W]$ is non-trivial. By 
Lemma~\ref{lem:sgood-rational-w}, it follows that there exists $w \in \cL_R \cap \SS(W)$ such that $W$ is $(s,w)$-semi-good. We want to apply Lemmas~\ref{lem:if-tube} and~\ref{lem:schonmann}, so let $Q \in \cP(W,w)$ satisfy
$$t(Q) = t_0(k,s,p)^{1/4d} \qquad \text{and} \qquad \tau(Q) = t(Q)^{2k} \cdot t_1(k,s,p).$$
Since $W$ is $s$-good and $(s,w)$-semi-good, it follows by Lemma~\ref{lem:if-tube} that\footnote{Moreover, the lemma implies the same bound for every translate of $Q$.}
$$\P_{2^k p}\big(I_W^\bullet(Q) \big) \to 1$$
as $p \to 0$. Observe that 
$$\diam(Q) = O\big( t(Q) + \tau(Q) \big) \ll t_0(k,s,p)^{2/3},$$
by~\eqref{eq:diamQ:small:enough}. Hence, by Lemma~\ref{lem:schonmann}, we have
$$\P_{2^k p}\big(I_W^\bullet(P) \big) \ge 1 - \exp\bigg( - \frac{t(P)}{t_0(k,s,p)} \bigg)$$
as required.
\end{proof}

\subsection{The proof of Theorem~\ref{thm:upper}}\label{final:proof:sec}

Combining Lemmas~\ref{lem:ihbase},~\ref{lem:ihbase:supercritical},~\ref{lem:ihb} and~\ref{lem:iha}, we immediately obtain the following proposition. 

\begin{prop}\label{prop:upper-ih}
$\ih(k,s)$ holds for all $1 \leq s \leq k \leq d$. \qed
\end{prop}

We are now in a position to complete the proof of Theorem~\ref{thm:upper}. The only remaining difficulty is that $r(\U) = r$ only implies (by Lemma~\ref{lem:Q-exists-rational-w}) that there exists $w \in \cL_R \cap \SS^{d-1}$ such that $\emptyset$ is $(r,w)$-semi-good; it may not be the case that $\emptyset$ is $r$-good. 

There are various ways to deal with this; for example, we could repeat the proof of Lemma~\ref{lem:IHb:likely}, growing the droplet until it fills the entire torus. We have chosen, however, a simpler (if somewhat less efficient) way of deducing the theorem from Proposition~\ref{prop:upper-ih}: for each site $x \in \Z_n^d$, we shall find an internally half-filled polytope on the line $x + \<w\>$, and then use Lemmas~\ref{lem:forwards-step-general} and~\ref{lem:forwards-step-probability} to grow in direction $w$ until we infect $x$.

\begin{proof}[Proof of Theorem~\ref{thm:upper}]
Let $1 \le r \le d$, and let $\U$ be a $d$-dimensional update family with $r(\U) = r$. By Lemma~\ref{lem:Q-exists-rational-w}, there exists $w \in \cL_R \cap \SS^{d-1}$ such that $\emptyset$ is $(r,w)$-semi-good. Let $n \in \N$ be sufficiently large, and set
\begin{equation}\label{eq:choice:of:p}
p := \bigg( \frac{1}{\log_{(r-1)} n} \bigg)^{\eps}
\end{equation}
for some sufficiently small constant $\eps = \eps(\U) > 0$. In order to avoid dealing with the geometry of the torus, let $A$ be a $2^d p$-random subset of $\{-3n,\ldots,4n\}^d$. We will show that 
$$\{1,\ldots,n\}^d \subset [A]_\U$$ 
with high probability, and then use this to deduce the claimed bound on $p_c( \Z_n^d,\U )$. 

Observe first that, by~\eqref{def:t0},~\eqref{eq:t:properties:a} and~\eqref{eq:choice:of:p}, we have
$$t_1(d,r,p) \le t_0(d,r,p) = \exp_{(r-1)}\big( p^{-\lambda(d)} \big) \le n^{1/2}$$
since $\eps$ is sufficiently small (in particular $\eps < 1 / \lambda(d)$). 
Now, for each $a \in \Z^d$, define 
$$P(a) := P(\emptyset,w;a,t,\tau),$$ 
where $\tau := t^{2d} \cdot t_1(d,r,p)$, and $t$ is chosen so that $\diam(P(a)) = n/4$. Note that $t \ge n^{1/6d}$, since $\tau = \Theta(n)$ and $t_1(d,r,p) \le n^{1/2}$. 

The following claim is an easy consequence of the induction hypothesis. 

\begin{claim}\label{claim:finalproof1}
With high probability, the event\/ $I_\emptyset^\circ\big( P(a) \big)$ holds for every\/ $a \in \{-n,\ldots,2n\}^d$. 
\end{claim}

\begin{clmproof}{claim:finalproof1}
By Proposition~\ref{prop:upper-ih}, $\ih(d,r)$ holds. Since $\emptyset$ is $(r,w)$-semi-good and $\tau = t^{2d} \cdot t_1(d,r,p)$, it follows, by Definition~\ref{def:ih}, that 
$$\P_{2^d p} \Big( I_\emptyset^\circ\big( P(a) \big) \Big) \ge 1 - e^{-t}$$
for each $a \in \{-n,\ldots,2n\}^d$. Hence, by the union bound, the probability (in $\P_{2^d p}$) that there exists $a \in \{-n,\ldots,2n\}^d$ such that $I_\emptyset^\circ\big( P(a) \big)$ fails to hold is at most
$$(3n+1)^d \cdot e^{-t} \le (3n+1)^d \cdot \exp\big( - n^{1/6d} \big) \to 0$$
as $n \to \infty$, as required. 
\end{clmproof}

Note that $a \in P(a)$, and therefore $P(a) \subset [-2n,3n]^d$ for every $a \in \{-n,\ldots,2n\}^d$. It therefore follows from Claim~\ref{claim:finalproof1} and Definition~\ref{def:int-half-filled} that, with high probability, for every $a \in \{-n,\ldots,2n\}^d$ there exists $Q(a) \in \cP(\emptyset,w;t)$ such that 
$$Q(a) \subset P(a) \qquad \text{and} \qquad Q(a) \cap \Z^d \subset \big[ P(a) \cap A \big]_{\U}.$$
Next we use Lemma~\ref{lem:forwards-step-probability} to show that, with high probability, every such polytope is forwards edge-filled by $A$.

\begin{claim}\label{claim:finalproof2}
With high probability, every polytope $Q \in \cP(\emptyset,w;t)$ with $Q \subset [-2n,3n]^d$ is forwards edge-filled by $A$.
\end{claim}

\begin{clmproof}{claim:finalproof2}
Since $\emptyset$ is $(r,w)$-semi-good, it follows by Lemmas~\ref{lem:forwards-step:supercritical} and~\ref{lem:forwards-step-probability} that for each $Q \in \cP(\emptyset,w;t)$, the probability that $Q$ is not forwards edge-filled by $A$ is at most
$$O(1) \cdot \exp\bigg( - \frac{t}{t_0(d-1,r-1,p)} \bigg).$$
Note that, by~\eqref{def:t0} and~\eqref{eq:choice:of:p}, we have
$$t_0(d-1,r-1,p) = n^{o(1)}.$$
Moreover, treating two polytopes as equivalent if the events that they are forwards edge-filled are identical,\footnote{That is, $Q_1$ and $Q_2$ are equivalent (for the purpose of Claim~\ref{claim:finalproof2}) if $Q_1$ is forwards edge-filled by $A$ if and only if $Q_2$ is forwards edge-filled by $A$.} there are at most $n^{O(1)}$ choices for the polytope $Q$ (cf. the proof of Lemma~\ref{lem:counting:polytopes}). Since $t \ge n^{1/6d}$, it follows that the expected number of $Q \in \cP(\emptyset,w;t)$ with $Q \subset [-2n,3n]^d$ that are not forwards edge-filled by $A$ is at most
$$n^{O(1)} \cdot \exp\big( - n^{1/6d - o(1)} \big) \to 0$$
as $n \to \infty$, so the result follows by Markov's inequality.
\end{clmproof}

We complete the proof with the following deterministic claim.

\begin{claim}\label{claim:finalproof3}
If\/ $I_\emptyset^\circ\big( P(a) \big)$ holds for every $a \in \{-n,\ldots,2n\}^d$, and every polytope $Q \in \cP(\emptyset,w;t)$ with $Q \subset [-2n,3n]^d$ is forwards edge-filled by $A$, then $\{1,\ldots,n\}^d \subset [A]_\U$. 
\end{claim}

\begin{clmproof}{claim:finalproof3}
For each $x \in \{1,\ldots,n\}^d$, choose a vertex $a \in \Z^d$ with $x - a \in \<w\>$ and within distance $O(1)$ of $x - (n/2) \cdot w$. The event $I_\emptyset^\circ(P(a))$ implies that there exists a polytope $Q_0 \in \cP(\emptyset,w;t)$ such that 
$$Q_0 \subset [P(a) \cap A]_\U.$$ 
Let $\G = (Q_i)_{i = 0}^\infty$ be the forwards growth sequence with seed $Q_0$, so 
$$Q_i = \fext(Q_{i-1})$$
for each $i \in \N$, and observe that $x \in Q_m$ for some $m \in \N$. Since $Q_i \in \cP(\emptyset,w;t)$ and $Q_i \subset [-2n,3n]^d$ for each $0 \le i \le m$, it follows from our assumption that $Q_i$ is forwards edge-filled by $A$. Thus, by Lemma~\ref{lem:forwards-step-general}, we have
$$B(Q_i) \cap \Z_n^d \subset \big[ B(Q_{i-1}) \cup A \big]_{\U},$$
and hence $Q_i \subset \big[ Q_0 \cup A \big]_{\U}$ for each $1 \le i \le m$. In particular, 
$$x \in Q_m \subset \big[ Q_0 \cup A \big]_\U.$$
Since $x \in \{1,\ldots,n\}^d$ was arbitrary, it follows that $\{1,\ldots,n\}^d \subset [A]_\U$.
\end{clmproof}

Claims~\ref{claim:finalproof1},~\ref{claim:finalproof2} and~\ref{claim:finalproof3} together imply that $\{1,\ldots,n\}^d \subset [A]_\U$ with high probability. Finally, to deduce the bound on $p_c( \Z_n^d,\U )$, observe that if $\phi \colon \Z^d \to \Z_n^d$ is the standard projection onto the torus, then for any set $A \subset \Z^d$, we have
$$\phi\big( [ A ]_\U \big) \subset \big[ \phi(A) \big]_\U.$$
Using this observation, we may couple the $\U$-bootstrap process on the torus, starting from a $p$-random subset of $\Z_n^d$, with the $\U$-bootstrap process on the infinite lattice $\Z^d$, starting from an $8^{-d} p$-random subset of $\{-3n,\ldots,4n\}^d$, and deduce that $[\phi(A) ]_\U = \Z_n^d$ with high probability, as required. By~\eqref{eq:choice:of:p}, and recalling that $n$ is sufficiently large, this completes the proof of Theorem~\ref{thm:upper}. 
\end{proof}

The proof of Theorem~\ref{thm:upper} given above immediately implies the following generalisation of Schonmann's theorem~\cite{Sch1} to arbitrary critical and supercritical models. 

\begin{theorem}\label{thm:not:subcritical}
Let\/ $\U$ be a $d$-dimensional update family. If\/ $\U$ is not subcritical, then 
$$p_c(\Z^d,\U) = 0.$$
\end{theorem}

\begin{proof}
We need to show, for each $p > 0$, that if $A$ is a $p$-random subset of $\Z^d$, then the set $[ A ]_\U$ contains the origin almost surely. To do so, set $A_n := A \cap \{-n,\ldots,n\}^d$ for each $n \in \N$, and observe that $r(\U) \le d$, by Lemma~\ref{lem:tri}, since $\U$ is not subcritical. Repeating the proof of Theorem~\ref{thm:upper}, we deduce from Claims~\ref{claim:finalproof1},~\ref{claim:finalproof2} and~\ref{claim:finalproof3} that
$$\Pr\big( \0 \in [ A_n ]_\U \big) \to 1$$
as $n \to \infty$. Since $n$ was arbitrary, it follows that $\0 \in [ A ]_\U$ almost surely, as required. 
\end{proof}

It also follows easily from the proof of Theorem~\ref{thm:upper} that if $\U$ is a $d$-dimensional update family with $r(\U) = r \le d$, then the expected infection time of the origin satisfies
\begin{equation}\label{eq:boot:hitting:time}
\Ex_p\big[ \min\{ t \ge 0 : \0 \in A_t \} \big] \le \exp_{(r-1)}( p^{-C} )
\end{equation}
for some $C = C(\U) > 0$, and all sufficiently small $p > 0$. Indeed, a careful examination of the proof shows that if $\eps(\U) \log_{(r)}(n) \ge \log (1/p)$, then the probability that the origin is uninfected at time $t \ge n^d$ is at most $\exp( - \sqrt{t} )$, and this implies~\eqref{eq:boot:hitting:time}.

\appendix

\section{Properties of canonical polytopes}\label{polytope:app}

This first appendix contains proofs of the basic properties of the polytopes $P(W,w)$ defined in Section~\ref{sec:polytopes}. We begin with proofs of Lemmas~\ref{lem:tube-is-tube} and~\ref{lem:tube-faces-to-P-faces}; these lemmas will then enable us to prove one of the key lemmas of Section~\ref{sec:polytopes}, Lemma~\ref{lem:tube-nbrs}.

First, let us recall the statement of Lemma~\ref{lem:tube-is-tube}.

\begin{lemma}\label{lem:tube-is-tube:app}
Let $W \subset \QQ$ and let $w \in \cL_R$. If $w \in W^\perp$, then
\begin{equation}\label{eq:tube-is-tube-emptyset:app}
P(W,w) = \bigcup_{\lambda \in [0,1]} \big( P(W) + \lambda w \big).
\end{equation}
\end{lemma}

\begin{proof}
We claim that
\begin{equation}\label{eq:tube-is-tube-emptyset:proof:app}
P(\emptyset,w) = \bigcup_{\lambda \in [0,1]} \big( P(\emptyset) + \lambda w \big).
\end{equation}
This will suffice to prove the lemma, since if we take the intersection of both sides of~\eqref{eq:tube-is-tube-emptyset:proof:app} with the set $\bigcap_{u \in W} \big\{ x \in \R^d : \< x,u \> = 1 \big\}$, then we obtain~\eqref{eq:tube-is-tube-emptyset:app}. Indeed, this follows from~\eqref{eq:P(W)} and~\eqref{eq:P(W,w)}, using the assumption that $w \in W^\perp$.

To prove~\eqref{eq:tube-is-tube-emptyset:proof:app}, we claim first that $P(\emptyset) + \lambda w \subset P(\emptyset,w)$ for each $\lambda \in [0,1]$. To see this, recall that if $x - \lambda w \in P(\emptyset)$, then $\<x - \lambda w, u\> \le 1$ for every $u \in \QQ$. Since $\lambda \ge 0$, it follows that $\< x,u \> \le 1$ if $\<u,w\> < 0$, and since $\lambda \le 1$, it also follows that $\< x - w, u \> \le 1$ if $\<u,w\> \ge 0$. By~\eqref{eq:P(emptyset,w)}, this implies that $x \in P(\emptyset,w)$, as claimed. 

It therefore remains to show that if $x \in P(\emptyset,w)$, then $x \in P(\emptyset) + \lambda w$ for some $\lambda \in [0,1]$. We shall show that if $0 \le \<x,w\> \le \|w\|^2$ then we can take 
$$\lambda = \frac{\<x,w\>}{\| w \|^2},$$
if $\<x,w\> < 0$ then we can take $\lambda = 0$, and that if $\<x,w\> > \|w\|^2$ then we can take $\lambda = 1$.

Suppose first that $\< x,w \> < 0$, and let $u \in \QQ$ be such that $x / \|x\| \in \Cell_\QQ(u)$. By Lemma~\ref{lem:quasi-cell}, it follows that $\< u,w \> \le 0$, and therefore $\< x,u \> \le 1$, since $x \in P(\emptyset,w)$ and recalling~\eqref{eq:P(emptyset,w)}. Note also that $\< x,u \> > 0$, because $\QQ$ intersects every open hemisphere of $\SS^{d-1}$, and observe that therefore $x / \< x,u \> \in P(\emptyset)$, by Lemma~\ref{lem:facets-to-cells}. Since $P(\emptyset)$ is convex and contains $\0$, and recalling that $\< x,u \> \leq 1$, it follows that $x \in P(\emptyset)$, as claimed.

Next, suppose that $\< x,w \> > \|w\|^2$, and let $u \in \QQ$ be such that $(x - w) / \|x-w\| \in \Cell_\QQ(u)$. The proof is now the same as in the previous case. Indeed, since $\<x - w,w\> > 0$, we have $\< u,w \> \ge 0$, by Lemma~\ref{lem:quasi-cell}, and therefore $\< x - w,u \> \le 1$, since $x \in P(\emptyset,w)$ and recalling~\eqref{eq:P(emptyset,w)}. Moreover $\< x - w,u \> > 0$, since $\QQ$ intersects every open hemisphere of $\SS^{d-1}$, and therefore $(x - w) / \< x - w,u \> \in P(\emptyset)$, by Lemma~\ref{lem:facets-to-cells}. Since $P(\emptyset)$ is convex and $\< x - w,u \> \le 1$, it follows that $x - w \in P(\emptyset)$, as claimed.

Finally, suppose that $0 \le \< x,w \> \le \|w\|^2$, and set $\lambda = \<x,w\> / \|w\|^2$. Let $z \in \{w\}^\perp$ be such that $x = \lambda w + z$, and note that $\lambda w \in P(\emptyset) + \lambda w$, so we may assume that $z \ne \0$. Let $u \in \QQ$ be such that $z / \|z\| \in \Cell_\QQ(u)$. Now, recall that $\QQ$ was chosen to satisfy the conclusion of Lemma~\ref{lem:quasi-new}, and note that therefore $\<u,w\> = 0$. By~\eqref{eq:P(emptyset,w)}, it follows that $\< u,z \> = \< x,u \> \le 1$, since $x \in P(\emptyset,w)$. Moreover, since $\< u,z \> > 0$ (again because $\QQ$ intersects every open hemisphere of $\SS^{d-1}$), it follows that $z / \< u,z \> \in P(\emptyset)$, by Lemma~\ref{lem:facets-to-cells}. Since $P(\emptyset)$ is convex, it follows that $z \in P(\emptyset)$. Since $x = \lambda w + z$, this implies that $x \in P(\emptyset) + \lambda w$, as required.
\end{proof}

Our next aim is to prove Lemma~\ref{lem:tube-faces-to-P-faces}. In order to do so, let us first define
\begin{equation}\label{eq:N*}
N_\QQ^\ast(W) := \big\{ u \in \QQ \setminus W \,:\, P(W \cup \{u\})
\neq \emptyset \big\}
\end{equation}
for each clique $W \subset \QQ$. We will need the following slight strengthening of Lemma~\ref{lem:P-nbrs}, which follows from the same proof. 

\begin{lemma}\label{lem:P-nbrs-N*}
Let $W \subset \QQ$ and suppose that $P(W) \neq \emptyset$. Then
\begin{equation}\label{eq:P-nbrs-N*}
P(W) = \bigcap_{u \in W} \big\{ x \in \R^d : \< x,u \> = 1 \big\} \cap
\bigcap_{u \in N_\QQ^\ast(W)} \big\{ x \in \R^d : \< x,u \> \leq 1
\big\}.
\end{equation}
\end{lemma}

\begin{proof}
The lemma follows from the proof of Lemma~\ref{lem:P-nbrs}, noting that if $z \in P(W \cup \{v\})$, then we have $v \in N_\QQ^\ast(W)$, by~\eqref{eq:N*}, while if $x \in P'(W)$ and $\< x,v \> > 1$, then $v \not\in N_\QQ^\ast(W)$. This contradiction implies that $x \in P(W)$, as required.
\end{proof}

We are ready to prove Lemma~\ref{lem:tube-faces-to-P-faces}. We shall actually prove the following  slightly more technical statement, which will be useful in the proof of Lemma~\ref{lem:spheres-equal-iff-tubes-equal:app}. 

\begin{lemma}\label{lem:tube-faces-to-P-faces:app}
Let $W \subset W'\subset \QQ$, with $P(W') \neq \emptyset$. Let $w \in \cL_R$, and suppose that $w \notin (W')^\perp$, and that either $P(W') = P(W)$, or $P(W',w) = P(W,w)$. Then
$$P(W,w) = P(W) + \delta(u,w) w$$
for each $u \in W'$ such that $\< u,w \> \ne 0$.
\end{lemma}

\begin{proof}
We claim first that if $\< u,w \> > 0$ for some $u \in W' \cup N_\QQ(W')$, then $\< u,w \> \ge 0$ for all $u \in W' \cup N_\QQ(W')$.\footnote{Recall the definition~\eqref{eq:QQW} of the set $N_\QQ(W)$.} Indeed, this follows from Lemma~\ref{lem:quasi-innerprod}, since $W'$ is a clique, by Lemma~\ref{lem:faces-to-cliques}, and $w \in \cL_R$. Similarly, if $\< u,w \> < 0$ for some $u \in W' \cup N_\QQ(W')$, then $\< u,w \> \le 0$ for all $u \in W' \cup N_\QQ(W')$. It follows, since $w \notin (W')^\perp$, that exactly one of these two possibilities holds, and thus, if we define
$$\delta' := \prod_{u \in W' \cup  N_\QQ(W')} \delta(u,w),$$
then we are required to show that
\begin{equation}\label{eq:stretchy-stp}
P(W,w) = P(W) + \delta' w.
\end{equation}
The following claim will be used several times.

\begin{claim}\label{clm:stretchy-delta'}
Let $u \in \QQ$ with $\< u,w \> \neq 0$. Suppose either that $u \in W' \cup N_\QQ(W')$, or that $P(W) = P(W')$ and $u \in N_\QQ^\ast(W)$. Then $\delta(u,w) = \delta'$.
\end{claim}

\begin{clmproof}{clm:stretchy-delta'}
If $u \in W' \cup N_\QQ(W')$ then this follows from the observations above. Indeed, if $\< u,w \> < 0$ then $\delta(u,w) = \delta' = 0$, and if $\< u,w \> > 0$ then we have $\< v,w \> \geq 0$ for all $v \in W' \cup N_\QQ(W')$, and therefore $\delta(u,w) = \delta' = 1$.

Let us therefore suppose  that $P(W) = P(W')$, and moreover that $u \in N_\QQ^\ast(W)$, and hence $P(W \cup \{u\}) \neq \emptyset$, by~\eqref{eq:N*}. Let $x \in P(W \cup \{u\})$ and note that $x / \|x\| \in \Cell_\QQ(u)$, by Lemma~\ref{lem:facets-to-cells}, and therefore $\< x,w \> \cdot \< u,w \> > 0$, by Lemma~\ref{lem:quasi-cell}. Now, recalling that $w \notin (W')^\perp$, let $v \in W'$ be such that $\<v,w\> \neq 0$, and observe that $\< x,w \> \cdot \< v,w \> > 0$, by Lemma~\ref{lem:quasi-cell}, since $x \in P(W \cup \{u\}) \subset P(W) = P(W')$, and so $x / \|x\| \in \Cell_\QQ(v)$, by Lemma~\ref{lem:facets-to-cells}. It follows that $\< u,w \> \cdot \< v,w \> > 0$, and hence $\delta(u,w) = \delta(v,w) = \delta'$, by the first part of the claim.
\end{clmproof}

To prove~\eqref{eq:stretchy-stp}, observe first that, by~\eqref{eq:P(W,w)} and Claim~\ref{clm:stretchy-delta'}, we have
\begin{equation}\label{eq:tube-faces-PWw1:app}
P(W,w) = P(\emptyset,w) \cap \bigcap_{u \in W} \big\{ x \in \R^d : \< x - \delta' w, u \> = 1 \big\},
\end{equation}
since $W \subset W'$, so for every $u \in W$ we have either $\<u,w\> = 0$ or $\delta(u,w) = \delta'$. Now, by Lemma~\ref{lem:tube-is-tube:app} applied with $W = \emptyset$, we have
$$P(W) + \delta' w \subset P(\emptyset) + \delta' w \subset P(\emptyset,w).$$
and moreover, by~\eqref{eq:P(W)}, we have
$$P(W) + \delta' w \subset \big\{ x \in \R^d : \< x - \delta' w, u \> = 1 \big\}$$
for every $u \in W$. By~\eqref{eq:tube-faces-PWw1:app}, it follows that $P(W) + \delta' w \subset P(W,w)$.

To show that $P(W,w) \subset P(W) + \delta' w$, we shall need to consider separately the cases $P(W') = P(W)$ and $P(W',w) = P(W,w)$. Suppose first that $P(W') = P(W) \ne \emptyset$, and observe that, by
Lemma~\ref{lem:P-nbrs-N*}, we have
$$P(W) + \delta' w = \bigcap_{u \in W} \big\{ x \in \R^d : \< x - \delta' w, u \> = 1 \big\} \cap \bigcap_{u \in N_\QQ^\ast(W)} \big\{ x \in \R^d : \< x - \delta' w, u \> \leq 1 \big\}.$$
Therefore, recalling~\eqref{eq:tube-faces-PWw1:app}, it suffices to show that
\begin{equation}\label{eq:tube-faces-P-faces-suff}
P(\emptyset,w) \subset \bigcap_{u \in N_\QQ^\ast(W)} \big\{ x \in \R^d
: \< x - \delta' w, u \> \le 1 \big\}.
\end{equation}
But we know from~\eqref{eq:P(emptyset,w)} that
\[
P(\emptyset,w) \subset \bigcap_{u \in N_\QQ^\ast(W)} \big\{ x \in \R^d
: \< x - \delta(u,w) w, u \> \le 1 \big\},
\]
so~\eqref{eq:tube-faces-P-faces-suff} follows from Claim~\ref{clm:stretchy-delta'}.

Suppose now, instead, that $P(W',w) = P(W,w)$. Then, by~\eqref{eq:P(W,w)}, we have
$$P(W,w) = P(\emptyset,w) \cap \bigcap_{u \in W'} \big\{ x \in \R^d : \< x - \delta' w, u \> = 1 \big\}.$$
Note that $P(W') + \delta' w \subset P(W) + \delta' w$, by~\eqref{eq:P(W)}, and that
$$P(W') + \delta' w = \bigcap_{u \in W'} \big\{ x \in \R^d : \< x - \delta' w, u \> = 1 \big\} \cap \bigcap_{u \in N_\QQ(W')} \big\{ x \in \R^d : \< x - \delta' w, u \> \leq 1 \big\},$$
by Lemma~\ref{lem:P-nbrs}, and since $P(W') \neq \emptyset$. Now observe that, by~\eqref{eq:P(emptyset,w)}, we have
\[
P(\emptyset,w) \subset \bigcap_{u \in N_\QQ(W')} \big\{ x \in \R^d :
\< x - \delta' w, u \> \le 1 \big\}
\]
since for each $u \in W' \cup N_\QQ(W')$ we have either $\<u,w\> = 0$ or $\delta(u,w) = \delta'$, by Claim~\ref{clm:stretchy-delta'}. It follows that $P(W,w) \subset P(W) + \delta' w$, as required.
\end{proof}

We can now prove an analogue of Lemma~\ref{lem:faces-to-cliques} for $P(W,w)$. Note, in particular, that we do not create new faces during stretching. 

\begin{lemma}\label{lem:tube-faces-to-cliques:app}
Let $W \subset \QQ$ and let $w \in \cL_R$. Then $P(W,w) \neq \emptyset$ if and only if $P(W) \neq \emptyset$, and if either holds then $W$ is a clique.
\end{lemma}

At first glance it might appear as if Lemma~\ref{lem:tube-faces-to-cliques:app} ought to follow as a simple consequence of Lemmas~\ref{lem:tube-is-tube:app} and~\ref{lem:tube-faces-to-P-faces:app}. However, the latter lemma requires prior knowledge that $P(W) \neq \emptyset$, so only allows us to deduce one direction.

\begin{proof}[Proof of Lemma~\ref{lem:tube-faces-to-cliques:app}]
Recall from Lemma~\ref{lem:faces-to-cliques} that if $P(W) \neq \emptyset$ then $W$ is a clique. It therefore suffices to show that $P(W,w) \neq \emptyset$ if and only if $P(W) \neq \emptyset$. Moreover, it follows from Lemma~\ref{lem:tube-is-tube:app} that if $w \in W^\perp$, then $P(W) \subset P(W,w)$, and from Lemma~\ref{lem:tube-faces-to-P-faces:app} (applied with $W' = W$) that if $w \notin W^\perp$, then $P(W,w) \in \{ P(W), P(W) + w \}$. In either case, $P(W) \neq \emptyset$ implies immediately that $P(W,w) \neq \emptyset$. 

The proof of the reverse implication is not quite so straightforward. Let $x' \in P(W,w)$, and set $x = x' + \lambda w$, where $\lambda$ is minimal\footnote{We shall not need this minimality condition until Claim~\ref{clm:tube-faces-u0:app}, but it is convenient to fix $x$ immediately.}  (noting that $P(W,w)$ is compact) such that $x \in P(W,w)$. The following claim is an easy consequence of Lemma~\ref{lem:tube-faces-to-P-faces:app}.

\begin{claim}\label{clm:tube-faces-u:app}
Let $u \in W$.
\begin{enumerate}
\item If $\< u,w \> < 0$ then $x \in P(\{u\})$.\smallskip
\item If $\< u,w \> > 0$ then $x - w \in P(\{u\})$.
\end{enumerate}
\end{claim}

\begin{clmproof}{clm:tube-faces-u:app}
Note that $u \in P(\{u\})$, and therefore, by Lemma~\ref{lem:tube-faces-to-P-faces:app} (applied with both sets equal to $\{u\}$), if $\< u,w \> < 0$ then $P(\{u\},w) = P(\{u\})$. Since $x \in P(W,w) \subset P(\{u\},w)$, it follows that $x \in P(\{u\})$, as claimed. The proof of part~$(b)$ is similar, except if $\< u,w \> > 0$, then by Lemma~\ref{lem:tube-faces-to-P-faces:app} we have $P(\{u\},w) = P(\{u\}) + w$. Since $x \in P(\{u\},w)$, it follows that $x - w \in P(\{u\})$, as claimed. 
\end{clmproof}

When $\< u,w \> = 0$, the situation is a little more complex; we shall need the minimality of $\lambda$, and we shall use Lemma~\ref{lem:tube-is-tube:app} instead of Lemma~\ref{lem:tube-faces-to-P-faces:app}.  

\begin{claim}\label{clm:tube-faces-u0:app}
Suppose $u \in W$ is such that $\< u,w \> = 0$. 
\begin{enumerate}
\item If $\< v,w \> \le 0$ for every $v \in W$, then $x \in P(\{u\})$.\smallskip
\item If $\< v,w \> > 0$ for some $v \in W$, then $x - w \in P(\{u\})$.
\end{enumerate}
\end{claim}

\begin{clmproof}{clm:tube-faces-u0:app}
Suppose first that $\< v,w \> \le 0$ for every $v \in W$. By our choice of $x$, we have $x - \eps w \notin P(W,w)$ for all $\eps > 0$. We claim that moreover $x - \eps w \notin P(\emptyset,w)$ for all $\eps > 0$. To see this, recall~\eqref{eq:P(W,w)}, and suppose first that there exists $v \in W \setminus \{w\}^\perp$. Since $\< v,w \> < 0$ and $x \in P(W,w)$, we have $\< x - \eps w, v \> > \< x,v \> = 1$ for every $\eps > 0$, and therefore, by~\eqref{eq:P(emptyset,w)}, we have $x - \eps w \notin P(\emptyset,w)$, as claimed. On the other hand, if $W \subset \{w\}^\perp$, then we have $\< x - \mu w,v \> = \< x,v \> = 1$ for every $v \in W$ and every $\mu \in \R$, and therefore $x - \eps w \in P(W,w)$ if and only if $x - \eps w \in P(\emptyset,w)$.

Now, by Lemma~\ref{lem:tube-is-tube:app}, we have 
\begin{equation}\label{eq:tube-is-tube:app}
P(\{u\},w) = \bigcup_{\lambda \in [0,1]} \big( P(\{u\}) + \lambda w \big).
\end{equation}
Since $P(W,w) \subset P(\{u\},w) \subset P(\emptyset,w)$ and $x \in P(W,w)$, it follows from the comments above that $x \in P(\{u\},w)$ but $x - \eps w \notin P(\{u\},w)$ for all $\eps > 0$. Hence $x \in P(\{u\})$, by~\eqref{eq:tube-is-tube:app}, as claimed.

Suppose now that $\< v,w \> > 0$ for some $v \in W$. In this case we claim that $x + \eps w \notin P(\emptyset,w)$ for all $\eps > 0$. Indeed, we have $\< x + \eps w - w, v \> > \< x - w, v \> = 1$ for every $\eps > 0$, since $x \in P(W,w)$, and so, by~\eqref{eq:P(emptyset,w)}, we have $x + \eps w \notin P(\emptyset,w)$. It follows, as before, that $x \in P(\{u\},w)$ and $x + \eps w \notin P(\{u\},w)$ for all $\eps > 0$, and hence, by~\eqref{eq:tube-is-tube:app}, we have $x \in P(\{u\}) + w$. 
\end{clmproof}

To deduce the lemma from the claims, let us first verify that we cannot have $u,v \in W$ such that $\< u,w \> < 0$ and $\< v,w \> > 0$. Indeed if we did, then by Claim~\ref{clm:tube-faces-u:app} and Lemma~\ref{lem:facets-to-cells}, we would have $x / \|x\| \in \Cell_\QQ(u)$ and $(x-w) / \|x-w\| \in \Cell_\QQ(v)$. It then follows, by Lemma~\ref{lem:quasi-cell}, that $\< x,w \> < 0$ and $\< x-w,w \> > 0$, and together these give a contradiction. Thus, either $\< u,w \> \leq 0$ for all $u \in W$ or $\< u,w \> \geq 0$ for all $u \in W$. 

Suppose first that $\< u,w \> \leq 0$ for all $u \in W$. Then by Claims~\ref{clm:tube-faces-u:app} and~\ref{clm:tube-faces-u0:app} we have $x \in P(\{u\})$ for all $u \in W$, and therefore $x \in P(W)$. In particular, it follows that $P(W) \neq \emptyset$. On the other hand, if $\< u,w \> \geq 0$ for all $u \in W$ (and not all are zero), then by Claims~\ref{clm:tube-faces-u:app} and~\ref{clm:tube-faces-u0:app} again we have $x - w \in P(\{u\})$ for all $u \in W$, and therefore $x - w \in P(W)$. Once again, this implies $P(W) \neq \emptyset$, and this completes the proof.
\end{proof}

Lemma~\ref{lem:tube-faces-to-cliques:app} easily implies our first aim of this section, which was Lemma~\ref{lem:tube-nbrs}. We restate that lemma now.

\begin{lemma}\label{lem:tube-nbrs:app}
Let $W \subset \QQ$, and let $w \in \cL_R$. If\/ $P(W,w) \neq \emptyset$, then $W$ is a clique, and
\begin{align}
& P(W,w) = \bigcap_{u \in W} \Big\{ x \in \R^d : \big\< x - \delta(u,w) w, \, u \big\> = 1 \Big\} \nonumber \\ 
& \hspace{5cm} \cap \bigcap_{u \in N_\QQ(W)} \Big\{ x \in \R^d : \big\< x - \delta(u,w) w, \, u \big\> \le 1 \Big\}.\label{eq:tube-nbrs:app}
\end{align}
\end{lemma}

\begin{proof} 
The deduction of this lemma from Lemma~\ref{lem:tube-faces-to-cliques:app} is essentially identical to that of Lemma~\ref{lem:P-nbrs} from Lemma~\ref{lem:faces-to-cliques}, but for the reader's convenience we shall spell out the details. Let us write $P'(W,w)$ for the right-hand side of~\eqref{eq:tube-nbrs:app}, and note that $P(W,w) \subset P'(W,w)$, by~\eqref{eq:P(W,w)}, since $N_\QQ(W) \subset \QQ$. To prove the lemma it is enough therefore to show that if $x \in P'(W,w)$ and $y \in P(W,w)$, then $x \in P(W,w)$.
 
To prove this, let $\lambda \ge 0$ be maximal such that $z := y + \lambda(x-y) \in P(W,w)$, and note first that if $\lambda \ge 1$ then $x \in P(W,w)$, since $P(W,w)$ is convex. On the other hand, if $\lambda < 1$ then, by~\eqref{eq:P(W,w)} and~\eqref{eq:tube-nbrs:app}, there exists $v \in \QQ \setminus (W \cup N_\QQ(W))$ such that 
\[
\< z - \delta(v,w)w, \, v \> = 1 \qquad \text{and} \qquad \< x - \delta(v,w)w, \, v \> > 1.
\]
Since $z \in P(W,w)$, it follows that $z \in P(W \cup \{v\},w)$, and therefore $W \cup \{v\}$ is a clique, by Lemma~\ref{lem:tube-faces-to-cliques:app}. But this is a contradiction, since $v \not\in W \cup N_\QQ(W)$. 
\end{proof}

Recall from Definition~\ref{def:curlyW} the set of `maximal' cliques
\begin{equation}\label{eq:W-reminder}
\W = \big\{ W \subset \QQ : \text{$W$ is a clique and $P(W') \neq P(W)$ for every $W \subsetneq W' \subset \QQ$} \big\}.
\end{equation}
The remainder of this section is primarily concerned with proving properties of this set. First, we need a preliminary lemma that relates faces of $P(W,w)$ to faces of $P(W)$.

\begin{lemma}\label{lem:spheres-equal-iff-tubes-equal:app}
Let $W \subset W' \subset \QQ$ be cliques, let $w \in \cL_R$, and suppose that $P(W) \neq \emptyset$. Then $P(W',w) = P(W,w)$ if and only if $P(W') = P(W)$.
\end{lemma}

\begin{proof}
Suppose first that $P(W') = P(W)$. If $w \in (W')^\perp \subset W^\perp$, then by Lemma~\ref{lem:tube-is-tube:app} we have
\[
P(W',w) = \bigcup_{\lambda \in [0,1]} \big( P(W') + \lambda w \big) = \bigcup_{\lambda \in [0,1]} \big( P(W) + \lambda w \big) = P(W,w),
\]
as required. If $w \notin (W')^\perp$, on the other hand, then we apply Lemma~\ref{lem:tube-faces-to-P-faces:app} twice, first with both sets equal to $W'$, and then to the pair $(W,W')$, to obtain 
\begin{equation}\label{eq:tube-faces-to-P-faces:application}
P(W',w) = P(W') + \delta' w = P(W) + \delta' w = P(W,w),
\end{equation}
where $\delta' = 1$ if  $\< u,w \> > 0$  for some $u \in W'$, and $\delta' = 0$ otherwise, as required.

Suppose now that $P(W',w) = P(W,w)$. If $w \in (W')^\perp \subset W^\perp$, then by Lemma~\ref{lem:tube-is-tube:app} we have
\begin{equation}\label{eq:tubes-equal:app}
\bigcup_{\lambda \in [0,1]} \big( P(W') + \lambda w \big) = P(W',w) = P(W,w) = \bigcup_{\lambda \in [0,1]} \big( P(W) + \lambda w \big).
\end{equation}
It follows that 
$$\max_{x \in P(W')} \< x,u \> = \max_{x \in P(W)} \< x,u \>$$
for every $u \in \SS^{d-1}$, and hence, since $P(W')$ and $P(W)$ are both compact, convex subsets of $\R^d$, they must be equal, as required. If $w \notin (W')^\perp$, on the other hand, then we wish to apply Lemma~\ref{lem:tube-faces-to-P-faces:app}, but to do so, we need to check that $P(W') \neq \emptyset$. To see that this holds, recall that $P(W) \neq \emptyset$, and observe that therefore $P(W',w) = P(W,w) \neq \emptyset$, by Lemma~\ref{lem:tube-faces-to-cliques:app}, and hence $P(W') \neq \emptyset$, again by Lemma~\ref{lem:tube-faces-to-cliques:app}. Thus, applying Lemma~\ref{lem:tube-faces-to-P-faces:app} twice, we obtain
\[
P(W') + \delta' w = P(W',w) = P(W,w) = P(W) + \delta' w,
\]
where $\delta'$ is as in~\eqref{eq:tube-faces-to-P-faces:application}. This completes the proof of the lemma.
\end{proof}

The first of the lemmas that we shall prove about the set of maximal cliques $\W$ is Lemma~\ref{lem:exists:WinW:trivial}, and it follows immediately from Lemma~\ref{lem:spheres-equal-iff-tubes-equal:app}.

\begin{lemma}\label{lem:exists:WinW:trivial:app}
Let $W \subset \QQ$ with $P(W) \ne \emptyset$. There exists $W' \in \W$ with 
$$W \subset W' \subset W \cup N_\QQ(W) \qquad \text{and} \qquad P(W,w) = P(W',w)$$
for every $w \in \cL_R$. 
\end{lemma}

\begin{proof}
Let $W' \supset W$ be maximal such that $P(W) = P(W')$. Note that $W'$ is a clique, by Lemma~\ref{lem:faces-to-cliques}, since $P(W) \ne \emptyset$, and therefore $W' \setminus W \subset N_\QQ(W)$. By Lemma~\ref{lem:spheres-equal-iff-tubes-equal:app}, it follows that $P(W,w) = P(W',w)$ for every $w \in \cL_R$, as required.
\end{proof}

Next we turn to Lemma~\ref{lem:clique-of-correct-dim}, which determines the dimension of a face. Recall that $\aff(X)$ denotes the affine span of a set $X \subset \R^d$.  

\begin{lemma}\label{lem:clique-of-correct-dim:app}
Let $W \in \W$ and $w \in \cL_R$. Then $P(W,w) \neq \emptyset$ and 
\[
\dim\big(\aff\big( P(W,w) \big) \big) = \dim(W^\perp).
\]
\end{lemma}

\begin{proof}
First let us note that $P(W) \ne \emptyset$ for every $W \in \W$. Indeed, this holds by~\eqref{eq:W-reminder}, because $P(\QQ) = \emptyset$ and $\QQ$ is not a clique. By Lemma~\ref{lem:tube-faces-to-cliques:app}, it follows that $P(W,w) \neq \emptyset$, and by Lemma~\ref{lem:spheres-equal-iff-tubes-equal:app}, that $P(W',w) \neq P(W,w)$ for every $W \subsetneq W' \subset \QQ$. 

Now, set $k := \dim\big(\aff\big( P(W,w) \big) \big)$, and recall from~\eqref{eq:P(W,w)} that $P(W,w)$ is contained in a translate of $W^\perp$, so $k \le \dim(W^\perp)$. Our task is therefore to show that $\dim(W^\perp) \le k$.

Recall from Lemma~\ref{lem:tube-nbrs:app} that 
\begin{align}
& P(W,w) = \bigcap_{u \in W} \big\{ x \in \R^d : \< x - \delta(u,w) w, \, u \> = 1 \big\} \nonumber\\ 
& \hspace{5cm} \cap \; \bigcap_{u \in N_\QQ(W)} \big\{ x \in \R^d : \< x - \delta(u,w) w, \, u \> \le 1 \big\},\label{eq:tube-nbrs:reminder:app}
\end{align}
set $U := \{ \pi(u,W^\perp) : u \in N_\QQ(W) \}$, and observe that $\< U \> = W^\perp$. Indeed, if this were not the case, then there would exist $v \in W^\perp$ such that $\< u,v \> = 0$ for all $u \in U$, and this would imply that $\<v\> \subset W^\perp \cap U^\perp$. But then $x + \<v\> \subset P(W,w)$ for all $x \in P(W,w)$, by~\eqref{eq:tube-nbrs:reminder:app}, and hence, since $P(W,w) \neq \emptyset$, it must be that $P(W,w)$ is unbounded. This is a contradiction, since~\eqref{eq:P(W,w)} implies that $P(W,w)$ is bounded, using the fact that $\QQ$ intersects every open hemisphere of $\SS^{d-1}$.

Next, set $F := \aff\big( P(W,w) \big) - a$ for some arbitrary $a \in \aff\big( P(W,w) \big)$, so $F$ is a subspace of $\R^d$ of dimension $k$. We claim that $U \subset F$, which (by the observations above) suffices to prove the lemma. To prove this, we claim first that, for each $u \in U$, the set
$$\aff\big( P(W,w) \big) \cap \big\{ x \in \R^d : \< x - \delta(u,w) w, \, u \> = 1 \big\}$$ 
is an affine subspace of $\R^d$ of dimension at most $k-1$. Indeed, if the dimension were $k$, then it would follow, using~\eqref{eq:P(W,w)}, that $P(W \cup \{u\},w) = P(W,w)$, and (by Lemma~\ref{lem:spheres-equal-iff-tubes-equal:app}) this would contradict our assumption that $W \in \W$. It follows that
\begin{equation}\label{eq:face-not-contained-in-sides}
P(W,w) \not\subset \bigcup_{u \in N_\QQ(W)} \big\{ x \in \R^d : \< x - \delta(u,w) w, \, u \> = 1 \big\},
\end{equation}
since a finite union of affine subspaces each of dimension at most $k-1$ cannot contain an affine subspace of dimension $k$, and since $P(W,w)$ is convex. 

Now, by~\eqref{eq:tube-nbrs:reminder:app} and~\eqref{eq:face-not-contained-in-sides} we may choose $x \in P(W,w)$ such that $\< x - \delta(u,w)w, u \> < 1$ for every $u \in N_\QQ(W)$. Observe that $x + \eps u' \in P(W,w)$ for every $u' \in U$ and sufficiently small $\eps > 0$. Indeed, we have $\< x + \eps u' - \delta(v,w)w, \, v \> \leq 1$ for all $v \in N_\QQ(W)$, and $\< x + \eps u' - \delta(u,w)w, v \> = 1$ for all $v \in W$, since $u' \in W^\perp$, so the assertion follows from~\eqref{eq:tube-nbrs:reminder:app}. It follows that $u' \in F$, and since $u'$ was arbitrary, we have $U \subset F$, as required.
\end{proof}

To finish this section, let us prove Lemma~\ref{lem:P(W,w,a,t,tau):def2}. 

\begin{lemma}\label{lem:P(W,w,a,t,tau):def2:app}
If $W \in \W$ and $P = P(W,w;a,t,\tau) \in \cP(W)$, then 
\begin{align*}
& P = \bigcap_{u \in W} \Big\{ x \in \R^d : \big\< x - a - \delta(u,w) \tau w, \, u \big\> = t \Big\} \nonumber \\ 
& \hspace{4cm} \cap \bigcap_{u \in N_\QQ(W)} \Big\{ x \in \R^d : \big\< x - a - \delta(u,w) \tau w, \, u \big\> \le t \Big\}. 
\end{align*}
\end{lemma}

\begin{proof}
First, note that $P(W,(\tau/t)w) \neq \emptyset$, by Lemma~\ref{lem:clique-of-correct-dim:app}. By Lemma~\ref{lem:tube-nbrs:app}, it follows that
\begin{align*}
& P(W, (\tau/t)w) = \bigcap_{u \in W} \Big\{ x \in \R^d : \big\< x - \delta(u,w) \tau w / t, \, u \big\> = 1 \Big\} \nonumber \\ 
& \hspace{5cm} \cap \bigcap_{u \in N_\QQ(W)} \Big\{ x \in \R^d : \big\< x - \delta(u,w) \tau w / t, \, u \big\> \le 1 \Big\}.
\end{align*}
Recalling from~\eqref{def:PWattau} that $P = a + t \cdot P\big( W, (\tau/t)w \big)$, the lemma follows.
\end{proof}

\section{The distance between faces of a polytope}\label{gamma:app}

In this appendix we prove several geometric properties of the canonical polytopes that were introduced in Section~\ref{sec:polytopes}. Each of these properties relates to the distance between the faces of these polytopes, and in particular to the constant $\gamma = \gamma(\QQ)$ defined in Definition~\ref{def:gamma}. We expect that these properties of polytopes are well-known, but we were unable to find references for them, and therefore provide the proofs for completeness. 

Two of the results proved in this section, Lemmas~\ref{lem:faces-far-away:app} and~\ref{lem:close-to-many-faces:app}, are used in Section~\ref{fundamental:sec}, while Lemmas~\ref{clm:faces-far-away} and~\ref{lem:faces-far-away:app2} are used in Appendix~\ref{cover:app}. The main step in the proof of Lemma~\ref{lem:faces-far-away:app} (which implies Lemma~\ref{lem:faces-far-away}) is the following variant of Lemma~\ref{lem:far:from:nonfaces}. 

\begin{lemma}\label{clm:faces-far-away}
Let $W \subset \QQ$, $w \in \cL_R$ and $v \in \QQ$. If $P(W \cup \{v\}) = \emptyset$, then
\begin{equation}\label{eq:faces-far-away:need2}
\big\< x - \delta(v,w) w, \, v \big\> \le 1 - \gamma
\end{equation}
for every $x \in P(W,w)$. 
\end{lemma}

\begin{proof}
We may assume that $P(W) \ne \emptyset$, otherwise there is nothing to prove. Moreover, if $w = \0$ then the claim follows from Lemma~\ref{lem:far:from:nonfaces}, so we may also assume that $w \ne \0$. We consider the cases $w \in W^\perp$ and $w \not\in W^\perp$ separately. 

\medskip
\noindent Case 1: $w \not\in W^\perp$.
\medskip

Let $u \in W$ be such that $\<u,w\> \ne 0$, and observe that
\begin{equation}\label{eq:faces-far-clm-x}
x - \delta(u,w) w \in P(W)
\end{equation}
by Lemma~\ref{lem:tube-faces-to-P-faces:app} (applied with $W = W'$). Since $P(W \cup \{v\}) = \emptyset$, it follows that 
\begin{equation}\label{eq:faces-far-away-1}
\big\< x - \delta(u,w) w, v \big\> \le 1 - \gamma,
\end{equation}
by Lemma~\ref{lem:far:from:nonfaces}. Now, observe that
$$\big( \delta(u,w) - \delta(v,w) \big) \< v,w \> \le 0,$$
since if $\< v,w \> \ge 0$ then $\delta(u,w) \le \delta(v,w)$, and if $\< v,w \> < 0$ then $\delta(u,w) \ge \delta(v,w)$. Combining this with~\eqref{eq:faces-far-away-1}, we obtain~\eqref{eq:faces-far-away:need2}.

\medskip
\noindent Case 2: $w \in W^\perp$.
\medskip

In this case the situation is slightly more complicated, because it may be that neither $x$ nor $x - w$ belongs to $P(W)$, a problem that we also encountered in the proof of Lemma~\ref{lem:tube-faces-to-cliques:app}. As in that proof, we deal with this difficulty by translating $x$ by some multiple of $w$. Thus, let us define
$$x' := x + \big( 2\delta(v,w) - 1 \big) \mu w,$$ 
where $\mu \ge 0$ is maximal such that $x' \in P(W,w)$; this is well-defined because $P(W,w)$ is compact. We claim that
\begin{equation}\label{eq:faces-far-away-x'}
x' - \delta(v,w) w \in P(W).
\end{equation}
To prove~\eqref{eq:faces-far-away-x'}, note that, since $w \ne \0$, we have
$$x' + \big( 2\delta(v,w) - 1 \big) \eps w \notin P(W,w)$$ 
for all $\eps > 0$, by our choice of $x'$. 
Since $x' \in P(W,w)$ and
$$P(W,w) = \bigcup_{\lambda \in [0,1]} \big( P(W) + \lambda w \big),$$
by Lemma~\ref{lem:tube-is-tube:app}, it follows that $x' \in P(W) + \delta(v,w) w$, as claimed.

Since $P(W \cup \{v\}) = \emptyset$, it follows, by Lemma~\ref{lem:far:from:nonfaces}, that  
\begin{equation}\label{eq:faces-far-away-2}
\big\< x' - \delta(v,w) w, \, v \big\> \le 1 - \gamma.
\end{equation}
To finish, we need to replace $x'$ by $x$ in~\eqref{eq:faces-far-away-2}. This is straightforward, because
$$\< x',v \> - \< x,v \> = \big( 2\delta(v,w) - 1 \big) \cdot \mu \cdot \< v,w \> = \mu \cdot | \< v,w \> | \ge 0,$$
since $\mu \ge 0$. 
This completes the proof of~\eqref{eq:faces-far-away:need2} when $w \in W^\perp$.
\end{proof}

We can now deduce the following lemma, which easily implies Lemma~\ref{lem:faces-far-away}, and which (as noted above) we shall use again in Appendix~\ref{cover:app}. 

\begin{lemma}\label{lem:faces-far-away:app2}
Let $W \subset W' \subset \QQ$, $w \in \cL_R$, $x \in P(W',w)$ and $y \in W^\perp$. If\/ $\| y \| \le \gamma$ and 
\begin{equation}\label{eq:faces-far-away:more:general:condition}
\< y,v \> \le 0 \quad \text{for every} \quad v \in \QQ \quad \text{such that} \quad P(W' \cup \{v\}) \ne \emptyset,
\end{equation}
then $x + y \in P(W,w)$. 
\end{lemma}

\begin{proof}
Observe first that, by Lemma~\ref{lem:tube-nbrs:app}, we have
$$\big\< x + y - \delta(v,w) w, \, v \big\> = 1$$ 
for all $v \in W$, since $x \in P(W',w) \subset P(W,w)$ and $y \in W^\perp$. To prove the lemma, it therefore suffices to show that 
\begin{equation}\label{eq:faces-far-away:need}
\big\< x + y - \delta(v,w) w, \, v \big\> \le 1
\end{equation}
for every $v \in N_\QQ(W)$, again by Lemma~\ref{lem:tube-nbrs:app}. 

If $P(W' \cup \{v\}) \ne \emptyset$, then $\< y,v \> \le 0$, by~\eqref{eq:faces-far-away:more:general:condition}, and so~\eqref{eq:faces-far-away:need} follows by Lemma~\ref{lem:tube-nbrs:app}, since $x \in P(W,w)$. On the other hand, if $P(W' \cup \{v\}) = \emptyset$ then, by Lemma~\ref{clm:faces-far-away}, 
\[
\big\< x - \delta(v,w) w, \, v \big\> \le 1 - \gamma
\]
for every $x \in P(W',w)$. Since $\|y\| \le \gamma$, it follows that~\eqref{eq:faces-far-away:need} holds, as required.
\end{proof}

As noted above, Lemma~\ref{lem:faces-far-away:app2} easily implies Lemma~\ref{lem:faces-far-away}, which is the case $W' = W \cup \{u\}$ of the following lemma. 

\begin{lemma}\label{lem:faces-far-away:app}
Let $W \subset W' \subset \QQ$, $w \in \cL_R$, $x \in P(W',w)$ and $y \in W^\perp$. If $\| y \| \le \gamma$ and 
\begin{equation}\label{eq:faces-far-away:condition:app}
\< y,v \> \le 0 \qquad \text{for every} \qquad v \in W' \cup N_\QQ(W').
\end{equation}
Then $x + y \in P(W,w)$. 
\end{lemma}

\begin{proof}
Observe that, by Lemma~\ref{lem:faces-to-cliques}, if $P(W' \cup \{v\}) \ne \emptyset$, then $v \in W' \cup N_\QQ(W')$, and therefore $\< y,v \> \le 0$. By Lemma~\ref{lem:faces-far-away:app2} it follows that $x + y \in P(W,w)$, as required.
\end{proof}

We next turn our attention to Lemma~\ref{lem:close-to-many-faces}. The proof of this lemma relies on the following standard fact about the distance between points and faces of polytopes. Let $\kappa = \kappa(\QQ) > 0$ be a sufficiently large constant depending on $\QQ$ and $\gamma$. 

\begin{lemma}\label{lem:innerprod-to-dist}
Let $W \subset \QQ$ and $u \in \QQ$ be such that $P(W') \neq \emptyset$, where $W' := W \cup \{u\}$. Then
\[
d\big( x, P(W') \big) \leq \kappa \cdot \big( 1 - \< x,u \> \big)
\]
for all\/ $x \in P(W)$. \qed
\end{lemma}

The proof of Lemma~\ref{lem:close-to-many-faces} will use Lemma~\ref{lem:innerprod-to-dist} directly, as well as the following consequence of that lemma.

\begin{lemma}\label{lem:dist-to-lower-dim-face}
Let $W, W' \subset \QQ$ and $u \in \QQ$ be such that $P(W \cup W' \cup \{u\}) \neq \emptyset$. Then 
\[
d\big( x, P(W \cup W' \cup \{u\}) \big) \leq 2\kappa \cdot \Big( d\big( x, P(W \cup W') \big) + d\big( x, P(W \cup \{u\}) \big) \Big)
\]
for all $x \in \R^d$.
\end{lemma}

\begin{proof}
Let $y \in P(W \cup W')$ be such that $d(x,y) = d\big(x,P(W \cup W')\big)$. Then
\begin{align}
d\big( x, P(W \cup W' \cup \{u\}) \big) &\leq d\big( x, P(W \cup W') \big) + d\big( y, P(W \cup W' \cup \{u\}) \big) \notag \\
&\leq d\big( x, P(W \cup W') \big) + \kappa \cdot d\big( y, H \big), \label{eq:dist-to-lower-1}
\end{align}
where $H := \big\{ z \in \R^d : \<z,u\> = 1 \big\}$, and the second inequality follows from Lemma~\ref{lem:innerprod-to-dist}. Continuing, we have
\[
d\big( y,H \big) \leq d\big( x,y \big) + d\big( x,H \big) \leq d\big( x,P(W \cup W') \big) + d\big( x,P(W \cup \{u\}) \big),
\]
where this time the second inequality follows because $d(x,y) = d\big(x,P(W \cup W')\big)$ and because $P(W \cup \{u\}) \subset H$, by~\eqref{eq:P(W)}. Combined with~\eqref{eq:dist-to-lower-1}, this completes the proof.
\end{proof}

We can now deduce Lemma~\ref{lem:close-to-many-faces}, which we restate  here for the reader's convenience. Recall from~\eqref{def:constants} that $\delta = \delta(\QQ) > 0$ is a sufficiently small constant; in particular, we will choose $\delta$ depending on $\gamma(\QQ)$ and $\kappa(\QQ)$. 


\begin{lemma}\label{lem:close-to-many-faces:app}
Let $W \in \W$ and $T \subset N_\QQ(W)$ be such that $P(W\cup\{u\}) \neq \emptyset$ for all $u \in T$. If there exists $x \in P(W)$ such that 
$$\< x,u \> \ge 1 - 2\delta$$
for every $u \in T$, then $W \cup T$ is a clique and $P(W \cup T) \neq \emptyset$.
\end{lemma}

\begin{proof} 
Let $x \in P(W)$ be such that $\< x,u \> \ge 1 - 2\delta$ for every $u \in T = \{u_1,\dots,u_\ell\}$. We shall prove by induction that $P(W \cup T_i) \neq \emptyset$ for each $i \in \{0,1,\dots,\ell\}$, where $T_i := \{u_1,\dots,u_i\}$; the result will then follow from the case $i = \ell$, together with Lemma~\ref{lem:faces-to-cliques}. 

In order to assist with the proof, we shall include in the induction hypothesis the additional assertion that
\begin{equation}\label{eq:close-to-many-hyp}
d\big( x, P(W \cup T_i) \big) \leq (2\kappa)^{i+1} \cdot i \cdot \delta.
\end{equation}
The base case $i=0$ is automatic for both parts, since $x \in P(W)$.

Let us suppose then that $P(W \cup T_i) \neq \emptyset$ and that~\eqref{eq:close-to-many-hyp} holds, for some $i < \ell$. Now, if $P(W \cup T_{i+1}) = \emptyset$, then
\begin{equation}\label{eq:close-to-many-assume}
D(W \cup T_i, u_{i+1}) \geq \gamma,
\end{equation}
by Definition~\ref{def:gamma}. On the other hand, we have
\begin{equation}\label{eq:close-to-many-dist}
d\big( x, P(W \cup \{u_{i+1}\}) \big) \le \kappa \cdot \big( 1 - \< x,u_{i+1} \> \big) \le 2\kappa\delta,
\end{equation}
by Lemma~\ref{lem:innerprod-to-dist}, since $P(W\cup\{u_{i+1}\}) \neq \emptyset$ and $x \in P(W)$, by assumption. 
Together with~\eqref{eq:close-to-many-hyp}, this implies that\footnote{For the first step we use Definition~\ref{def:gamma} and the fact that $P(W \cup \{u_{i+1}\}) \subset \big\{ y : \<y,u_{i+1}\> = 1 \big\}$.}
\[
D(W \cup T_i, u_{i+1}) \leq d\big( P(W \cup T_i), P(W \cup \{u_{i+1}\}) \big) \leq ((2\kappa)^{i+1} \cdot i + 2\kappa) \cdot \delta < \gamma,
\]
since $\delta$ is sufficiently small (in terms of $\kappa$ and $\gamma$). This contradicts~\eqref{eq:close-to-many-assume}, and hence proves that $P(W \cup T_{i+1}) \neq \emptyset$.

To complete the proof, we must prove the induction step for~\eqref{eq:close-to-many-hyp}. Applying Lemma~\ref{lem:dist-to-lower-dim-face} (with $W' = T_i$ and $u = u_{i+1}$, and using the fact that that $P(W \cup T_{i+1}) \neq \emptyset$), we have
\[
d\big( x, P(W \cup T_{i+1}) \big) \le 2\kappa \cdot \Big( d\big( x, P(W \cup T_i) \big) + d\big( x, P(W \cup \{ u_{i+1} \} ) \big) \Big).
\]
Combining this with~\eqref{eq:close-to-many-hyp} and~\eqref{eq:close-to-many-dist}, we obtain
\[
d\big( x, P(W \cup T_{i+1}) \big) \leq 2\kappa \cdot \big( (2\kappa)^{i+1} \cdot i \cdot \delta + 2\kappa\delta \big) \leq (2\kappa)^{i+2} \cdot (i+1) \cdot \delta,
\]
completing the induction. As noted above, the case $i = \ell$ of the induction hypothesis implies that $P(W \cup T) \neq \emptyset$, and hence $W \cup T$ is a clique, as required.
\end{proof}

\section{Interiors and extensions}\label{buffers:app}

This appendix contains proofs of the lemmas stated in Section~\ref{buffers:sec}. 

\subsection{The interior and the closed interior of a polytope}

In this subsection we prove Lemmas~\ref{lem:interiors-are-canonical},~\ref{lem:clint:Delta} and~\ref{lem:x-notin-Delta:Wv}. Before proving the first of these, Lemma~\ref{lem:interiors-are-canonical}, we need to establish two simple facts about interiors and closed interiors.

Recall from~\eqref{def:intP} and Definition~\ref{def:closed:interior} the definitions of the interior $\interior(P)$ and the closed interior $\clint(P)$ of a polytope $P \in \cP(W)$. 
The first of our two facts says that $\clint(P)$ and $P$ are contained in the same translate of $W^\perp$. This observation will also be useful in Section~\ref{extension:app}. 

\begin{lemma}\label{lem:Ds:same:shift:app}
Let $W \in \W$ and $P = P(W,w;a,t,\tau) \in \cP(W)$, let $y \in P(W)$ and $\eps \in \R$, and suppose that $Q = P(W,w;a+\eps y,t-\eps,\tau) \in \cP(W)$. Then
\begin{equation}\label{eq:Ds:same:shift:app}
P - Q \subset W^\perp.
\end{equation}
\end{lemma}

\begin{proof}
It follows from Lemma~\ref{lem:P(W,w,a,t,tau):def2:app} that $P$ and $Q$ are contained in translates of $W^\perp$, so to prove~\eqref{eq:Ds:same:shift:app} we need to show that $\< x,u \> = \<x',u\>$ for every $x \in P$, $x' \in Q$, and $u \in W$. To see this, observe that, again by Lemma~\ref{lem:P(W,w,a,t,tau):def2:app},  
\[
\big\< x' - (a + \eps y) - \delta(u,w) \tau w, \, u \big\> = t - \eps = \big\< x - a - \delta(u,w) \tau w, \, u \big\> - \eps.
\]
Since $y \in P(W)$ implies that $\<y,u\> = 1$, it follows that $\< x,u \> = \<x',u\>$, as claimed.
\end{proof}

\begin{lemma}\label{lem:interior:def2}
If\/ $W \in \W$ and $P = P(W,w;a,t,\tau) \in \cP(W)$, then 
\begin{align}
& \interior(P) = \bigcap_{u \in W} \Big\{ x \in \R^d : \big\< x - a - \delta(u,w) \tau w, \, u \big\> = t \Big\} \nonumber\\ 
& \hspace{4cm} \cap \; \bigcap_{u \in N_\QQ(W)} \Big\{ x \in \R^d : \big\< x - a -  \delta(u,w) \tau w, \, u \big\> < t \Big\}.\label{eq:tube-nbrs:interior:app}
\end{align}
Moreover, 
\begin{align}
& \interior\big( P(W,w) \big) = \bigcap_{u \in W} \Big\{ x \in \R^d : \big\< x - \delta(u,w) w, \, u \big\> = 1 \Big\} \nonumber\\ 
& \hspace{5cm} \cap \; \bigcap_{u \in N_\QQ(W)} \Big\{ x \in \R^d : \big\< x -  \delta(u,w) w, \, u \big\> < 1 \Big\}.\label{eq:canonical:interior:app}
\end{align}

\end{lemma}

\begin{proof}
Suppose first that $x \in P$, and that $x$ is not in the right-hand side of~\eqref{eq:tube-nbrs:interior:app}. It follows, by Lemma~\ref{lem:P(W,w,a,t,tau):def2:app}, that $\big\< x - a -  \delta(u,w) \tau w, \, u \big\> = t$ for some $u \in N_\QQ(W)$. Set $W' = W \cup \{u\}$, and note that $N_\QQ(W') \subset N_\QQ(W)$. By Lemma~\ref{lem:P(W,w,a,t,tau):def2:app}, and recalling that $x \in P$, it follows that $x \in \Delta(P,W')$, and hence $x \not\in \interior(P)$, by~\eqref{def:intP}. 

Now suppose that $x$ is in the right-hand side of~\eqref{eq:tube-nbrs:interior:app}, so $x \in P$, by Lemma~\ref{lem:P(W,w,a,t,tau):def2:app}. Moreover, $x \not\in \Delta(P,W \cup \{u\})$ for each $u \in N_\QQ(W)$, since $\big\< x - a -  \delta(u,w) \tau w, \, u \big\> < t$, and thus $x \in \interior(P)$, as required. To deduce~\eqref{eq:canonical:interior:app}, recall that 
$$\interior\big( P(W,w) \big) = t^{-1} \big( \interior\big( t \cdot P(W,w) \big) \big)$$ 
for any $t > C$. 
\end{proof}

We can now deduce the key properties of $\clint(P)$ in Lemma~\ref{lem:interiors-are-canonical}. 

\begin{lemma}\label{lem:interiors-are-canonical:app}
Let $W \in \W$ and $P \in \cP(W)$. Then\/ $\clint(P) \in \cP(W)$, 
\begin{equation}\label{eq:interior:droplets:equal:app}
\clint(P) \subset \interior(P) \qquad \text{and} \qquad \clint(P) \cap \Z^d = \interior(P) \cap \Z^d.
\end{equation}
\end{lemma}

\begin{proof} 
Let $P = P(W,w;a,t,\tau)$ and $\clint(P) = P(W,w;a+\eps y,t-\eps,\tau)$, where $y \in \interior\big( P(W) \big)$ and $\eps > 0$. Observe first that $\clint(P) \in \cP(W)$ if $\eps$ is sufficiently small, by Definition~\ref{def:cPW}, since $t > C$. Note also that, by Lemma~\ref{lem:Ds:same:shift:app}, we have
\begin{equation}\label{eq:Ds:same:shift:app2}
P - \clint(P) \subset W^\perp.
\end{equation}

Now, let $x \in \clint(P)$ and observe that
$$\big\< x - (a + \eps y) - \delta(u,w) \tau w, \, u \big\> \leq t - \eps$$
for every $u \in N_\QQ(W)$, by Lemma~\ref{lem:P(W,w,a,t,tau):def2:app}. Since $y \in \interior\big(P(W)\big)$, we have $\< y,u \> < 1$, by Lemma~\ref{lem:interior:def2}, so 
\begin{equation}\label{eq:need:for:DinD:1:app}
\big\< x - a - \delta(u,w) \tau w, \, u \big\> < t.
\end{equation}
for every $u \in N_\QQ(W)$. Since $P$ and $\clint(P)$ are contained in the same translate of $W^\perp$, by~\eqref{eq:Ds:same:shift:app2}, it follows from Lemma~\ref{lem:interior:def2} that $x \in \interior(P)$, as required.

It remains to show that $D = D'$, where $D := \interior(P) \cap \Z^d$ and $D' := \clint(P) \cap \Z^d$. One inclusion is clear, since we have $\clint(P) \subset \interior(P)$, and therefore $D' \subset D$ by taking the intersection with $\Z^d$. To show that $D \subset D'$, note first that~\eqref{eq:need:for:DinD:1:app} holds for every $x \in D$ and $u \in N_\QQ(W)$, by Lemma~\ref{lem:interior:def2} and since $D \subset \interior(P)$. Since $D$ is finite, 
it follows that 
\[
t' := \max_{x \in D} \max_{u \in N_\QQ(W)} \big\< x - a - \delta(u,w) \tau w, \, u \big\> < t.
\]
Moreover, $y \in \interior\big( P(W) \big)$ implies that $\< y,u \> < 1$ for every $u \in N_\QQ(W)$, by Lemma~\ref{lem:interior:def2}, and in particular,
$$c := \max_{u \in N_\QQ(W)} \big( 1 - \<y,u\> \big) > 0.$$
(One might have expected `$\min$' instead of `$\max$' here, but the maximum is intended.)
Hence, if $0 < \eps < (t - t')/c$, then for every $x \in D$ and $u \in N_\QQ(W)$, 
\[
\big\< x - (a + \eps y) - \delta(u,w) \tau w, \, u \big\> \leq t' - \eps \< y,u \> \le t - \eps,
\]
since $\eps\big(1 - \<y,u\>\big) \leq t - t'$. Recalling that $P$ and $\clint(P)$ are contained in the same translate of $W^\perp$, by~\eqref{eq:Ds:same:shift:app2}, it follows by Lemma~\ref{lem:P(W,w,a,t,tau):def2:app} that $x \in \clint(P)$. Finally, since $x \in D \subset \Z^d$, we obtain $x \in D'$, as required.
\end{proof}

\begin{remark}\label{rem:clint-generality}
All that is used about $\clint(P)$ in the proof of Lemma~\ref{lem:interiors-are-canonical:app} is that it is equal to $P(W,w;a+\eps y,t-\eps,\tau)$ for some $y \in \interior(P)$ and some $\eps > 0$ that is sufficiently small (uniformly in $y$). Thus, we may replace $\clint(P)$ by $P(W,w;a+\eps y,t-\eps,\tau)$ (for any $y \in \interior(P)$ and sufficiently small $\eps > 0$) throughout the statement of the lemma.
\end{remark}

We move on to the proofs of Lemmas~\ref{lem:clint:Delta} and~\ref{lem:x-notin-Delta:Wv}. Both lemmas relate to properties of the $W'$-shifted closed interior of a polytope $P \in \cP(W)$, which was defined in Definition~\ref{def:shifted:interior}. We divide the proof of Lemma~\ref{lem:clint:Delta} into two; the first part is a straightforward consequence of Lemma~\ref{lem:interiors-are-canonical:app} and Remark~\ref{rem:clint-generality}.

\begin{lemma}\label{lem:clint:Delta:app2}
Let $W \in \W$ and $P \in \cP(W)$. Then $\clint(P \to W') \in \cP(W)$, 
$$\Delta\big( \clint(P \to W'),W' \big) \subset \interior\big( \Delta(P,W') \big),$$
and
$$\Delta\big( \clint(P \to W'),W' \big) \cap \Z^d = \clint \big( \Delta(P,W') \big) \cap \Z^d$$
for every $W \subset W' \in \W$.
\end{lemma}

\begin{proof}
Let $P = P(W,w;a,t,\tau)$, so that $\Delta(P,W') = P(W',w;a,t,\tau)$, and
$$\Delta\big( \clint(P \to W'),W' \big) = P\big( W',w; a+\eps y, t-\eps,\tau \big)$$
for some $y \in \interior\big( P(W') \big)$ and sufficiently small $\eps > 0$. Note first that $\clint(P \to W') \in \cP(W)$, by Definitions~\ref{def:cPW} and~\ref{def:shifted:interior}, and since $t > C$ and $\eps$ is sufficiently small. 

Now, by Lemma~\ref{lem:interiors-are-canonical:app} and Remark~\ref{rem:clint-generality}, we have 
$$P\big( W',w; a+\eps y, t-\eps,\tau \big) \subset \interior\big( P( W',w; a,t,\tau ) \big)$$
and 
$$P\big( W',w; a+\eps y, t-\eps,\tau \big) \cap \Z^d = \clint\big( P( W',w; a,t,\tau ) \big) \cap \Z^d$$
for any $y \in \interior\big( P(W') \big)$ and sufficiently small $\eps > 0$, as required.
\end{proof}

We next show that $\clint(P \to W')$ is contained in $P$; the proof is similar to that of Lemma~\ref{lem:interiors-are-canonical:app}.

\begin{lemma}\label{lem:clint:Delta:app1}
Let $W \in \W$ and $P \in \cP(W)$. Then 
$$\clint(P \to W') \subset P$$
for every $W \subset W' \in \W$.
\end{lemma}

\begin{proof}
Let $P = P(W,w;a,t,\tau)$ and $\clint(P \to W') = P(W,w;a+\eps y,t-\eps,\tau) \in \cP(W)$ (by Lemma~\ref{lem:clint:Delta:app2}), where $y \in \interior\big( P(W') \big) \subset P(W)$ and $\eps > 0$. By Lemma~\ref{lem:Ds:same:shift:app}, we have
\begin{equation}\label{eq:same:shift:app3}
P - \clint(P \to W') \subset W^\perp.
\end{equation}
Now, let $x \in \clint(P \to W')$ and observe that
$$\big\< x - (a + \eps y) - \delta(u,w) \tau w, \, u \big\> \leq t - \eps$$
for every $u \in N_\QQ(W)$, by Lemma~\ref{lem:P(W,w,a,t,tau):def2:app}. Since $y \in P(\emptyset)$, we have $\< y,u \> \le 1$, by~\eqref{def:P}, so 
$$\big\< x - a - \delta(u,w) \tau w, \, u \big\> \le t.$$
Since $P$ and $\clint(P \to W')$ are contained in the same translate of $W^\perp$, by~\eqref{eq:same:shift:app3}, it follows from Lemma~\ref{lem:P(W,w,a,t,tau):def2:app} that $x \in P$, which completes the proof.
\end{proof}

Next we prove Lemma~\ref{lem:x-notin-Delta:Wv}, which allows us to control which of the faces of $P$ have non-empty intersection with $\clint(P \to W')$.

\begin{lemma}\label{lem:x-notin-Delta:Wv:app}
Let $W \in \W$ and $P \in \cP(W)$, and let $W \subset W' \in \W$. If $x \in \clint(P \to W')$, then
$$x \notin \Delta\big( P, W \cup \{v\} \big)$$
for every $v \in \QQ \setminus W'$.
\end{lemma}

\begin{proof}
Let $y \in \interior\big( P(W') \big)$ and $\eps > 0$ be such that
\[
P = P(W,w;a,t,\tau) \qquad \text{and} \qquad \clint(P \to W') = P(W,w; a+\eps y,t-\eps,\tau).
\]
We claim that $\<y,v\> < 1$ for every $v \in \QQ \setminus W'$. To see this, note first that $y \in P(W') \subset P(\emptyset)$, so by~\eqref{def:P} and~\eqref{eq:P(W)} we have $\<y,v\> \le 1$. If $\<y,v\> = 1$, then $y \in P(W' \cup \{v\})$, so we must have $v \in N_\QQ(W')$ by Lemma~\ref{lem:faces-to-cliques}. But $y \in \interior\big( P(W') \big)$, so by Lemma~\ref{lem:interior:def2} we have $\<y,v\> < 1$ for each $v \in N_\QQ(W')$. Hence $\<y,v\> < 1$ for every $v \in \QQ \setminus W'$.

Now let $v \in \QQ \setminus W'$, and note that, since $x \in \clint(P \to W')$, we have
$$\big\< x - (a + \eps y) - \delta(u,w) \tau w, \, v \big\> \leq t - \eps,$$
by~Lemma~\ref{lem:P(W,w,a,t,tau):def2:app}. Thus, since $\<y,v\> < 1$, it follows that
$$\big\< x - a - \delta(u,w) \tau w, \, v \big\> < t,$$
and hence, by Lemma~\ref{lem:P(W,w,a,t,tau):def2:app}, we have $x \not\in \Delta(P,W \cup \{v\})$. 
\end{proof}

\subsection{Forwards and sideways faces}

Our next task is to prove Lemmas~\ref{lem:Delta:forwards} and~\ref{lem:max:face:interior}. We first prove the latter lemma, and then use it to deduce the former. 

\begin{lemma}\label{lem:max:face:interior:app} 
Let $W \subset \QQ$ with $P(W) \ne \emptyset$, let $P \in \cP(W)$, and let $x \in P$. If $W' \subset \QQ$ is maximal such that 
$$W \subset W' \in \W \qquad \text{and} \qquad x \in \Delta(P,W'),$$ 
then $W' \setminus W \subset N_\QQ(W)$ and $x \in \interior\big( \Delta(P,W') \big)$.
\end{lemma}

\begin{proof}
Note that $W'$ is a clique (since $W' \in \W$), and therefore $W' \setminus W \subset N_\QQ(W)$. (Alternatively, this follows by Lemma~\ref{lem:faces-to-cliques}, since $x \in \Delta(P,W')$.) 

Now, if $x \in \Delta(P,W') \setminus \interior\big( \Delta(P,W') \big)$, then $x \in \Delta(P,W' \cup \{u\})$ for some $u \in N_\QQ(W')$, by~\eqref{def:Delta:face} and~\eqref{def:intP}. By Lemma~\ref{lem:exists:WinW:trivial:app}, it follows that  
$$x \in \Delta(P,W' \cup \{u\}) = \Delta(P,W'') \qquad \text{for some} \qquad W' \cup \{u\} \subset W'' \in \W.$$
Since $W'' \in \W$ and $W' \subsetneq W''$, this contradicts the maximality of $W'$.
\end{proof}

We can now deduce Lemma~\ref{lem:Delta:forwards}. Recall the definition of the forwards boundary $\Delta(P)$ of a polytope $P \in \cP(W)$ (see Definition~\ref{def:forwards:boundary}). 

\begin{lemma}\label{lem:Delta:forwards:app}
Let $W \in \W$ and $P \in \cP(W)$. Then
$$\Delta(P) = \bigcup_{u \in N_\QQ(W,w)} \Delta( P, W \cup \{u\}).$$
\end{lemma}

\begin{proof}
If $x \in \Delta(P)$ then $x \in \interior\big( \Delta(P,W') \big)$ for some $W' \in \W^\to(P)$. Recalling~\eqref{def:WtoP}, let $u \in (W' \setminus W) \cap N_\QQ(W,w)$, and observe that, by~\eqref{def:Delta:face}, we have\footnote{Observe that if $W \subset W'$, then $\Delta(P,W') \subset \Delta(P,W)$.} 
$$x \in \Delta(P,W') \subset \Delta(P,W \cup \{u\}).$$ 
On the other hand, if $x \in \Delta(P,W \cup \{u\})$ for some $u \in N_\QQ(W,w)$, then let 
$$W \cup \{u\} \subset W' \in \W$$ 
be maximal such that $x \in \Delta(P,W')$. By Lemma~\ref{lem:max:face:interior}, we have $x \in \interior\big( \Delta(P,W') \big)$. Since $u \in W' \cap N_\QQ(W,w)$ and $W' \in \W$, it follows that $W' \in \W^\to(P)$, as required. 
\end{proof}

To finish this subsection, let us prove Lemma~\ref{lem:forwardsfaces:tau}.

\begin{lemma}\label{lem:forwardsfaces:tau:app}
Let $W \in \W$ and $P \in \cP(W)$. If\/ $W' \in \W^\to(P)$, then there exists $Q \in \cP(W')$ with $\tau(Q) = 0$ such that $Q = \Delta(P,W')$ (as subsets of $\R^d$).
\end{lemma}

\begin{proof}
Let $P = P(W,w;a,t,\tau)$, and define $P_0 := P(W,w;a,t,0)$. Recall that $W'$ intersects $N_\QQ(W,w)$, by the definition of $\W^\to(P)$ in~\eqref{def:WtoP}, and therefore $\<u,w\> > 0$ for some $u \in W'$. Hence, by Lemma~\ref{lem:tube-faces-to-P-faces},~\eqref{def:PWattau} and~\eqref{def:Delta:face}, 
it follows that $\Delta(P,W') = Q$, where
\[
Q = \Delta(P_0,W') + \tau w = P(W',w;a+\tau w,t,0),
\]
so $\tau(Q) = 0$, as claimed.
\end{proof}

\subsection{The extension of a polytope}\label{extension:app}

Our next task is to prove Lemma~\ref{lem:exterior}, our two key properties of the extension $\ext(P)$ of a polytope $P \in \cP(W)$ (see Definition~\ref{def:ext:P}). First, however, we will show that $\ext(P)$ is well-defined. 

\begin{lemma}\label{lem:ext:exists}
Let $W \in \W$ with $\dim(W^\perp) \ne 0$, let $P = P(W,w;a,t,\tau) \in \cP(W)$  with $P \cap \Z^d \ne \emptyset$, and let $y \in \interior\big( P(W) \big)$. There exists $\eps > 0$ such that 
$$P\big( W,w; a - \eps y, t + \eps, \tau \big) \cap \Z^d \ne P \cap \Z^d.$$
\end{lemma}

\begin{proof}
Set $P(\eps) := P\big( W,w; a - \eps y, t + \eps, \tau \big)$ for each $\eps > 0$, and observe that, by Lemma~\ref{lem:Ds:same:shift:app} (applied to $-\eps$), we have
$$P - P(\eps) \subset W^\perp.$$
Now, since $P$ is bounded and $\L(W) = W^\perp \cap \Z^d$ is infinite (since $\dim(W^\perp) \ne 0$), it follows that there exists $x \in \Z^d \setminus P$ with $x - P \subset W^\perp$. We claim that if $\eps$ is sufficiently large, then $x \in P(\eps)$. By Lemma~\ref{lem:P(W,w,a,t,tau):def2:app}, and since $x - P(\eps) \subset W^\perp$, it will suffice to show that if $u \in N_\QQ(W)$ then
$$\big\< x - a + \eps y - \delta(u,w) \tau w, \, u \big\> \le t + \eps,$$
which follows (for $\eps$ large) because $\<y,u\> < 1$, by Lemma~\ref{lem:interior:def2}, since $y \in \interior\big( P(W) \big)$. 
\end{proof}
 
We will next prove the two properties claimed in Lemma~\ref{lem:exterior}. We shall do so in two separate lemmas; in the first of these, we show that $P \subset \ext(P)$. The proof is similar to those of Lemmas~\ref{lem:interiors-are-canonical:app} and~\ref{lem:clint:Delta:app1}.

\begin{lemma}\label{lem:exterior:app}
Let $W \in \W$ and $P \in \cP(W)$ with $P \cap \Z^d \ne \emptyset$. Then $P \subset \ext(P)$.
\end{lemma}

\begin{proof}
Let 
$$P = P(W,w;a,t,\tau) \qquad \text{and} \qquad \ext(P) = P(W,w;a-\eps y,t+\eps,\tau),$$ 
where $y \in P(W)$ and $\eps > 0$. By Lemma~\ref{lem:Ds:same:shift:app} (applied to $-\eps$), we have
$$P - \ext(P) \subset W^\perp.$$
Let $x \in P$, and observe that
$$\big\< x - a - \delta(u,w) \tau w, \, u \big\> \le t$$
for every $u \in N_\QQ(W)$, by Lemma~\ref{lem:P(W,w,a,t,tau):def2:app}. Since $y \in P(W)$, we have $\< y,u \> \le 1$, by~\eqref{eq:P(W)}, so
$$\big\< x - (a - \eps y) - \delta(u,w) \tau w, \, u \big\> \le t + \eps$$
for every $u \in N_\QQ(W)$. Since $P$ and $\ext(P)$ are contained in the same translate of $W^\perp$, it follows by Lemma~\ref{lem:P(W,w,a,t,tau):def2:app} that $x \in \ext(P)$, as required.
\end{proof}

Next, we prove the second property in Lemma~\ref{lem:exterior}, which states that every lattice point of $\ext(P) \setminus P$ is contained in one of the faces of $P$.  

\begin{lemma}\label{lem:exterior:faces:app}
Let $W \in \W$ and $P \in \cP(W)$, and set $P' := \ext(P)$. Then
$$\interior(P') \cap \Z^d \subset P.$$
\end{lemma}

\begin{proof}
We claim that if $x \in (P' \setminus P) \cap \Z^d$, then 
$$x \in \Delta(P',W \cup \{u\})$$ 
for some $u \in N_\QQ(W)$. By the definition~\eqref{def:intP} of the interior, this will suffice. 

Let $P = P(W,w;a,t,\tau)$ and $P' = P(W,w;a-\eps' y,t+\eps',\tau)$, where $y \in \interior\big( P(W) \big)$ and $\eps' > 0$, and for each $\eps \geq 0$, set
$$P(\eps) := P(W,w;a-\eps y,t+\eps,\tau).$$
Note that, by Lemma~\ref{lem:Ds:same:shift:app}, 
$$P - P(\eps) \subset W^\perp$$
for all $\eps \geq 0$. By the minimality of $\eps'$, we have $P \cap \Z^d = P(\eps) \cap \Z^d$ for every $0 \leq \eps < \eps'$. Moreover, $\<y,u\> < 1$ for every $u \in N_\QQ(W)$, by Lemma~\ref{lem:interior:def2}, since $y \in \interior\big( P(W) \big)$. Thus, if we had
\[
\big\< x - (a - \eps' y) - \delta(u,w) \tau w, \, u \big\> < t + \eps'
\]
for all $u \in N_\QQ(W)$, then we would have $x \in P(\eps) \cap \Z^d$ for some $0 \leq \eps < \eps'$, by Lemma~\ref{lem:P(W,w,a,t,tau):def2:app}, which would be a contradiction because $x \notin P \cap \Z^d$. Therefore,
$$\big\< x - (a - \eps' y) - \delta(u,w) \tau w, \, u \big\> = t + \eps'$$
for some $u \in N_\QQ(W)$, and hence $x \in \Delta(P',W \cup \{u\})$, as claimed.
\end{proof}

\subsection{The forwards extension and retraction of a polytope}\label{retraction:app}

Next we present proofs of Lemmas~\ref{lem:fext} and~\ref{lem:fret}. Before doing so, let us prove a simple observation about the elements of $\Delta(P)$. Recall from Definition~\ref{def:buffer} that if $x \in P = P(W,w;a,t,\tau)$, then 
$$\tau_P(x) := \inf\big\{ \tau^* \ge 0 : x \in P\big( W,w;a,t,\tau^* \big) \big\}.$$

\begin{lemma}\label{lem:maxtau:inDelta}
Let $W \in \W$ and $P = P(W,w;a,t,\tau) \in \cP(W)$, with $w \in W^\perp$ and $\tau > 0$. Then 
$$\Delta(P) = \big\{ x \in P : \tau_P(x) = \tau \big\}.$$ 
\end{lemma}

\begin{proof}
Let $x \in P$. By Lemma~\ref{lem:P(W,w,a,t,tau):def2:app}, and since $w \in W^\perp$ and $\tau > 0$, we have $\tau_P(x) = \tau$ if and only if
$$\big\< x - a - \delta(u,w) \tau w, \, u \big\> = t$$
for some $u \in N_\QQ(W)$ with $\delta(u,w) = 1$ and $\<w,u\> \ne 0$. Since $w \in W^\perp$, and again using Lemma~\ref{lem:P(W,w,a,t,tau):def2:app}, this is equivalent to 
$$x \in \Delta(P,W \cup \{u\})$$
for some $u \in N_\QQ(W,w)$. By Lemma~\ref{lem:Delta:forwards:app}, this is true if and only if $x \in \Delta(P)$, as required.
\end{proof}

Turning to Lemma~\ref{lem:fext}, recall from Definition~\ref{def:fext:P} the forwards extension $\fext(P)$ and forwards retraction $\fret(P)$ of a polytope $P \in \cP(W)$.

\begin{lemma}\label{lem:fext:app}
Let $W \in \W$ and $P \in \cP(W)$, with $w(P) \in W^\perp$. Then
\begin{equation}\label{eq:fext-app}
\fext(P) \cap \Z^d \subset P \cup \Delta\big( \fext(P) \big)
\end{equation}
and
\begin{equation}\label{eq:fret-app}
P \cap \Z^d \subset \fret(P) \cup \Delta(P).
\end{equation}
\end{lemma}

\begin{proof}
We first claim that $\fext(P)$ is well-defined, i.e., that 
$$P\big( W,w;a,t,\tau' \big) \cap \Z^d \ne P\big( W,w;a,t,\tau \big) \cap \Z^d$$
for some $\tau' > \tau$. This follows because $w \in W^\perp$ and $t > C$, using Lemma~\ref{lem:tube-is-tube}. Indeed, there are infinitely many points $x \in \L(W) = W^\perp \cap \Z^d$ such that $x \in P(W,w;a,t,\tau) + \lambda w$ for some $\lambda > 0$, since $C = C(\QQ)$ is sufficiently large. Since $P$ is bounded, the claim follows.

To prove~\eqref{eq:fext-app}, set $P' := \fext(P)$ and let $x \in (P' \setminus P) \cap \Z^d$. By Definition~\ref{def:fext:P}, $\tau(P') > \tau(P)$ is minimal such that $P' \cap \Z^d \ne P \cap \Z^d$. It follows that $\tau_{P'}(x) = \tau(P')$, and therefore $x \in \Delta(P')$, by Lemma~\ref{lem:maxtau:inDelta}, proving~\eqref{eq:fext-app}.

Similarly, to prove~\eqref{eq:fret-app}, observe that if $x \in P \cap \Z^d$ and $x \not\in \fret(P)$, then by Definition~\ref{def:fext:P} we must have $\tau_P(x) = \tau(P) > 0$, since $\tau(P) - \tau(P')$ is sufficiently small. By Lemma~\ref{lem:maxtau:inDelta}, it follows that $x \in \Delta(P)$, as required.
\end{proof}

The proof of Lemma~\ref{lem:fret} is also straightforward.

\begin{lemma}\label{lem:fret:app}
Let $W \in \W$ and $P \in \cP(W)$, with $w(P) \in W^\perp$. If $\tau(P) > 0$ and $\Delta(P) \cap \Z^d \ne \emptyset$, then  
$$\fext\big( \fret(P) \big) = P.$$
\end{lemma}

\begin{proof} 
If $\Delta(P) \cap \Z^d \ne \emptyset$ and $\tau(P) > 0$, then by Lemma~\ref{lem:maxtau:inDelta} there exists $x \in P \cap \Z^d$ with $\tau_P(x) = \tau(P)$. Note also that $\fret(P) \subset P$, by Lemma~\ref{lem:tube-is-tube:app}, since $w \in W^\perp$ and $\tau' < \tau$. Since $\tau(P) - \tau\big( \fret(P) \big)$ is sufficiently small, it follows that $\tau(P)$ is minimal such that $\tau(P) > \tau\big( \fret(P) \big)$ and $P \cap \Z^d \ne \fret(P) \cap \Z^d$, so $\fext\big( \fret(P) \big) = P$, as required.
\end{proof}

\subsection{Growth sequences}

Next, we prove the lemmas stated in Section~\ref{sec:growth:seq}. Recall from Definition~\ref{def:forwards:growth:sequence} that the forwards growth sequence $\G$ with seed $Q$, where $Q \in \cP(W,w)$ for some $W \in \W$ and $w \in \cL_R \cap \SS(W)$, is defined by setting $Q_0 := Q$ and $Q_j := \fext(Q_{j-1})$ for each $j \ge 1$. We will first prove Lemma~\ref{lem:growing:linearly:forwards}. 

\begin{lemma}\label{lem:growing:linearly:forwards:app}
There exists a constant $\xi = \xi(\QQ) > 0$ such that the following holds. Let $W \in \W$ and $w \in \cL_R \cap \SS(W)$, and let $\G$ be the forwards growth sequence with seed $Q \in \cP(W,w)$. Then
\begin{equation}\label{eq:growing:linearly:forwards:app}
\tau(Q_j) \ge \tau(Q_0) + \xi \cdot j
\end{equation}
for all $j \ge 1/\xi$. 
\end{lemma}

\begin{proof}
By Definition~\ref{def:fext:P}, for each $j \ge 0$ there exists $x \in \big( Q_{j+1} \setminus Q_j \big) \cap \Z^d$. It follows that $x \in \Delta(Q_{j+1})$, by Lemma~\ref{lem:fext:app}, and therefore
\begin{equation}\label{eq:growing:forwards:face}
x \in \Delta( Q_{j+1}, W \cup \{u\})
\end{equation}
for some $u \in N_\QQ(W,w)$, by Lemma~\ref{lem:Delta:forwards:app}. Note that $\big( x + (W \cup \{u\})^\perp \big) \cap \Z^d$ contains a copy of the lattice $\L(W \cup \{u\})$, since $x \in \Z^d$, and that $x \not\in \Delta(Q_i)$ for every $i \ne j + 1$, by Lemma~\ref{lem:maxtau:inDelta}. It follows that~\eqref{eq:growing:forwards:face} can hold at most
$$O\big( \tau(Q_j) - \tau(Q_0) + 1 \big)$$
times in the first $j$ steps for each $u \in N_\QQ(W,w)$, where the implicit constant depends on $W \cup \{u\}$ and $w$. Since $w \in \cL_R$, it follows that there exists $\xi = \xi(\QQ,R) > 0$ such that~\eqref{eq:growing:linearly:forwards:app} holds for all $j \ge 1/\xi$, as required.
\end{proof}

We will next prove Lemma~\ref{lem:poly-seq-eats-end}.

\begin{lemma}\label{lem:poly-seq-eats-end:app}
Let $W \in \W$, and let $P,Q \in \cP(W)$ be such that
$$Q \subset P, \qquad a(P) = a(Q), \qquad t(P) = t(Q) \qquad \text{and} \qquad w(P) = w(Q) \in W^\perp.$$
Let $\G = (Q_j)_{j \ge 0}$ be the forwards growth sequence with seed $Q$, and let $m$ be maximal such that $Q_m \subset P$. Then
\[
P \cap \Z^d \subset Q_m.
\]
\end{lemma}

\begin{proof}
Note that $m$ is finite, by Lemma~\ref{lem:growing:linearly:forwards:app}. Let $Q_i = P(W,w;a,t,\tau_i)$ for each $i \ge 0$, and let $\tau = \tau(P)$, so
$$\tau_m \le \tau < \tau_{m+1},$$ 
by Definition~\ref{def:fext:P}. 
We claim that $\tau_P(x) \le \tau_m$ for every $x \in P \cap \Z^d$, which implies that $x \in Q_m$, as required. Indeed, if $\tau_P(x) > \tau_m$ then $\tau_P(x) \ge \tau_{m+1}$, by Definition~\ref{def:fext:P}. But then $x \not\in P$, which is a contradiction. This proves the lemma.
\end{proof}

We move on now to (general) growth sequences, which were defined in Definition~\ref{def:growth:sequence}. Our next proof is of Lemma~\ref{lem:growth-seq-eats-sides}. 

\begin{lemma}\label{lem:growth-seq-eats-sides:app}
Let $W \in \W$, $w \in \cL_R \cap \SS(W)$ and $t > C$, and let $Q \in \cP(W,w)$ with $t(Q) \le t$.  Let $\G$ be a growth sequence with seed $Q$, and let $m$ be maximal such that $t\big( Q^{(m)}_0 \big) \le t$. Then
$$P \cap \Z^d \subset Q^{(m)}_*$$
for some $P \in \cP(W,w;t)$.
\end{lemma}

\begin{proof}
Let $Q^{(m)}_* = P(W,w;a,t',\tau)$ and, recalling Definition~\ref{def:ext:P}, let $y \in \interior\big( P(W) \big)$ and $\eps' > 0$ be such that $Q^{(m+1)}_0 = P(W,w;a-\eps' y,t'+\eps',\tau)$. For each $\eps \in \R$, define
$$Q(\eps) := P(W,w;a-\eps y,t'+\eps,\tau),$$ 
and observe that $Q(\eps) \cap \Z^d = Q^{(m)}_* \cap \Z^d$ for every $0 \le \eps < \eps'$, by Definition~\ref{def:ext:P}. 

Now, by the maximality of $m$, we have 
$$t' = t\big( Q^{(m)}_0 \big) \le t < t\big( Q^{(m+1)}_0 \big) = t' + \eps',$$
and hence there exists $0 \le \eps  < \eps'$ such that $t = t' + \eps$. It follows that 
$$Q(\eps) \in \cP(W,w;t) \qquad \text{and} \qquad Q(\eps) \cap \Z^d \subset Q^{(m)}_*,$$
as required. 
\end{proof}

Our next task is to prove Lemma~\ref{lem:growing:linearly:sideways}. Recall that a polytope $P$ is a grower if either 
$$t\big( \ext(P) \big) \ge t(P) + \xi,$$
or there exists $W' \in \W^\perp(P)$ such that 
$$\Delta\big( \ext(P), W' \big) \cap \Z^d \ne \emptyset,$$
and that a growth sequence $\G$ is happy if $Q^{(i)}_*$ is a grower for every $i \in \N$. 

\begin{lemma}\label{lem:growing:linearly:sideways:app}
Let $W \in \W$ and $w \in \cL_R \cap \SS(W)$, and let $\G$ be a happy growth sequence with seed $Q \in \cP(W,w)$. Then
\begin{equation}\label{eq:speed:of:growing:sideways:app}
t(Q^{(i)}_0) \ge t(Q) + \xi' \cdot i
\end{equation}
for every $i \ge 1/\xi'$. 
\end{lemma}

\begin{proof}
Consider, for each $W' \in \W^\perp(P)$, the set of $i \in \N$ such that
$$\Delta\big( Q^{(i+1)}_0, W' \big) \cap \Z^d \ne \emptyset.$$
As in the proof of Lemma~\ref{lem:growing:linearly:forwards:app}, for each $m \in \N$ this set contains at most
$$O\big( t(Q^{(m)}_0) - t(Q) + 1 \big)$$
elements less than $m$, where the implicit constant depends only on $W'$. Since $\G$ is happy and $W' \subset \QQ$, it follows that
$$t( Q^{(i+1)}_0 ) \ge t( Q^{(i)}_0 ) + \xi$$
for all but $O\big( t(Q^{(m)}_0) - t(Q) \big)$ elements of $[m]$, where the implicit constant depends on $\QQ$. We therefore obtain 
$$t(Q^{(m)}_0) - t(Q) \ge \Big( m - O\big( t(Q^{(m)}_0) - t(Q) \big) \Big) \cdot \xi,$$
which implies~\eqref{eq:speed:of:growing:sideways:app} for some $\xi' > 0$.  
\end{proof}

Finally, let us discuss how to construct a happy growth sequence. Observe first that, by Lemma~\ref{lem:growing:linearly:forwards}, every $1/\xi$ consecutive members of a forwards growth sequence contains at least one polytope $Q_j$ such that 
$$\tau(Q_{j+1}) \ge \tau(Q_j) + \xi.$$
Using this fact,  we may adjust both the `forwards' and `backwards' ends of our polytope $Q^{(i)}_*$ so that there is no lattice point $x \in \Z^d$ such that  
$$x \not\in Q^{(i)}_* \qquad \text{but} \qquad x \pm \xi w \in Q^{(i)}_*.$$
However, when we extend $Q^{(i)}_*$ (and therefore increase $t$), we might  (at least, in theory) nevertheless find a lattice point arbitrarily close to $Q^{(i)}_*$ in one of the new forwards (or backwards) faces of $Q^{(i+1)}_0$. This is not really a problem, however, since it implies the existence of a copy of some lattice $\L(W \cup \{u\})$ very close to $Q^{(i)}_*$, and we can choose $Q^{(i)}_*$ to avoid this, via a modification of the proof of Lemma~\ref{lem:growing:linearly:forwards}. 

In order to make the observations above precise, recall from~\eqref{def:Pminus} the definition of the polytope $P^-$. The following lemma allows us to construct happy growth sequences. 

\begin{lemma}\label{lem:finding:a:grower:app}
Let $W \in \W$ and $w \in \cL_R \cap \SS(W)$, and let $Q \in \cP(W,w)$. Then there exists $Q' \in \cP(W,w)$, with
$$t(Q) = t(Q'), \qquad |\tau(Q) - \tau(Q')| \le 1 \qquad \text{and} \qquad a(Q) - a(Q') = \mu w$$
for some $0 \le \mu \le 1$, such that $Q'$ is a grower. 
\end{lemma}

\begin{proof}
In order to show that $Q'$ is a grower, we need to choose $a(Q')$ so that $Q'$ either intersects or is sufficiently far from every lattice $x + \L(W \cup \{u\})$, where $x \in \Z^d$ and $u \in N_\QQ(W)$ with $\<u,w\> < 0$, and then choose $\tau(Q')$ so that the same holds for those $u \in N_\QQ(W)$ with $\<u,w\> > 0$. 

In order to choose $a(Q')$, observe that for each $u \in N_\QQ(W)$ with $\<u,w\> < 0$, there are $O(1)$ copies of $\L(W \cup \{u\}) \subset \Z^d$, that intersect $Q - w$ but not $Q$, where the implicit constant depends only on $\QQ$ and $R$. Since $\xi = \xi(\QQ,R)$ is sufficiently small, it follows that there exists $\mu \in [0,1]$ such that if
$$Q_0 := Q - \mu w = P(W,w;a,t,\tau),$$
then 
$$\Delta\big( Q_0(\eps), W \cup \{u\} \big) \cap \Z^d = \emptyset$$ 
for every $0 < \eps < \xi$ and every $u \in N_\QQ(W)$ with $\<u,w\> < 0$, where 
$$Q_0(\eps) := P\big( W,w; a - \eps y, t + \eps, \tau \big),$$
for each $\eps \in \R$ and some (arbitrary) $y \in \interior\big( P(W) \big)$ (cf.~Definition~\ref{def:ext:P}).

Similarly, for each $u \in N_\QQ(W)$ with $\<u,w\> > 0$, there are a bounded number of copies of $\L(W \cup \{u\}) \subset \Z^d$, that intersect $Q_0 + w$ but not $Q_0$, and therefore there exists $\tau \le \tau' \le \tau + 1$ such that if
$$Q' := P(W,w;a,t,\tau'),$$
then 
$$\Delta\big( Q'(\eps), W \cup \{u\} \big) \cap \Z^d = \emptyset$$ 
for every $0 < \eps < \xi$ and every $u \in N_\QQ(W)$ with $\<u,w\> > 0$, where 
$$Q'(\eps) := P\big( W,w; a - \eps y, t + \eps, \tau' \big),$$
for each $\eps \in \R$. It follows that either 
$$t\big( \ext(Q') \big) \ge t(Q') + \xi,$$
or there exists $W' \in \W^\perp(P)$ such that 
$$\Delta\big( \ext(Q'), W' \big) \cap \Z^d \ne \emptyset,$$
as required. 
\end{proof}

\subsection{Buffers}\label{app:buffers:sec}

In this section we prove the four lemmas of Section~\ref{subsec:buffers}, namely Lemmas~\ref{lem:newbuffer:contains:oldbuffer}--\ref{lem:sideways:event}. Recall from Definition~\ref{def:buffer} that
$$B(P) := \big\{ x \in P : \tau_P(x) > \tau - C \big\}.$$
We begin with Lemma~\ref{lem:newbuffer:contains:oldbuffer}.

\begin{lemma}\label{lem:newbuffer:contains:oldbuffer:app}
Let $W \in \W$, $w \in \cL_R \cap \SS(W)$ and $P \in \cP(W,w)$. Then
$$\bigcup_{y \in \Delta(P)} \big\{ x \in P : \|x - y\| \le R_0 \big\} \subset B(P).$$
\end{lemma}

\begin{proof}
Let $x \in P = P(W,w;a,t,\tau)$, and suppose that there exists $u \in N_\QQ(W,w)$ and $y \in \Delta( P, W \cup \{u\})$ such that $d(x,y) \le R_0$. By Lemma~\ref{lem:P(W,w,a,t,tau):def2:app}, and since $\< y-x,u \> \leq \| x-y \| \leq R_0$, we have 
$$\big\< x - a - \tau w, \, u \big\> \ge \big\< y - a - \tau w, \, u \big\> - R_0 = t - R_0,$$
using the fact that $\delta(u,w) = 1$ for all $u \in N_\QQ(W,w)$. Since $C = C(\QQ,R)$ is sufficiently large and $\<w,u\> > 0$, it follows that\footnote{Note that here we need $C > R_0 \cdot \<w,u\>^{-1}$ for all $u \in \QQ$ and $w \in \cL_R \cap \SS^{d-1}$ with $\<w,u\> > 0$.}
$$\big\< x - a - (\tau - C) w, \, u \big\> \ge t - R_0 + C \< w,u \> > t,$$
and therefore, by Lemma~\ref{lem:P(W,w,a,t,tau):def2:app}, $x \not\in P(W,w;a,t,\tau - C)$. But by Lemma~\ref{lem:Delta:forwards:app}, every $y \in \Delta(P)$ is an element of $\Delta(P,W\cup\{u\})$ for some $u \in N_\QQ(W,w)$, so we are done.
\end{proof}

Next we prove Lemmas~\ref{lem:forwards:event} and~\ref{lem:fext:fret:buffers}. Both straightforward consequences of the definitions and Lemma~\ref{lem:fext:app}.

\begin{lemma}\label{lem:forwards:event:app}
Let $W \in \W$, $w \in \cL_R \cap \SS(W)$ and $P \in \cP(W,w)$, and set $P' := \fext(P)$. Then
\begin{equation}\label{eq:forwards:event:app}
B(P') \cap \Z^d \subset B(P) \cup \Delta(P').
\end{equation}
\end{lemma}

\begin{proof}
Let $x \in B(P') \cap \Z^d$, and note that, since $B(P') \subset P'$ and
$$P' \cap \Z^d \subset P \cup \Delta(P')$$
by Lemma~\ref{lem:fext:app}, we may assume that $x \in P$. We want to show that $x \in B(P)$, so suppose instead (since $x \in P$) that $\tau_P(x) \le \tau(P) - C$. This implies that 
$$x \in P\big( W,w;a,t,\tau(P) - C \big) \subset P\big( W,w;a,t,\tau(P') - C \big),$$
where $P = P(W,w;a,t,\tau(P))$, and so $\tau_{P'}(x) \le \tau(P') - C$. But this means that $x \not\in B(P')$, which is a contradiction. 
\end{proof}

\begin{lemma}\label{lem:fext:fret:buffers:app}
Let $W \in \W$, $w \in \cL_R \cap \SS(W)$ and $P \in \cP(W,w)$, and set $P' := \fret(P)$. Then
$$B(P) \cap \Z^d \subset B(P') \cup \Delta(P).$$
\end{lemma}

\begin{proof}
The proof is identical to that of Lemma~\ref{lem:forwards:event:app}, swapping the roles of $P$ and $P'$.  
\end{proof}

Our next task is to show that the forward half $F(P) = P \cap ( P + \tau w / 2 )$ of a polytope $P = P(W,w;a,t,\tau)$ satisfies
$$F(P) = P\big( W,w; a + \tau w/2,t,\tau/2 \big).$$
We shall prove a slightly more general statement, since this more general form will also be useful in Section~\ref{sec:deterministic}.

\begin{lemma}\label{lem:Pzero:properties1}
Let $W \in \W$ and $w \in \cL_R \cap \SS(W)$, and let $P = P(W,w;a,t,\tau) \in \cP(W,w)$ and $c \in [0,1]$. Then 
$$P \cap \big( P + c \tau w \big) = P\big( W,w; a + c\tau w,t,(1-c)\tau \big) \in \cP(W,w).$$
\end{lemma}

\begin{proof}
Observe that, by~\eqref{def:PWattau} and Lemma~\ref{lem:tube-is-tube}, and since $w \in W^\perp$, we have
$$P = a + t \cdot \bigcup_{\lambda \in [0,1]} \big( P(W) + \lambda (\tau/t) w \big),$$
and therefore 
$$P + c \tau w = a + c \tau w + t \cdot \bigcup_{\lambda \in [0,1]} \big( P(W) + \lambda (\tau/t) w \big).$$
It follows that
$$P \cap \big( P + c \tau w \big) = a + c \tau w + t \cdot \bigcup_{\lambda \in [0,1-c]} \big( P(W) + \lambda (\tau/t) w \big),$$
and hence, again by Lemma~\ref{lem:tube-is-tube}, we have 
$$P \cap \big( P + c \tau w \big) = P\big( W,w; a + c\tau w,t,(1-c)\tau \big) \in \cP(W,w),$$
as claimed.
\end{proof}

Finally, we need to prove Lemma~\ref{lem:sideways:event}.

\begin{lemma}\label{lem:sideways:event:app}
Let\/ $W \in \W$,\/ $w \in \cL_R \cap \SS(W)$ and\/ $P \in \cP(W,w)$, and set\/ $P' := \ext(P)$ and\/ $P'' := \fret(P')$. If\/ $\tau(P) \ge 5C$, then
\begin{equation}\label{eq:sideways:event:app}
B(P'') \cap \Z^d \subset P \cup \bigcup_{W' \in \W^\perp(P)} \interior\big( \Delta(F,W') \big), 
\end{equation}
where $F := F(P')$. 
\end{lemma}

This lemma is not quite so straightforward, and requires some preliminary results. The first of these is the following consequence of Lemmas~\ref{lem:quasi-cell},~\ref{lem:facets-to-cells} and~\ref{lem:tube-is-tube}. 

\begin{lemma}\label{lem:forwardsfaces:far}
Let $W \in \W$ and $w \in \cL_R \cap \SS(W)$, and let $P = P(W,w;a,t,\tau) \in \cP(W)$. If $x \in \Delta(P)$, then 
$$\<x-a,w\> > \tau.$$  
\end{lemma}

\begin{proof}
By Lemma~\ref{lem:Delta:forwards:app}, there exists $u \in N_\QQ(W,w)$ such that $x \in \Delta(P,W \cup \{u\})$, and hence
\begin{equation}\label{eq:forwardsfaces-far-x}
\< x - a - \tau w, u \> = t,
\end{equation}
by Lemma~\ref{lem:P(W,w,a,t,tau):def2:app} and since $\delta(u,w) = 1$. Also, by Lemma~\ref{lem:quasi-cell}, and since $\< u,w \> > 0$, we have $\< v,w \> > 0$ for all $v \in \Cell_\QQ(u)$. 

Suppose first that $\tau = 0$, and set $y := t^{-1}(x-a)$, so $y \in P(W \cup \{u\})$. Since $y \in P(\emptyset)$ and $\< y,u \> = 1$, it follows by Lemma~\ref{lem:facets-to-cells} that $y / \| y \| \in \Cell_\QQ(u)$. As noted above, this implies that $\< y,w \> > 0$, and hence $\<x-a,w\> > 0$, completing the proof when $\tau = 0$.

For general $\tau \geq 0$, we have $\< u,w \> > 0$ for all $u \in N_\QQ(W,w)$, and therefore
\[
\Delta(P) = \bigcup_{u \in N_\QQ(W,w)} \Delta( P, W \cup \{u\}) = \bigcup_{u \in N_\QQ(W,w)} \big( \Delta( P_0, W \cup \{u\}) + \tau w \big) = \Delta(P_0) + \tau w,
\]
where $P_0 := P(W,w;a,t,0)$, the first and third equalities by Lemma~\ref{lem:Delta:forwards:app}, and the second by Lemmas~\ref{lem:tube-faces-to-P-faces} and~\ref{lem:P(W,w,a,t,tau):def2:app}. Since $\< x'-a, w \> > 0$ for each $x' \in \Delta(P_0)$, by the case $\tau = 0$, it follows that $\< x-a, w \> > \tau$ for each $x \in \Delta(P)$, as required.
\end{proof}

For the rest of this section, let us fix $W \in \W$ and $w \in \cL_R \cap \SS(W)$, and polytopes $P \in \cP(W,w)$, $P' := \ext(P) = P(W,w;a,t,\tau)$, $P'' := \fret(P')$, and 
$$F := F(P') = P' \cap \big( P' + \tau w / 2 \big)$$ 
as in the statement of Lemma~\ref{lem:sideways:event:app}. Note that
\begin{equation}\label{eq:F-rep}
F = P\big( W,w; a + \tau w/2,t,\tau/2 \big) \in \cP(W,w),
\end{equation}
by Lemma~\ref{lem:Pzero:properties1}. We shall need the following property of $F$. 

\begin{lemma}\label{lem:Pzero:properties2}
If $\tau(P) \ge 5C$, then  
$$\big\{ x \in \interior\big( \Delta(P',W') \big) : \tau_{P'}(x) > \tau(P') - 2C \big\} \subset \interior\big( \Delta(F,W') \big)$$
for every $W' \in \W^\perp(P)$.
\end{lemma}

\begin{proof}
Let $x \in \interior\big( \Delta(P',W') \big)$ be such that $\tau_{P'}(x) > \tau(P') - 2C$, which is the same as $\tau_{P'}(x) > \tau - 2C$ since $P' = P(W,w;a,t,\tau)$. We must show that $x \in \interior\big( \Delta(F,W') \big)$, where
$$\Delta(F,W') = P\big( W',w; a + \tau w/2,t,\tau/2 \big),$$
by~\eqref{eq:F-rep}.

Note that $W' \subset \{w\}^\perp$, since $W' \in \W^\perp(P)$. By Lemma~\ref{lem:P(W,w,a,t,tau):def2:app}, it follows that 
$$\big\< x - (a  + \tau w / 2) - \delta(u,w) \tau w / 2, \, u \big\> = \big\< x - a - \delta(u,w) \tau w, \, u \big\> = t$$
for every $u \in W'$, since $x \in \Delta(P',W') = P(W',w;a,t,\tau)$ and $\< u,w \> = 0$. Therefore, by Lemma~\ref{lem:interior:def2}, we only need to show that 
\begin{equation}\label{eq:Pzero:properties:finalclaim}
\big\< x - a  - \tau w / 2 - \delta(u,w) \tau w / 2, \, u \big\> < t
\end{equation}
for every $u \in N_\QQ(W')$. Since $x \in \interior\big( \Delta(P',W') \big)$, we have
$$\big\< x - a  - \delta(u,w) \tau w, \, u \big\> < t$$
for every $u \in N_\QQ(W')$, which is what we want if $\delta(u,w) = 1$. It therefore only remains to prove~\eqref{eq:Pzero:properties:finalclaim} when $\delta(u,w) = 0$; that is, 
when $u \in N_\QQ(W',-w)$.

To do so, we claim first that 
\begin{equation}\label{eq:Pzero:provingthe:finalclaim}
\big\< x - a - \tau_{P'}(x)w, u \big\> \leq t
\end{equation}
for all $u \in N_\QQ(W')$. This follows by Lemma~\ref{lem:P(W,w,a,t,tau):def2:app}, since
\[
x \in \Delta\big( P(W',w;a,t,\tau_{P'}(x) ) \big) \subset P(W',w;a,t,0) + \tau_{P'}(x)w,
\]
the first inclusion holding by Lemma~\ref{lem:maxtau:inDelta} and since $\tau_{P'}(x) > 0$, and the second by Lemma~\ref{lem:tube-is-tube} and~\eqref{def:PWattau}, 
and since $w \in W'^\perp$.

Since $\<u,w\> < 0$ for each $u \in N_\QQ(W',-w)$, to deduce~\eqref{eq:Pzero:properties:finalclaim} from~\eqref{eq:Pzero:provingthe:finalclaim}, we need to show that $\tau_{P'}(x) > \tau/2$. This follows because
$$\tau_{P'}(x) > \tau - 2C > \frac{\tau}{2},$$
by our assumptions, and since $\tau(P) \ge 5C$ implies that $\tau > 4C$. This completes the proof of~\eqref{eq:Pzero:properties:finalclaim}, and hence of the lemma.
\end{proof}

We can now deduce our final buffer lemma.

\begin{proof}[Proof of Lemma~\ref{lem:sideways:event:app}]
Let $x \in B(P'') \cap \Z^d$. We must show that
\[
x \in P \cup \bigcup_{W' \in \W^\perp(P)} \interior\big( \Delta(F,W') \big).
\]
By Lemma~\ref{lem:exterior:faces:app}, we have $\interior(P') \cap \Z^d \subset P$, so it follows that either $x \in P$, in which case we are done, or $x \in P' \setminus \interior(P')$, since $B(P'') \subset P'' \subset P'$. By the definition~\eqref{def:intP} of $\interior(P')$, we may therefore assume that $x \in \Delta\big( P',W \cup \{u\} \big)$ for some $u \in N_\QQ(W)$, and it follows, by Lemma~\ref{lem:max:face:interior:app}, that 
$$x \in \interior\big( \Delta(P',W') \big)$$ 
for some $W \subsetneq W' \in \W$. 

We claim that $W' \in \W^\perp(P)$. To show this, note first that
$$\fret(P') \cap \Delta(P') = \emptyset$$
by Lemma~\ref{lem:maxtau:inDelta}, so $W' \not\in \W^\to(P)$, by Definition~\ref{def:forwards:boundary}. We therefore need to show that
$$W' \cap N_\QQ(W,-w) = \emptyset.$$
Thus, let $P' = P(W,w;a,t,\tau)$, and observe that 
$$\tau_{P'}(x) = \tau_{P''}(x) > \tau(P'') - C > 0,$$
since $x \in B(P'')$ and $\tau(P'') \ge \tau(P) - C > C$. It follows that $x \not\in P(W,w;a,t,0)$, and hence $W' \in \W^\perp(P)$, as claimed. 

Finally, by Lemma~\ref{lem:Pzero:properties2}, if $x \in \interior\big( \Delta(P',W') \big)$ for some $W' \in \W^\perp(P)$, and 
$$\tau_{P'}(x) = \tau_{P''}(x) > \tau(P'') - C > \tau(P') - 2C,$$ 
then $x \in \interior\big( \Delta(F,W') \big)$, as required.
\end{proof}

To conclude this appendix we prove the following lemma in preparation for Appendix~\ref{cover:app}. The lemma may be viewed as a discrete variant of Lemma~\ref{lem:tube-is-tube}.

\begin{lemma}\label{lem:tube-is-discrete-tube}
There exists $\gamma' > 0$ depending only on $\QQ$ such that the following holds. If $W \subset \QQ$ and $w \in \cL_R \cap \SS(W)$, then
\begin{equation}\label{eq:tube-is-discrete-tube}
P(W,w) = \bigcup_{\lambda \in \Omega} \big( P(W) + \lambda w \big)
\end{equation}
for some finite set $\Omega \subset [0,1]$ such that $|x - y| \ge \gamma'$ for every $x,y \in \Omega$ with $x \ne y$. 
\end{lemma}

\begin{proof}
Note first that, by Lemma~\ref{lem:tube-is-tube:app}, the right-hand side of~\eqref{eq:tube-is-discrete-tube} is contained in the left-hand side for every $\Omega \subset [0,1]$. For the reverse inclusion, we choose $\Omega$ by adding elements to the set $\{0,1\}$ one by one, until every point of the interval $[0,1]$ lies within distance $\gamma'$ of some point of $\Omega$. Let $x \in P(W,w)$ and note that, by Lemma~\ref{lem:tube-is-tube:app}, there exists $\lambda_x \in [0,1]$ such that $x \in P(W) + \lambda_x w$. We claim that 
\begin{equation}\label{eq:x:in:left:or:right}
x \in \big( P(W) + a w \big) \cup \big( P(W) + b w \big)
\end{equation}
where $a = \max\{ \lambda \in \Omega : \lambda \le \lambda_x \}$ and $b = \min\{ \lambda \in \Omega : \lambda \ge \lambda_x \}$. 

To prove~\eqref{eq:x:in:left:or:right}, note first that $b - a \le 2\gamma'$, by our choice of $\Omega$. Now, if~\eqref{eq:x:in:left:or:right} fails to hold then fix $a < a' \le \lambda_x \le b' < b$ with $a'$ minimal and $b'$ maximal such that 
$$x \in \big( P(W) + a' w \big) \cap \big( P(W) + b' w \big),$$
noting that $a'$ and $b'$ exist because $x \in P(W) + \lambda_x w$. Observe (cf.~Lemmas~\ref{lem:Delta:forwards:app} and~\ref{lem:maxtau:inDelta}) that there exist $u,v \in \QQ$, with $\<u,w\> > 0 > \<v,w\>$, such that 
$$x \in \big( P(W \cup \{u\}) + a' w \big) \cap \big( P(W \cup \{v\}) + b' w \big).$$
In particular, we have $\<x - aw, u\> > \<x - a'w, u\> = 1$.

Now, since $\<u,w\> > 0 > \<v,w\>$, it follows from Lemma~\ref{lem:quasi-innerprod} that $uv \notin E\big(\Vor(\QQ)\big)$. Therefore, by Lemma~\ref{lem:faces-to-cliques}, we have $P(W \cup \{u,v\}) = \emptyset$, and hence, by Lemma~\ref{lem:far:from:nonfaces}, 
$$\<x - b'w, u\> \le 1 - \gamma,$$ 
since $x - b' w \in P(W \cup \{v\})$. Combining this with the inequality $\<x - aw, u\> > 1$, we deduce that $(b'-a) \<w,u\> > \gamma$, which contradicts the fact that $b' - a < b - a \le 2\gamma'$ if the constant $\gamma' = \gamma'(\QQ)$ is chosen sufficiently small.
\end{proof}

\section{Some technical details from Sections~\ref{sec:deterministic} and~\ref{proof:sec}}\label{sideways:app}

In this short appendix we provide some technical details that were omitted from the final two sections. First, in Section~\ref{app:sideways:claims}, we will prove two claims from Section~\ref{sec:deterministic}; then, in Section~\ref{app:t:properties}, we will prove some inequalities involving the functions defined in Section~\ref{proof:sec}. 

\subsection{Some inclusions involving polytopes}\label{app:sideways:claims}

Our first task is to prove Claim~\ref{claim:Qplus:minusQ}, from the proof of Lemma~\ref{lem:sideways-step-part1}. We will deduce the claim from the following lemma. 

\begin{lemma}\label{lem:Qplus:minusQ}
Let $W \in \W$ and $w \in \cL_R \cap \SS(W)$, and let $Q_1 = P(W,w;a,t,\tau) \in \cP(W)$. If
$$Q_2 = P\big( W,w;a-3\tau'w,t,\tau + 3\tau' \big)$$
for some $\tau' \ge \sqrt{C} \cdot t$, then there exists a polytope $Q' \in \cP( W,w; t, \tau')$ with 
$$Q' \subset Q_2 \setminus Q_1.$$
\end{lemma}

\begin{proof}
Since $w \in W^\perp$, it follows from Lemma~\ref{lem:tube-is-tube} and~\eqref{def:PWattau} that
$$Q_1 = \bigcup_{0 \le \lambda \le \tau} \big( Q_0 + \lambda w \big) \qquad \text{and} \qquad Q_2 = \bigcup_{-3\tau' \le \lambda \le \tau} \big( Q_0 + \lambda w \big)$$
for some $Q_0 \in \cP(W,w;t,0)$. We claim that 
\begin{equation}\label{def:Qprime:in:QplusminusQ}
Q' := \bigcup_{-2\tau' \le \lambda \le -\tau'} \big( Q_0 + \lambda w \big) \subset Q_2 \setminus Q_1,
\end{equation}
which will suffice to prove the lemma, since $Q' \in \cP(W,w;t,\tau')$. To prove the inclusion in~\eqref{def:Qprime:in:QplusminusQ}, note that $Q' \subset Q_2$ is immediate, and that $Q' \cap Q_1 = \emptyset$ follows from
$$\max_{x \in Q'} \< x,w\> = \max_{x \in Q_0} \< x,w\> - \tau' < \min_{x \in Q_0} \< x,w\> = \min_{x \in Q_1} \< x,w\>,$$
where the inequality holds since $\tau' \ge \sqrt{C} \cdot t$, and $C = C(\QQ)$ is sufficiently large.
\end{proof}

Recall from Section~\ref{sec:growth-sideways} that we fixed $W \in \W_k$ and $w \in \cL_R \cap \SS(W)$, and polytopes 
$$P \in \cP(W,w), \qquad P' = \ext(P) \qquad \text{and} \qquad P'' = \fret(P')$$
as in the statement of Lemma~\ref{lem:sideways-step-part2}, so $P'$ is long (meaning $\tau' \ge Ct(P')$, where $\tau' := \tau(P')$) and sideways edge-filled by $A$ (see Definition~\ref{def:sideways-edge-filled}). Recall also from~\eqref{def:Pjs} that 
\begin{equation}\label{app:def:Pjs}
P_j := P' \cap \big( P' + 4(k - j) \alpha \tau' \cdot w \big)
\end{equation}
where $\alpha = 1/8d$, from the statement of Lemma~\ref{lem:sideways-step-part1} that $0 \le j \le k - 1$, that $W' \in \W^\perp(P)$ with $\dim(W'^\perp) = j$, and that $Q = \interior\big( \Delta(P_{j+1},W') \big)$, and from the proof of Lemma~\ref{lem:sideways-step-part1} that $Q^* = \clint\big( \Delta(P',W') \big)$ and $t^* = t(Q^*)$. 

To prove Claim~\ref{claim:Qplus:minusQ} (restated below), we apply Lemma~\ref{lem:Qplus:minusQ} to two suitably chosen polytopes that are both contained in $Q^*$. 

\begin{claim}\label{claim:Qplus:minusQ:app}
There exists a polytope $Q' \in \cP\big( W',w; t^*, \alpha \tau' \big)$ with 
$$Q' \subset \big( Q^* \cap Q \big) \setminus P_j.$$ 
\end{claim}

\begin{proof}[Proof of Claim~\ref{claim:Qplus:minusQ:app}]
Let 
$$P' = P(W,w;a',t,\tau') \qquad \text{and} \qquad Q^* =  P(W',w;a^*,t^*,\tau'),$$
and, recalling~\eqref{app:def:Pjs}, set $c = \big( 4(k - j) - 1 \big) \alpha \in [0,1]$ and define 
$$Q_1 = P(W',w; a,t^*,\tau) \qquad \text{and} \qquad Q_2 =  P(W',w;a - 3\alpha \tau' w,t^*,\tau + 3\alpha \tau'),$$
where $a = a^* + c \tau' w$ and $\tau = (1-c)\tau'$. Note that $Q_1 \subset Q_2 \subset Q^*$, by Lemma~\ref{lem:Pzero:properties1}, and since $c \ge 3\alpha$. Note that $\alpha \tau' = \tau' / 8d \ge \sqrt{C} \cdot t \ge \sqrt{C} \cdot t^*$, since $\tau' \ge Ct$, and that $w \in W'^\perp$, since $W' \in \W^\perp(P)$ and recalling~\eqref{def:Wperp}. It therefore follows from Lemma~\ref{lem:Qplus:minusQ} that there exists a polytope $Q' \in \cP\big( W',w; t^*, \alpha \tau' \big)$ with 
$$Q' \subset Q_2 \setminus Q_1.$$
It therefore only remains to show that 
$$Q_2 \setminus Q_1 \subset \big( Q^* \cap Q \big) \setminus P_j.$$
We have already observed that $Q_2 \subset Q^*$, so our next task is to show that
\begin{equation}\label{eq:intersection:of:interiors}
Q_2 \subset Q = \interior\big( \Delta(P_{j+1},W') \big) = \interior\big( \Delta(P',W') \big) \cap \Big( \interior\big( \Delta(P',W') \big) + c' \tau' w \Big),
\end{equation}
where $c' = 4(k - j - 1) \alpha \in [0,1]$. To do so, recall that $w \in W'^\perp$, and observe that 
$$Q_2 = Q^* \cap \big( Q^* + c' \tau' w \big),$$
by Lemma~\ref{lem:Pzero:properties1}, and since $c - c' = 3\alpha$. Since $Q^* \subset \interior\big( \Delta(P',W') \big)$, by Lemma~\ref{lem:interiors-are-canonical}, we obtain~\eqref{eq:intersection:of:interiors}. Finally, we are required to show that 
$$\big( Q_2 \setminus Q_1 \big) \cap P_j = \emptyset.$$
As in the proof of Lemma~\ref{lem:Qplus:minusQ}, we do so by considering the inner product with $w$. To be precise, observe that
$$\max_{x \in Q_2 \setminus Q_1} \< x,w\> < \min_{x \in Q_1} \< x,w\> + \frac{\alpha \tau'}{2} < \min_{x \in P_j} \< x,w\>,$$
where the first inequality holds since $Q_1 = Q_2 \cap ( Q_2 + 3\alpha \tau' w)$, by Lemma~\ref{lem:Pzero:properties1}, the second by the definition~\eqref{app:def:Pjs} of $P_j$ and Definition~\ref{def:closed:interior}, and where for both inequalities we used the bound $\tau' \ge \sqrt{C} \cdot t$, and the fact that $C$ is a sufficiently large constant. 
\end{proof}

Our next task is to prove Claim~\ref{claim:extensions:cover:tube}, which is also part of the proof of Lemma~\ref{lem:sideways-step-part1}. Recall that $\G$ is the forwards growth sequence whose seed\footnote{Here $W' \in \W^\perp(P)$, $w = w(P)$ and $t^* = t(Q^*)$ are all as in Claim~\ref{claim:Qplus:minusQ:app}.} $Q_0 \in \cP(W',w;t^*)$ satisfies $Q_0 \subset Q'$, where $Q' \in \cP\big( W',w; t^*, \alpha \tau' \big)$ is the polytope with $Q' \subset \big( Q^* \cap Q \big) \setminus P_j$ whose existence is guaranteed by Claim~\ref{claim:Qplus:minusQ:app}, and that $m \in \N$ is maximal such that $Q_m \subset Q^*$. 

\begin{claim}\label{claim:extensions:cover:tube:app}
\begin{equation}\label{eq:claim:extensions:cover:tube:app}
\interior\big( \Delta(P_j,W') \big) \cap \Z^d \subset Q_m.
\end{equation}
\end{claim}

\begin{clmproof}{claim:extensions:cover:tube:app}
Since $Q_0 \in \cP(W',w;t^*)$ and $Q_0 \subset Q^*$, by Lemma~\ref{lem:Pzero:properties1} it follows that there exists $c \ge 0$ such that the polytope
$$Q'' := Q^* \cap \big( Q^* + c w \big)$$ 
satisfies $a(Q'') = a(Q_0)$ and $Q_0 \subset Q''$. Our plan is to apply Lemma~\ref{lem:poly-seq-eats-end} to the polytopes $Q_0$ and $Q''$, so note first that
$$t(Q'') = t(Q_0) \quad \text{and} \quad w(Q'') = w(Q_0) \in W^\perp,$$
by Lemma~\ref{lem:Pzero:properties1} and our assumption that $w \in \SS(W)$. Observe also that $m$ is maximal such that $Q_m \subset Q''$. Indeed, since $Q'' \subset Q^*$ we have $Q_i \not\subset Q''$ for all $i > m$, and if $x \in Q_m \subset Q^*$ then $x \in Q^* + cw$, and hence $x \in Q''$, by Definitions~\ref{def:fext:P} and~\ref{def:forwards:growth:sequence}, and since $Q_0 \subset Q^* + cw$ and $c \ge 0$. By Lemma~\ref{lem:poly-seq-eats-end}, it follows that
$$Q'' \cap \Z^d = Q^* \cap \big( Q^* + c w \big) \cap \Z^d \subset Q_m.$$

To deduce~\eqref{eq:claim:extensions:cover:tube:app}, we are required to show that 
\begin{equation}\label{eq:claim:extensions:final:inclusion}
\interior\big( \Delta(P_j,W') \big) \cap \Z^d \subset Q^* \cap \big( Q^* + c w \big).
\end{equation}
To prove~\eqref{eq:claim:extensions:final:inclusion}, it will suffice to show that $c \le 4(k - j) \alpha \tau'$, since this will imply that $P_j \subset P' \cap \big( P' + c w \big)$, by~\eqref{app:def:Pjs} and Lemma~\ref{lem:Pzero:properties1}, and moreover that
$$\Delta(P_j,W') \subset \Delta(P',W') \cap \big( \Delta(P',W') + c w \big).$$
since $w \in W'^\perp$. Recalling that $Q^* = \clint\big( \Delta(P',W') \big)$, it will then follow by Lemma~\ref{lem:interiors-are-canonical} that~\eqref{eq:claim:extensions:final:inclusion} holds, as required. The claimed bound on $c$ follows easily from the fact that $Q_0 \subset P' \setminus P_j$ (cf.~the proof of Claim~\ref{claim:Qplus:minusQ:app}), so this completes the proof of~\eqref{eq:claim:extensions:cover:tube:app}. 
\end{clmproof}

\subsection{Bounding the functions in Section~\ref{proof:sec}}\label{app:t:properties}

In this section we will prove the various simple inequalities stated in Observation~\ref{obs:t:properties}. 

\begin{obs}\label{obs:t:properties:app}
Let\/ $1 \le s \le k \le d$, and let\/ $p > 0$ be sufficiently small. Then 
\begin{equation}\label{eq:t:properties:app}
t_0(k,s,p) \ge t_1(k,s,p)^{8d}.
\end{equation}
If $s^* := \min\{s, k-1\}$, then
\begin{equation}\label{eq:t:properties:cpp}
t_0(k,s,p) \ge t_0(k-1,s^*,p)^{8d},
\end{equation}
and if $s \ge 2$, then
\begin{equation}\label{eq:t:properties:bpp}
t_0(k-1,s-1,p) \ge 2 \cdot \log t_1(k-1,s^*,p).
\end{equation}
Moreover, $t_1(1,1,p) \ge t_0(1,1,p) \ge p^{-2R_0}$.
\end{obs}

\begin{proof}
Each inequality follows easily from~\eqref{def:t0} and~\eqref{def:t1}, the definitions of $t_0(k,s,p)$ and $t_1(k,s,p)$. Indeed, for~\eqref{eq:t:properties:app} observe that
$$\log t_0(k,1,p) = \lambda(k) \log(1/p) \ge 8d \cdot C^{4d} \cdot \log(1/p) = 8d \cdot \log t_1(k,1,p),$$
since $\lambda(k) = (8d)^k \cdot C^{3d}$, and that if $s \ge 2$, then  
$$\log t_0(k,s,p) \ge t_0(k-1,s-1,p)^{8d} \gg t_0(k-1,s-1,p)^{7d} \log(1/p) \gg \log t_1(k,1,p),$$
since $\lambda(k) = 8d \cdot \lambda(k-1)$.\footnote{Note that when $s \ge 3$ we only need $\lambda(k) > \lambda(k-1)$ here.}  
For~\eqref{eq:t:properties:cpp}, note that 
$$t_0(k,k,p) \ge t_0(k-1,k-1,p)^{8d},$$
which proves the case $s = k$, and that 
$$t_0(k,s,p) \ge t_0(k-1,s,p)^{8d},$$
which proves the case $s < k$. Finally, for~\eqref{eq:t:properties:bpp}, note that if $k \ge 2$ then 
$$t_0(k-1,k-1,p) \gg  t_0(k-2,k-2,p)^{7d} \log(1/p) \gg \log t_1(k-1,k-1,p),$$
which proves the case $s = k$, and that if $2 \le s < k$ then
$$t_0(k-1,s-1,p) \gg t_0(k-2,s-1,p)^{7d} \log(1/p) \gg \log t_1(k-1,s,p),$$
since $\lambda(k-1) = 8d \cdot \lambda(k-2)$. Finally, for the last part, note that 
$$t_0(1,1,p) = p^{-8d C^{3d}} \ge p^{- C^{3d}} = t_1(1,1,p),$$
and recall that $C$ was chosen to be sufficiently large.
\end{proof}

\section{Perfectly covering a polytope with smaller polytopes}\label{cover:app}

In this final appendix we prove Lemma~\ref{lem:cover:exists}, which says that there exists a bounded degree perfect cover of a large polytope $P$ by a small polytope $Q$. Recall from~\eqref{def:GC:edges} that, given a cover of $P$ by copies of $Q$, we define the graph $G_\C$ on vertex set $\C$ to have edge set
\[
E(G_\C) = \big\{ Q_1Q_2 : d(Q_1,Q_2) \le 2R_0 \big\}.
\]

\begin{lemma}\label{lem:cover:exists:app}
There exists a constant $\Delta > 0$ depending only on $\QQ$ such that the following holds. Let $W \in \W$ and $w \in \cL_R \cap \SS^{d-1}$, and let $P,Q \in \cP(W,w)$ satisfy
\begin{equation}\label{eq:PQ:condition:app}
\diam(Q) \le \frac{t(P)}{\Delta}.
\end{equation}
Then there exists a perfect cover $\C$ of $P$ by copies of $Q$ such that $\Delta( G_\C ) \le \Delta$. 
\end{lemma}

The proof is by induction on $k = \dim(W^\perp)$, and the main challenge will be to (perfectly) cover the points within distance $O(1)$ of the boundary of $P$. Roughly speaking, we cover each face of $P$ (using the induction hypothesis; see Lemma~\ref{lem:cover:exists:IH}), and show that the union of these covers is contained in $P$ and covers all points sufficiently close to the faces. We then repeat this process $O(1)$ times (except using slightly smaller polytopes). Once we have done this, it will then be straightforward to complete the covering. 

The main complication in the proof is that we must choose a suitable induction hypothesis so that our perfect covers of the faces also cover all points close to the faces. In order to do so, we shall use the following definition (cf.~Definition~\ref{def:closed:interior}).

\begin{definition}\label{def:eps:interior}
For each $W \in \W$ and $\eps > 0$, define the \emph{$\eps$-interior} of a polytope $P = P(W,w;a,t,\tau) \in \cP(W)$ to be
$$\interior^{\eps}(P) = P(W,w;a + \eps y,t - \eps,\tau)$$
for some $y \in \interior\big( P(W) \big)$. 
\end{definition}

We remark that our bound on $\Delta$ will depend on the choices of $y$ in Definition~\ref{def:eps:interior}. More precisely, let us fix, for each $W \in \W$, a vector $y = y(W) \in \interior( P(W) )$ to use when defining the $\eps$-interiors of polytopes in $\cP(W)$, and define
$$\sigma_u(W) := 1 - \<y(W),u\>$$
for each $u \in N_\QQ(W)$. Note that $\sigma_u(W) > 0$ for each $W \in \W$ and $u \in N_\QQ(W)$, by Lemma~\ref{lem:interior:def2}, so 
$$\sigma := \min_{W \in \W} \min_{u \in N_\QQ(W)} \sigma_u(W) > 0.$$
The following description of $\interior^{\eps}(P)$ is an immediate consequence of Lemma~\ref{lem:P(W,w,a,t,tau):def2:app}, using the fact that $y \in \interior( P(W) )$ (cf.~Lemmas~\ref{lem:Ds:same:shift:app} and~\ref{lem:interior:def2}). 

\begin{lemma}\label{lem:eps:interior:props:app}
Let $W \in \W$ and $P = P(W,w;a,t,\tau) \in \cP(W)$, and let $\eps > 0$. Then 
\begin{multline}
\interior^{\eps}(P) = \bigcap_{u \in W} \Big\{ x \in \R^d : \big\< x - a - \delta(u,w) \tau w, \, u \big\> = t \Big\} \\ \cap \bigcap_{u \in N_\QQ(W)} \Big\{ x \in \R^d : \big\< x - a - \delta(u,w) \tau w, \, u \big\> \le t - \eps \cdot \sigma_u(W) \Big\}.
\end{multline}
\end{lemma}

In particular, by Lemma~\ref{lem:interior:def2}, it follows that
$$\interior^{\eps}(P) \subset \interior(P).$$
In order to state our induction hypothesis, we need to fix three sequences of constants, which are chosen as follows. First, choose
\begin{equation}\label{def:epsilons:sequence}
1 \gg \eps(0) \gg \eps'(0) \gg \eps(1) \gg \eps'(1) \gg \cdots \gg \eps(d) \gg \eps'(d) > 0,
\end{equation}
chosen from left to right, and depending on $\QQ$ and $\sigma$, and then 
\begin{equation}\label{def:cks:sequence}
1 \ll c(0) \ll c(1) \ll \cdots \ll c(d),
\end{equation}
again chosen from left to right, and depending on the sequences $\eps(k)$ and $\eps'(k)$. 

We also need the following notion of maximum degree, which is more convenient for the induction step than $\Delta( G_\C )$. Given a finite collection $\C$ of copies of $Q$, define
$$\hat{\Delta}(\C) := \max_{\ell \ge 1} \frac{1}{\ell^d} \max_{a \in \R^d} \Big| \Big\{ Q' \in \C : d\big( a+Q,Q' \big) \le \ell \Big\} \Big|;$$
thus, in the inner maximization, we choose a translate $a + Q$ with the maximal number of elements of $\C$ lying within distance $\ell$ of that translate. Observe that $\hat{\Delta}(\C) \leq |\C|$ (in particular, it is finite); that $\Delta( G_\C ) \le (2R_0)^d \cdot \hat{\Delta}(\C)$ (by taking $\ell = 2R_0$ and using the definition of $G_\C$ from~\eqref{def:GC:edges}), and that $\hat{\Delta}(\C_1 \cup \C_2) \le \hat{\Delta}(\C_1) + \hat{\Delta}(\C_2)$ (by performing both maximizations separately for $\C_1$ and $\C_2$). The following observation will be used in the proof Lemma~\ref{lem:cover:exists:IH} and is the main motivation for the definition of $\hat{\Delta}$.

\begin{obs}\label{obs:Deltahat-intQ-to-Q}
Let $\C$ be a cover of $P$ by copies of $\interior^{\eps'(k)}(Q)$, for some $0 \leq k \leq d$, and let $\C'$ be the cover of $P$ obtained from $\C$ by replacing the copies of $\interior^{\eps'(k)}(Q)$ by the corresponding copies of $Q$. Then
\[
\hat{\Delta}(\C') \leq 2^d \cdot \hat{\Delta}(\C).
\]
\end{obs}

\begin{proof}
This is just a consequence of the fact that if $d(a + Q, Q') \leq \ell$, then
\[
d\big( a + \interior^{\eps'(k)}(Q), \, \interior^{\eps'(k)}(Q') \big) \leq \ell + O\big( \eps'(k) \big) \leq 2\ell
\]
since $\ell \geq 1$ and $\eps'(k)$ was chosen sufficiently small.
\end{proof}

To reduce the number of parameters we need to deal with during the induction, we shall assume that $a(P) = \0$ and $\tau(P) = 0$, and use Lemma~\ref{lem:tube-is-discrete-tube} to deduce the general version from this case. We also set $a(Q) = \0$ and rescale by $t(Q)$, so that $Q = P(W, w)$ for some $W \in \W$ and $w \in \cL_R$. Recall from Lemma~\ref{lem:tube-faces-to-P-faces:app} that if $w \not\in W^\perp$ then $P(W,w)$ is just a translate of $P(W)$, and if $w \in W^\perp$ then
\begin{equation}\label{eq:diamQ:bounds}
\diam(Q) = \Theta\big( 1 + \|w\| \big).
\end{equation}
We will prove the following statement by induction on $k = \dim(W^\perp)$.

\begin{lemma}\label{lem:cover:exists:IH}
Let $1 \le k \le d$, $W \in \W_k$, $w \in \cL_R$ and $t \ge c(k) \big( 1 + \|w\| \big)$, and set
\begin{equation}\label{eq:cover:exists:PandQ}
P = t \cdot P(W) \qquad \text{and} \qquad Q = P(W,w).
\end{equation}
There exists a perfect cover\/ $\C$ of\/ $P$ with copies of\/ $Q$ such that\/ $\hat{\Delta}(\C) \le c(k)$, and every point of the $\eps(k)$-interior of $P$ is contained in the $\eps'(k)$-interior of some member of\/ $\C$.
\end{lemma}

The proof of Lemma~\ref{lem:cover:exists:IH} will take up most of the rest of this section; in the lemmas below, we will assume that $1 \le k \le d$, $W \in \W_k$, $w \in \cL_R$ and $t \ge c(k) ( 1 + \|w\| )$ are fixed, set $P = t \cdot P(W)$, and assume that $Q$ is an arbitrary translate of $P(W,w)$ with $\0 \in Q$ (since we are only interested in copies of $Q$, this does not affect the statement). Note that $\aff(P)$ and $\aff(Q)$ are both translates of $W^\perp$, by Lemma~\ref{lem:clique-of-correct-dim}, and hence $\aff(Q) = W^\perp$. The first step is to construct a covering of the $c(0)$-interior of $P$ with $\eps(0)$-interiors. 

\begin{lemma}\label{lem:covering:the:middle}
There exists a collection\/ $\C'$ of copies of\/ $Q$, each contained in $P$, such that $\hat{\Delta}(\C') \le c(0)$, and such that every point of the $c(0)$-interior of $P$ is contained in the $\eps(0)$-interior of some member of $\C'$. 
\end{lemma}

The proof of this lemma is not too difficult, but the details require a little care. The rough idea is simply to choose a suitable lattice and take one copy of $Q$ for each lattice point. When $w \in W^\perp$, we will need the following simple lemma. 

\begin{lemma}\label{lem:covering:tube-fits-near-faces}
Suppose that $w \in W^\perp$, and that $a \in \{w\}^\perp$ satisfies $a + Q \subset \aff(P)$ and $\<a,u\> \le t - 2\sqrt{d}$ for all\/ $u \in N_\QQ(W) \cap \{w\}^\perp$. Then 
\[
a + Q \subset P.
\]
\end{lemma}

\begin{proof}
Since $a + Q \subset \aff(P)$, and recalling Lemma~\ref{lem:P-nbrs}, we only need to prove that $\<x,v\> \le t$ for all $x \in a + Q$ and $v \in N_\QQ(W)$. Suppose first that $v \in \{w\}^\perp$, and note that we have $\<y,v\> \le 2\sqrt{d}$ for every $y \in Q$, by the definition of $P(W,w)$, and since $\{e_1,\ldots,e_d\} \subset \QQ$ and $\0 \in Q$. By our assumption that $\<a,v\> \le t - 2\sqrt{d}$, it follows that
\[
\<x,v\> = \<x-a,v\> + \<a,v\> \le t,
\]
since $x - a \in Q$ and $v \in N_\QQ(W) \cap \{w\}^\perp$. 

To deal with the remaining case, when $v \in N_\QQ(W) \setminus \{w\}^\perp$, we first claim that if $x \in P(\emptyset) \cap \{w\}^\perp$ is such that $\<x,u\> = 1$ for some $u \in \QQ$, then $\<u,w\> = 0$. Indeed,  by Lemma~\ref{lem:facets-to-cells} we have $x/\|x\| \in \Cell_\QQ(u)$, and the claim then follows by Lemma~\ref{lem:quasi-new}. Since $P(\emptyset) \cap \{w\}^\perp$ is compact, $v \in \QQ \setminus \{w\}^\perp$ and $w \in \cL_R$, it follows that there exists $\eps > 0$, depending only on $\QQ$, such that $\<x,v\> \leq 1 - \eps$ for all $x \in P(\emptyset) \cap \{w\}^\perp$.

We next claim that $a / t \in P(W) \subset P(\emptyset)$. To see this, note first that $a \in \aff(P)$, since $\aff(Q) = W^\perp$ and $\aff(P)$ is a translate of $W^\perp$, and therefore $\<a,u\> = t$ for all $u \in W$, by~\eqref{eq:P(W)}. By Lemma~\ref{lem:P-nbrs}, it will therefore suffice to show that $\<a,u\> \le t$ for all $u \in N_\QQ(W)$. Choose $u \in N_\QQ(W)$ with $\<a,u\>$ maximal, and let $\lambda > 0$ be such that $\< \lambda a / t, u\> = 1$. By our choice of $u$, we have $\< \lambda a / t, u'\> \le 1$ for every $u' \in N_\QQ(W)$, and thus $\lambda a / t \in P(W) \subset P(\emptyset)$. Now, by Lemma~\ref{lem:facets-to-cells}, it follows that $a / \|a\| \in \Cell_\QQ(u)$, and hence, by Lemma~\ref{lem:quasi-new} and since $a \in \{w\}^\perp$ and $w \in \cL_R$, we deduce that $\<u,w\> = 0$. But we have $\<a,u\> \le t$ for all $u \in N_\QQ(W) \cap \{w\}^\perp$, by assumption, so $a / t \in P(W)$, as claimed.

Combining the observations above, and recalling that $a \in \{w\}^\perp$, we deduce that $\<a,v\> \leq (1-\eps)t$ for some $\eps > 0$ depending only on $\QQ$. Now, observe that if $y \in Q$, then $\|y\| \le 2\sqrt{d} + \|w\|$, by Lemma~\ref{lem:tube-is-tube} and since $\{e_1,\ldots,e_d\} \subset \QQ$ and $\0 \in Q$. It follows that if $x \in a + Q$, then
\[
\<x,v\> = \<a,v\> + \<x-a,v\> \leq (1-\eps)t + 2\sqrt{d} + \|w\| \le t
\]
provided $t \geq (2\sqrt{d} + \|w\|)/\eps$, which holds since $c(k)$ was chosen sufficiently large.
\end{proof}

We can now construct our covering of the $c(0)$-interior of $P$.

\begin{proof}[Proof of Lemma~\ref{lem:covering:the:middle}]
We will construct our covering by defining a suitable lattice, and then considering the set of translates $x + Q$ such that $x$ is in the lattice, and $x + Q \subset P$. We divide the proof into two cases, depending on whether or not $w \in W^\perp$. 

\medskip

\noindent Case 1: $w \not\in W^\perp$.
\medskip

Note that $Q$ is a translate of $P(W)$, by Lemma~\ref{lem:tube-faces-to-P-faces}, and recall that $\aff(Q) = W^\perp$ and that $\aff(P)$ is a translate of $W^\perp$. Let $b_1, \dots, b_k$ be an orthonormal basis for $W^\perp$, and consider the lattice
\[
\eps(0) \cdot \big( b_1\Z + \dots + b_k\Z \big).
\]
Let $\Lambda_0$ be a copy of this lattice embedded in $\aff(P)$, and let
$$\Lambda := \big\{ x \in \Lambda_0 \,:\, x + Q \subset P \big\} \qquad \text{and} \qquad \C' := \big\{ x + Q \,:\, x \in \Lambda \big\}.$$ 
Note that every member of $\C'$ is contained in $P$, by construction. The bound on $\hat{\Delta}(\C')$ also follows easily from the definition, since there are at most $O\big(\ell/\eps(0)\big)^k$ lattice points in a ball of radius $\ell + 2 \cdot \diam(Q)$, and $c(0)$ is allowed to depend on $\eps(0)$. 

It remains to show that every point of the $c(0)$-interior of $P$ is contained in the $\eps(0)$-interior of some member of $\C'$. To do so, let $y \in P$, and suppose that $y \not\in \interior^{\eps(0)}(Q')$ for every $Q' \in \C'$. Let $B$ be the $k$-dimensional $\ell_\infty$-ball (with respect to the basis $b_1,\dots,b_k$) of radius $\eps(0)$ embedded in $W^\perp$, and observe that, since $\eps(0)$ is sufficiently small and $Q$ is an arbitrary translate of $P(W)$ containing the origin, we may assume that
$$B \subset \interior^{\eps(0)}(Q),$$
and therefore 
$$y \in \aff(P) = \bigcup_{x \in \Lambda_0} \big( x + B \big) = \bigcup_{x \in \Lambda_0} \big( x + \interior^{\eps(0)}(Q) \big).$$
Since $y \not\in \interior^{\eps(0)}(Q')$ for every $Q' \in \C'$, it follows that $y \in x + Q \not\subset P$ for some $x \in \Lambda_0$. Since $\diam(Q) \le 2\sqrt{d}$, and recalling that $c(0)$ is sufficiently large, we deduce that $y$ is not in the $c(0)$-interior of $P$, as required.

\medskip
\noindent Case 2: $w \in W^\perp$.
\medskip

In this case we need to slightly tweak the argument above. Set $\hat{w} := w / \|w\|$, and let $b_1,\dots,b_{k-1},\hat{w}$ be an orthonormal basis for $W^\perp$. Let $\Lambda_0$ be an embedding of the lattice 
$$\eps(0) \cdot \big( b_1\Z + \dots + b_{k-1}\Z \big)$$ 
in $\aff(P) \cap \{w\}^\perp$, and for each $x \in \Lambda_0$, define 
\[
I(x) := \big\{ \mu \in \R : x + \mu \hat{w} + Q \subset P \big\}.
\]
Note that $I(x)$ is either empty, or a closed and bounded interval, and let $J(x) \subset I(x)$ be a finite subset containing the endpoints of $I(x)$, and containing either one or two points of each sub-interval of $I(x)$ of length $\eps(0) + \|w\|$. Define
\[
\Lambda := \big\{ x + \mu \hat{w} \,:\, x \in \Lambda_0 \,\text{ and }\, \mu \in J(x) \big\} \qquad \text{and} \qquad \C' := \{ x + Q : x \in \Lambda \},
\]
and note that every member of $\C'$ is contained in $P$, since $J(x) \subset I(x)$. The claimed bound on $\hat{\Delta}(\C')$ also follows easily, since $J(x)$ contains at most two points of each sub-interval of $I(x)$ of length $\eps(0) + \|w\|$, and there are therefore at most $O\big(\ell/\eps(0)\big)^k$ members of $\C'$ within distance $\ell$ of a given copy of $Q$. 

It remains to show that the $c(0)$-interior of $P$ is contained in the set 
\[
C := \bigcup_{Q' \in \C'} \interior^{\eps(0)}(Q').
\]
To do so, let $y \in P \setminus C$, and observe that, by the definition of $\Lambda_0$, there exist $x \in \Lambda_0$, $x' \in W^\perp \cap \{w\}^\perp$ and $\lambda \in \R$ such that 
$$y = x + x' + \lambda\hat{w} \qquad \text{and} \qquad \|x'\|_\infty \leq \eps(0)$$
(with respect to the basis $b_1,\dots,b_{k-1}$). Suppose first that $x + Q \not\subset P$. Noting that $x \in \{w\}^\perp$ and $x + Q \subset \aff(P)$, it follows by Lemma~\ref{lem:covering:tube-fits-near-faces} that $\<x,u\> > t - 2\sqrt{d}$ for some $u \in N_\QQ(W) \cap \{w\}^\perp$, and hence $\<y,u\> > t - 3\sqrt{d}$. But this implies that $y \not\in \interior^{c(0)}(P)$, by Lemma~\ref{lem:eps:interior:props:app} and since $c(0)$ was chosen to be sufficiently large depending on $\sigma$. 

On the other hand, if $x + Q \subset P$, then $I(x)$ is non-empty. Let $B_0$ be the $k$-dimensional $\ell_\infty$-ball (with respect to the basis $b_1,\dots,b_{k-1},\hat{w}$) of radius $\eps(0)$ embedded in $W^\perp$, and let
\[
B := \bigcup_{\mu\in[0,1]} \big( B_0 + \mu w \big)
\]
(so that $B$ is a hyperrectangle). As in Case~1, and using Lemma~\ref{lem:tube-is-tube}, we may assume that $B \subset \interior^{\eps(0)}(Q)$, and therefore that $x + B \subset C$ for every $x \in \Lambda$. Moreover, 
\begin{equation}\label{eq:covering:J-is-cts}
\bigcup_{\mu \in I(x)} \big(x + \mu\hat{w} + B\big) = \bigcup_{\mu \in J(x)} \big(x + \mu\hat{w} + B\big) \subset C
\end{equation}
for every $x \in \Lambda_0$, since $J(x)$ contains at least one element in each subinterval of $I(x)$ of length $\eps(0) + \|w\|$, which is also the length of $B$ in direction $w$. Observe that it also follows from~\eqref{eq:covering:J-is-cts} that $\lambda \notin I(x)$, since $y \in x + \lambda\hat{w} + B$, but $y \notin C$. Let $I(x) = [\lambda_0,\lambda_1]$, and without loss of generality let us assume that $\lambda > \lambda_1$.

Now, note that there exists $z \in Q_1 := x + \lambda_1\hat{w} + Q$ such that $z$ is in the forwards boundary of $P$. In particular, there exists $v \in N_\QQ(W)$ with $\<v,w\> > 0$ such that
\[
\< z,v \> = t.
\]
Let $y' = y - \lambda'\hat{w}$, where $\lambda' \ge 0$ is minimal such that $y' \in Q_1$ (note that such a $\lambda'$ exists since $y \in x + \lambda\hat{w} + B$ and $x + \lambda_1\hat{w} + B \subset Q_1$, and since $\lambda > \lambda_1$). Observe that both $y'$~and $z$ are in the copy of $P(W)$ that contains the forwards boundary of $Q_1$, and therefore $d(y',z) \le 2\sqrt{d}$. Since $\<v,w\> > 0$ and $\lambda' \ge 0$, it follows that
\[
\< y,v \> \ge \< y',v \> \ge \< z,v \> - 2\sqrt{d} = t - 2\sqrt{d}, 
\]
and hence $y$ is not in the $c(0)$-interior of $P$, as required.
\end{proof}

It remains to cover the points of $P$ that are within a bounded distance of one of the faces. In order to do so, we will use the induction hypothesis to cover each of the faces of $P$ with copies of the corresponding face of $Q$, and then repeat this a constant number of times. The next lemma shows that the corresponding copies of $Q$ are contained in $P$. 

\pagebreak

\begin{lemma}\label{lem:faces:subsets:are:subsets}
Let\/ $W \subset W' \in \W$. If\/ $a + \Delta(Q,W') \subset \Delta(P,W')$, then $a + Q \subset P$.
\end{lemma}

\begin{proof}
We need to check that $x \in a + Q$ satisfies the various equations and inequalities that define $P$. Observe first that $a + Q$ and $P$ are in the same translate of $W^\perp$, since $a + \Delta(Q,W') \subset \Delta(P,W')$, so by Lemma~\ref{lem:P(W,w,a,t,tau):def2:app} it suffices to show that 
\begin{equation}\label{eq:faces:subsets:check}
\< x, u \> \le t
\end{equation}
for every $u \in N_\QQ(W)$. To prove~\eqref{eq:faces:subsets:check} when $u \in W'$, let $y \in a + \Delta(Q,W') \subset \Delta(P,W')$ and observe that
$$\<x,u\> \le \<y,u\> = t.$$
So we may assume that $u \notin W'$. If $P(W' \cup \{u\})$ is non-empty, then let $y \in a + \Delta(Q,W' \cup \{u\})$, and observe that $y \in \Delta(P,W') \subset P$, so again
$$\<x,u\> \le \<y,u\> \le t.$$
Finally, if $P(W' \cup \{u\})$ is empty, then by Lemma~\ref{lem:far:from:nonfaces} we have 
$$\<x',u\> \le 1 - \gamma$$ 
for every $x' \in P(W')$, and therefore
$$\<y,u\> \le t - \gamma t$$
for every $y \in a + \Delta(Q,W') \subset \Delta(P,W')$. Since $x,y \in a + Q$, and recalling~\eqref{eq:diamQ:bounds}, and that $t \ge c(k) \cdot \big( 1 + \|w\| \big)$ and $c(k)$ is sufficiently large, it follows that 
$$\< x - y, u \> \le \gamma \cdot c(k) \cdot \big( 1 + \|w\| \big) \le \gamma t,$$ 
and therefore we obtain~\eqref{eq:faces:subsets:check} for all $u \in N_\QQ(W)$, as required.
\end{proof}

We also need the following slight variant of the lemma above, which follows from almost the same proof. 

\begin{lemma}\label{lem:interior:faces:implies:subset}
Let $\eps \ge \eps' \ge 0$, set $P' := \interior^{\eps}(P)$ and $Q' := \interior^{\eps'}(Q)$, and let $W \subset W' \in \W$. If $a + \Delta(Q',W') \subset \Delta(P',W')$, then $a + Q \subset P$.
\end{lemma}

\begin{proof}
We repeat the proof of Lemma~\ref{lem:faces:subsets:are:subsets}, except we use Lemma~\ref{lem:eps:interior:props:app} in order to deduce~\eqref{eq:faces:subsets:check}. Indeed, given $x \in a + Q$ and $u \in N_\QQ(W)$, if $u \in W'$ then 
$$\<x,u\> \le \<y,u\> + \eps' \sigma_u(W) = t + (\eps' - \eps) \sigma_u(W) \le t$$
for each $y \in a + \Delta(Q',W') \subset \Delta(P',W')$, and similarly if $P(W' \cup \{u\})$ is non-empty, then 
$$\<x,u\> \le \<y,u\> + \eps' \sigma_u(W) \le t + (\eps' - \eps) \sigma_u(W) \le t$$
for each $y \in a + \Delta(Q',W' \cup \{u\}) \subset P'$, in both cases by Lemma~\ref{lem:eps:interior:props:app} and since $\eps \ge \eps' \ge 0$. The rest of the proof is the same. 
\end{proof}

The next ingredient that we need for our proof of Lemma~\ref{lem:cover:exists:IH} says that if $x \in P$ is not in the $\eps(k)$-interior of $P$ then it is covered by one of the families obtained by the induction hypothesis applied to the faces. To show this, we choose a point $y$ on the boundary of $P$ that is sufficiently close to $x$, and that is contained in the $\eps(\ell)$-interior of $\Delta(P,W')$ for some $W' \in \W_\ell$. By the induction hypothesis, such a point $y$ is contained in the $\eps'(\ell)$-interior of the $W'$-face of some member of our cover. We shall show that this polytope also contains $x$.  

The point $y = y(x)$ will be chosen by the following process. We define sequences 
$$(x_0, x_1,\ldots,x_m), \qquad k_0 > k_1 > \cdots > k_m \ge 0 \qquad \text{and} \qquad W_0 \subset W_1 \subset \cdots \subset W_m,$$
with $x_0 := x$, $k_0 := k$ and $W_0 := W$, 
such that 
$$W_i \in \W_{k_i}, \qquad x_i \in \Delta(P,W_i), \qquad \text{and} \qquad d(x_{i-1},x_i) \le \kappa \cdot \eps(k_{i-1}) \cdot \sigma_u(W_{i-1})$$
for each $i \in [m]$, where $\kappa = \kappa(\QQ)$ is defined just before Lemma~\ref{lem:innerprod-to-dist}, and such that either $k_m = 0$ or $x_m$ is in the $\eps(k_m)$-interior of $\Delta(P,W_m)$. 

To do so, let $i \ge 0$ and suppose that we have already constructed $W_i \in \W_{k_i}$ and $x_i \in \Delta(P,W_i)$. If either $k_i = 0$ or $x_i$ is in the $\eps(k_i)$-interior of $\Delta(P,W_i)$, then set $m := i$ and stop; otherwise, recalling Lemma~\ref{lem:eps:interior:props:app}, let $u \in N_\QQ(W_i)$ be such that
\begin{equation}\label{eq:E-ip-large-ind}
\< x_i, u \> > t - \eps(k_i) \cdot \sigma_u(W_i).
\end{equation}
Observe that $P(W_i \cup \{u\})$ is non-empty, by Lemma~\ref{lem:far:from:nonfaces}, since $x_i \in \Delta(P,W_i)$ and $\eps(k_i)$ is sufficiently small. Define $x_{i+1}$ to be the nearest element of $\Delta(P,W_i \cup \{u\})$ to $x_i$, and note that
\begin{equation}\label{eq:E-ip-to-dist}
d(x_i,x_{i+1}) \le \kappa \cdot \eps(k_i) \cdot \sigma_u(W_i)
\end{equation}
by~\eqref{eq:E-ip-large-ind} and Lemma~\ref{lem:innerprod-to-dist}. Now, by Lemma~\ref{lem:exists:WinW:trivial}, there exists $W_i \cup \{u\} \subset W_{i+1} \in \W$ with $P(W_{i+1}) = P(W_i \cup \{u\})$. Set $k_{i+1} := \dim(W_{i+1}^\perp)$, and observe that $k_{i+1} < k_i$, by Lemma~\ref{lem:faces:lower:dim}; that $W_{i+1} \in \W_{k_{i+1}}$; and that $x_{i+1} \in \Delta(P,W_{i+1})$. Since the $k_i$ are strictly decreasing, the process eventually stops. 

For each $x \in P$, define $y(x) := x_m$, and observe that 
\begin{equation}\label{eq:distance:x:to:xi}
d(x,y) \le \kappa \cdot \sum_{i = 0}^{m-1} \eps(k_i) \cdot \sigma_u(W_i) \ll \eps'(k_m),
\end{equation}
by~\eqref{eq:E-ip-to-dist} and~\eqref{def:epsilons:sequence}, and since $k_0 > \cdots > k_m$ and $\kappa$ depends only on $\QQ$. We also define $W'(x) := W_m$ and $\ell(x) := k_m$, so
\begin{equation}\label{eq:cover:y:properties}
y(x) \in \Delta(P,W'(x)) \qquad \text{and} \qquad W'(x) \in \W_{\ell(x)},
\end{equation}
and either $\ell(x) = 0$ or $y(x)$ is in the $\eps(\ell(x))$-interior of $\Delta(P,W'(x))$. 

We shall use the induction hypothesis to cover $y(x)$ using the $W'$-face of $Q$, and then apply the next two lemmas to deduce that the corresponding copies of $Q$ cover $x$. The first deals with the case $\ell(x) = 0$. 


\begin{lemma}\label{lem:covering:near:the:edges:0}
Let $x \in P$, and suppose that $\ell(x) = 0$, so $y = y(x)$ is a vertex of $P$. If $a \in \R^d$ is such that $a + \Delta(Q,W') = \{y\}$, where $W' = W'(x)$, then $x \in a + Q$.  
\end{lemma}

\begin{proof}
We shall use Lemma~\ref{lem:faces-far-away:app2} to show that 
\begin{equation}\label{eq:covering:near:the:edges:need}
x - a = (y - a) + (x - y) \in P(W,w) = Q.
\end{equation}
To do so, observe first that $y - a \in P(W',w)$, since $\Delta(Q,W') = \{y - a\}$ by assumption and $P(W',w) = \Delta(Q,W')$, and observe also that $\|x - y\| \le \eps'(k) \le \gamma$, by~\eqref{eq:distance:x:to:xi} and since $\eps'(k)$ is sufficiently small. Next, note that since $x \in P$ and $y \in \Delta(P,W')$, we have $x - y \in W^\perp$ and $\<x - y, v\> \le 0$ for every $v \in W'$. Moreover, since $W' \in \W_0$, we have $P(W' \cup \{v\}) = \emptyset$ for every $v \not\in W'$. Hence, by Lemma~\ref{lem:faces-far-away:app2}, we obtain~\eqref{eq:covering:near:the:edges:need}, as required.
\end{proof}

The next lemma deals with the case $\ell(x) \ge 1$. 

\begin{lemma}\label{lem:covering:near:the:edges:1}
Let $x \in P$, and suppose that $y = y(x)$ is contained in the $\eps'(\ell)$-interior of $a + \Delta(Q,W')$, where $\ell = \ell(x) \ge 1$ and $W' = W'(x)$. Then $x \in a + Q$.  
\end{lemma}

\begin{proof}
By Lemma~\ref{lem:eps:interior:props:app}, since $y$ is contained in the $\eps'(\ell)$-interior of $a + \Delta(Q,W')$, we have
\begin{equation}\label{eq:E-near-edges-y}
\big\< y - a - \delta(u,w) w, \, u \big\> = 1
\end{equation}
for every $u \in W'$, and 
$$\big\< y - a - \delta(u,w) w, \, u \big\> \le 1 - \eps'(\ell) \cdot \sigma_u(W')$$ 
for every $u \in N_\QQ(W')$, and therefore 
\begin{equation}\label{eq:covering:near:the:edges:need2}
\big\< x - a - \delta(u,w) w, \, u \big\> \le 1
\end{equation}
for every $u \in N_\QQ(W')$, by~\eqref{eq:distance:x:to:xi}. Moreover, since $x \in P$ and $y \in \Delta(P,W')$, we have $x - y \in W^\perp$ and hence
\[
\big\< x - a - \delta(u,w) w, \, u \big\> = 1
\]
for all $u \in W \subset W'$, by~\eqref{eq:E-near-edges-y}, and we also have $\<x-y,u\> \le 0$ for every $u \in W'$, so~\eqref{eq:covering:near:the:edges:need2} holds for each such $u$, again by~\eqref{eq:E-near-edges-y}. To deduce that $x \in a + Q$, it only remains to show that~\eqref{eq:covering:near:the:edges:need2} holds for each $u \in N_\QQ(W) \setminus \big( W' \cup N_\QQ(W') \big)$. 

To do so, observe that $P(W' \cup \{u\}) = \emptyset$ for each such $u$, and therefore 
$$\big\< y - a - \delta(u,w) w, \, u \big\> \le 1 - \gamma$$ 
by Lemma~\ref{clm:faces-far-away}, since $y \in a + P(W',w)$. Recalling~\eqref{eq:distance:x:to:xi} and that $\eps'(\ell) \le \gamma$, it follows that~\eqref{eq:covering:near:the:edges:need2} holds, as required.
\end{proof}

Finally, we will use the following simple observation, which follows from Definition~\ref{def:eps:interior}. 

\begin{obs}\label{obs:rescaling}
If $P = t \cdot P(W)$ and $\mu > 0$, then 
$$\mu \cdot \interior^\eps(P) = \interior^{\mu \eps}(\mu \cdot P).$$
\end{obs}

We are finally ready to prove Lemma~\ref{lem:cover:exists:IH}. The plan is to use the induction hypothesis to obtain perfect covers of the faces of $P$ that also cover the points outside the $\eps(k)$-interior of $P$. We then iterate this process a bounded number of times, except using copies of the $\eps'(k)$-interior of $Q$. Finally, we use Lemma~\ref{lem:covering:the:middle} to cover the remaining points of $P$.  

\begin{proof}[Proof of Lemma~\ref{lem:cover:exists:IH}]
The proof is by induction on $k = \dim(W^\perp)$. Suppose first that $k = 1$, and note that in this case $P$ and $Q$ are lines, by Lemma~\ref{lem:clique-of-correct-dim}, and that the length of $P$ is $\Theta(t)$ and the length of $Q$ is either $\Theta(1 + \|w\|)$ (if $w \in W^\perp$) or $\Theta(1)$ (otherwise). Since $t \ge c(1) \cdot (1 + \|w\|)$ and $0 < \eps'(1) \ll \eps(1) \ll 1$, it is trivial to construct a perfect cover $\C$ of $P$ with copies of $Q$ such that $\hat{\Delta}(\C) \le 3$, and such that every point of the $\eps(1)$-interior of $P$ is contained in the $\eps'(1)$-interior of some member of $\C$.

So let $k \ge 2$, and assume that the lemma holds for all smaller values of $k$. First we construct a family $\C_0$ of copies of $Q \subset P$ that cover $P \setminus \interior^{\eps(k)}(P)$ such that all elements of $\C_0$ are contained in $P$ and satisfying $\hat{\Delta}(\C_0) \ll c(k)$. We construct $\C_0$ by applying the induction hypothesis to each face of $P$; in doing so, it will only be necessary to assume that $t \ge c(k-1) \cdot \big( 1 + \|w\| \big)$. 

Let $W \subsetneq W' \in \W$, set $k' := \dim(W'^\perp)$, and observe that $k' < k$, by Lemma~\ref{lem:faces:lower:dim}, that 
$$\Delta(P,W') = t \cdot P(W') \qquad \text{and} \qquad \Delta(Q,W') = P(W', w),$$
and that $t \ge c(k') \cdot (1 + \|w\|)$, by~\eqref{def:cks:sequence}. If $k' \ge 1$, then it follows by the induction hypothesis that there exists a perfect cover $\C(W')$ of $\Delta(P,W')$ with copies of $\Delta(Q,W')$ such that $\hat{\Delta}(\C(W')) \le c(k')$, and every point of the $\eps(k')$-interior of $\Delta(P,W')$ is contained in the $\eps'(k')$-interior of some member of $\C(W')$. If $k' = 0$, on the other hand, then $\Delta(P,W') = \{z\}$ for some vertex $z$ of $P$, by Lemma~\ref{lem:clique-of-correct-dim}, and in this case we set $\C(W') := \{ a + \Delta(Q,W') \}$, where $a + \Delta(Q,W') = \{z\}$. Now, define
$$\C_0 := \bigcup_{W \subsetneq W' \in \W} \big\{ a + Q : a + \Delta(Q,W') \in \C(W') \big\},$$
and observe that, by Lemma~\ref{lem:faces:subsets:are:subsets}, every member of $\C_0$ is contained in $P$. Moreover, 
$$\hat{\Delta}(\C_0) \le 2^{|N_\QQ(W)|} c(k - 1) \ll c(k),$$
by~\eqref{def:cks:sequence}. The following claim will therefore complete this first stage of the proof. 

\begin{claim}
If $x \in P \setminus \interior^{\eps(k)}(P)$, then $x$ is contained in some member of $\C_0$. 
\end{claim}

\begin{proof}
Let $y = y(x)$ and $W' = W'(x)$, and suppose first that $\ell := \ell(x) = 0$. Recall from~\eqref{eq:cover:y:properties} that $y(x) \in \Delta(P,W')$, and observe that therefore $a + \Delta(Q,W') = \{y\}$, where $\C(W') = \{a + \Delta(Q,W')\}$. By Lemma~\ref{lem:covering:near:the:edges:0}, it follows that $x \in a + Q \in \C_0$, as claimed.

If $\ell \ge 1$, on the other hand, then $y$ is in the $\eps(\ell)$-interior of $\Delta(P,W')$, by construction. It follows, by our choice of $\C(W')$, that $y$ is contained in the $\eps'(\ell)$-interior of some $a + \Delta(Q,W') \in \C(W')$. By Lemma~\ref{lem:covering:near:the:edges:1}, it follows that $x \in a + Q \in \C_0$, as required.
\end{proof}

It remains to cover $\interior^{\eps(k)}(P)$. However, our task here is different to before, because of the condition in the statement of the lemma that every point of the $\eps(k)$-interior of $P$ is contained in the $\eps'(k)$-interior of some member of $\C$. In order to construct a perfect cover $\C$ satisfying this constraint, we iterate the above process a bounded number of times, except now we use the $\eps'(k)$-interior of $Q$ to cover the points of $\interior^{\eps(k)}(P)$. We do this by showing that the argument above can be applied (in exactly the same way) to suitable interiors of $P$ and $Q$. To be precise, define
$$Q' := s \cdot P(W,w'),$$
where $s := 1 - \eps'(k)$ and $w' = s^{-1} \cdot w \in \cL_R$. Observe that $Q'$ is a translate\footnote{Indeed, $Q = P(W,w;\0,1,1)$, so $\interior^{\eps'(k)} (Q) = P(W,w;\eps'(k) y,s,1) = \eps'(k) y + s \cdot P(W,w')$.} of $\interior^{\eps'(k)} (Q)$.

 For each $i \in \N$, define $\phi(i) := \big( 1 + (i-1) s \big) \cdot \eps(k)$, and set $m := 2c(1) / \eps(k)$. Now, for each $i \in [m]$, define 
$$P_i := \interior^{\phi(i)} (P) = a_i + t_i \cdot P(W)$$
for some $a_i \in \R^d$, where $t_i = t - \phi(i)$. Observe that 
$$t_i \ge t - 2c(1) \ge c(k) \cdot (1 + \|w\|) - 2c(1) \gg c(k-1) \cdot (1 + \|w'\|),$$
for every $i \in [m]$, by~\eqref{def:cks:sequence} and since $k \ge 2$.

We would like to obtain covers of each $P_i \setminus P_{i+1}$ using copies of $Q'$, or equivalently with copies of $\interior^{\eps'(k)} (Q)$, and then note that the corresponding copies of $Q$ are still contained in $P$ by Lemma~\ref{lem:interior:faces:implies:subset}. However, in order to keep our induction hypothesis as simple as possible, we have only stated it in the case where the set of which we take copies has the form $P(W',w')$, for some $W' \in \W$ and $w' \in \cL_R$. We therefore rescale $P_i \setminus P_{i+1}$ and $Q'$ so that we can apply the induction hypothesis, and then revert the rescaling to obtain the desired cover. Thus, let us apply the argument above to the pair
$$P_i' := s^{-1} P_i \qquad \text{and} \qquad s^{-1} \cdot Q' = P(W,w')$$
to obtain a cover $\C'_i$ of $P_i' \setminus \interior^{\eps(k)}(P_i')$ with copies of $s^{-1} \cdot Q'$ such that $\hat{\Delta}(\C'_i) \ll c(k)$. Now, by Observation~\ref{obs:rescaling}, we have 
$$P_{i+1} = \interior^{\phi(i+1) - \phi(i)} (P_i) = \interior^{s \cdot \eps(k)} (P_i) = s \cdot \interior^{\eps(k)} (P_i'),$$
and therefore, rescaling by a factor of $s$, we obtain a cover $\C_i''$ of $P_i \setminus P_{i+1}$ with copies of $\interior^{\eps'(k)} (Q)$ such that $\hat{\Delta}(\C_i'') \ll c(k)$ (since $s < 1$), as desired. Finally, we replace each copy of $\interior^{\eps'(k)} (Q)$ in $\C_i''$ with the corresponding copy of $Q$ to obtain a cover $\C_i$ of $P_i \setminus P_{i+1}$. Note that this cover satisfies $\hat{\Delta}(\C_i) \leq 2^d \cdot \hat{\Delta}(\C_i'') \ll c(k)$ by Observation~\ref{obs:Deltahat-intQ-to-Q}. Moreover, since $\eps'(k) \ll \eps(k) \leq \phi(i)$, it follows from Lemma~\ref{lem:interior:faces:implies:subset} (applied with $W' = W$) that each copy of $Q$ in $\C_i$ is contained in $P$.

To complete the proof, let $\C'$ be the family constructed in Lemma~\ref{lem:covering:the:middle}, and define
$$\C := \C' \cup \bigcup_{i = 0}^m \C_i.$$
Observe that, by our choice of $m$, we have
$$\interior^{\eps(k)} (P) \setminus \interior^{c(1)} (P) \subset P_1 \setminus P_m,$$
and that, for each $i \in [m]$, the $\eps'(k)$-interiors of the members of $\C_i$ cover $P_i \setminus P_{i+1}$. Recall also, from Lemma~\ref{lem:covering:the:middle}, that the $\eps'(k)$-interiors (which contain the $\eps(0)$-interiors) of the members of $\C'$ cover $\interior^{c(1)} (P)$. Since $\C_0$ covers $P \setminus \interior^{\eps(k)} (P)$, it follows that $\C$ is a cover of $P$ with copies of $Q$, and every point of the $\eps(k)$-interior of $P$ is contained in the $\eps'(k)$-interior of some member of $\C$. Furthermore, each copy of $Q$ in $\C$ is contained in $P$ by our observations above, so $\C$ is a perfect cover of $P$. Finally, since $\C$ is the union of $\C'$ with a bounded number of families with such that $\hat{\Delta}(\C_i) \ll c(k)$, it follows that $\hat{\Delta}(\C) \le c(k)$, as required. This completes the induction step, and hence also the proof of the lemma. 
\end{proof}

Finally, let us deduce Lemma~\ref{lem:cover:exists:app} from Lemma~\ref{lem:cover:exists:IH}. For this we shall use Lemma~\ref{lem:tube-is-discrete-tube}.

\begin{proof}[Proof of Lemma~\ref{lem:cover:exists:app}]
Let $\Delta$ be a sufficiently large constant, depending on $\QQ$; we will assume in particular that $\Delta \gg c(d)$. Let $W \in \W$ and $w \in \cL_R \cap \SS^{d-1}$, set $k := \dim(W^\perp)$, and let $P = P(W,w;a,t,\tau)$ and $Q = P(W,w;a',t',\tau')$ satisfy~\eqref{eq:PQ:condition:app}. 

Recalling~\eqref{def:PWattau}, set $w' := (\tau' / t') w$ and define
$$P' := (t/t') \cdot P(W) \qquad \text{and} \qquad Q' := P(W,w').$$
By~\eqref{eq:PQ:condition:app} and~\eqref{eq:diamQ:bounds}, and noting that $\diam(Q) = t' \cdot \diam(Q')$, we have
$$\frac{t}{t'} \ge \Delta \cdot \diam(Q') \ge c(d) \cdot \big( 1 + \|w'\| \big),$$
since since $\Delta \gg c(d)$. By Lemma~\ref{lem:cover:exists:IH}, 
it follows that there exists a perfect cover $\C'$ of $P'$ with copies of $Q'$ such that $\hat{\Delta}(\C') \le c(k)$. Rescaling by a factor of $t' > C$, we obtain a perfect cover $\C$ of $t \cdot P(W)$ with copies of $Q$ such that $\hat{\Delta}(\C) \le c(k)$.

Now, since $\Delta( G_\C ) \le (2R_0)^d \cdot \hat{\Delta}(\C)$, and $\Delta$ was chosen sufficiently large, we are done if $\tau = 0$. Moreover, by Lemma~\ref{lem:tube-faces-to-P-faces:app}, if $w \not\in W^\perp$ then $P(W,w)$ is just a translate of $P(W)$, and therefore $P$ is a translate of $t \cdot P(W)$. We are therefore also done in this case. 

Finally, if $w \in W^\perp$ and $\tau > 0$, then we need to apply Lemma~\ref{lem:tube-is-discrete-tube}, which provides us with a perfect cover $\C''$ of $P$ with copies of $t \cdot P(W)$. Moreover, for each member of $\C''$ and each $\ell \ge 1$, there are $O(\ell)$ other members of $\C''$ within distance $\ell$. Since we may perfectly cover each copy of $t \cdot P(W)$ with copies of $Q$, as above, it follows that there exists a perfect cover $\C$ of $P$ with copies of $Q$ such that $\Delta( G_\C ) \le \Delta$, as required.
\end{proof}

\section*{Acknowledgements}

The authors would like to thank Hugo Duminil-Copin for a very useful discussion at the outset of this project, and for many other interesting conversations over the years. We are also grateful to the anonymous referee for reading the proof extremely carefully, and for a large number of very helpful comments which improved the presentation.

\bibliographystyle{amsplain}
\bibliography{bprefs}

\providecommand{\noopsort}[1]{}
\providecommand{\bysame}{\leavevmode\hbox to3em{\hrulefill}\thinspace}
\providecommand{\MR}{\relax\ifhmode\unskip\space\fi MR }
\providecommand{\MRhref}[2]{%
  \href{http://www.ams.org/mathscinet-getitem?mr=#1}{#2}
}
\providecommand{\href}[2]{#2}
\begin{thebibliography}{10}

\bibitem{Adler}
J.~Adler, \emph{Bootstrap percolation}, Phys. A \textbf{171} (1991), 453--470.

\bibitem{ADE}
J.~Adler, A.~{\noopsort{Enter}}{van Enter}, and J.A.M.S. Duarte,
  \emph{Finite-size effects for some bootstrap percolation models}, J. Stat.
  Phys. \textbf{60} (1990), no.~3, 323--332.

\bibitem{AL}
M.~Aizenman and J.L. Lebowitz, \emph{Metastability effects in bootstrap
  percolation}, J. Phys. A \textbf{21} (1988), no.~19, 3801--3813.

\bibitem{And93}
E.~Andjel, \emph{Characteristic exponents for two-dimensional bootstrap
  percolation}, Ann. Probab. \textbf{21} (1993), no.~2, 926--935.

\bibitem{AMS}
E.~Andjel, T.~Mountford, and R.~Schonmann, \emph{Equivalence of exponential
  decay rates for bootstrap percolation like cellular automata}, Ann. Inst.
  Henri Poincar\'e Probab. Stat. \textbf{31} (1995), no.~1, 13--25.

\bibitem{ALBB}
F.~Arceri, F.P. Landes, L.~Berthier, and G.~Biroli, \emph{A statistical
  mechanics perspective on glasses and aging}, Encyclopedia of Complexity and
  Systems Science (R.A. Meyers, ed.), Springer, Berlin, Heidelberg, 2022.

\bibitem{BBMSsub}
P.~Balister, B.~Bollob\'as, R.~Morris, and P.~Smith, \emph{Subcritical monotone
  cellular automata}, Random Structures Algorithms, to appear.

\bibitem{BBMSuncomp}
\bysame, \emph{Uncomputability of critical probabilities for monotone cellular
  automata}, in preparation.

\bibitem{BBMSlower}
\bysame, \emph{Universality for monotone cellular automata}, submitted,
  arXiv:2203.13806.

\bibitem{BBPS}
P.~Balister, B.~Bollob\'as, M.~Przykucki, and P.~Smith, \emph{Subcritical
  $\mathcal{U}$-bootstrap percolation models have non-trivial phase
  transitions}, Trans. Amer. Math. Soc. \textbf{368} (2016), 7385--7411.

\bibitem{BBDM}
J.~Balogh, B.~Bollob\'as, H.~Duminil-Copin, and R.~Morris, \emph{The sharp
  threshold for bootstrap percolation in all dimensions}, Trans. Amer. Math.
  Soc. \textbf{364} (2012), no.~5, 2667--2701.

\bibitem{BBM3d}
J.~Balogh, B.~Bollob\'as, and R.~Morris, \emph{Bootstrap percolation in three
  dimensions}, Ann. Probab. \textbf{37} (2009), no.~4, 1329--1380.

\bibitem{Blanq22}
D.~Blanquicett, \emph{The $d$-dimensional bootstrap percolation models with
  threshold at least double exponential}, preprint, arXiv:2201.09029.

\bibitem{Blanq19}
\bysame, \emph{Anisotropic bootstrap percolation in three dimensions}, Ann.
  Probab. \textbf{48} (2020), no.~5, 2591--2614.

\bibitem{BDMSDuarte}
B.~Bollob\'as, H.~Duminil-Copin, R.~Morris, and P.~Smith, \emph{The sharp
  threshold for the {D}uarte model}, Ann. Probab. \textbf{45} (2017), no.~6B,
  4222--4272.

\bibitem{BDMS}
\bysame, \emph{Universality of two-dimensional critical cellular automata},
  Proc. Lond. Math. Soc. \textbf{126} (2023), no.~2, 620--703.

\bibitem{BSU}
B.~Bollob\'as, P.~Smith, and A.~Uzzell, \emph{Monotone cellular automata in a
  random environment}, Combin. Probab. Comput. \textbf{24} (2015), no.~4,
  687--722.

\bibitem{BH16}
A.~Bovier and F.~{\noopsort{Hollander}}{den Hollander}, \emph{Metastability: A
  {P}otential-{T}heoretic {A}pproach}, Grundlehren der mathematischen
  Wissenschaften, vol. 351, Springer, 2015.

\bibitem{CMRT}
N.~Cancrini, F.~Martinelli, C.~Roberto, and C.~Toninelli, \emph{Kinetically
  constrained spin models}, Probab. Theory Related Fields \textbf{140} (2008),
  459--504.

\bibitem{CMRT2}
\bysame, \emph{Kinetically constrained models}, New Trends in Mathematical
  Physics (V.~Sidoravicius, ed.), Springer, 2009, pp.~741--752.

\bibitem{CC}
R.~Cerf and E.~Cirillo, \emph{Finite size scaling in three-dimensional
  bootstrap percolation}, Ann. Probab. \textbf{27} (1999), no.~4, 1837--1850.

\bibitem{CM}
R.~Cerf and F.~Manzo, \emph{The threshold regime of finite volume bootstrap
  percolation}, Stochastic Process. Appl. \textbf{101} (2002), no.~1, 69--82.

\bibitem{CMnuc2}
\bysame, \emph{Nucleation and growth for the {Ising} model in $d$ dimensions at
  very low temperatures}, Ann. Probab. \textbf{41} (2013), no.~6, 3697--3785.

\bibitem{CLR}
J.~Chalupa, P.L. Leath, and G.R. Reich, \emph{Bootstrap percolation on a
  {B}ethe lattice}, J. Phys. C \textbf{12} (1979), no.~1, L31--L35.

\bibitem{CFM}
P.~Chleboun, A.~Faggionato, and F.~Martinelli, \emph{Time scale separation and
  dynamic heterogeneity in the low temperature {E}ast model}, Comm. Math. Phys.
  \textbf{328} (2014), 955--993.

\bibitem{DS01}
P.G. Debenedetti and F.H. Stillinger, \emph{Supercooled liquids and the glass
  transition}, Nature \textbf{410} (2001), 259--267.

\bibitem{DS1}
P.~Dehghanpour and R.~Schonmann, \emph{Metropolis dynamics relaxation via
  nucleation and growth}, Comm. Math. Phys. \textbf{188} (1997), 89--119.

\bibitem{Hugo20}
H.~Duminil-Copin, \emph{Lectures on the {I}sing and {P}otts models on the
  hypercubic lattice}, Random {G}raphs, {P}hase {T}ransitions, and the
  {G}aussian {F}ree {F}ield, PIMS-CRM Summer School in Probability, Springer,
  2020.

\bibitem{DE}
H.~Duminil-Copin and A.~{\noopsort{Enter}}{van Enter}, \emph{Sharp
  metastability threshold for an an\-iso\-tropic bootstrap percolation model},
  Ann. Probab. \textbf{41} (2013), no.~3A, 1218--1242.

\bibitem{DEH}
H.~Duminil-Copin, A.~{\noopsort{Enter}}{van Enter}, and T.~Hulshof,
  \emph{Higher order corrections for anisotropic bootstrap percolation},
  Probab. Theory Related Fields \textbf{172} (2018), 191--243.

\bibitem{DHar}
H.~Duminil-Copin and I.~Hartarsky, \emph{Sharp metastability transition for
  two-dimensional bootstrap percolation with symmetric isotropic threshold
  rules}, preprint, arXiv:2303.13920.

\bibitem{DH}
H.~Duminil-Copin and A.~Holroyd, \emph{Finite volume bootstrap percolation with
  threshold rules on $\mathbb{Z}^2$: balanced case}, Unpublished manuscript,
  2012.

\bibitem{vE}
A.~{\noopsort{Enter}}{van Enter}, \emph{Proof of {S}traley's argument for
  bootstrap percolation}, J. Stat. Phys. \textbf{48} (1987), 943--945.

\bibitem{vEF}
A.~{\noopsort{Enter}}{van Enter} and A.~Fey, \emph{Metastability thresholds for
  anisotropic bootstrap percolation in three dimensions}, J. Stat. Phys.
  \textbf{147} (2012), 97--112.

\bibitem{FSS}
L.~Fontes, R.~Schonmann, and V.~Sidoravicius, \emph{Stretched exponential
  fixation in stochastic {I}sing models at zero temperature}, Comm. Math. Phys.
  \textbf{228} (2002), 495--518.

\bibitem{FA}
G.H. Fredrickson and H.C. Andersen, \emph{Kinetic {I}sing model of the glass
  transition}, Phys. Rev. Lett. \textbf{53} (1984), 1244--1247.

\bibitem{FV17}
S.~Friedli and Y.~Velenik, \emph{{Statistical Mechanics of Lattice Systems: A
  Concrete Mathematical Introduction}}, Cambridge, 2017.

\bibitem{GST}
J.P. Garrahan, P.~Sollich, and C.~Toninelli, \emph{Kinetically constrained
  models}, Dynamical heterogeneities in glasses, colloids, and granular media
  (L.~Berthier, G.~Biroli, J.-P. Bouchaud, L.~Cipelletti, and W.~van Saarloos,
  eds.), International series of monographs in physics, Oxford, 2011,
  pp.~341--369.

\bibitem{GG93}
J.~Gravner and D.~Griffeath, \emph{Threshold growth dynamics}, Trans. Amer.
  Math. Soc. \textbf{340} (1993), no.~2, 837--870.

\bibitem{GG99}
\bysame, \emph{Scaling laws for a class of critical cellular automaton growth
  rules}, Proceedings of the Erd{\Horig{o}}s Center Workshop on Random Walks,
  1999, pp.~167--188.

\bibitem{Har}
I.~Hartarsky, \emph{Refined universality for critical {KCM}: upper bounds},
  preprint, arXiv:2104.02329.

\bibitem{HarMa}
I.~Hartarsky and L.~Mar\^ech\'e, \emph{Refined universality for critical {KCM}:
  lower bounds}, Combin. Probab. Comput. \textbf{31} (2022), no.~5, 879--906.

\bibitem{HMT2}
I.~Hartarsky, L.~Mar\^ech\'e, and C.~Toninelli, \emph{Universality for critical
  {KCM}: infinite number of stable directions}, Probab. Theory Related Fields
  \textbf{178} (2020), 289--326.

\bibitem{HMT1}
I.~Hartarsky, F.~Martinelli, and C.~Toninelli, \emph{Universality for critical
  {KCM}: finite number of stable directions}, Ann. Probab. \textbf{49} (2021),
  no.~5, 2141--2174.

\bibitem{HMT3}
\bysame, \emph{Sharp threshold for the {FA-2f} kinetically constrained model},
  Prob. Theory Related Fields \textbf{185} (2023), no.~3, 993--1037.

\bibitem{HM}
I.~Hartarsky and R.~Morris, \emph{The second term for two-neighbour bootstrap
  percolation in two dimensions}, Trans. Amer. Math. Soc. \textbf{372} (2019),
  6465--6505.

\bibitem{HS22}
I.~Hartarsky and R.~Szab\'o, \emph{Subcritical bootstrap percolation via {T}oom
  contours}, Electron. Commun. Probab. \textbf{27} (2022), 13pp.

\bibitem{H04}
F.~{\noopsort{Hollander}}{den Hollander}, \emph{Metastability under stochastic
  dynamics}, Stochastic Process. Appl. \textbf{114} (2004), no.~1, 1--26.

\bibitem{Hol}
A.E. Holroyd, \emph{Sharp metastability threshold for two-dimensional bootstrap
  percolation}, Probab. Theory Related Fields \textbf{125} (2003), no.~2,
  195--224.

\bibitem{HMod}
\bysame, \emph{The metastability threshold for modified bootstrap percolation
  in $d$ dimensions}, Electron. J. Probab. \textbf{11} (2006), 418--433.

\bibitem{KO}
R.~Koteck\'y and E.~Olivieri, \emph{Droplet dynamics for asymmetric {I}sing
  model}, J. Stat. Phys. \textbf{70} (1993), 1121--1148.

\bibitem{MaMaT}
L.~Mar\^ech\'e, F.~Martinelli, and C.~Toninelli, \emph{Exact asymptotics for
  {D}uarte and supercritical rooted kinetically constrained models}, Ann.
  Probab. \textbf{48} (2020), 317--342.

\bibitem{M99}
F.~Martinelli, \emph{Lectures on {G}lauber {D}ynamics for {D}iscrete {S}pin
  {M}odels}, Lectures on Probability Theory and Statistics (P.~Bernard, ed.),
  Lecture Notes in Mathematics, vol. 1717, Springer, 1999.

\bibitem{MMT}
F.~Martinelli, R.~Morris, and C.~Toninelli, \emph{Universality results for
  kinetically constrained spin models in two dimensions}, Comm. Math. Phys.
  \textbf{369} (2019), 761--809.

\bibitem{MOS91}
F.~Martinelli, E.~Olivieri, and E.~Scoppola, \emph{On the {S}wendsen and {W}ang
  dynamics. {II}: {C}ritical droplets and homogeneous nucleation at low
  temperature}, J. Stat. Phys. \textbf{62} (1991), 135--159.

\bibitem{MT}
F.~Martinelli and C.~Toninelli, \emph{Towards a universality picture for the
  relaxation to equilibrium of kinetically constrained models}, Ann. Probab.
  \textbf{47} (2019), 324--361.

\bibitem{MGlauber}
R.~Morris, \emph{Zero-temperature {G}lauber dynamics on $\mathbb{Z}^d$},
  Probab. Theory Related Fields \textbf{149} (2011), no.~3, 417--434.

\bibitem{M17}
\bysame, \emph{Bootstrap percolation and other automata}, European J. Combin.
  \textbf{66} (2017), 250--263.

\bibitem{Mount92}
T.S. Mountford, \emph{Rates for the probability of large cubes being
  non-internally spanned in modified bootstrap percolation}, Probab. Theory
  Related Fields \textbf{93} (1992), no.~2, 159--167.

\bibitem{Mount}
\bysame, \emph{Critical length for semi-oriented bootstrap percolation},
  Stochastic Process. Appl. \textbf{56} (1995), 185--205.

\bibitem{NS91}
E.~Neves and R.~Schonmann, \emph{Critical droplets and metastability for a
  {G}lauber dynamics at very low temperatures}, Comm. Math. Phys. \textbf{137}
  (1991), 209--230.

\bibitem{NS92}
\bysame, \emph{Behavior of droplets for a class of {G}lauber dynamics at very
  low temperature}, Probab. Theory Related Fields \textbf{91} (1992), 331--354.

\bibitem{RS03}
F.~Ritort and P.~Sollich, \emph{Glassy dynamics of kinetically constrained
  models}, Adv. Phys. \textbf{52} (2003), no.~4, 219--342.

\bibitem{Sch3}
R.~Schonmann, \emph{Critical points of two-dimensional bootstrap
  percolation-like cellular automata}, J. Stat. Phys. \textbf{58} (1990),
  no.~5, 1239--1244.

\bibitem{Sch91}
\bysame, \emph{The pattern of escape from metastability of a stochastic {I}sing
  model}, Comm. Math. Phys. \textbf{147} (1991), 231--240.

\bibitem{Sch1}
\bysame, \emph{On the behavior of some cellular automata related to bootstrap
  percolation}, Ann. Probab. \textbf{20} (1992), no.~1, 174--193.

\bibitem{Sch4}
\bysame, \emph{Slow droplet-driven relaxation of stochastic {I}sing models in
  the vicinity of the phase coexistence region}, Comm. Math. Phys. \textbf{161}
  (1994), 1--49.

\bibitem{Sch98}
\bysame, \emph{Metastability and the {I}sing model}, Proceedings of the
  International Congress of Mathematicians, Berlin 1998 (G.~Fischer and
  U.~Rehmann, eds.), Doc. Math., Extra Vol. ICM III, 1998, pp.~173--181.

\bibitem{SS98}
R.~Schonmann and S.~Shlosman, \emph{Wulff droplets and the metastable
  relaxation of kinetic {I}sing models}, Comm. Math. Phys. \textbf{194} (1998),
  no.~2, 389--462.

\bibitem{S16}
E.~Steinitz, \emph{Bedingt konvergente {Reihen} und konvexe {Systeme}}, J.
  Reine Angew. Math. \textbf{143} (1913), 128--176.

\bibitem{T22}
C.~Toninelli, \emph{Interacting particle systems with kinetic constraints},
  {SpringerBriefs in Mathematical Physics}, in preparation.

\end{thebibliography}

\end{document}